\documentclass{amsart}

\usepackage{amssymb}
\usepackage{amsmath}
\usepackage{amsthm}
\usepackage{amsfonts}

\usepackage{a4wide}
\usepackage{longtable}
\usepackage{multirow}
\usepackage{xcolor}

\newtheorem{lemma}{Lemma}
\newtheorem{proposition}{Proposition}
\newtheorem{theorem}{Theorem}
\newtheorem{corollary}{Corollary}

\theoremstyle{definition}
\newtheorem{definition}{Definition}
\newtheorem{example}{Example}

\newtheorem{remark}{Remark}

\numberwithin{equation}{section}

\usepackage{hyperref}
\hypersetup{
    colorlinks=true, 
    citecolor=blue,  
    urlcolor=blue,   
}

\begin{document}

\title[Cyclic cubic number fields]{Group theoretic approach to cyclic cubic fields}

\author{Siham Aouissi}
\address{\emph{Algebraic theories and applications research team (ATA). Ecole Normale Sup\`erieure of Moulay Ismail University (ENS-UMI).} Address: ENS, B.P. 3104, Toulal, Mekn\`es, Morocco}
\email{s.aouissi@umi.ac.ma}

\author{Daniel C. Mayer}
\address{Naglergasse 53\\8010 Graz\\Austria}
\email{algebraic.number.theory@algebra.at}
\urladdr{http://www.algebra.at}

\thanks{Research of second author supported by the Austrian Science Fund (FWF): projects J0497-PHY, P26008-N25, and by the Research Executive Agency of the European Union (EUREA)}



\date{Wednesday, 27 December 2023}

\begin{abstract}
Let \((k_\mu)_{\mu=1}^4\) be a quartet
of cyclic cubic number fields
sharing a common conductor \(c=pqr\)
divisible by exactly three prime(power)s \(p,q,r\).
For those components of the quartet whose \(3\)-class group
\(\mathrm{Cl}_3(k_\mu)\simeq(\mathbb{Z}/3\mathbb{Z})^2\)
is elementary bicyclic,
the automorphism group
\(\mathfrak{M}=\mathrm{Gal}(\mathrm{F}_3^2(k_\mu)/k_\mu)\)
of the maximal metabelian unramified \(3\)-extension of \(k_\mu\)
is determined by conditions
for cubic residue symbols between \(p,q,r\)
and for ambiguous principal ideals in subfields
of the common absolute \(3\)-genus field \(k^\ast\)
of all \(k_\mu\).
With the aid of the relation rank \(d_2(\mathfrak{M})\),
it is decided
whether \(\mathfrak{M}\) coincides with the Galois group
\(\mathfrak{G}=\mathrm{Gal}(\mathrm{F}_3^\infty(k_\mu)/k_\mu)\)
of the maximal unramified pro-\(3\)-extension of \(k_\mu\).
\end{abstract}

\maketitle

\hypersetup{linkcolor=blue}
\tableofcontents


\section{Introduction}
\label{s:Intro}

\noindent
Let \(k\) be a cyclic cubic number field,
that is, an abelian extension of the rational number field \(\mathbb{Q}\)
with degree \(\lbrack k:\mathbb{Q}\rbrack=3\) and
some positive integer conductor \(c>1\) (see \S\
\ref{ss:Multiplicity}).
In 1973, Georges Gras
\cite{Gr1973}
determined the rank \(\varrho=\varrho(k)\)
of the \(3\)-class group \(\mathrm{Cl}_3(k)\)
in dependence on the number \(t\) of
prime(power) divisors \(q_1,\ldots,q_t\) of \(c\)
and on the cubic residue symbols
\(\left(\frac{q_i}{q_j}\right)_3\) for \(i\ne j\).
For mutual cubic residues,
\(\left(\frac{q_i}{q_j}\right)_3=\left(\frac{q_j}{q_i}\right)_3=1\),
we write
\(q_i\leftrightarrow q_j\),
otherwise
\(q_i\not\leftrightarrow q_j\).

It turned out that
\(\varrho=0\) for \(t=1\),
and \(\varrho=1\) if \(t=2\) and \(q_1\not\leftrightarrow q_2\).
So in the former case,
the maximal unramified pro-\(3\)-extension \(\mathrm{F}_3^\infty(k)\) of \(k\)
is the base field \(k\) itself,
and in the latter case,
it is the Hilbert \(3\)-class field \(\mathrm{F}_3^1(k)\) of \(k\),
in fact, \(\lbrack\mathrm{F}_3^1(k):k\rbrack=3\),
since \(\varrho=1\) iff
\(\mathrm{Cl}_3(k)\simeq \mathbb{Z}/3\mathbb{Z}\) is elementary cyclic.
If \(t=2\) and \(q_1\leftrightarrow q_2\), then \(\varrho=2\),
\(\mathrm{Cl}_3(k)\) is bicyclic, but may be non-elementary (singular).

In 1995, Ayadi
\cite{Ay1995}
proved that there are only two possibilities for the Galois group
\(\mathfrak{G}=\mathrm{Gal}(\mathrm{F}_3^\infty(k)/k)\)
of the \(3\)-class field tower of \(k\) with length \(\ell_3(k)\),
when \(t=2\), \(q_1\leftrightarrow q_2\), and 
\(\mathrm{Cl}_3(k)\simeq(\mathbb{Z}/3\mathbb{Z})^2\)
is elementary bicyclic (regular), namely,
in the notation of
\cite{BEO2005},
either \(\mathfrak{G}\simeq\mathrm{SmallGroup}(9,2)\simeq(\mathbb{Z}/3\mathbb{Z})^2\) is abelian
or \(\mathfrak{G}\simeq\mathrm{SmallGroup}(27,4)\) is the extra special \(3\)-group with exponent \(9\).

The impact of \(t\) and \(q_1,\ldots,q_t\) on the tower group \(\mathfrak{G}\)
and its metabelianization \(\mathfrak{M}=\mathfrak{G}/\mathfrak{G}^{\prime\prime}\),
i.e., the group \(\mathfrak{M}=\mathrm{Gal}(\mathrm{F}_3^2(k)/k)\)
of the second Hilbert \(3\)-class field \(\mathrm{F}_3^2(k)\) of \(k\),
is shown in Table
\ref{tbl:Impact}.


\renewcommand{\arraystretch}{1.1}

\begin{table}[ht]
\caption{Known and unknown impact of \(t\) and \(q_1,\ldots,q_t\) on \(\varrho(k)\) and \(\mathfrak{M}\), \(\mathfrak{G}\)}
\label{tbl:Impact}
\begin{center}
{\small
\begin{tabular}{|c|c||c|c|c|c|c|}
\hline
 \(t\)   & conditions                                                 & \(\varrho(k)\)       & \(\mathrm{F}_3^\infty(k)\) & \(\mathfrak{M}\)         & \(\mathfrak{G}\) & \(\ell_3(k)\) \\
\hline
 \(t=1\) &                                                            & \(\varrho=0\)        & \(=k\)                     & \(=1\)                         & \(=1\)            & \(=0\) \\
 \(t=2\) & \(q_1\not\leftrightarrow q_2\)                             & \(\varrho=1\)        & \(=\mathrm{F}_3^1(k)\)     & \(=\mathbb{Z}/3\mathbb{Z}\)    & \(=\mathfrak{M}\) & \(=1\) \\
\hline
 \(t=2\) & \(q_1\leftrightarrow q_2\), \(\mathrm{Cl}_3(k)\) elem.     & \(\varrho=2\)        & \(=\mathrm{F}_3^1(k)\)     & \(=\mathrm{SmallGroup}(9,2)\)  & \(=\mathfrak{M}\) & \(=1\) \\
         &                                                            & or                   & \(=\mathrm{F}_3^2(k)\)     & \(=\mathrm{SmallGroup}(27,4)\) & \(=\mathfrak{M}\) & \(=2\) \\
\hline
 \(t=2\) & \(q_1\leftrightarrow q_2\), \(\mathrm{Cl}_3(k)\) non-elem. & \(\varrho=2\)        & \(\ge\mathrm{F}_3^2(k)\)   & ? & ? & \(\ge 2\) \\
 \(t=3\) &                                                            & \(2\le\varrho\le 4\) & \(\ge\mathrm{F}_3^1(k)\)   & ? & ? & \(\ge 1\) \\
\hline
\end{tabular}
}
\end{center}
\end{table}


However, according to Gras
\cite{Gr1973},
\(\varrho=2\) is also possible for \(t=3\), and,
according to Ayadi
\cite{Ay1995},
\(\varrho=2\) iff
\(\mathrm{Cl}_3(k)\simeq(\mathbb{Z}/3\mathbb{Z})^2\) is elementary bicyclic,
when \(t=3\).

For this situation \(t=3\), \(c=pqr\), \(\varrho=2\), \(\mathrm{Cl}_3(k)\simeq(\mathbb{Z}/3\mathbb{Z})^2\),
the present article identifies the Galois group
\(\mathfrak{M}=\mathrm{Gal}(\mathrm{F}_3^2(k)/k)\)
in dependence on the cubic residue symbols between \(p,q,r\).
The crucial techniques are based on the lucky coincidence
that the four unramified cyclic extensions of degree
\(\lbrack E_i:k\rbrack=3\), \(1\le i\le 4\),
can always be found among the \(13\) bicyclic bicubic subfields
\(B_1,\ldots,B_{13}\)
of the absolute \(3\)-genus-field \(k^\ast\) of \(k\),
for which Parry
\cite{Pa1990}
has established a useful class number relation
and a structure theory of the unit group.
With the aid of
the relation rank \(d_2(\mathfrak{M})\le 4\) or \(d_2(\mathfrak{M})\ge 5\),
it is decided
whether \(\mathfrak{M}\) coincides with the tower group
\(\mathfrak{G}\) or not.

The examination of cyclic cubic fields \(k\) with
\(\varrho=3\) and elementary tricyclic
\(\mathrm{Cl}_3(k)\simeq(\mathbb{Z}/3\mathbb{Z})^3\)
is reserved for a future paper,
since among the \(13\) unramified cyclic extensions of degree
\(\lbrack E_i:k\rbrack=3\), \(1\le i\le 13\),
only four are bicyclic bicubic,
and the remaining nine \(E_i\) arise in three triplets
of pairwise isomorphic non-Galois nonic fields.
Similarly, \(\mathrm{Cl}_3(k)\) non-elementary for
\(t=2\) and \(p\leftrightarrow q\).

The present work illuminates Ayadi's doctoral thesis
\cite{Ay1995}
from the perspective of group theory,
and completely clarifies the question mark \lq\lq ?\rq\rq\ 
for the group \(\mathfrak{M}\) in the
last row of Table
\ref{tbl:Impact},
partially also the \lq\lq ?\rq\rq\ for the group \(\mathfrak{G}\),
provided that \(\mathrm{Cl}_3(k)\simeq(\mathbb{Z}/3\mathbb{Z})^2\)
is elementary bicyclic. 


\section{Construction of cyclic fields of odd prime degree}
\label{s:Construction}

\subsection{Multiplicity of conductors and discriminants}
\label{ss:Multiplicity}

\noindent
For a fixed odd prime number \(\ell\ge 3\),
let \(k\) be a \textit{cyclic number field} of degree \(\ell\),
that is,
\(k/\mathbb{Q}\) is a Galois extension of degree
\(\lbrack k:\mathbb{Q}\rbrack=\ell\)
with absolute automorphism group 
\(\mathrm{Gal}(k/\mathbb{Q})=\langle\sigma\mid\sigma^\ell=1\rangle\).
According to the \textbf{Theorem of Kronecker, Weber and Hilbert}
on abelian extensions of the rational number field \(\mathbb{Q}\),
the \textit{conductor} \(c\) of \(k\)
is the smallest positive integer such that
\(k=k_c\) is contained in the cyclotomic field
\(K=\mathbb{Q}(\zeta_c)\),
where \(\zeta_c=\exp(2\pi\sqrt{-1}/c)\) denotes a primitive \(c\)-th root of unity,
more precisely, in the \(\ell\)-\textit{ray class field modulo} \(c\) of \(\mathbb{Q}\),
denoted by \(\mathrm{F}_{\ell,c}(\mathbb{Q})\),
which lies in the maximal real subfield \(K^+=\mathbb{Q}(\zeta_c+\zeta_c^{-1})\)
of \(K=\mathbb{Q}(\zeta_c)\).

\begin{theorem}
\label{thm:Multiplicity}
The \textbf{conductor} of a cyclic field of odd prime degree \(\ell\)
has the shape
\(c=\ell^e\cdot q_1\cdots q_\tau\),
where
\(e\in\lbrace 0,2\rbrace\)
and the \(q_i\) are pairwise distinct prime numbers
\(q_i\equiv +1\,(\mathrm{mod}\,\ell)\),
for \(1\le i\le\tau\).
The \textbf{discriminant} of \(k=k_c\) is the perfect \((\ell-1)\)-th power
\(d_k=c^{\ell-1}\),
and the number of rational primes which are (totally) ramified in \(k\)
is given by
\begin{equation}
\label{eqn:Ramification}
t:=
\begin{cases}
\tau   & \text{ if } e=0\ (\ell \text{ is unramified in } k),\\
\tau+1 & \text{ if } e=2\ (\ell \text{ is ramified in } k).
\end{cases}
\end{equation}
In the last case, we formally put \(q_{\tau+1}:=\ell^2\).
The number of non-isomorphic cyclic number fields
\(k_{c,1},\ldots,k_{c,m}\)
of degree \(\ell\), sharing the common conductor \(c\),
is given by the \textbf{multiplicity formula}
\begin{equation}
\label{eqn:Multiplicity}
m=m(c)=(\ell-1)^{t-1}.
\end{equation}
\end{theorem}

\begin{proof}
See \cite[p. 831]{Ma1992}.
\end{proof}


\subsection{Construction as ray class fields}
\label{ss:RayClassFields}

\noindent
For the construction of all cyclic number fields \(k=k_c\) of degree \(\ell\)
with ascending conductors \(b\le c\le B\)
between an assigned lower bound \(b\) and upper bound \(B\)
by means of the computational algebra system Magma
\cite{MAGMA2023},
the class field theoretic routines by Fieker
\cite{Fi2001}
can be used without the need of preparing a list of
suitable generating polynomials of \(\ell\)-th degree.
The big advantage of this technique is that
the cyclic number fields of degree \(\ell\) are produced as a \textit{multiplet}
\((k_{c,1},\ldots,k_{c,m})\)
of pairwise non-isomorphic fields sharing the common conductor \(c\)
with \textit{multiplicity} \(m\in\lbrace 1,\ell-1,(\ell-1)^2,(\ell-1)^3,\ldots\rbrace\)
in dependence on the number \(t\in\lbrace 1,2,3,4,\ldots\rbrace\)
of primes dividing the conductor \(c\),
according to Formula
\eqref{eqn:Multiplicity}.
Our algorithms for the construction,
and statistics of \(\ell\)-class groups,
have been presented in
\cite[Alg. 1--3, pp. 4--7, Tbl. 1.1--1.6, pp. 7--11]{Ma2022}.
From now on, let \(\ell=3\),
for the remainder of this article.


\section{Arithmetic of cyclic cubic fields}
\label{s:Arithmetic}

\noindent
Generally, \(t\) denotes the number of prime divisors
of the conductor \(c\) of a cyclic cubic number field \(k\),
and \(\varrho(k)=\varrho_3(k)\) denotes the rank
\(\dim_{\mathbb{F}_3}(\mathrm{Cl}_3(k)/\mathrm{Cl}_3(k)^3)\)
of the \(3\)-class group
\(\mathrm{Cl}_3(k)=\mathrm{Syl}_3\mathrm{Cl}(k)\).
In formulas concerning principal factors (\S\
\ref{ss:AmbiguousIdeals}),
the prime power conductor \(3^2\) must be replaced by \(3\).


\subsection{Rank of 3-class groups}
\label{ss:3ClassRank}

\noindent
Since the rank \(\varrho_3(k)\) of the \(3\)-class group \(\mathrm{Cl}_3(k)\)
of a cyclic cubic field \(k\)
depends on the mutual cubic residue conditions
between the prime(power) divisors \(q_1,\ldots,q_t\) of the conductor \(c\),
Gras
\cite[pp. 21--22]{Gr1973}
has introduced directed graphs with \(t\) vertices \(q_1,\ldots,q_t\)
whose directed edges \(q_i\to q_j\) describe values of cubic residue symbols.
We use a simplified notation of these graphs, fitting in a single line,
but occasionally requiring the repetition of a vertex.


\begin{definition}
\label{dfn:CubicResidues}
Let \(\zeta_3\) be a fixed primitive third root of unity.
For each pair \((q_i,q_j)\) with \(1\le i\ne j\le t\),
the value of the \textit{cubic residue symbol}
\(\left(\frac{q_i}{q_j}\right)_3=\zeta_3^{a_{ij}}\)
is determined uniquely by the integer \(a_{ij}\in\lbrace -1,0,1\rbrace\).
Let a \textit{directed edge} \(q_i\to q_j\) be defined if and only if \(\left(\frac{q_i}{q_j}\right)_3=1\),
that is, \(q_i\) is a cubic residue modulo \(q_j\) (and thus \(a_{ij}=0\)).
The \textbf{combined cubic residue symbol} \(\lbrack q_1,\ldots,q_t\rbrack_3:=\)

\begin{equation}
\label{eqn:CubicResidues}
\Biggl\lbrace q_i\to q_j\Biggm\vert i\ne j,\left(\frac{q_i}{q_j}\right)_3=1\Biggr\rbrace
\bigcup
\Biggl\lbrace q_i\Biggm\vert (\forall j\ne i)\,\left(\frac{q_i}{q_j}\right)_3\ne 1,\left(\frac{q_j}{q_i}\right)_3\ne 1 \Biggr\rbrace
\end{equation}

\noindent
where the subscripts \(i\) and \(j\) run from \(1\) to \(t\),
is defined as the union of the set of all \textit{directed edges}
which occur in the graph associated with \(q_1,\ldots,q_t\) in the sense of Gras,
and the set of all \textit{isolated vertices}.
For \(t=3\),
we additionally need the invariant \(\delta:=a_{12}a_{23}a_{31}-a_{13}a_{32}a_{21}\)
in order to distinguish two subcases of the case with three isolated vertices.
\end{definition}


\begin{theorem}
\label{thm:3ClassRank}
(Rank Distribution, \textbf{G. Gras, 1973},
\cite[Prop. VI.5, pp. 21--22]{Gr1973}.) \\
Let \(k\) be a cyclic cubic field of conductor \(c=q_1\cdots q_t\) with \(1\le t\le 3\).
We indicate \textbf{mutual cubic residues} simply by writing \(q_1\leftrightarrow q_2\) instead of \(q_1\to q_2\to q_1\).

\begin{itemize}

\item
If \(t=1\), then \(m=1\), \(k\) forms a singlet, \(\lbrack q_1\rbrack_3=\lbrace q_1\rbrace\), and \(\varrho(k)=0\).

\item
If \(t=2\), then \(m=2\), \(k\) is member of a doublet \((k_1,k_2)\), and there arise two possibilities.

\begin{enumerate}
\item
\((\varrho(k_1),\varrho(k_2))=(1,1)\), if
\begin{equation}
\label{eqn:Rank1}
\lbrack q_1,q_2\rbrack_3=
\begin{cases}
\lbrace q_1,q_2\rbrace & \textbf{, Graph } \mathbf{1}, \text{ or} \\
\lbrace q_i\to q_j\rbrace & \textbf{, Graph } \mathbf{2}, \text{ with } i\ne j.
\end{cases}
\end{equation}
\item
\((\varrho(k_1),\varrho(k_2))=(2,2)\), if
\begin{equation}
\label{eqn:Rank2}
\lbrack q_1,q_2\rbrack_3=\lbrace q_1\leftrightarrow q_2\rbrace \textbf{, Graph }\mathbf{3}.
\end{equation}
\end{enumerate}

\item
If \(t=3\), then \(m=4\), \(k\) is member of a quartet \((k_1,\ldots,k_4)\), and there arise five cases.

\begin{enumerate}
\item
\((\varrho(k_1),\ldots,\varrho(k_4))=(2,2,2,2)\), called \textbf{Category} \(\mathbf{III}\), if
\begin{equation}
\label{eqn:Category3}
\lbrack q_1,q_2,q_3\rbrack_3=
\begin{cases}
\lbrace q_1,q_2,q_3;\delta\not\equiv 0\,(\mathrm{mod}\,3)\rbrace & \textbf{, Graph }\mathbf{1},\text{ or} \\
\lbrace q_i\rightarrow q_j;q_l\rbrace & \textbf{, Graph }\mathbf{2},\text{ or} \\
\lbrace q_i\rightarrow q_j\rightarrow q_l\rbrace & \textbf{, Graph }\mathbf{3},\text{ or} \\
\lbrace q_i\rightarrow q_j\rightarrow q_l\rightarrow q_i\rbrace & \textbf{, Graph }\mathbf{4},\text{ or} \\
\lbrace q_i\leftrightarrow q_j;q_l\rbrace & \textbf{, Graph }\mathbf{5},\text{ or} \\
\lbrace q_i\leftrightarrow q_j\rightarrow q_l\rbrace & \textbf{, Graph }\mathbf{6},\text{ or} \\
\lbrace q_i\leftrightarrow q_j\leftarrow q_l\rbrace & \textbf{, Graph }\mathbf{7},\text{ or} \\
\lbrace q_l\rightarrow q_i\leftrightarrow q_j\leftarrow q_l\rbrace & \textbf{, Graph }\mathbf{8},\text{ or} \\
\lbrace q_l\rightarrow q_i\leftrightarrow q_j\rightarrow q_l\rbrace & \textbf{, Graph }\mathbf{9} \\
\end{cases}
\end{equation}
with \(i,j,l\) pairwise distinct.
\item
\((\varrho(k_1),\ldots,\varrho(k_4))=(3,2,2,2)\), called \textbf{Category} \(\mathbf{I}\), if
\begin{equation}
\label{eqn:Category1}
\lbrack q_1,q_2,q_3\rbrack_3=
\begin{cases}
\lbrace q_1,q_2,q_3;\delta\equiv 0\,(\mathrm{mod}\,3)\rbrace & \textbf{, Graph }\mathbf{1},\text{ or} \\
\lbrace q_i\leftarrow q_j\rightarrow q_l\rbrace & \textbf{, Graph }\mathbf{2}
\end{cases}
\end{equation}
with \(i,j,l\) pairwise distinct.
\item
\((\varrho(k_1),\ldots,\varrho(k_4))=(3,3,2,2)\), called \textbf{Category} \(\mathbf{II}\), if
\begin{equation}
\label{eqn:Category2}
\lbrack q_1,q_2,q_3\rbrack_3=
\begin{cases}
\lbrace q_i\rightarrow q_j\leftarrow q_l\rbrace & \textbf{, Graph }\mathbf{1},\text{ or} \\
\lbrace q_i\rightarrow q_j\leftarrow q_l\rightarrow q_i\rbrace & \textbf{, Graph }\mathbf{2}
\end{cases}
\end{equation}
with \(i,j,l\) pairwise distinct.
\item
\((\varrho(k_1),\ldots,\varrho(k_4))=(3,3,3,3)\), called \textbf{Category} \(\mathbf{IV}\), if
\begin{equation}
\label{eqn:Rank3}
\lbrack q_1,q_2,q_3\rbrack_3=
\begin{cases}
\lbrace q_i\leftarrow q_j\leftrightarrow q_l\rightarrow q_i\rbrace & \textbf{, Graph }\mathbf{1},\text{ or} \\
\lbrace q_i\leftrightarrow q_j\leftrightarrow q_l\rbrace & \textbf{, Graph }\mathbf{2},\text{ or} \\
\lbrace q_i\leftrightarrow q_j\leftrightarrow q_l\rightarrow q_i\rbrace & \textbf{, Graph }\mathbf{3}
\end{cases}
\end{equation}
with \(i,j,l\) pairwise distinct.
\item
\((\varrho(k_1),\ldots,\varrho(k_4))=(4,4,4,4)\), called \textbf{Category} \(\mathbf{V}\), if
\begin{equation}
\label{eqn:Rank4}
 \lbrack q_1,q_2,q_3\rbrack_3=\lbrace q_1\leftrightarrow q_2\leftrightarrow q_3\leftrightarrow q_1\rbrace.
\end{equation}
\end{enumerate}

\end{itemize}

\end{theorem}

\begin{proof}
See
\cite[Prp. VI.5, pp. 21--22]{Gr1973}.
Multiplicities \(m\in\lbrace 1,2,4\rbrace\) are taken from Theorem
\ref{thm:Multiplicity}.
\end{proof}


\begin{remark}
\label{rmk:3ClassRank}
Ayadi introduced categories in
\cite[pp. 45--47]{Ay1995}.
He investigated the cases
\(t=2\), Formula
\eqref{eqn:Rank2},
and
\(t=3\), Formulae
\eqref{eqn:Category3},
\eqref{eqn:Category1},
\eqref{eqn:Category2},
in Theorem
\ref{thm:3ClassRank}.
For \(t=3\), he denoted
the nine subcases of Formula
\eqref{eqn:Category3}
by Graph 1,2,3,4,5,6,7,8,9 of Category III,
the two subcases of Formula
\eqref{eqn:Category1}
by Graph 1,2 of Category I,
and the two subcases of Formula
\eqref{eqn:Category2}
by Graph 1,2 of Category II.
For the Categories I and II,
Ayadi did \textbf{not} study the fields with
\(3\)-class rank \(\varrho_3(k_{\mu})=3\),
\(1\le\mu\le 4\).
Our algorithms for the classification
by categories and graphs,
and their statistics,
have been presented in
\cite[Alg. 4--5, Tbl. 2.1, pp. 15--19]{Ma2022}.
\end{remark}


For \(t=3\),
we also write briefly
\(p=q_1\), \(q=q_2\) and \(r=q_3\)
for the prime(power)s dividing the conductor \(c=pqr\).
Graph \(1\) of Category \(\mathrm{I}\)
with symbol
\(\lbrack p,q,r\rbrack_3=\lbrace p,q,r;\delta\equiv 0\,(\mathrm{mod}\,3)\rbrace\)
and Graph \(1\) of Category \(\mathrm{III}\)
with symbol
\(\lbrack p,q,r\rbrack_3=\lbrace p,q,r;\delta\not\equiv 0\,(\mathrm{mod}\,3)\rbrace\)
are the only two situations
without any trivial cubic residue conditions between \(p,q,r\).
We show the impact of the \(\delta\)-invariant.

\begin{lemma}
\label{lem:Delta}
Consider three cubic residue symbols for products of two primes,
\(\left(\frac{qr}{p}\right)_3\),
\(\left(\frac{pr}{q}\right)_3\),
\(\left(\frac{pq}{r}\right)_3\)
with respect to triviality, i.e., being equal to \(1\).

If \(\delta\equiv 0\),
then zero or two of the symbols are trivial.

If \(\delta\not\equiv 0\),
then one or three of the symbols are trivial.
\end{lemma}

\begin{proof}
For each of the two triplets
\((a_{12},a_{23},a_{31})\) and \((a_{32},a_{13},a_{21})\)
of exponents in Definition
\ref{dfn:CubicResidues},
there are \(2^3=8\) combinatorial possibilities.
The product of the components is
\(+1\) if zero or two components are negative, and it is
\(-1\) if one or three components are negative.

For Graph \(\mathrm{I}.1\) with \(\delta\equiv 0\),
triplets with equal product must be combined.
Consequently, for each choice of a fixed first triplet,
one of the four admissible second triplets
(namely \((a_{32},a_{13},a_{21})=(a_{12},a_{23},a_{31})\))
produces no trivial symbol,
and three of the second triplets produce two trivial symbols each.

For Graph \(\mathrm{III}.1\) with \(\delta\not\equiv 0\),
triplets with distinct product must be combined.
Consequently, for each choice of a fixed first triplet,
one of the four admissible second triplets
(namely \((a_{32},a_{13},a_{21})=-(a_{12},a_{23},a_{31})\))
produces three trivial symbols,
and three of the second triplets produce a single trivial symbol each.
\end{proof}


\subsection{Ambiguous principal ideals}
\label{ss:AmbiguousIdeals}

\noindent
The number of \textit{primitive ambiguous ideals} of a cyclic cubic field \(k\),
which are invariant under \(\mathrm{Gal}(k/\mathbb{Q})=\langle\sigma\rangle\),
increases with the number \(t\) of prime factors of the conductor \(c\).
According to Hilbert's Theorem \(93\), the number is given by
\begin{equation}
\label{eqn:AmbiguousIdeals}
\#\left(\mathcal{I}_k^{\langle\sigma\rangle}/\mathcal{I}_{\mathbb{Q}}\right)=3^t.
\end{equation}
However, the number of \textit{primitive ambiguous} \textbf{principal} \textit{ideals} of \(k\)
is a fixed invariant of all cyclic cubic fields, regardless of the number \(t\).

\begin{theorem}
\label{thm:AmbiguousPrincipalIdeals}
The number of ambiguous principal ideals of any cyclic cubic field \(k\) is given by
\begin{equation}
\label{eqn:AmbiguousPrincipalIdeals}
\#\left(\mathcal{P}_k^{\langle\sigma\rangle}/\mathcal{P}_{\mathbb{Q}}\right)=3.
\end{equation}
\end{theorem}

\begin{proof}
The well-known theorem on the \textbf{Herbrand quotient}
of the unit group \(U_k\) of \(k\)
as a Galois module over the group \(\mathrm{Gal}(k/\mathbb{Q})=\langle\sigma\rangle\),
which can be expressed by abstract cohomology groups
\(\#\mathrm{H}^{-1}(\langle\sigma\rangle,U_k)/\#\hat{\mathrm{H}}^{0}(\langle\sigma\rangle,U_k)=\lbrack k:\mathbb{Q}\rbrack\),
can also be stated more ostensively as
\(\#\left(\mathcal{P}_k^{\langle\sigma\rangle}/\mathcal{P}_{\mathbb{Q}}\right)
=\#\left(E_{k/\mathbb{Q}}/U_k^{1-\sigma}\right)
=\lbrack k:\mathbb{Q}\rbrack\cdot\#\left(U_{\mathbb{Q}}/N_{k/\mathbb{Q}}(U_k)\right)=3\),
since the unit norm index is given by \(\left(U_{\mathbb{Q}}:N_{k/\mathbb{Q}}(U_k)\right)=1\).
Here, \(E_{k/\mathbb{Q}}=
\lbrace\varepsilon\in U_k\mid N_{k/\mathbb{Q}}(\varepsilon)=1\rbrace\)
are the relative units.
\end{proof}

\noindent
Consequently, if we speak about a \textit{non-trivial primitive ambiguous principal ideal} of \(k\),
then we either mean \((\alpha)=\alpha\mathcal{O}_k\) or \((\alpha^2/b)=(\alpha^2/b)\mathcal{O}_k\),
where \(\mathcal{P}_k^{\langle\sigma\rangle}/\mathcal{P}_{\mathbb{Q}}=\lbrace 1,(\alpha),(\alpha^2/b)\rbrace\).
The norms of these two elements are divisors of the square \(c^2=q_1^2\cdots q_t^2\)
of the conductor \(c\) of \(k\),
where \(q_t\) must be replaced by \(3\) if \(q_t=9\).
When \(N_{k/\mathbb{Q}}(\alpha)=a\cdot b^2\)
with square free coprime integers \(a,b\), then
\(N_{k/\mathbb{Q}}(\alpha^2/b)=a^2\cdot b^4/b^3=a^2\cdot b\).
It follows that both norms are cube free integers.


\begin{definition}
\label{dfn:PrincipalFactor}
The minimum of the two norms
of non-trivial primitive \textit{ambiguous principal ideals}
\((\alpha),(\alpha^2/b)\)
of a cyclic cubic field \(k\) is called the \textbf{principal factor}
(of the discriminant \(d_k=c^2\)) of the field \(k\), denoted by
\(A(k):=\min\lbrace a\cdot b^2,a^2\cdot b\rbrace\), that is
\begin{equation}
\label{eqn:PrincipalFactor}
A(k)=
\begin{cases}
a\cdot b^2 & \text{ if } b<a, \\
a^2\cdot b & \text{ if } a<b.
\end{cases}
\end{equation}
\end{definition}

\noindent
Ayadi 
\cite[Rem. 2.6, p. 18]{Ay1995},
\cite{AAI2001}
speaks about the \textit{Parry constant} or \textit{Parry invariant} \(A(k)\) of \(k\),
and Derhem
\cite{Dh2002}
calls \(A(k)=N_{k/\mathbb{Q}}(R)\)
with \(R=1+\varepsilon+\varepsilon^{1+\sigma}\),
\(\varepsilon=R^{1-\sigma}\),
the \textit{Kummer resolvent} of \(k\),
when \(U_k=\langle -1,\varepsilon,\varepsilon^\sigma\rangle\)
as a \(\langle\sigma\rangle\)-module
is generated by \(-1\) and the fundamental unit \(\varepsilon\).
However, the concept of \textit{principal factors}
has been coined much earlier by Barrucand and Cohn
\cite{BC1971}.
Our algorithm for the determination
of principal factors
has been presented in
\cite[Alg. 6, pp. 20--21]{Ma2022}.


\begin{theorem}
\label{thm:TwoPrimeCond}
(Principal factor criterion, \textbf{Ayadi, 1995},
\cite[Thm. 3.3, p. 37]{Ay1995}.) \\
Let \(c\) be a conductor divisible by two primes, \(t=2\), such that
\(\mathrm{Cl}_3(k_{c,\mu})\simeq (3,3)\)
for both cyclic cubic fields
\(k_{c,\mu}\), \(1\le\mu\le 2\),
with conductor \(c\).
Denote by \(\mathcal{P}\) the number of prime divisors of the norm
\(A(k)=\mathrm{N}_{k/\mathbb{Q}}(\alpha)\)
of a non-trivial primitive ambiguous principal ideal \((\alpha)\),
i.e. a \textbf{principal factor},
of any of the two fields \(k=k_{c,\mu}\).
Then \(\mathcal{P}\in\lbrace 1,2\rbrace\), \\
and the second \(3\)-class group
\(\mathfrak{M}=\mathrm{Gal}(\mathrm{F}_3^2(k)/k)\)
of both fields \(k=k_{c,\mu}\) is given by
\begin{equation}
\label{eqn:TwoPrimeCondAmbiguous}
\mathfrak{M}\simeq
\begin{cases}
\langle 9,2\rangle \text{ with capitulation type } \mathrm{a}.1,\ \varkappa(k)=(0000),  & \text{ if } \mathcal{P}=2, \\
\langle 27,4\rangle \text{ with capitulation type } \mathrm{A}.1,\ \varkappa(k)=(1111), & \text{ if } \mathcal{P}=1.
\end{cases}
\end{equation}
The length of the Hilbert \(3\)-class field tower is
\(\ell_3(k)=1\) with \(\mathrm{F}_3^\infty(k)=\mathrm{F}_3^1(k)\) if \(\mathcal{P}=2\), and
\(\ell_3(k)=2\) with \(\mathrm{F}_3^\infty(k)=\mathrm{F}_3^2(k)\) if \(\mathcal{P}=1\).
In both cases, \(\mathfrak{G}=\mathrm{Gal}(\mathrm{F}_3^\infty(k)/k)=\mathfrak{M}\).
\end{theorem}

\begin{proof}
See
\cite[Prp. 3.6, p. 32, Thm. 3.1, p. 34, Thm. 3.3, p. 37]{Ay1995}
and 
\cite[pp. 31--33]{Ma2022}.
\end{proof}

The first example \(c=19\cdot 1129=21\,451\)
for \(\mathfrak{M}\simeq\langle 27,4\rangle\)
is due to Scholz and Taussky
\cite[pp. 209--210]{Ta1932}.
It was misprinted as \(19\cdot 1429=27\,151\) in
\cite[p. 383]{Ta1952}.
Systematic tables have been presented at
\texttt{http://www.algebra.at/ResearchFrontier2013ThreeByThree.htm}
in \S\S\ 1.1--1.2.

Concerning the \(3\)-capitulation types \(\mathrm{a}.1\) and \(\mathrm{A}.1\),
viewed as transfer kernel types (TKT),
and the related concept of transfer target types (TTT),
i.e., abelian type invariants (ATI),
see
\cite{Ma2012}.


\section{Unramified extensions of cyclic cubic fields}
\label{s:UnramifiedExtensions}

\noindent
In this crucial section,
we first introduce the absolute \(3\)-genus field \(k^\ast\) (\S\
\ref{ss:GenusField})
of a cyclic cubic number field \(k\).
Then we show that the bicyclic bicubic subfields \(B<k^\ast\)
constitute unramified cyclic cubic relative extensions \(B/k\)
of a cyclic cubic number field \(k\) with \(t=3\).
Finally,
using the unramified cyclic cubic relative extensions \(E/k\)
as capitulation targets (\S\
\ref{ss:CapitulationTargets}),
we define the capitulation kernels (\S\
\ref{ss:CapitulationKernels})
of a cyclic cubic number field \(k\)
with non-trivial \(3\)-class group \(\mathrm{Cl}_3(k)\).


\subsection{The absolute 3-genus field}
\label{ss:GenusField}

\noindent
The \textit{absolute \(3\)-genus field} \(k^\ast=(k/\mathbb{Q})^\ast\)
of a cyclic cubic field \(k\) is
the maximal unramified \(3\)-extension \(k^\ast/k\)
with abelian absolute Galois group \(\mathrm{Gal}(k^\ast/\mathbb{Q})\).
If the conductor \(c=q_1\cdots q_t\) of \(k=k_c\) has \(t\) prime divisors,
then \(k^\ast\) is the compositum of the multiplet \((k_{c,1},\ldots,k_{c,m})\)
of all cyclic cubic fields sharing the common conductor \(c\),
where \(m=m(c)=2^{t-1}\), according to the multiplicity formula
\eqref{eqn:Multiplicity}.
The absolute Galois group \(\mathrm{Gal}(k^\ast/\mathbb{Q})\)
is the elementary abelian \(3\)-group \((\mathbb{Z}/3\mathbb{Z})^t\).
In particular, if \(t=1\), \(c=q_1\),
then \(k^\ast=k\) is the cyclic cubic field itself,
and if \(t=2\), \(c=q_1q_2\), 
then \(k^\ast=k_{c,1}\cdot k_{c,2}\) is a bicyclic bicubic field
with conductor \(c\) and discriminant
\begin{equation}
\label{eqn:Discriminant}
d(k^\ast)=d(k_{q_1})\cdot d(k_{q_2})\cdot d(k_{c,1})\cdot d(k_{c,2})=
q_1^2\cdot q_2^2\cdot(q_1q_2)^2\cdot(q_1q_2)^2=c^6.
\end{equation}


In 1990, Parry
\cite{Pa1990}
investigated the arithmetic of a general \textit{bicyclic bicubic field}
\(B/\mathbb{Q}\)
with conductor \(c=q_1\cdots q_t\), \(t\ge 2\),
and four cyclic cubic subfields \(k_1,\ldots,k_4\).
In particular, he determined the \textit{class number relation}
in terms of the \textit{index \(I\) of subfield units} of \(B\).

\begin{theorem}
\label{thm:Parry}
Let \(M:=(e_{i,j})\) be the \((4\times t)\)-\textbf{matrix of} integer exponents
in the following representation of the \textbf{principal factors}
\(A(k_i)=\prod_{j=1}^t\,q_j^{e_{i,j}}\), for \(1\le i\le 4\).
Then:
\begin{enumerate}
\item
The Galois group \(\mathrm{Gal}(B/\mathbb{Q})\simeq (3,3)\) is elementary bicyclic.
\item
The index \(I:=(U:V)\) of the subgroup \(V:=\langle U_1,\ldots,U_4\rangle\)
generated by the unit groups \(U_i:=U_{k_i}\), \(1\le i\le 4\),
in the unit group \(U:=U_B\) is bounded by \(I=3^e\), \(0\le e\le 3\).
\item
The \textbf{class number} of \(B\) satisfies the following \textbf{relation}:
\begin{equation}
\label{eqn:ParryFormula}
h(B)
=\frac{I}{3^5}\cdot\prod_{i=1}^4\,h(k_i)
=\frac{(U:V)}{243}\cdot h(k_1)\cdot h(k_2)\cdot h(k_3)\cdot h(k_4),
\end{equation}
where \(I\)
denotes the abovementioned \textbf{index of subfield units} of \(B\).
\item
\(3\nmid h(B)\) if and only if \(c=pq\), i.e. \(t=2\),
and \(p\not\leftrightarrow q\) are \textbf{not} mutual cubic residues,
i.e., the graph of \(p,q\) is either Graph \(1\) or Graph \(2\).
If \(3\nmid h(B)\), then \(I=27\).
\item
In dependence on the rank \(2\le r_M:=\mathrm{rank}(M)\le 4\) of the matrix \(M\),
the \textbf{index} \(I\) takes the following values:
\begin{equation}
\label{eqn:UnitIndex}
I=(U:V)=
\begin{cases}
1  & \text{ if } r_M=4, \\
3  & \text{ if } r_M=3, \\
9 \text{ or } 27  & \text{ if } r_M=2.
\end{cases}
\end{equation}
\end{enumerate}
\end{theorem}

\begin{proof}
For the \textit{class number relation}, see Parry
\cite[Prp. 7, p. 496, Thm. 9, p. 497]{Pa1990}.
Generally, 
the \textit{index of subfield units},
\(I\), is a divisor of \(27=3^3\)
\cite[Lem. 11, p. 500, Thm. 13, p. 501]{Pa1990}.
See also Ayadi
\cite[Prop. 2.7.(2) and Prop. 2.8, p. 20]{Ay1995}.
Note that \(p\not\leftrightarrow q\) implies
\(\prod_{i=1}^4\,h(k_i)=9\).
\end{proof}


\begin{corollary}
\label{cor:ClassNumber}
Let \(t=3\) and \(B\) be a bicyclic bicubic field with conductor \(c=pqr\)
such that there are \textbf{no mutual cubic residues} among \(p,q,r\). Then:
\begin{enumerate}
\item
For all \(1\le j\le 4\), \(h_3(B_j)=\frac{(U_j:V_j)}{3^2}h_3(k_j)\).
\item
For all \(5\le j\le 10\), \(h_3(B_j)=\frac{(U_j:V_j)}{3^4}h_3(k_i)h_3(k_\ell)\),
where \(1\le i,\ell\le 4\), \(i\ne\ell\), and \(k_i\), \(k_\ell\)
are the two components of the quartet which are contained in \(B_j\).
\end{enumerate}
\end{corollary}

\begin{proof}
By
\eqref{eqn:Three9Four},
the first statement is valid,
since \(h_3(k)=3\) for the six subfields \(k\) with \(t=2\).
By 
\eqref{eqn:Three9Ten},
the second statement holds,
since \(h_3(k)=1\) for the three subfields \(k\) with \(t=1\).
\end{proof}


For a cyclic cubic field \(k\) with \(t=2\), \(c=pq\),
the \(3\)-class numbers of the \(3\)-genus field \(k^\ast\),
which is bicyclic bicubic,
and of its four cyclic cubic subfields
can be summarized as follows.

\begin{theorem}
\label{thm:TwoPrimeGenus}
Let \(k^\ast=k_{p}\cdot k_{q}\cdot k_{c,1}\cdot k_{c,2}\)
be the genus field
of the two cyclic cubic fields \(k_{c,1}\) and \(k_{c,2}\)
with conductor \(c=pq\).
Denote the \(3\)-valuations of the class numbers
\(h^\ast\), \(h_1\), \(h_2\), \(h_3\), \(h_4\)
of \(k^\ast\), \(k_{p}\), \(k_{q}\), \(k_{c,1}\), \(k_{c,2}\), respectively,
by \(v^\ast\), \(v_1\), \(v_2\), \(v_3\), \(v_4\).
Then \(v_1=v_2=0\), and
\begin{equation}
\label{eqn:TwoPrimeGenus}
v^\ast
\begin{cases}
=0,\ v_3=v_4=1,\ I=27, & \text{ if } p\not\leftrightarrow q, \\
=1, & \text{ if } p\leftrightarrow q,\ v_3=v_4=2,\ I=9, \\
=2, & \text{ if } p\leftrightarrow q,\ v_3=v_4=2,\ I=27, \\
\ge 3, & \text{ if } p\leftrightarrow q,\ v_3\ge 3,\ v_4\ge 3,\ I\ge 9.
\end{cases}
\end{equation}
\end{theorem}

\begin{proof}
According to Theorem
\ref{thm:3ClassRank},
we generally have \(v_1=v_2=0\),
\(v_3\ge 1\), \(v_4\ge 1\) if \(p\not\leftrightarrow q\),
and \(v_3\ge 2\), \(v_4\ge 2\) if \(p\leftrightarrow q\).
Now, the claim is a consequence of Formula
\eqref{eqn:ParryFormula},
which yields
\[v^\ast=v_3(h^\ast)=v_3(I)-5+\sum_{i=1}^4\,v_3(h_i)=
v_3(I)-5+v_1+v_2+v_3+v_4=v_3(I)-5+v_3+v_4.\]
The combination of
\cite[Thm. 9, p. 497]{Pa1990}
and
\cite[Cor. 1, p. 498]{Pa1990}
shows that \(v^\ast=0\) if and only if \(p\not\leftrightarrow q\),
and \(v^\ast=0\) implies \(v_3(I)=3\),
whence necessarily \(v_3=v_4=1\).
However, if \(p\leftrightarrow q\), then
\(v_3=2\) is equivalent with \(v_4=2\), according to
\cite[Thm. 4.1, p. 472]{AAI2001}.
\end{proof}

\begin{remark}
\label{rmk:TwoPrimeGenus}
For \(v_3=v_4=2\), we have \(\mathrm{Cl}_3(k_{pq,\mu})\simeq (3,3)\).
The smallest occurrences of \(v_3=v_4=3\) are the conductors
\(7\cdot 673=4\,711\)
(\lq\lq Eau de Cologne\rq\rq, \textbf{singular} with \(\mathrm{Cl}_3(k^\ast)\simeq (3,3,3)\)) and
\(7\cdot 769=5\,383\) (\textbf{super-singular} with \(\mathrm{Cl}_3(k^\ast)\simeq (9,3,3)\))
both with \(\mathrm{Cl}_3(k_{pq,\mu})\simeq (9,3)\),
for \(\mu\in\lbrace 1,2\rbrace\).
\end{remark}


\noindent
For a cyclic cubic field \(k\) with \(t=3\) and conductor \(c=q_1q_2q_3\),
the \(3\)-genus field \(k^\ast\) contains \(13\) bicyclic bicubic subfields.
Three of them are the \textit{sub genus fields}
\(B_i:=(k_{f_{i-10}})^\ast\), \(11\le i\le 13\),
of the cyclic cubic fields with conductors
\(f_1=q_1q_2\), \(f_2=q_1q_3\), \(f_3=q_2q_3\),
respectively.
In the numerical tables of
\cite{Ma2022},
we always start with the leading three sub genus fields
\(B_i\), \(11\le i\le 13\),
separated by a semicolon from the trailing ten remaining bicyclic bicubic subfields,
when we give a family of invariants for these \(13\) subfields \(B_1,\ldots,B_{13}\),
\begin{equation}
\label{eqn:SubGenusFields}
\text{in particular, }
\left\lbrack\mathrm{Cl}_3{B_i}\right\rbrack_{1\le i\le 13}:=
\lbrack \mathrm{Cl}_3(B_{11}),\ldots,\mathrm{Cl}_3(B_{13});
\mathrm{Cl}_3(B_1),\ldots,\mathrm{Cl}_3(B_{10})\rbrack.
\end{equation}


\subsection{Capitulation kernels}
\label{ss:CapitulationKernels}

\noindent
We recall the connection between the size of the capitulation kernel \(\ker(T_{E/k})\)
and the unit norm index \((U_k:\mathrm{N}_{E/k}(U_E))\)
of an unramified cyclic cubic extension \(E/k\) of a cyclic cubic field \(k\).
Here, \(T_{E/k}:\,\mathrm{Cl}_3(k)\to\mathrm{Cl}_3(E)\),
\(\mathfrak{a}\mathcal{P}_k\mapsto(\mathfrak{a}\mathcal{O}_E)\mathcal{P}_E\),
denotes the extension homomorphism or \textbf{transfer} of \(3\)-classes from \(k\) to \(E\).

\begin{theorem}
\label{thm:Capitulation}
The order of the \(3\)-capitulation kernel
or \textbf{transfer kernel}
of \(E/k\) is given by
\begin{equation}
\label{eqn:CapitulationKernel}
\#\ker(T_{E/k})=
\begin{cases}
3, \\
9, \\
27,
\end{cases}
\text{ if and only if }\quad
(U_k:\mathrm{N}_{E/k}(U_E))=
\begin{cases}
1, \\
3, \\
9.
\end{cases}
\end{equation}
\end{theorem}

\begin{proof}
According to the \textbf{Herbrand Theorem} on the cohomology of the unit group \(U_E\)
as a Galois module with respect to \(G=\mathrm{Gal}(E/k)\),
we have the relation
\(\#\ker(T_{E/k})=\lbrack E:k\rbrack\cdot(U_k:\mathrm{N}_{E/k}(U_E))\),
since \(\ker(T_{E/k})\simeq H^1(G,U_E)\)
when \(E/k\) is unramified of odd prime degree \(\lbrack E:k\rbrack=3\)
and \(U_k/\mathrm{N}_{E/k}(U_E)\simeq\hat{H}^0(G,U_E)\).
The cyclic cubic base field \(k\) has signature \((r_1,r_2)=(3,0)\)
and torsionfree Dirichlet unit rank \(r=r_1+r_2-1=3+0-1=2\).
Thus, there are three possibilities for the unit norm index
\((U_k:\mathrm{N}_{E/k}(U_E))\in\lbrace 1,3,9\rbrace\).
\end{proof}

\begin{remark}
\label{rmk:Capitulation}
When \(k\) is a cyclic cubic field
with \(3\)-class group \(O:=\mathrm{Cl}_3(k)\) of elementary tricyclic type \((3,3,3)\),
viewed as a vector space of dimension \(3\) over the finite field \(\mathbb{F}_3\), then
\(\#\ker(T_{E/k})=3\) if and only if \(\ker(T_{E/k})=L_i\)
is a \textbf{line} for some \(1\le i\le 13\),
\(\#\ker(T_{E/k})=9\) if and only if \(\ker(T_{E/k})=P_i\)
is a \textbf{plane} for some \(1\le i\le 13\),
and
\(\#\ker(T_{E/k})=27\) if and only if \(\ker(T_{E/k})=O\)
is the \textbf{entire vector space} over \(\mathbb{F}_3\).
Details are reserved for a future paper.
Our algorithms for the determination
of the capitulation kernels
for \(\mathrm{Cl}_3(k)\) of type \((3,3)\) and \((3,3,3)\)
have been presented in
\cite[Alg. 8--9, pp. 26--30]{Ma2022}.
\end{remark}


\noindent
In our theorems on cyclic cubic fields with \(t=3\)
belonging to the various graphs of each category,
we shall frequently find particular statements
which relate several similar capitulation types.

\begin{definition}
\label{dfn:PartOrd}
Let \(G\) be a \(3\)-group
with generator rank \(d_1(G)=2\)
and elementary bicyclic commutator quotient
\(G/G^\prime\simeq (3,3)\).
By \(T_{G,H_i}:\,G/G^\prime\to H_i/H_i^\prime\)
we denote the transfers
from \(G\) to the four maximal normal subgroups
\(H_i\), \(1\le i\le 4\).
Then the set of all \textit{ordered transfer kernel types}
\(\varkappa=(\varkappa_i)_{1\le i\le 4}\)
with \(\varkappa_i:=\ker(T_{G,H_i})\)
is endowed with a \textit{partial order} relation
\(\varkappa\le\varkappa^\prime\) by 
\((\forall\,1\le i\le 4)\) \(\varkappa_i\le\varkappa^\prime_i\).
The order is strict, \(\varkappa<\varkappa^\prime\),
when \(\varkappa\le\varkappa^\prime\) and
\((\exists\,1\le j\le 4)\) \(\varkappa_j<\varkappa^\prime_j\).
\end{definition}

\noindent
The possibilities for a strict order are rather limited,
since a transfer kernel is either cyclic of order \(3\)
(\textit{partial} --- by Hilbert's Theorem \(94\), it cannot be trivial)
or bicyclic of type \((3,3)\) (\textit{total}).
As usual, we abbreviate
\(\varkappa_i=j\) if \((\exists\,1\le j\le 4)\) \(\ker(T_{G,H_i})=H_j/G^\prime\),
and \(\varkappa_i=0\) if \(\ker(T_{G,H_i})=G/G^\prime\),
for fixed \(1\le i\le 4\).
So, \(\varkappa<\varkappa^\prime\) \(\Longleftrightarrow\)
\((\exists\,1\le j,i\le 4)\) \(\varkappa_j=H_i/G^\prime<G/G^\prime=\varkappa^\prime_j\).
The arithmetical application of this group theoretic Definition
\ref{dfn:PartOrd}
is given in the following definition.

\begin{definition}
\label{dfn:mTKT}
Let \(K\) be an algebraic number field with
elementary bicyclic \(3\)-class goup
\(\mathrm{Cl}_3(K)\simeq (3,3)\).
Then \(K\) has four unramified cyclic cubic relative extensions
\(E_i/K\), \(1\le i\le 4\),
and corresponding class extension homomorphisms
\(T_{E_i/K}:\,\mathrm{Cl}_3(K)\to\mathrm{Cl}_3(E_i)\). 
Let \(\mathfrak{M}:=\mathrm{Gal}(\mathrm{F}_3^2(K)/K)\)
be the Galois group of the second Hilbert \(3\)-class field of \(K\),
that is, the maximal metabelian unramified \(3\)-extension of \(K\).
Then \(\varkappa(K):=\varkappa(\mathfrak{M})\) is called
the \textbf{minimal transfer kernel type} (mTKT) of \(K\),
if \(\varkappa(K)\le\varkappa^\prime(K)\),
for any other possible capitulation type \(\varkappa^\prime(K)\).
\end{definition}


\subsection{Capitulation targets}
\label{ss:CapitulationTargets}

\noindent
The precise
constitution of the lattice of all subfields
of the absolute \(3\)-genus field \(k^\ast\)
of a cyclic cubic field \(k=k_{pqr}\)
with \(t=3\) and conductor \(c=pqr\)
is as follows.

\begin{theorem}
\label{thm:Three}
The genus field \(k^\ast\) of \(k\) contains
\(13\) cyclic cubic fields,
\begin{equation}
\label{eqn:Three3}
\begin{aligned}
& k_{p,1},k_{q,1},k_{r,1},k_{pq,1},k_{pq,2},k_{pr,1},k_{pr,2},k_{qr,1},k_{qr,2},k_{pqr,1},k_{pqr,2},k_{pqr,3},k_{pqr,4}, \text{ briefly } \\
& k_p,k_q,k_r,k_{pq},\tilde{k}_{pq},k_{pr},\tilde{k}_{pr},k_{qr},\tilde{k}_{qr},k_{1},k_{2},k_{3},k_{4}.
\end{aligned}
\end{equation}
The composita \(L:=k_{pq}k_{pr}k_{qr}\) and \(\tilde{L}:=\tilde{k}_{pq}\tilde{k}_{pr}\tilde{k}_{qr}\) satisfy
the \textbf{skew balance of degrees} \\
\(\lbrack L:\mathbb{Q}\rbrack\cdot\lbrack\tilde{L}:\mathbb{Q}\rbrack=243\)
with \(\lbrack L:\mathbb{Q}\rbrack=9\) and \(\lbrack\tilde{L}:\mathbb{Q}\rbrack=27\),
or vice versa. \\
\textbf{Alert:} Always in the sequel, the \textbf{normalization} \(\lbrack L:\mathbb{Q}\rbrack=9\) is assumed. \\
The genus field \(k^\ast\) of \(k\) contains
\(13\) bicyclic bicubic fields,
\begin{equation}
\label{eqn:Three9Four}
\begin{aligned}
4 \textbf{ single capitulation targets } \quad
B_1    &:= k_{pq}k_{pr} = k_1k_{pq}k_{pr}k_{qr}, \\
B_2    &:= \tilde{k}_{pr}\tilde{k}_{qr} = k_2k_{pq}\tilde{k}_{pr}\tilde{k}_{qr}, \\
B_3    &:= \tilde{k}_{pq}\tilde{k}_{pr} = k_3\tilde{k}_{pq}\tilde{k}_{pr}k_{qr}, \\
B_4    &:= \tilde{k}_{pq}\tilde{k}_{qr} = k_4\tilde{k}_{pq}k_{pr}\tilde{k}_{qr},
\end{aligned}
\end{equation}
\begin{equation}
\label{eqn:Three9Ten}
\begin{aligned}
6 \textbf{ double capitulation targets } \quad
B_5    &:= k_p\tilde{k}_{qr} = k_1k_3k_p\tilde{k}_{qr}, \\
B_6    &:= k_q\tilde{k}_{pr} = k_1k_4k_q\tilde{k}_{pr}, \\
B_7    &:= k_r\tilde{k}_{pq} = k_1k_2k_r\tilde{k}_{pq}, \\
B_8    &:= k_pk_{qr} = k_2k_4k_pk_{qr}, \\
B_9    &:= k_qk_{pr} = k_2k_3k_qk_{pr}, \\
B_{10} &:= k_rk_{pq} = k_3k_4k_rk_{pq}, 
\end{aligned}
\end{equation}
\begin{equation}
\label{eqn:Three9Thirteen}
\begin{aligned}
\text{and } 3 \textbf{ sub genus fields } \quad
B_{11} &:= k_{pq}\tilde{k}_{pq} = k_pk_qk_{pq}\tilde{k}_{pq}, \\
B_{12} &:= k_{pr}\tilde{k}_{pr} = k_pk_rk_{pr}\tilde{k}_{pr}, \\
B_{13} &:= k_{qr}\tilde{k}_{qr} = k_qk_rk_{qr}\tilde{k}_{qr},
\end{aligned}
\end{equation}
provided that \(k_{pq,}k_{pr},k_{qr}\) are normalized.
The conductor of \(B_1,\ldots,B_{10}\) is \(c=pqr\),
the conductor of \(B_{11}\) is \(f_1=pq\),
the conductor of \(B_{12}\) is \(f_2=pr\), and
the conductor of \(B_{13}\) is \(f_3=qr\).
\end{theorem}

\begin{proof}
See
\cite[Prop. 4.1, p. 40, Lem. 4.1, p. 42]{Ay1995}.
The \textbf{short form} suffices for construction.
\end{proof}

\noindent
The algorithm for the determination
of bicyclic bicubic fields
has been presented in
\cite[Alg. 7, pp. 24--26]{Ma2022},
but \(B_5,\ldots,B_{10}\) should be defined as in Formula
\eqref{eqn:Three9Ten}
(short form without \(k_1,\ldots,k_4\)).


\begin{corollary}
\label{cor:Three}
The \textbf{capitulation targets}, i.e. the
unramified cyclic cubic relative extensions 
of \(k_1\), respectively \(k_2\),
respectively \(k_3\), respectively \(k_4\),
among the absolutely bicyclic bicubic subfields
of the \(3\)-genus field \(k^\ast=k_pk_qk_r\)
are \(B_1,B_5,B_6,B_7\),
respectively \(B_2,B_7,B_8,B_9\),
respectively \(B_3,B_5,B_9,B_{10}\),
respectively \(B_4,B_6,B_8,B_{10}\).
In particular,
\(B_7\) is common to both, \(k_1\) and \(k_2\),
\(B_5\) is common to \(k_1\) and \(k_3\),
\(B_6\) is common to \(k_1\) and \(k_4\),
\(B_9\) is common to \(k_2\) and \(k_3\),
\(B_8\) is common to \(k_2\) and \(k_4\), and
\(B_{10}\) is common to \(k_3\) and \(k_4\).
\end{corollary}

\begin{proof}
This follows immediately from Theorem
\ref{thm:Three},
Equations
\eqref{eqn:Three9Four}
and
\eqref{eqn:Three9Ten}.
\end{proof}


\begin{proposition}
\label{prp:Bicub}
If there exists \(1\le j\le 10\) such that
\(h_3(B_j)=3\),
then \(h_3(B_\ell)=3\), for all \(1\le\ell\le 10\), and
\(h_3(k_i)=9\), for all \(1\le i\le 4\).

The \(3\)-class number of \(B_j\), \(1\le j\le 10\),
satisfies the \textbf{tame} condition \(h_3(B_j)=(U_j:V_j)\)
if and only if
for each cyclic cubic subfield \(k\) of \(B_j\)
the Hilbert \(3\)-class field \(\mathrm{F}_3^1(k)\) of \(k\) coincides with
the genus field \(k^\ast\) of \(k\).
Otherwise the \textbf{wild} condition \(h_3(B_j)>(U_j:V_j)\) holds.

If there exists \(1\le j\le 10\) such that
\(h_3(B_j)>(U_j:V_j)\),
then \(9\mid h_3(B_\ell)\), for all \(1\le\ell\le 10\).
\end{proposition}

\begin{proof}
The condition is trivial for the subfields \(k\) with \(t=1\),
since \(h_3(k)=\lbrack\mathrm{F}_3^1(k):k\rbrack=\lbrack k^\ast:k\rbrack=1\)
is satisfied anyway.
However, the subfields \(k\) with \(t=2\) must have the \(3\)-class number
\(h_3(k)=\lbrack\mathrm{F}_3^1(k):k\rbrack=\lbrack k^\ast:k\rbrack=3\),
in particular, the prime divisors of the conductor are not mutual cubic residues,
and the subfields \(k\) with \(t=3\) must have \(3\)-class number
\(h_3(k)=\lbrack\mathrm{F}_3^1(k):k\rbrack=\lbrack k^\ast:k\rbrack=9\),
that is, they cannot have \(3\)-class rank \(\varrho(k)\ge 3\).
For details see
\cite[pp. 47--48, i.p. Prop. 4.5]{Ay1995}.
\end{proof}


\noindent
Let \(t=3\) and \(k_\mu\), \(1\le\mu\le 4\), be one of the four
cyclic cubic number fields
sharing the common conductor \(c=pqr\), and
suppose \(B_j\), \(1\le j\le 10\),
is one of the ten bicyclic bicubic subfields
of the absolute \(3\)-genus field \(k^\ast\) of \(k_\mu\)
such that \(B_j/k_\mu\) is an unramified cyclic extension of degree \(3\).
We denote by \(U_j\) the unit group of \(B_j\),
by \(V_j\) the subgroup generated by all subfield units,
by \(r_j\) the rank of the principal factor matrix \(M_j\) of \(B_j\),
and by \(A=(a_{\iota\lambda})\) the right upper triangular \((8\times 8)\)-matrix such that
\((\gamma_1^3,\ldots,\gamma_8^3)=(\varepsilon_1,\ldots,\varepsilon_8)\cdot A\)
(in the sense of exponentiation),
for a suitable torsion free basis \((\gamma_1,\ldots,\gamma_8)\) of \(U_j\)
and a canonical basis \((\varepsilon_1,\ldots,\varepsilon_8)\) of \(V_j\),
according to
\cite[pp. 497--503]{Pa1990}
and
\cite[pp. 19--22]{Ay1995}.

For several times,
Ayadi
\cite{Ay1995}
alludes to the following fact:
the \textit{minimal subfield unit index}
\((U_j:V_j)=3\) for the matrix rank \(r_j=3\)
of \(B_j\) corresponds to
the \textit{maximal unit norm index} \((U(k_\mu):N_{B_j/k_\mu}(U_j))=3\),
associated with a \textit{total transfer kernel}
\(\#\ker(T_{B_j/k_\mu})=9\)
of \(B_j/k_\mu\).
Since he does not give a prove,
we summarize all related issues in a lemma.

\begin{lemma}
\label{lem:UnitIndices}
The following statements are equivalent, row by row:
\begin{equation}
\label{eqn:UnitIndices}
\begin{aligned}
(U_j:V_j)=3 & \ \Longleftrightarrow\  a_{77}=3,\ a_{88}=3,\ a_{66}=1 \ \Longrightarrow\  
(U(k_\mu):N_{B_j/k_\mu}(U_j))=3, \\
(U_j:V_j)=9 & \ \Longleftrightarrow\  a_{77}=3,\ a_{88}=1,\ a_{66}=1 \ \Longrightarrow\ 
(U(k_\mu):N_{B_j/k_\mu}(U_j))=3, \\
(U_j:V_j)=27 & \ \Longleftrightarrow\  a_{77}=1,\ a_{88}=1,\ a_{66}=1 \ \Longleftrightarrow\ 
(U(k_\mu):N_{B_j/k_\mu}(U_j))=1.
\end{aligned}
\end{equation}
\end{lemma}

\begin{proof}
According to Theorem
\ref{thm:Parry},
\(r_j=3\) \(\Longleftrightarrow\) \((U_j:V_j)=3\), and
\(r_j=2\) \(\Longleftrightarrow\) \((U_j:V_j)\in\lbrace 9,27\rbrace\). \\
Now,
\(a_{77}=1\) implies
\(\gamma_7^3=(\prod_{\iota=1}^6\,\varepsilon_{\iota}^{a_{\iota 7}})\cdot\varepsilon_7\),
\(N_{B_j/k_\mu}(\gamma_7)=\pm\varepsilon_7\),
\quad \((U(k_\mu):N_{B_j/k_\mu}(U_j))=1\), \\
but
\(a_{77}=3\) implies
\(\gamma_7^3=(\prod_{\iota=1}^6\,\varepsilon_{\iota}^{a_{\iota 7}})\cdot\varepsilon_7^3\),
\(N_{B_j/k_\mu}(\gamma_7)=\pm\varepsilon_7^3\),
\quad \((U(k_\mu):N_{B_j/k_\mu}(U_j))=3\).

Finally, Theorem
\ref{thm:Capitulation}
on the \textbf{Herbrand quotient} of \(U_j\)
shows the cardinality of the transfer kernel,
\(\#\ker{T_{B_j/k_\mu}}=\lbrack k_\mu:\mathbb{Q}\rbrack\cdot (U(k_\mu):N_{B_j/k_\mu}(U_j))
=3\cdot (U(k_\mu):N_{B_j/k_\mu}(U_j))\).
\end{proof}


\begin{proposition}
\label{prp:Bicyc}
Let \(\ell\) be an odd prime,
and suppose that \(B=K\cdot L\) is a bicyclic field of degree \(\ell^2\),
compositum of two cyclic fields \(K\) and \(L\) of degree \(\ell\).
If \(p\) is a prime number which ramifies in both, \(K\) and \(L\),
i.e., \(p\mathcal{O}_K=\mathfrak{p}_1^\ell\) and  \(p\mathcal{O}_L=\mathfrak{p}_2^\ell\),
then the extension ideals
\(\mathfrak{p}_1\mathcal{O}_B=\mathfrak{p}_2\mathcal{O}_B\) coincide.
\end{proposition}

\begin{proof}
If the decomposition invariants of \(p\) in \(B\) are
\((e,f,g)=(\ell,1,\ell)\), resp. \((\ell,\ell,1)\), resp. \((\ell^2,1,1)\),
then those of \(\mathfrak{p}_1\) and \(\mathfrak{p}_2\) in \(B\)
must be identical
\((e,f,g)=(1,1,\ell)\), resp. \((1,\ell,1)\), resp. \((\ell,1,1)\),
and unique prime decomposition enforces
\(\mathfrak{p}_1\mathcal{O}_B=\mathfrak{p}_2\mathcal{O}_B\).
\end{proof}


\begin{corollary}
\label{cor:Capitulation}
Let \(\mu\in\lbrace 1,2,3,4\rbrace\) and
\(p\mathcal{O}_{k_\mu}=\mathfrak{p}^3\),
\(q\mathcal{O}_{k_\mu}=\mathfrak{q}^3\),
\(r\mathcal{O}_{k_\mu}=\mathfrak{r}^3\).
Then the following \textbf{capitulation laws} for ideal classes hold
independently of the combined cubic residue symbol
\(\lbrack p,q,r\rbrack_3\).
\begin{enumerate}
\item
\(\lbrack\mathfrak{p}\rbrack\) capitulates
in \(B_5/k_\mu\), for \(\mu=1,3\), and
in \(B_8/k_\mu\), for \(\mu=2,4\).
\item
\(\lbrack\mathfrak{q}\rbrack\) capitulates
in \(B_6/k_\mu\), for \(\mu=1,4\), and
in \(B_9/k_\mu\), for \(\mu=2,3\).
\item
\(\lbrack\mathfrak{r}\rbrack\) capitulates
in \(B_7/k_\mu\), for \(\mu=1,2\), and
in \(B_{10}/k_\mu\), for \(\mu=3,4\).
\end{enumerate}
\end{corollary}

\begin{proof}
We show that \(\lbrack\mathfrak{p}\rbrack\in\mathrm{Cl}_3(k_1)\) capitulates in \(B_5\). 
Everythig else is proved in the same way, always using Proposition
\ref{prp:Bicyc}
with \(\ell=3\).
The bicyclic bicubic field \(B_5=k_1k_3k_p\tilde{k}_{qr}\)
is compositum of the cyclic cubic fields \(k_p\) and \(k_1\).
Since the conductor of \(k_p\) is \(p\),
the principal factor \(A(k_p)=p\) is determined uniquely,
and \(p\mathcal{O}_{k_p}=\mathfrak{p}_0^3\) is totally ramified,
whence \(\mathfrak{p}_0=\alpha\mathcal{O}_{k_p}\)
with \(\alpha\in k_p^\times\) is necessarily a principal ideal.
Since the conductor of \(k_1\) is \(c=pqr\),
the prime \(p\mathcal{O}_{k_1}=\mathfrak{p}^3\) is also totally ramified,
and Proposition
\ref{prp:Bicyc}
asserts that \(\mathfrak{p}\mathcal{O}_{B_5}=\mathfrak{p}_0\mathcal{O}_{B_5}\),
which is the principal ideal \(\alpha\mathcal{O}_{B_5}\).
Thus the class \(\lbrack\mathfrak{p}\rbrack\) capitulates in \(B_5\).
\end{proof}


\begin{proposition}
\label{prp:Principal2}
If \(\left(\frac{p}{q}\right)_3=1\) but \(\left(\frac{q}{p}\right)_3\ne 1\),
then \(\mathrm{Cl}_3(k_{pq})\simeq (3)\), \(\mathrm{Cl}_3(\tilde{k}_{pq})\simeq (3)\),
and \textbf{two principal factors} are given by
\(A(k_{pq})=p\), \(A(\tilde{k}_{pq})=p\).
\end{proposition}

\begin{proof}
If \(p\rightarrow q\),
then \(p\) splits in \(k_q\),
\(p\mathcal{O}_{k_q}=\wp_1\wp_2\wp_3\),
and \(\mathrm{Cl}_3(k_{pq})\simeq (3)\),
according to Georges Gras
\cite{Gr1973}.
The Hilbert \(3\)-class field
\(\mathrm{F}_3^1(k_{pq})\) of \(k_{pq}\)
with \(\lbrack\mathrm{F}_3^1(k_{pq}):k_{pq}\rbrack=3\)
coincides with the absolute \(3\)-genus field
\(k^\ast=k_p\cdot k_q=k_{pq}\cdot\tilde{k}_{pq}\) of the doublet
\((k_{pq},\tilde{k}_{pq})\)
with \(\lbrack k^\ast:\mathbb{Q}\rbrack=9\)
and \(\lbrack k^\ast:k_{pq}\rbrack=3\).

Since the conductor of \(k_{pq}\) is \(pq\),
\(p\mathcal{O}_{k_{pq}}=\mathfrak{p}^3\) is ramified in \(k_{pq}\),
but \(k^\ast\) is unramified over \(k_{pq}\),
and the decomposition invariants of \(p\) in \(k^\ast\) are 
\((e,f,g)=(3,1,3)\),
those of \(\mathfrak{p}\) in \(k^\ast=\mathrm{F}_3^1(k_{pq})\) are
\((e,f,g)=(1,1,3)\),
i.e. \(\mathfrak{p}\) splits completely in \(\mathrm{F}_3^1(k_{pq})\),

By the decomposition law of the Hilbert \(3\)-class field,
\(\mathfrak{p}=\alpha\mathcal{O}_{k_{pq}}\)
is principal with \(\alpha\in k_{pq}^\times\).
Therefore the unique principal factor of \(k_{pq}\) is \(A(k_{pq})=p\).
The same reasoning is true for \(\tilde{k}_{pq}\).
\end{proof}


\begin{proposition}
\label{prp:Principal3}
Let \(\mu\in\lbrace 1,2,3,4\rbrace\),
such that \(\mathrm{Cl}_3(k_\mu)\simeq (3,3)\). \\
If \(\left(\frac{p}{q}\right)_3=1\) and \(\left(\frac{p}{r}\right)_3=1\),
then \textbf{the principal factor} of \(k_\mu\) is
\(A(k_\mu)=p\).
\end{proposition}

\begin{proof}
Since \(\mathrm{Cl}_3(k_\mu)\simeq (3,3)\),
the Hilbert \(3\)-class field
\(\mathrm{F}_3^1(k_\mu)\) of \(k_\mu\)
with \(\lbrack\mathrm{F}_3^1(k_\mu):k_\mu\rbrack=9\)
coincides with the absolute \(3\)-genus field
\(k^\ast=k_p\cdot k_q\cdot k_r\) of the quartet
\((k_1,\ldots,k_4)\)
with \(\lbrack k^\ast:\mathbb{Q}\rbrack=27\)
and \(\lbrack k^\ast:k_\mu\rbrack=9\).

Since the conductor of \(k_\mu\) is \(c=pqr\),
\(p\mathcal{O}_{k_\mu}=\mathfrak{p}^3\) is ramified in \(k_\mu\),
but \(k^\ast\) is unramified over \(k_\mu\).

If \(q\leftarrow p\rightarrow r\) is universally repelling,
then \(p\) splits in \(k_q\) and in \(k_r\),
and the decomposition invariants of \(p\) in \(k^\ast\) are 
\((e,f,g)=(3,1,9)\),
those of \(\mathfrak{p}\) in \(k^\ast=\mathrm{F}_3^1(k_\mu)\) are
\((e,f,g)=(1,1,9)\),
i.e. \(\mathfrak{p}\) splits completely in \(\mathrm{F}_3^1(k_\mu)\),
and the decomposition law of the Hilbert \(3\)-class field
implies that \(\mathfrak{p}=\alpha\mathcal{O}_{k_\mu}\)
is principal with \(\alpha\in k_\mu^\times\).
Therefore the unique principal factor of \(k_\mu\) is \(A(k_\mu)=p\).
\end{proof}


\section{Finite 3-groups of type (3,3)}
\label{s:GroupTheory}

\noindent
In the following tables,
we list those invariants of finite \(3\)-groups \(G\)
with elementary bicyclic commutator quotient \(G/G^\prime\simeq (3,3)\)
which qualify metabelian groups \(\mathfrak{M}\) as
second \(3\)-class groups \(\mathrm{Gal}(\mathrm{F}_3^2(k)/k)\)
and non-metabelian groups \(\mathfrak{G}\) as
\(3\)-class field tower groups \(\mathrm{Gal}(\mathrm{F}_3^\infty(k)/k)\)
of cyclic cubic number fields \(k\).
The process of searching for suitable groups in descendant trees
with the strategy of pattern recognition
\cite{Ma2020}
is governed by the \textit{Artin pattern}
\(\mathrm{AP}=(\alpha,\varkappa)\)
\cite[p. 27]{Ma2015c},
where \(\alpha=\alpha_1\), respectively \(\varkappa=\varkappa_1\), denotes the first layer
of the transfer target type (TTT), respectively transfer kernel type (TKT).
Additionally, we give the top layer \(\alpha_2\) of the TTT,
which consists of the abelian quotient invariants
of the commutator subgroup \(\mathfrak{M}^\prime\),
corresponding to the \(3\)-class group of the first Hilbert \(3\)-class field
\(\mathrm{F}_3^1(k)\) of \(k\).
The \textit{nuclear rank} \(\nu\) is responsible for the search complexity.
The \(p\)-multiplicator rank \(\mu\) of a group \(G\) is precisely its \textit{relation rank}
\(d_2(G)=\dim_{\mathbb{F}_3}\mathrm{H}^2(G,\mathbb{F}_3)\),
which decides whether \(G\) is admissible as \(\mathrm{Gal}(\mathrm{F}_3^\infty(k)/k)\),
according to the Shafarevich Theorem
\cite{Sh1964},
\cite{Ma2015c}.
In the case of cyclic cubic fields \(k\), it is limited by the \textit{Shafarevich bound}
\(\mu\le\varrho+r+\theta\),
where
\(\varrho=d_1(G)=\dim_{\mathbb{F}_3}\mathrm{H}^1(G,\mathbb{F}_3)\)
denotes the \textit{generator rank} of \(G\),
which coincides with the \(3\)-class rank \(\varrho\) of \(k\),
\(r=r_1+r_2-1=2\) is the torsion free Dirichlet unit rank of the field \(k\) with signature \((r_1,r_2)=(3,0)\),
and \(\theta=0\) indicates the absence of a (complex) primitive third root of unity
in the totally real field \(k\).
Finally, \(\pi(\mathfrak{M})=\mathfrak{M}/\gamma_c(\mathfrak{M})\) denotes
the parent of \(\mathfrak{M}\),
that is the quotient by the last non-trivial lower central
with \(c=\mathrm{cl}(\mathfrak{M})\).


\begin{theorem}
\label{thm:ClassGroup33}
Let \(k\) be a cyclic cubic number field
with elementary bicyclic \(3\)-class group \(\mathrm{Cl}_3(k)\simeq (3,3)\).
Denote by \(\mathfrak{M}=\mathrm{Gal}(\mathrm{F}_3^2(k)/k)\)
the second \(3\)-class group of \(k\),
and by \(\mathfrak{G}=\mathrm{Gal}(\mathrm{F}_3^\infty(k)/k)\)
the \(3\)-class field tower group of \(k\).
Then, the Artin pattern \((\alpha,\varkappa)\) of \(k\)
identifies the groups \(\mathfrak{M}\) and \(\mathfrak{G}\),
and determines the length \(\ell_3(k)\) of the \(3\)-class field tower of \(k\),
according to the following \textbf{deterministic laws}.
(See the associated descendant tree \(\mathcal{T}^1\langle 9,2\rangle\) in
\cite[Fig. 6.1, p. 44]{Ma2022}.)
\begin{enumerate}
\item
If \(\alpha=\lbrack 1,1,1,1\rbrack\), \(\varkappa=(0000)\) (type \(\mathrm{a}.1\)), then \(\mathfrak{G}\simeq\langle 9,2\rangle\) and \(\ell_3(k)=1\).
\item
If \(\alpha\sim\lbrack 11,2,2,2\rbrack\), \(\varkappa\sim (1111)\) (type \(\mathrm{A}.1\)), then \(\mathfrak{G}\simeq\langle 27,4\rangle\).
\item
If \(\alpha\sim\lbrack 111,11,11,11\rbrack\), \(\varkappa\sim (2000)\) (type \(\mathrm{a}.3^\ast\)), then \(\mathfrak{G}\simeq\langle 81,7\rangle\).
\item
If \(\alpha\sim\lbrack 21,11,11,11\rbrack\), \(\varkappa\sim (2000)\) (type \(\mathrm{a}.3\)), then \(\mathfrak{G}\simeq\langle 81,8\rangle\).
\item
If \(\alpha\sim\lbrack 21,11,11,11\rbrack\), \(\varkappa\sim (1000)\) (type \(\mathrm{a}.2\)), then \(\mathfrak{G}\simeq\langle 81,10\rangle\).
\item
If \(\alpha\sim\lbrack 22,11,11,11\rbrack\), \(\varkappa\sim (2000)\) (type \(\mathrm{a}.3\)), then \(\mathfrak{G}\simeq\langle 243,25\rangle\).
\item
If \(\alpha\sim\lbrack 22,11,11,11\rbrack\), \(\varkappa\sim (1000)\) (type \(\mathrm{a}.2\)), then \(\mathfrak{G}\simeq\langle 243,27\rangle\).
\end{enumerate} 
Except for the abelian tower in item (1), the tower is metabelian with \(\ell_3(k)=2\).
\end{theorem}

\begin{proof}
Generally, a cyclic cubic field \(k\)
has signature \((r_1,r_2)=(3,0)\), torsion free unit rank \(r=r_1+r_2-1=2\),
does not contain primitive third roots of unity,
and thus possesses the maximal admissible relation rank
\(d_2\le d_1+r=4\) for the group \(\mathfrak{G}\),
when its \(3\)-class rank, i.e. the generator rank of \(\mathfrak{G}\), is \(d_1=\varrho=2\).
Consequently, \(\ell_3(k)\ge 3\) in the case of \(d_2(\mathfrak{M})\ge 5\).

For item (1), we have
\(\mathfrak{M}=\mathrm{Gal}(\mathrm{F}_3^2(k)/k)\simeq\langle 9,2\rangle\simeq (3,3)\simeq\mathrm{Cl}_3(k)\simeq\mathrm{Gal}(\mathrm{F}_3^1(k)/k)\),
whence \(\ell_3(k)=1\).
We always identify groups according to
\cite{BEO2005}
and
\cite{GNO2006}.

For item (2) to item (7),
the group \(\mathfrak{M}\) is of maximal class (coclass \(\mathrm{cc}(\mathfrak{M})=1\)),
and thus coincides with \(\mathfrak{G}\), whence \(\ell_3(k)=2\).

In each case, the Artin pattern \((\alpha,\varkappa)\) identifies \(\mathfrak{M}=\mathfrak{G}\) uniquely,
and the relation ranks are
\(d_2\langle 9,2\rangle=3\), 
\(d_2\langle 27,4\rangle=2\), 
\(d_2\langle 81,7\rangle=3\), 
\(d_2\langle 81,8\rangle=3\), 
\(d_2\langle 81,10\rangle=3\), 
\(d_2\langle 243,25\rangle=3\),
\(d_2\langle 243,27\rangle=3\),
each of them less than \(4\).
\end{proof}


\begin{corollary}
\label{cor:ClassGroup33}
Under the assumptions of Theorem
\ref{thm:ClassGroup33},
the abelian type invariants \(\alpha_2\) of the \(3\)-class group
\(\mathrm{Cl}_3(\mathrm{F}_3^1(k))\)
of the first Hilbert \(3\)-class field of \(k\)
are required for the unambiguous identification
of the following groups \(\mathfrak{G}\) respectively \(\mathfrak{M}\). 
(See the associated descendant tree \(\mathcal{T}^2\langle 729,40\rangle\) in
\cite[Fig. 6.2, p. 45]{Ma2022}.) \\
If \(\alpha\sim\lbrack 21,11,11,11\rbrack\), \(\varkappa=(0000)\), \(\mathrm{a}.1\), then
\(\mathfrak{G}\simeq
\begin{cases}
\langle 81,9\rangle & \text{ for } \alpha_2=\lbrack 11\rbrack, \\
\langle 243,28..30\rangle & \text{ for } \alpha_2\sim\lbrack 21\rbrack.
\end{cases}\) \\
If \(\alpha\sim\lbrack 21,21,111,111\rbrack\), \(\varkappa\sim (0043)\), \(\mathrm{b}.10\), then
\(\mathfrak{M}\simeq
\begin{cases}
\langle 729,34..36\rangle & \text{ for } \alpha_2=\lbrack 1111\rbrack, \\
\langle 729,37..39\rangle & \text{ for } \alpha_2\sim\lbrack 211\rbrack.
\end{cases}\) 
\end{corollary}

\begin{proof}
The Artin pattern \((\alpha,\varkappa)\) of \(k\) alone is not able
to identify the groups \(\mathfrak{M}\) and \(\mathfrak{G}\) unambiguously.
Ascione
\cite{AHL1977}
uses the notation
\(\langle 729,34\rangle=H\),
\(\langle 729,35\rangle=I\),
\(\langle 729,37\rangle=A\),
\(\langle 729,38\rangle=C\).
\end{proof}


\renewcommand{\arraystretch}{1.1}

\begin{table}[hb]
\caption{Invariants of metabelian \(3\)-groups \(\mathfrak{M}\) with \(\mathfrak{M}/\mathfrak{M}^\prime\simeq (3,3)\)}
\label{tbl:Metabelian33}
\begin{center}
{\tiny
\begin{tabular}{|c|r|l|c|c|c|r|r|c|}
\hline
\(\mathfrak{M}\)                     & cc    & Type & \(\varkappa\) & \(\alpha\)        & \(\alpha_2\) & \(\nu\) & \(\mu\) & \(\pi(\mathfrak{M})\) \\
\hline
\(\langle 9,2\rangle\)               & \(1\) & a.1  & \((0000)\)    & \(1,1,1,1\)       & \(0\)        &   \(3\) &   \(3\) & --- \\
\(\langle 27,4\rangle\)              & \(1\) & A.1  & \((1111)\)    & \(11,2,2,2\)      & \(1\)        &   \(0\) &   \(2\) & \(\langle 9,2\rangle\) \\
\(\langle 81,7\rangle\)              & \(1\) & a.3\({}^\ast\)  & \((2000)\)    & \(111,11,11,11\)  & \(11\)       &   \(0\) &   \(3\) & \(\langle 27,3\rangle\) \\
\(\langle 81,8\rangle\)              & \(1\) & a.3  & \((2000)\)    & \(21,11,11,11\)   & \(11\)       &   \(0\) &   \(3\) & \(\langle 27,3\rangle\) \\
\(\langle 81,9\rangle\)              & \(1\) & a.1  & \((0000)\)    & \(21,11,11,11\)   & \(11\)       &   \(1\) &   \(4\) & \(\langle 27,3\rangle\) \\
\(\langle 81,10\rangle\)             & \(1\) & a.2  & \((1000)\)    & \(21,11,11,11\)   & \(11\)       &   \(0\) &   \(3\) & \(\langle 27,3\rangle\) \\
\(\langle 243,25\rangle\)            & \(1\) & a.3  & \((2000)\)    & \(22,11,11,11\)   & \(21\)       &   \(0\) &   \(3\) & \(\langle 81,9\rangle\) \\
\(\langle 243,27\rangle\)            & \(1\) & a.2  & \((1000)\)    & \(22,11,11,11\)   & \(21\)       &   \(0\) &   \(3\) & \(\langle 81,9\rangle\) \\
\(\langle 243,28..30\rangle\)        & \(1\) & a.1  & \((0000)\)    & \(21,11,11,11\)   & \(21\)       &   \(0\) &   \(3\) & \(\langle 81,9\rangle\) \\
\hline
\(\langle 243,3\rangle\)             & \(2\) & b.10 & \((0043)\)    & \(21,21,111,111\) & \(111\)      &   \(2\) &   \(4\) & \(\langle 27,3\rangle\) \\
\(\langle 729,34\rangle=H\)          & \(2\) & b.10 & \((0043)\)    & \(21,21,111,111\) & \(1111\)     &   \(2\) &   \(5\) & \(\langle 243,3\rangle\) \\
\(\langle 729,35\rangle=I\)          & \(2\) & b.10 & \((0043)\)    & \(21,21,111,111\) & \(1111\)     &   \(1\) &   \(4\) & \(\langle 243,3\rangle\) \\
\(\langle 729,37\rangle=A\)          & \(2\) & b.10 & \((0043)\)    & \(21,21,111,111\) & \(211\)      &   \(2\) &   \(5\) & \(\langle 243,3\rangle\) \\
\(\langle 729,38\rangle=C\)          & \(2\) & b.10 & \((0043)\)    & \(21,21,111,111\) & \(211\)      &   \(1\) &   \(4\) & \(\langle 243,3\rangle\) \\
\(\langle 729,40\rangle=B\)          & \(2\) & b.10 & \((0043)\)    & \(22,21,111,111\) & \(211\)      &   \(2\) &   \(5\) & \(\langle 243,3\rangle\) \\
\(\langle 729,41\rangle=D\)          & \(2\) & d.19 & \((4043)\)    & \(22,21,111,111\) & \(211\)      &   \(1\) &   \(4\) & \(\langle 243,3\rangle\) \\
\(\langle 729,42\rangle\)            & \(2\) & d.23 & \((1043)\)    & \(22,21,111,111\) & \(211\)      &   \(0\) &   \(3\) & \(\langle 243,3\rangle\) \\
\(\langle 729,43\rangle\)            & \(2\) & d.25 & \((2043)\)    & \(22,21,111,111\) & \(211\)      &   \(0\) &   \(3\) & \(\langle 243,3\rangle\) \\
\(\langle 2187,248\vert 249\rangle\) & \(2\) & d.19 & \((4043)\)    & \(32,21,111,111\) & \(221\)      &   \(0\) &   \(4\) & \(\langle 729,40\rangle\) \\
\(\langle 2187,250\rangle\)          & \(2\) & d.23 & \((1043)\)    & \(32,21,111,111\) & \(221\)      &   \(0\) &   \(4\) & \(\langle 729,40\rangle\) \\
\(\langle 2187,251\vert 252\rangle\) & \(2\) & d.25 & \((2043)\)    & \(32,21,111,111\) & \(221\)      &   \(0\) &   \(4\) & \(\langle 729,40\rangle\) \\
\(\langle 2187,253\rangle\)          & \(2\) & b.10 & \((0043)\)    & \(22,21,111,111\) & \(221\)      &   \(1\) &   \(5\) & \(\langle 729,40\rangle\) \\
\(\langle 6561,1989\rangle\)         & \(2\) & d.19 & \((4043)\)    & \(33,21,111,111\) & \(321\)      &   \(0\) &   \(4\) & \(\langle 2187,247\rangle\) \\
\hline
\(\langle 243,8\rangle\)             & \(2\) & c.21 & \((0231)\)    & \(21,21,21,21\)   & \(111\)      &   \(1\) &   \(3\) & \(\langle 27,3\rangle\) \\
\(\langle 729,52\rangle=S\)          & \(2\) & G.16 & \((4231)\)    & \(22,21,21,21\)   & \(211\)      &   \(1\) &   \(3\) & \(\langle 243,8\rangle\) \\
\(\langle 729,54\rangle=U\)          & \(2\) & c.21 & \((0231)\)    & \(22,21,21,21\)   & \(211\)      &   \(2\) &   \(4\) & \(\langle 243,8\rangle\) \\
\(\langle 2187,301\vert 305\rangle\) & \(2\) & G.16 & \((4231)\)    & \(32,21,21,21\)   & \(221\)      &   \(1\) &   \(4\) & \(\langle 729,54\rangle\) \\
\(\langle 2187,303\rangle\)          & \(2\) & c.21 & \((0231)\)    & \(32,21,21,21\)   & \(221\)      &   \(1\) &   \(4\) & \(\langle 729,54\rangle\) \\
\hline
\(\langle 2187,64\rangle=P_7\)       & \(3\) & b.10 & \((0043)\)    & \(22,22,111,111\) & \(2111\)     &   \(4\) &   \(6\) & \(\langle 243,3\rangle\) \\
\(\langle 2187,65\vert 67\rangle\)   & \(3\) & H.4  & \((3343)\)    & \(22,22,111,111\) & \(2111\)     &   \(3\) &   \(5\) & \(\langle 243,3\rangle\) \\
\(\langle 2187,66\vert 73\rangle\)   & \(3\) & F.11 & \((1143)\)    & \(22,22,111,111\) & \(2111\)     &   \(2\) &   \(4\) & \(\langle 243,3\rangle\) \\
\(\langle 2187,69\rangle\)           & \(3\) & G.16 & \((1243)\)    & \(22,22,111,111\) & \(2111\)     &   \(2\) &   \(4\) & \(\langle 243,3\rangle\) \\
\(\langle 2187,71\rangle\)           & \(3\) & G.19 & \((2143)\)    & \(22,22,111,111\) & \(2111\)     &   \(2\) &   \(4\) & \(\langle 243,3\rangle\) \\
\(\langle 6561,676\vert 677\rangle\) & \(3\) & d.19 & \((4043)\)    & \(32,22,111,111\) & \(2211\)     &   \(0\) &   \(5\) & \(\langle 2187,64\rangle\) \\
\(\langle 6561,678\rangle\)          & \(3\) & d.23 & \((1043)\)    & \(32,22,111,111\) & \(2211\)     &   \(0\) &   \(5\) & \(\langle 2187,64\rangle\) \\
\(\langle 6561,679\vert 680\rangle\) & \(3\) & d.25 & \((2043)\)    & \(32,22,111,111\) & \(2211\)     &   \(0\) &   \(5\) & \(\langle 2187,64\rangle\) \\
\(\langle 6561,693..698\rangle\)          & \(3\) & b.10 & \((0043)\)    & \(22,22,111,111\) & \(2211\)     &   \(0\) &   \(5\) & \(\langle 2187,64\rangle\) \\
\hline
\(P_7-\#2;34\vert 35\)               & \(4\) & H.4  & \((3343)\)    & \(32,32,111,111\) & \(2221\)     &   \(1\) &   \(5\) & \(\langle 2187,64\rangle\) \\
\hline
\end{tabular}
}
\end{center}
\end{table}


In Table
\ref{tbl:Metabelian33},
we begin with metabelian groups \(\mathfrak{M}\) of generator rank \(d_1(\mathfrak{M})=2\).
The Shafarevich bound
\cite[Thm. 5.1, p. 28]{Ma2015c}
is given by \(\mu\le\varrho+r+\theta=2+2+0=4\).
For order \(6561\) see
\cite{MAGMA6561}.


\renewcommand{\arraystretch}{1.1}

\begin{table}[hb]
\caption{Invariants of non-metabelian \(3\)-groups \(\mathfrak{G}\) with \(\mathfrak{G}/\mathfrak{G}^\prime\simeq (3,3)\)}
\label{tbl:NonMetabelian33}
\begin{center}
{\scriptsize
\begin{tabular}{|c|r|l|c|c|c|r|r|c|}
\hline
\(\mathfrak{G}\)                      & cc    & Type & \(\varkappa\) & \(\alpha\)        & \(\alpha_2\) & \(\nu\) & \(\mu\) & \(\mathfrak{G}/\mathfrak{G}^{\prime\prime}\) \\
\hline
\(\langle 2187,263..265\rangle\) & \(2\) & d.19 & \((4043)\)    & \(22,21,111,111\) & \(211\)      &   \(0\) &   \(3\) & \(\langle 729,41\rangle\) \\
\(\langle 2187,307\vert 308\rangle\)  & \(2\) & c.21 & \((0231)\)    & \(22,21,21,21\)   & \(211\)      &   \(0\) &   \(3\) & \(\langle 729,54\rangle\) \\
\hline
\(\langle 6561,619\vert 623\rangle\)  & \(3\) & G.16 & \((4231)\)    & \(32,21,21,21\)   & \(221\)      &   \(1\) &   \(3\) & \(\langle 2187,301\vert 305\rangle\) \\
\hline
\end{tabular}
}
\end{center}
\end{table}


Capital letters for \(\mathfrak{M}\) are due to Ascione
\cite{AHL1977}.
For the metabelian groups \(\mathfrak{M}\) with non-trivial cover \(\mathrm{cov}(\mathfrak{M})\)
\cite[p. 30]{Ma2015c},
we need non-metabelian groups \(\mathfrak{G}\) in the cover, which are given in Table
\ref{tbl:NonMetabelian33},
where we begin with groups \(\mathfrak{G}\) of generator rank \(d_1(\mathfrak{G})=2\).
For \(d_1(\mathfrak{G})=3\), we refer to a forthcoming paper.
Instead of the parent \(\pi(\mathfrak{G})\),
we give the metabelianization \(\mathfrak{G}/\mathfrak{G}^{\prime\prime}\).


\section{Categories I and II}
\label{s:Cat1And2}

\noindent
Common feature of these two categories
is the inhomogeneity of \(3\)-class ranks
of the four components in the quartet \((k_\mu)_{\mu=1}^4\)
sharing the conductor \(c=pqr\).
In the present article,
we restrict ourselves to \(3\), respectively \(2\), components
with elementary bicyclic \(3\)-class group
\(\mathrm{Cl}_3(k_\mu)\simeq (3,3)\),
for Category \(\mathrm{I}\), respectively Category \(\mathrm{II}\),
and we postpone elementary tricyclic \(\mathrm{Cl}_3(k_\mu)\simeq (3,3,3)\)
to a future paper.
All computations for examples were performed with Magma
\cite{BCP1997,BCFS2023,MAGMA2023}.


\begin{definition}
\label{dfn:SingularQuartet}
According to the \(3\)-class numbers \(h_3(k_\mu)\),
a quartet \((k_\mu)_{\mu=1}^4\) of cyclic cubic fields
with common conductor \(c=pqr\)
belonging to Category \(\mathrm{I}\) or \(\mathrm{II}\)
is called
\begin{equation}
\label{eqn:SingularQuartet}
\begin{cases}
\textbf{regular}        & \\
\textbf{singular}       & \\
\textbf{super-singular} & \\
\end{cases}
\text{ if } \max\lbrace h_3(k_\mu)\mid 1\le\mu\le 4\rbrace
\begin{cases}
=27, & \\
=81, & \\
\ge 243. & \\
\end{cases}
\end{equation}
\end{definition}

\noindent
In a regular, respectively singular, respectively super-singular, quartet,
there occurs a \(3\)-class group \(\mathrm{Cl}_3(k_\mu)\simeq (3,3,3)\),
respectively \(\mathrm{Cl}_3(k_\mu)\simeq (9,3,3)\),
respectively \(\mathrm{Cl}_3(k_\mu)\simeq (9,9,3)\),
for some \(1\le\mu\le 4\).


\subsection{Category I, Graph 1}
\label{ss:Cat1Gph1}

\noindent
Let \((k_1,\ldots,k_4)\) be a quartet of cyclic cubic number fields
sharing the common conductor \(c=pqr\),
belonging to Graph \(1\) of Category \(\mathrm{I}\)
with combined cubic residue symbol
\(\lbrack p,q,r\rbrack_3=\lbrace p,q,r;\delta\equiv 0\,(\mathrm{mod}\,3)\rbrace\).


Since there are no trivial cubic residue symbols
among the three prime(power) divisors \(p,q,r\)
of the conductor \(c=pqr\),
the principal factors of the subfields with \(t=2\)
of the absolute genus field \(k^\ast\)
must be divisible by both relevant primes,
and we can use the general approach
\begin{equation}
\label{eqn:ElEmEnCat1Gph1}
\begin{aligned}
A(k_{pq})=p^\ell q, & \quad A(\tilde{k}_{pq})=p^{-\ell}q, \\
A(k_{pr})=p^m r, & \quad A(\tilde{k}_{pr})=p^{-m}r, \text{ and } \\
A(k_{qr})=q^nr, & \quad A(\tilde{k}_{qr})=q^{-n}r,
\end{aligned}
\end{equation}
with \(\ell,m,n\in\lbrace -1,1\rbrace\),
identifying \(-1\equiv 2\,(\mathrm{mod}\,3)\),
since it is easier to manage:
\(\ell^2=m^2=n^2=1\).

\begin{lemma}
\label{lem:ElEmEn}
The product \(\ell\cdot m\cdot n=-1\) is negative
(that is, either one or three among \(\ell,m,n\) are negative)
if and only if
the compositum \(L=k_{pq}k_{pr}k_{qr}\) satisfies
the normalization \(\lbrack L:\mathbb{Q}\rbrack=9\):
\begin{equation}
\label{eqn:ElEmEn}
\begin{aligned}
\ell\cdot m\cdot n=-1 \quad & \Longleftrightarrow \quad \lbrack L:\mathbb{Q}\rbrack=9, \\
\ell\cdot m\cdot n=+1 \quad & \Longleftrightarrow \quad \lbrack L:\mathbb{Q}\rbrack=27.
\end{aligned}
\end{equation}
\end{lemma}

\begin{proof}
By Theorem
\ref{thm:Three},
the fields \(L\) and \(\tilde{L}=\tilde{k}_{pq}\tilde{k}_{pr}\tilde{k}_{qr}\),
satisfy a skew balance of their degrees \(\in\lbrace 9,27\rbrace\) in the product
\(\lbrack L:\mathbb{Q}\rbrack\cdot\lbrack\tilde{L}:\mathbb{Q}\rbrack=243\).

Suppose \(lmn=+1\).
Then we produce a contradiction by the assumption that
\(\lbrack L:\mathbb{Q}\rbrack=9\) and
\(\lbrack\tilde{L}:\mathbb{Q}\rbrack=27\).
We define the compositum \(K:=\tilde{k}_{pq}\tilde{k}_{pr}\) of degree \(9\).
Then \(K\) contains one of the fields \(k_\mu\), \(\mu=1,\ldots,4\),
and either \(\tilde{k}_{qr}\) or \(k_{qr}\).
In the former case, \(K=\tilde{L}\) would have degree \(27\).
So \(K=\tilde{k}_{pq}\tilde{k}_{pr}k_{qr}\),
and we calculate the following sub-determinants
of the principal factor matrices \(M_L\) and \(M_K\),
with respect to the fields with \(t=2\) only (ignoring the field with \(t=3\)): \\
\(
\left\vert
\begin{matrix}
\ell & 1 & 0 \\
m & 0 & 1 \\
0 & n & 1
\end{matrix}
\right\vert
=-\ell n-m=0
\)
\(\quad\Longleftrightarrow\quad\)
\(
\left\vert
\begin{matrix}
-\ell & 1 & 0 \\
-m & 0 & 1 \\
0 & n & 1
\end{matrix}
\right\vert
=\ell n+m=0
\)
\(\quad\Longleftrightarrow\quad\)
\(\ell n=-m\). \\
However, \(lmn=+1\) implies \(\ell n=m\) and thus rank \(3\) of \(M_L\) and \(M_K\).
By
\eqref{eqn:UnitIndex},
this gives indices of subfield units
\((U_L:V_L)=3\) and \((U_K:V_K)=3\).
At least one among \(L\) and \(K\), say \(X\),
does not contain the critical field \(k_\mu\) with \(\varrho_3(k_\mu)=3\),
whence it is tame with \(h_3(X)=(U_X:V_X)=3\),
in contradiction to \(9\mid h_3(X)\), by Proposition
\ref{prp:Bicub}.
Thus we must have \(\lbrack L:\mathbb{Q}\rbrack=27\).

With nearly identical arguments,
it is easy to show that
\(lmn=-1\) implies \(\lbrack L:\mathbb{Q}\rbrack=9\).
\end{proof}


\begin{lemma}
\label{lem:Cat1Gph1}
\textbf{(3-class ranks of components for \(\mathrm{I}.1\).)}
Without loss of generality,
precisely three components \(k_2\), \(k_3\), \(k_4\) of the quartet
have elementary bicyclic \(3\)-class groups
\(\mathrm{Cl}_3(k_\mu)\simeq (3,3)\), \(2\le\mu\le 4\),
whereas the single remaining component \(k_1\)
has \(3\)-class rank \(\varrho_3(k_1)=3\).
In dependence on the 
\textbf{decisive principal factors} in Equation
\eqref{eqn:ElEmEnCat1Gph1},
the principal factors of \(k_\mu\) are
\begin{equation}
\label{eqn:PFCat1Gph1}
\begin{aligned}
A(k_2)=pq^2r,\quad A(k_3)=pqr,\quad A(k_4)=pqr^2 \quad & \text{ if } (\ell,m,n)=(1,1,2), \\
A(k_2)=p^2qr,\quad A(k_3)=pqr^2,\quad A(k_4)=pqr \quad & \text{ if } (\ell,m,n)=(1,2,1), \\
A(k_2)=pqr,\quad A(k_3)=pq^2r,\quad A(k_4)=p^2qr \quad & \text{ if } (\ell,m,n)=(2,1,1), \\
A(k_2)=pqr^2,\quad A(k_3)=p^2qr,\quad A(k_4)=pq^2r \quad & \text{ if } (\ell,m,n)=(2,2,2).
\end{aligned}
\end{equation}
The \textbf{tame} condition \(9\mid h_3(B_j)=(U_j:V_j)\in\lbrace 9,27\rbrace\)
with \(r_j=2\) is satisfied
for \(j\in\lbrace 2,3,4,8,9,10\rbrace\).

A further 
\textbf{decisive principal factor}
\(A(k_1)=p^{e_1}q^{e_2}r^{e_3}\)
and the associated invariant counter
\(\mathcal{D}:=\#\lbrace 1\le i\le 3\mid e_i\ne 0\rbrace\)
admit several conclusions for
\textbf{wild} ranks:
\begin{equation}
\label{eqn:RanksCat1Gph1}
r_5=r_6=r_7=3 \quad \text{ iff } \quad \mathcal{D}=2 \quad
\text{ iff } \quad A(k_1) \text{ has precisely two prime divisors}.
\end{equation}
\end{lemma}

\begin{proof}
According to
\cite[Prop. 4.4, pp. 43--44]{Ay1995},
the required condition to distinguish the unique component \(k_1\)
with \(3\)-class rank \(\varrho(k_1)=3\)
in the quartet \((k_\mu)_{\mu=1}^4\)
is the set of decomposition invariants \((e,f,g)=(3,1,3)\)
simultaneously for \(p,q,r\)
in the bicyclic bicubic field \(B_1=k_1k_{pq}k_{pr}k_{qr}\),
that is, \\
\(p\) splits in \(k_{qr}\), and thus also in \(B_j\) for \(j\in\lbrace 1,3,8\rbrace\), \\
\(q\) splits in \(k_{pr}\), and thus also in \(B_j\) for \(j\in\lbrace 1,4,9\rbrace\), \\
\(r\) splits in \(k_{pq}\), and thus also in \(B_j\) for \(j\in\lbrace 1,2,10\rbrace\).

Then, exactly the six fields
\(B_2=k_2k_{pq}\tilde{k}_{pr}\tilde{k}_{qr}\),
\(B_3=k_3\tilde{k}_{pq}\tilde{k}_{pr}k_{qr}\),
\(B_4=k_4\tilde{k}_{pq}k_{pr}\tilde{k}_{qr}\),
\(B_8=k_2k_4k_pk_{qr}\),
\(B_9=k_2k_3k_qk_{pr}\),
\(B_{10}=k_3k_4k_rk_{pq}\)
do not contain \(k_1\),
and satisfy the \textit{tame} relation
\(9\mid h_3(B_j)=(U_j:V_j)\in\lbrace 9,27\rbrace\)
with ranks \(r_j=2\) for \(j=2,3,4,8,9,10\),
by Proposition
\ref{prp:Bicub}.

This fact can be exploited
for each tame bicyclic bicubic field \(B_j\),
by calculating the rank \(r_j\) with row operations
on the associated principal factor matrix \(M_j\)
and drawing conclusions for the exponents \(x_\mu,y_\mu,z_\mu\)
in the approach
\(A(k_\mu)=p^{x_\mu}q^{y_\mu}r^{z_\mu}\), \(1\le\mu\le 4\): \\
\(M_8=
\begin{pmatrix}
x_2 & y_2 & z_2 \\
x_4 & y_4 & z_4 \\
1 & 0 & 0 \\
0 & n & 1 \\
\end{pmatrix}\),
\(M_9=
\begin{pmatrix}
x_2 & y_2 & z_2 \\
x_3 & y_3 & z_3 \\
0 & 1 & 0 \\
m & 0 & 1 \\
\end{pmatrix}\),
\(M_{10}=
\begin{pmatrix}
x_3 & y_3 & z_3 \\
x_4 & y_4 & z_4 \\
0 & 0 & 1 \\
\ell & 1 & 0 \\
\end{pmatrix}\).

For \(B_8=k_2k_4k_pk_{qr}\),
\(M_8\)
leads to decisive pivot elements \(z_2-ny_2\) and \(z_4-ny_4\)
in the last column,
for \(B_9=k_2k_3k_qk_{pr}\), 
\(M_9\)
leads to \(z_2-mx_2\) and \(z_3-mx_3\)
in the last column,
and for \(B_{10}=k_3k_4k_rk_{pq}\),
\(M_{10}\)
leads to \(y_3-\ell x_3\) and \(y_4-\ell x_4\)
in the middle column.

So \(r_8=r_9=r_{10}=2\) implies
\(ny_2\equiv z_2\),
\(ny_4\equiv z_4\),
\(z_2\equiv mx_2\),
\(z_3\equiv mx_3\),
\(\ell x_3\equiv y_3\),
\(\ell x_4\equiv y_4\).
Or, in combined form,
\(mx_2\equiv ny_2\equiv z_2\),
\(mx_3\equiv -ny_3\equiv z_3\),
\(-mx_4\equiv ny_4\equiv z_4\).
This yields \eqref{eqn:PFCat1Gph1}.

Additionally, we use the remaining three tame ranks for \\
\(M_2=
\begin{pmatrix}
x_2 & y_2 & z_2 \\
\ell & 1 & 0 \\
-m & 0 & 1 \\
0 & -n & 1 \\
\end{pmatrix}\),
\(M_3=
\begin{pmatrix}
x_3 & y_3 & z_3 \\
-\ell & 1 & 0 \\
-m & 0 & 1 \\
0 & n & 1 \\
\end{pmatrix}\),
\(M_4=
\begin{pmatrix}
x_4 & y_4 & z_4 \\
-\ell & 1 & 0 \\
m & 0 & 1 \\
0 & -n & 1 \\
\end{pmatrix}\).

For \(B_2=k_2k_{pq}\tilde{k}_{pr}\tilde{k}_{qr}\),
\(M_2\)
leads to the decisive pivot elements
\(z_2+mx_2+ny_2\), \(\ell m+n\)
in the last column,
for \(B_3=k_3\tilde{k}_{pq}\tilde{k}_{pr}k_{qr}\), 
\(M_3\)
leads to
\(z_3+mx_3-ny_3\), \(-\ell m-n\)
in the last column,
and
for \(B_4=k_4\tilde{k}_{pq}k_{pr}\tilde{k}_{qr}\),
\(M_4\)
leads to
\(z_4-mx_4+ny_4\), \(\ell m+n\)
in the last column.
So \(r_2=r_3=r_4=2\) implies
\(mx_2+ny_2\equiv -z_2\),
\(mx_3-ny_3\equiv -z_3\),
\(-mx_4+ny_4\equiv -z_4\),
since the other pivot elements vanish a priori,
\(\ell m+n=0\), i.e. \(\ell m=-n\),
because \(\ell mn=-1\) and \(n^2=1\)
in Lemma
\ref{lem:ElEmEn}.
The congruences follow already
from those for \(r_8=r_9=r_{10}=2\).

For each \textit{wild} bicyclic bicubic field \(B_j\),
\(j\in\lbrace 1,5,6,7\rbrace\),
the rank \(r_j\) is now calculated with row operations
on the associated principal factor matrix \(M_j\): \\
\(M_1=
\begin{pmatrix}
x_1 & y_1 & z_1 \\
\ell & 1 & 0 \\
m & 0 & 1 \\
0 & n & 1 \\
\end{pmatrix}\),
\(M_5=
\begin{pmatrix}
x_1 & y_1 & z_1 \\
x_3 & y_3 & z_3 \\
1 & 0 & 0 \\
0 & -n & 1 \\
\end{pmatrix}\),
\(M_6=
\begin{pmatrix}
x_1 & y_1 & z_1 \\
x_4 & y_4 & z_4 \\
0 & 1 & 0 \\
-m & 0 & 1 \\
\end{pmatrix}\),
\(M_7=
\begin{pmatrix}
x_1 & y_1 & z_1 \\
x_2 & y_2 & z_2 \\
0 & 0 & 1 \\
-\ell & 1 & 0 \\
\end{pmatrix}\).

For \(B_1=k_1k_{pq}k_{pr}k_{qr}\),
\(M_1\)
leads to the decisive pivot element
\(z_1-mx_1-ny_1\),
since \(-\ell m-n\equiv 0\).
So, \(r_1=2\) implies \(z_1\equiv mx_1+ny_1\).
For \(B_5=k_1k_3k_p\tilde{k}_{qr}\),
\(M_5\)
leads to
\(z_1+ny_1\), \(z_3+ny_3\)
in the last column.
So, \(r_5=3\) iff 
either \(-z_1\not\equiv ny_1\)
or \(-z_3\not\equiv ny_3\) modulo \(3\).
For \(B_6=k_1k_4k_q\tilde{k}_{pr}\), 
\(M_6\)
leads to
\(z_1+mx_1\), \(z_4+mx_4\).
So, \(r_6=3\) iff
either \(-z_1\not\equiv mx_1\)
or \(-z_4\not\equiv mx_4\).
For \(B_7=k_1k_2k_r\tilde{k}_{pq}\), 
\(M_7\)
leads to \(y_1+\ell x_1\), \(y_2+\ell x_2\)
in the middle column.
So, \(r_7=3\) iff
either \(-y_1\not\equiv\ell x_1\)
or \(-y_2\not\equiv\ell x_2\).
For each of these three ranks,
the \textit{second} condition can \textit{never} be satisfied.

Since at most one of the exponents \(x_1,y_1,z_1\)
may vanish, the new congruences immediately lead to
\eqref{eqn:RanksCat1Gph1}.
For instance,
\(z_1=0\) \(\Rightarrow\)
\(-z_1=0\not\equiv ny_1\), \(-z_1=0\not\equiv mx_1\) \(\Rightarrow\)
\(r_5=r_6=3\);
but \(0=z_1\equiv mx_1+ny_1\) also implies
\(mx_1\equiv -ny_1\), \(mnx_1\equiv -y_1\), \(\ell x_1\equiv y_1\)
and thus \(r_7=3\).
Conversely, suppose \(\mathcal{D}=3\).
If \(-z_1\not\equiv ny_1\), then \(z_1\equiv ny_1\),
and \(z_1\equiv mx_1+ny_1\) implies \(mx_1\equiv 0\),
and thus the contradiction \(x_1=0\).
\end{proof}


\begin{proposition}
\label{prp:Cat1Gph1}
\textbf{(Sub-triplet with \(3\)-rank two for \(\mathrm{I}.1\).)}
For fixed \(\mu\in\lbrace 2,3,4\rbrace\),
let \(\mathfrak{p},\mathfrak{q},\mathfrak{r}\) be the prime ideals of \(k_\mu\)
over \(p,q,r\), that is
\(p\mathcal{O}_{k_\mu}=\mathfrak{p}^3\),
\(q\mathcal{O}_{k_\mu}=\mathfrak{q}^3\),
\(r\mathcal{O}_{k_\mu}=\mathfrak{r}^3\),
then the \(3\)-class group of \(k_\mu\) is generated by any two among
the non-trivial classes \(\lbrack\mathfrak{p}\rbrack,\lbrack\mathfrak{q}\rbrack,\lbrack\mathfrak{r}\rbrack\),
that is,
\begin{equation}
\label{eqn:GenCat1Gph1}
\mathrm{Cl}_3(k_\mu)=
\langle\lbrack\mathfrak{p}\rbrack,\lbrack\mathfrak{q}\rbrack\rangle=
\langle\lbrack\mathfrak{p}\rbrack,\lbrack\mathfrak{r}\rbrack\rangle=
\langle\lbrack\mathfrak{q}\rbrack,\lbrack\mathfrak{r}\rbrack\rangle
\simeq (3,3).
\end{equation}
The unramified cyclic cubic relative extensions of \(k_\mu\)
are among the absolutely bicyclic bicubic subfields \(B_i\), \(1\le i\le 10\),
of the common genus field \(k^\ast\)
of the four components of the quartet \((k_1,\ldots,k_4)\).
The unique \(B_9>k_\mu\), \(\mu\in\lbrace 2,3\rbrace\),
has norm class group
\(N_{B_9/k_\mu}(\mathrm{Cl}_3(B_9))=\langle\lbrack\mathfrak{q}\rbrack\rangle\),
and potential \textbf{fixed point} transfer kernel
\[
\ker(T_{B_9/k_\mu})\ge\langle\lbrack\mathfrak{q}\rbrack\rangle.
\]
The unique \(B_{10}>k_\mu\), \(\mu\in\lbrace 3,4\rbrace\),
has norm class group
\(N_{B_{10}/k_\mu}(\mathrm{Cl}_3(B_{10}))=\langle\lbrack\mathfrak{r}\rbrack\rangle\),
and potential \textbf{fixed point} transfer kernel
\[
\ker(T_{B_{10}/k_\mu})\ge\langle\lbrack\mathfrak{r}\rbrack\rangle.
\]
The unique \(B_8>k_\mu\), \(\mu\in\lbrace 2,4\rbrace\),
has norm class group
\(N_{B_8/k_\mu}(\mathrm{Cl}_3(B_8))=\langle\lbrack\mathfrak{qr}^n\rbrack\rangle\),
and potential \textbf{fixed point} transfer kernel
\[
\ker(T_{B_8/k_\mu})\ge\langle\lbrack\mathfrak{q}^n\mathfrak{r}\rbrack\rangle.
\]
The remaining \(B_i>k_\mu\), \(i\in\lbrace 5,6,7\rbrace\),
more precisely, \(i=7\) for \(\mu=2\),
\(i=5\) for \(\mu=3\),
and \(i=6\) for \(\mu=4\),
have norm class group 
\(\langle\lbrack\mathfrak{q}^2\mathfrak{r}^{-n}\rbrack\rangle\)
and a hidden or explicit \textbf{transposition} transfer kernel,
with respect to the corresponding \(\mu\).
\end{proposition}

\begin{proof}
As mentioned in the proof of Lemma
\ref{lem:Cat1Gph1},
\(q\) splits in \(B_9\),
\(r\) splits in \(B_{10}\), and
\(p\) splits in \(B_8\),
where \(\lbrack\mathfrak{p}\rbrack=\lbrack\mathfrak{qr}^n\rbrack\),
according to 
\eqref{eqn:PFCat1Gph1},
independently of \(n\in\lbrace 1,2\rbrace\). \\
Now we use Corollary
\ref{cor:Capitulation}
and Proposition
\ref{prp:Bicyc}. \\
Since \(\mathfrak{q}\) is principal in \(k_q\),
\(\lbrack\mathfrak{q}\rbrack\) capitulates in
\(B_6=k_1k_4k_q\tilde{k}_{pr}\) and
\(B_9=k_2k_3k_qk_{pr}\). \\
Since \(\mathfrak{r}\) is principal in \(k_r\),
\(\lbrack\mathfrak{r}\rbrack\) capitulates in
\(B_7=k_1k_2k_r\tilde{k}_{pq}\) and
\(B_{10}=k_3k_4k_rk_{pq}\). \\
Since \(\mathfrak{q}^n\mathfrak{r}\) is principal in \(k_{qr}\),
\(\lbrack\mathfrak{q}^n\mathfrak{r}\rbrack\) capitulates in
\(B_3=k_3\tilde{k}_{pq}\tilde{k}_{pr}k_{qr}\) and
\(B_8=k_2k_4k_pk_{qr}\).
\end{proof}


\noindent
In terms of capitulation targets in Corollary
\ref{cor:Three},
Proposition
\ref{prp:Cat1Gph1}
and parts of its proof are now summarized in Table
\ref{tbl:UniCat1Gph1}
for the minimal transfer kernel type (mTKT) and \(n=2\),
with transposition in \textbf{boldface} font.
This essential new perspective admits
progress beyond Ayadi's work
\cite{Ay1995}.


\renewcommand{\arraystretch}{1.1}

\begin{table}[ht]
\caption{Norm class groups and minimal transfer kernels with \(n=2\) for Graph I.1}
\label{tbl:UniCat1Gph1}
\begin{center}
{\normalsize
\begin{tabular}{|c||c|c|c|c||c|c|c|c||c|c|c|c|}
\hline
 Base         & \multicolumn{4}{c||}{\(k_2\)} & \multicolumn{4}{c||}{\(k_3\)} & \multicolumn{4}{c|}{\(k_4\)} \\
\hline
 Ext          & \(B_2\) & \(B_7\) & \(B_8\) & \(B_9\) & \(B_3\) & \(B_5\) & \(B_9\) & \(B_{10}\) & \(B_4\) & \(B_6\) & \(B_8\) & \(B_{10}\) \\
\hline
 NCG          & \(\mathfrak{r}\)    & \(\mathfrak{qr}\)   & \(\mathfrak{qr}^2\) & \(\mathfrak{q}\)    
              & \(\mathfrak{qr}\)   & \(\mathfrak{qr}^2\) & \(\mathfrak{q}\)    & \(\mathfrak{r}\) & \(\mathfrak{q}\)  & \(\mathfrak{qr}\) & \(\mathfrak{qr}^2\) & \(\mathfrak{r}\) \\
 TK           & \(\mathfrak{qr}\)   & \(\mathfrak{r}\)    & \(\mathfrak{qr}^2\) & \(\mathfrak{q}\) 
              & \(\mathfrak{qr}^2\) & \(\mathfrak{qr}\)   & \(\mathfrak{q}\)    & \(\mathfrak{r}\) & \(\mathfrak{qr}\) & \(\mathfrak{q}\)  & \(\mathfrak{qr}^2\) & \(\mathfrak{r}\) \\
\(\varkappa\) &  \(\mathbf{2}\) & \(\mathbf{1}\) & \(3\) & \(4\) &  \(\mathbf{2}\) & \(\mathbf{1}\) & \(3\) & \(4\) & \(\mathbf{2}\) & \(\mathbf{1}\) & \(3\) & \(4\) \\
\hline
\end{tabular}
}
\end{center}
\end{table}


\begin{theorem}
\label{thm:Cat1Gph1}
\textbf{(Second \(3\)-class group for \(\mathrm{I}.1\).)}
Let \((k_1,\ldots,k_4)\) be a quartet of cyclic cubic number fields
sharing the common conductor \(c=pqr\)
belonging to Graph \(1\) of Category \(\mathrm{I}\), that is,
\(\lbrack p,q,r\rbrack_3=\lbrace p,q,r;\delta\equiv 0\,(\mathrm{mod}\,3)\rbrace\).
Without loss of generality, suppose that 
\(\mathrm{Cl}_3(k_\mu)\simeq (3,3)\), for \(\mu=2,3,4\),
and \(\varrho_3(k_1)=3\).

Then the \textbf{minimal transfer kernel type} (mTKT) of \(k_\mu\), \(2\le\mu\le 4\),
is \(\varkappa_0=(4231)\), type \(\mathrm{G}.16\),
and other possible capitulation types in ascending partial order
\(\varkappa_0<\varkappa<\varkappa^\prime,\varkappa^{\prime\prime}<\varkappa^{\prime\prime\prime}\)
are
\(\varkappa=(0231)\), type \(\mathrm{c}.21\),
\(\varkappa^\prime=(0001)\), type \(\mathrm{a}.3\),
\(\varkappa^{\prime\prime}=(0200)\), type \(\mathrm{a}.2\), and
\(\varkappa^{\prime\prime\prime}=(0000)\), type \(\mathrm{a}.1\).

In order to identify the second \(3\)-class group
\(\mathfrak{M}=\mathrm{Gal}(\mathrm{F}_3^2(k_\mu)/k_\mu)\), \(2\le\mu\le 4\),
let the \textbf{principal factor} of \(k_1\) be
\(A(k_1)=p^{e_1}q^{e_2}r^{e_3}\),
and define \(\mathcal{D}:=\#\lbrace 1\le i\le 3\mid e_i\ne 0\rbrace\).
In the \textbf{regular} situation where
\(\mathrm{Cl}_3(k_1)\simeq (3,3,3)\)
is elementary tricyclic, we have
\begin{equation}
\label{eqn:RegCat1Gph1}
\mathfrak{M}\simeq
\begin{cases}
\langle 81,8\rangle,\ \alpha=\lbrack 11,11,11,21\rbrack,\ \varkappa=(0001) & \text{ once if } \mathcal{D}=2,\ \mathcal{N}=1, \\
\langle 81,10\rangle^2,\ \alpha=\lbrack 11,21,11,11\rbrack,\ \varkappa=(0200) & \text{ twice if } \mathcal{D}=2,\ \mathcal{N}=1, \\
\langle 243,25\rangle,\ \alpha=\lbrack 11,11,11,22\rbrack,\ \varkappa=(0001) & \text{ once if } \mathcal{D}=3,\ \mathcal{N}=1, \\
\langle 243,28..30\rangle^2,\ \alpha=\lbrack 21,11,11,11\rbrack,\ \varkappa=(0000) & \text{ twice if } \mathcal{D}=3,\ \mathcal{N}=0,
\end{cases}
\end{equation}
where \(\mathcal{N}:=\#\lbrace 1\le j\le 10\mid k_\mu<B_j,\ I_j=27\rbrace\).
In the \textbf{(super-)singular} situation where
\(81\mid h(k_1)\) and \(\mathrm{Cl}_3(k_1)\)
is non-elementary tricyclic, we have \(\mathfrak{M}\simeq\)
\begin{equation}
\label{eqn:SngCat1Gph1}
\begin{cases}
\langle 243,8\rangle^3,\ \alpha=\lbrack 21,21,21,21\rbrack,\ \varkappa=(0231) & \text{ if } h_3(k_1)=81,\ \mathcal{D}=2,\ \mathcal{N}=3, \\
\langle 729,54\rangle^3,\ \alpha=\lbrack 22,21,21,21\rbrack,\ \varkappa=(0231) & \text{ if } h_3(k_1)=81,\ \mathcal{D}=3,\ \mathcal{N}=3, \\
\langle 2187,301\vert 305\rangle^3,\ \alpha=\lbrack 32,21,21,21\rbrack,\ \varkappa=(4231) & \text{ if }
h_3(k_1)=81,\ \mathcal{D}=3,\ \mathcal{N}=4, \\
\langle 2187,303\rangle^3,\ \alpha=\lbrack 32,21,21,21\rbrack,\ \varkappa=(0231) & \text{ if } 
h_3(k_1)=243,\ \mathcal{D}=3,\ \mathcal{N}=3.
\end{cases}
\end{equation}
With exception of the last three rows, the \(3\)-class field tower has the group
\(\mathfrak{G}=\mathrm{Gal}(\mathrm{F}_3^\infty(k_\mu)/k_\mu)\simeq\mathfrak{M}\),
\(\ell_3(k_\mu)=2\),
since \(d_2(\mathfrak{M})\le 4\).
For the last three rows,
the tower length is \(2\le\ell_3(k_\mu)\le 3\)
\cite{Ma2015c}.
(See the associated descendant tree \(\mathcal{T}^2\langle 243,8\rangle\) in
\cite[Fig. 6.4, p. 63]{Ma2022}.)
\end{theorem}

\begin{proof}
In the \textit{non-uniform regular} situations,
we have \(r_j=2\), \(h_3(B_j)=I_j\in\lbrace 9,27\rbrace\) for the tame bicyclic bicubic fields \(j\in\lbrace 2,3,4,8,9,10\rbrace\).
Now we use Lemma
\ref{lem:UnitIndices}
and Lemma
\ref{lem:Cat1Gph1}.

If \(\mathcal{D}=3\),
then all tame indices of subfield units \(I_j=9\) are minimal,
and the ranks of wild bicyclic bicubic fields are \(r_j=2\) for \(j=5,6,7\),
but non-uniform indices
two times \(I_5=I_6=9\), i.e. \(\mathcal{N}=0\),
and one time \(I_7=27\), i.e. \(\mathcal{N}=1\),
corresponding to total capitulation twice and
non-fixed point capitulation once
(due to a hidden transposition). 
According to Theorem
\ref{thm:ClassGroup33},
the common \(\alpha_2=(21)\),
and Corollary
\ref{cor:ClassGroup33},
the Artin pattern
\(\alpha=\lbrack 21,11,11,11\rbrack\) and \(\varkappa=(0000)\)
determines three possible groups
\(\langle 243,28..30\rangle\),
and \(\alpha=\lbrack 11,11,11,22\rbrack\), \(\varkappa=(0001)\)
uniquely leads to \(\langle 243,25\rangle\).

If \(\mathcal{D}=2\),
then tame indices of subfield units \(I_j\) are non-uniform,
two times \(I_j=9\), for \(j=2,4,9,10\),
and one time \(I_j=27\), for \(j=3,8\),
the latter corresponding to 
fixed point capitulation twice, \(j=8\) over \(\mu=2,4\), and
non-fixed point capitulation once, \(j=3\) over \(\mu=3\). 
So \(\mathcal{N}=1\),
since the ranks of wild bicyclic bicubic fields are \(r_j=3\)
with uniform index \(I_j=3\) for \(j=5,6,7\),
corresponding to a total capitulation.
The Artin pattern
\(\alpha=\lbrack 11,21,11,11\rbrack\), \(\varkappa=(0200)\)
uniquely determines the group
\(\langle 81,10\rangle\),
and \(\alpha=\lbrack 11,11,11,21\rbrack\), \(\varkappa=(0001)\)
uniquely leads to \(\langle 81,8\rangle\).

In the \textit{uniform singular} situation with TKT \(\mathrm{c}.21\), \(\varkappa=(0231)\),
\(\mathcal{N}=3\), the ATI decide about the group:
\(\alpha=\lbrack 21,21,21,21\rbrack\) uniquely identifies \(\langle 243,8\rangle\),
\(\alpha=\lbrack 22,21,21,21\rbrack\) leads to \(\langle 729,54\rangle\),
and in the super-singular situation,
\(\alpha=\lbrack 32,21,21,22\rbrack\) leads to \(\langle 2187,303\rangle\).
In contrast,
for TKT \(\mathrm{G}.16\), \(\varkappa=(4231)\),
\(\mathcal{N}=4\), the ATI \(\alpha=\lbrack 32,21,21,21\rbrack\)
lead to \(\langle 2187,301\vert 305\rangle\).

The regular groups are of maximal class, which guarantees length \(\ell_3(k_\mu)=2\) of the tower.
The annihilator ideal of \(\langle 243,8\rangle\) is \(\mathfrak{L}\),
which enforces \(\ell_3(k_\mu)=2\), according to Scholz and Taussky
\cite{SoTa1934}.
The (super-)singular groups \(\langle 729,54\rangle\) and \(\langle 2187,303\rangle\) have non-metabelian descendants.
Although they satisfy the bound \(d_2(\mathfrak{G})\le 4\) for the relation rank,
a tower with three stages could only be excluded
by means of computationally expensive invariants \(\alpha^{(2)}\) of second order.
\end{proof}


\begin{corollary}
\label{cor:UniCat1Gph1}
\textbf{(Non-uniformity of the sub-triplet for \(\mathrm{I}.1\).)}
Only two components of the sub-triplet with \(3\)-rank two
share a common capitulation type
\(\varkappa(k_\lambda)\sim\varkappa(k_\mu)\),
common abelian type invariants
\(\alpha(k_\lambda)\sim\alpha(k_\mu)\),
and a common second \(3\)-class group
\(\mathrm{Gal}(\mathrm{F}_3^2(k_\lambda)/k_\lambda)\simeq
\mathrm{Gal}(\mathrm{F}_3^2(k_\mu)/k_\mu)\).
The invariants of the third component \(k_\nu\)
\textbf{differ} in the \textbf{regular} situation \(\mathrm{Cl}_3(k_1)\simeq (3,3,3)\),
however, they \textbf{agree} in the \textbf{(super-)singular} situation \(81\mid h(k_1)\).
Here, \(\lbrace\lambda,\mu,\nu\rbrace=\lbrace 2,3,4\rbrace\).
\end{corollary}

\begin{proof}
This is an immediate consequence of Theorem
\ref{thm:Cat1Gph1}.
\end{proof}


\begin{example}
\label{exm:Cat1Gph1}
Examples
\ref{exm:Cat1Gph1}--\ref{exm:Cat3Gph9}
are supplemented by
\cite[Tbl. 6.4--21, pp. 49--67]{Ma2022}.
The prototypes for Graph \(\mathrm{I}.1\), that is,
the minimal conductors for each scenario in Theorem
\ref{thm:Cat1Gph1},
are as follows.

There are \textbf{regular} cases:
\(c=4\,977\) with symbol
\(\lbrace 9,7,79\rbrace\)
and, non-uniformly, \(\mathfrak{G}=\mathfrak{M}=\langle 243,25\rangle\),
\(\langle 243,28..30\rangle^2\) (Corollary
\ref{cor:ClassGroup33});
\(c=11\,349\) with symbol
\(\lbrace 9,13,97\rbrace\)
and, non-uniformly, \(\mathfrak{G}=\mathfrak{M}=\langle 81,8\rangle\),
\(\langle 81,10\rangle^2\).

Further,
\textbf{singular} cases:
\(c=28\,791\) with symbol
\(\lbrace 9,7,457\rbrace\)
and \(\mathfrak{M}=\langle 729,54\rangle^3\);
\(c=38\,727\) with symbol
\(\lbrace 9,13,331\rbrace\)
and \(\mathfrak{G}=\mathfrak{M}=\langle 243,8\rangle^3\);
and, with \textbf{considerable statistic delay}, there occurred
\(c=417\,807\) with ordinal number \(189\), symbol
\(\lbrace 9,13,3571\rbrace\)
and \(\mathfrak{M}=\langle 2187,301\vert 305\rangle^3\).

And \textbf{super-singular} cases:
\(c=67\,347\) with symbol
\(\lbrace 9,7,1069\rbrace\)
and \(\mathfrak{M}=\langle 2187,303\rangle^3\);
\(c=436\,267\) with symbol
\(\lbrace 13,37,907\rbrace\)
and \(\mathfrak{M}=(\langle 6561,2050\rangle-\#1;3\vert 5)^3\).
\end{example}


In Table
\ref{tbl:ProtoCat1Gph1},
we summarize the prototypes of graph \(\mathrm{I}.1\).
Data comprises
ordinal number No.,
conductor \(c\) of \(k\),
combined cubic residue symbol \(\lbrack p,q,r\rbrack_3\),
regularity, resp. (super-)singularity, expressed by
\(3\)-valuation \(v=v_3(\#\mathrm{Cl}(k_1))\)
of class number of critical field \(k_1\),
critical exponents \(x,y,z\) in principal factor \(A(k_1)=p^xq^yr^z\)
and \(\ell,m,n\) in \(A(k_{pq})=p^\ell q\), \(A(k_{pr})=p^mr\), \(A(k_{qr})=q^nr\),
capitulation type of \(k\),
second \(3\)-class group
\(\mathfrak{M}=\mathrm{Gal}(\mathrm{F}_3^2(k)/k)\) of \(k\),
and length \(\ell_3(k)\) of \(3\)-class field tower of \(k\).
We put \(R:=\langle 6561,2050\rangle\) for abbreviation.

\renewcommand{\arraystretch}{1.1}

\begin{table}[ht]
\caption{Prototypes for Graph I.1}
\label{tbl:ProtoCat1Gph1}
\begin{center}
{\scriptsize
\begin{tabular}{|r|r||c|c|c|c|c|c|}
\hline
 No. & \(c\)        & \(p,q,r\)    & \(v\) & \(x,y,z;\ell,m,n\) & capitulation type             & \(\mathfrak{M}\)                                           & \(\ell_3(k)\) \\
\hline
   1 &   \(4\,977\) & \(9,7,79\)   & \(3\) & \(2,1,1;1,1,2\) & \(\mathrm{a}.3,\mathrm{a}.1\) & \(\langle 243,25\rangle,\langle 243,28..30\rangle^2\)      & \(=2\) \\
   3 &  \(11\,349\) & \(9,13,97\)  & \(3\) & \(0,1,1;2,1,1\) & \(\mathrm{a}.3,\mathrm{a}.2\) & \(\langle 81,8\rangle,\langle 81,10\rangle^2\)             & \(=2\) \\
  10 &  \(28\,791\) & \(9,7,457\)  & \(4\) & \(2,1,1;1,1,2\) & \(\mathrm{c}.21\)             & \(\langle 729,54\rangle^3\)                                & \(\ge 2\) \\
  14 &  \(38\,727\) & \(9,13,331\) & \(4\) & \(1,0,1;2,1,1\) & \(\mathrm{c}.21\)             & \(\langle 243,8\rangle^3\)                                 & \(=2\) \\
  27 &  \(67\,347\) & \(9,7,1069\) & \(5\) & \(2,1,1;1,1,2\) & \(\mathrm{c}.21\)             & \(\langle 2187,303\rangle^3\)                              & \(\ge 2\) \\
 189 & \(417\,807\) & \(9,13,3571\)& \(4\) & \(2,2,1;2,2,2\) & \(\mathrm{G}.16\)             & \(\langle 2187,301\vert 305\rangle^3\)                     & \(\ge 2\) \\
 198 & \(436\,267\) & \(13,37,907\)& \(6\) & \(1,1,1;2,2,2\) & \(\mathrm{G}.16\)             & \((R-\#1;3\vert 5)^3\)                     & \(\ge 2\) \\
\hline
\end{tabular}
}
\end{center}
\end{table}


\subsection{Category I, Graph 2}
\label{ss:Cat1Gph2}

\noindent
Let \((k_1,\ldots,k_4)\) be a quartet of cyclic cubic number fields
sharing the common conductor \(c=pqr\),
belonging to Graph \(2\) of Category \(\mathrm{I}\)
with combined cubic residue symbol
\(\lbrack p,q,r\rbrack_3=\lbrace q\leftarrow p\rightarrow r\rbrace\).


\begin{lemma}
\label{lem:Cat1Gph2}
\textbf{(3-class ranks of components for \(\mathrm{I}.2\).)}
Under the normalizing assumptions that
\(q\) splits in \(k_{pr}\) and
\(r\) splits in \(k_{pq}\),
precisely the three components \(k_2\), \(k_3\), \(k_4\) of the quartet
have elementary bicyclic \(3\)-class group
\(\mathrm{Cl}_3(k_\mu)\simeq (3,3)\), \(\mu=2,3,4\),
of rank \(2\), whereas the remaining component
has \(3\)-class rank \(\varrho_3(k_1)=3\).
Thus, the \textbf{tame} condition
\(9\mid h_3(B_j)=(U_j:V_j)\in\lbrace 9,27\rbrace\), \(r_j=2\),
is satisfied
for the bicyclic bicubic fields \(B_j\) with \(j\in\lbrace 2,3,4,8,9,10\rbrace\).
\end{lemma}

\begin{proof}
\(p\) is universally repelling
\(\lbrace q\leftarrow p\rightarrow r\rbrace\).
Since \(p\rightarrow r\), \(p\) splits in \(k_r\).
Since \(q\leftarrow p\), \(p\) splits in \(k_q\).
Thus \(p\) also splits in \(k_{qr}\) and  \(\tilde{k}_{qr}\).
By the normalizing assumptions that
\(q\) splits in \(k_{pr}\) and \(r\) splits in \(k_{pq}\),
the primes \(p,q,r\) share the common decomposition type
\((e,f,g)=(3,1,3)\) in the bicyclic bicubic field
\(B_1=k_1k_{pq}k_{pr}k_{qr}\),
which implies that \(\varrho_3(k_1)=3\),
according to
\cite[Prop. 4.4, pp. 43--44]{Ay1995}.
Finally, none among
\(B_2=k_2k_{pq}\tilde{k}_{pr}\tilde{k}_{qr}\),
\(B_3=k_3\tilde{k}_{pq}\tilde{k}_{pr}k_{qr}\),
\(B_4=k_4\tilde{k}_{pq}k_{pr}\tilde{k}_{qr}\),
\(B_8=k_2k_4k_pk_{qr}\),
\(B_9=k_2k_3k_qk_{pr}\),
\(B_{10}=k_3k_4k_rk_{pq}\)
contains \(k_1\).
\end{proof}


\begin{proposition}
\label{prp:Cat1Gph2}
\textbf{(Sub-triplet with \(3\)-rank two for \(\mathrm{I}.2\).)}
For fixed \(\mu\in\lbrace 2,3,4\rbrace\),
let \(\mathfrak{p},\mathfrak{q},\mathfrak{r}\) be the prime ideals of \(k_\mu\)
over \(p,q,r\), that is
\(p\mathcal{O}_{k_\mu}=\mathfrak{p}^3\),
\(q\mathcal{O}_{k_\mu}=\mathfrak{q}^3\),
\(r\mathcal{O}_{k_\mu}=\mathfrak{r}^3\),
then the \(3\)-class group of \(k_\mu\) is generated by
the non-trivial classes \(\lbrack\mathfrak{q}\rbrack,\lbrack\mathfrak{r}\rbrack\),
that is,
\begin{equation}
\label{eqn:GenCat1Gph2}
\mathrm{Cl}_3(k_\mu)=
\langle\lbrack\mathfrak{q}\rbrack,\lbrack\mathfrak{r}\rbrack\rangle\simeq (3,3).
\end{equation}

The unramified cyclic cubic relative extensions of \(k_\mu\)
are among the absolutely bicyclic bicubic fields \(B_i\), \(1\le i\le 10\).

In terms of \textbf{decisive principal factors}
\(A(k_1)=p^xq^yr^z\), \(x,y,z\in\lbrace 0,1,2\rbrace\), and
\(A(k_{qr})=qr^n\), \(n\in\lbrace 1,2\rbrace\),
the ranks of principal factor matrices of \textbf{wild} bicyclic bicubic fields are
\begin{equation}
\label{eqn:Cat1Gph2Ranks}
r_5=3 \text{ iff } -z\not\equiv ny\,(\mathrm{mod}\,3),\
r_6=3 \text{ iff } z\ne 0 \text{ iff } r\mid A(k_1),\
r_7=3 \text{ iff } y\ne 0 \text{ iff } q\mid A(k_1).
\end{equation}

The field \(B_2=k_2k_{pq}\tilde{k}_{pr}\tilde{k}_{qr}\)
has norm class group
\(N_{B_2/k_2}(\mathrm{Cl}_3(B_2))=\langle\lbrack\mathfrak{r}\rbrack\rangle\),
and transfer kernel
\[
\ker(T_{B_2/k_2})\ge\langle\lbrack\mathfrak{q}^2\mathfrak{r}^n\rbrack\rangle.
\]

The field \(B_4=k_4\tilde{k}_{pq}k_{pr}\tilde{k}_{qr}\)
has norm class group
\(N_{B_4/k_4}(\mathrm{Cl}_3(B_2))=\langle\lbrack\mathfrak{q}\rbrack\rangle\),
and transfer kernel
\[
\ker(T_{B_4/k_4})\ge\langle\lbrack\mathfrak{q}^2\mathfrak{r}^n\rbrack\rangle.
\]

The field \(B_9=k_2k_3k_qk_{pr}\), which contains \(k_2\) and \(k_3\),
has norm class group
\(N_{B_9/k_\mu}(\mathrm{Cl}_3(B_9))=\langle\lbrack\mathfrak{q}\rbrack\rangle\),
for \(\mu=2,3\), and possible \textbf{fixed point} transfer kernel
\begin{equation}
\label{eqn:Fix1Cat1Gph2}
\ker(T_{B_9/k_\mu})\ge\langle\lbrack\mathfrak{q}\rbrack\rangle.
\end{equation}

The field \(B_{10}=k_3k_4k_rk_{pq}\), which contains \(k_3\) and \(k_4\),
has norm class group
\(N_{B_{10}/k_\mu}(\mathrm{Cl}_3(B_{10}))=\langle\lbrack\mathfrak{r}\rbrack\rangle\),
for \(\mu=3,4\), and possible \textbf{fixed point} transfer kernel
\begin{equation}
\label{eqn:Fix2Cat1Gph2}
\ker(T_{B_{10}/k_\mu})\ge\langle\lbrack\mathfrak{r}\rbrack\rangle.
\end{equation}

The remaining two \(B_i>k_\mu\), \(i\in\lbrace 3,5,6,7,8\rbrace\),
more precisely, \(i\in\lbrace 7,8\rbrace\) for \(\mu=2\),
\(i\in\lbrace 3,5\rbrace\) for \(\mu=3\),
and \(i\in\lbrace 6,8\rbrace\) for \(\mu=4\),
have norm class group 
\(\langle\lbrack\mathfrak{qr}\rbrack\rangle\) respectively
\(\langle\lbrack\mathfrak{qr}^2\rbrack\rangle\). Among them,
the tame extensions \(B_i>k_\mu\) with
either \(i=\mu=3\) or \(i=8\), \(\mu=2,4\),
have \textbf{partial} transfer kernel
\begin{equation}
\label{eqn:PrtCat1Gph2}
\ker(T_{B_i/k_\mu})=\langle\lbrack\mathfrak{qr}\rbrack\rangle
\end{equation}
of order \(3\),
giving rise to either a \textbf{transposition}
or a \textbf{fixed point}.
\end{proposition}

\begin{proof}
Since \(p\rightarrow r\),
two principal factors are \(A(k_{pr})=A(\tilde{k}_{pr})=p\); and
since \(q\leftarrow p\),
two further principal factors are \(A(k_{pq})=A(\tilde{k}_{pq})=p\),
by Proposition
\ref{prp:Principal2}.
Since \(p\) is universally repelling
\(\lbrace q\leftarrow p\rightarrow r\rbrace\),
three further principal factors are \(A(k_\mu)=p\) for \(2\le\mu\le 4\),
by Proposition
\ref{prp:Principal3}.

Thus, \(\lbrack\mathfrak{p}\rbrack=1\) is trivial, and
the non-trivial classes \(\lbrack\mathfrak{q}\rbrack,\lbrack\mathfrak{r}\rbrack\)
generate \(\mathrm{Cl}_3(k_\mu)=
\langle\lbrack\mathfrak{q}\rbrack,\lbrack\mathfrak{r}\rbrack\rangle\simeq (3,3)\).

Since \(q\) splits in \(k_{pr}\),
it also splits in 
\(B_4=k_4\tilde{k}_{pq}k_{pr}\tilde{k}_{qr}\),
\(B_9=k_2k_3k_qk_{pr}\).

Since \(r\) splits in \(k_{pq}\),
it also splits in
\(B_2=k_2k_{pq}\tilde{k}_{pr}\tilde{k}_{qr}\),
\(B_{10}=k_3k_4k_rk_{pq}\).

Since the \textit{tame} condition \(9\mid h_3(B_j)=(U_j:V_j)\) is satisfied
for \(j\in\lbrace 2,3,4,8,9,10\rbrace\),
the rank of the corresponding principal factor matrix \(M_j\)
must be \(r_j=2\).
This can also be verified directly and
has no further consequences.

We propose the principal factors \(A(k_1)=p^{x}q^{y}r^{z}\) and
and \(A(k_{qr})=qr^n\), \(A(\tilde{k}_{qr})=q^2r^n\)
with \(n\in\lbrace 1,2\rbrace\).
For each \textit{wild} bicyclic bicubic field \(B_j\),
\(j\in\lbrace 5,6,7\rbrace\),
the rank \(r_j\) is now calculated with row operations
on the associated principal factor matrices \(M_j\): \\
\(M_5=
\begin{pmatrix}
x & y & z \\
1 & 0 & 0 \\
1 & 0 & 0 \\
0 & 2 & n \\
\end{pmatrix}\),
\(M_6=
\begin{pmatrix}
x & y & z \\
1 & 0 & 0 \\
0 & 1 & 0 \\
1 & 0 & 0 \\
\end{pmatrix}\),
\(M_7=
\begin{pmatrix}
x & y & z \\
1 & 0 & 0 \\
0 & 0 & 1 \\
1 & 0 & 0 \\
\end{pmatrix}\).

For \(B_5=k_1k_3k_p\tilde{k}_{qr}\),
\(M_5\)
leads to the decisive pivot element
\(z+ny\)
in the last column.
So, \(r_5=3\) iff \(-z\not\equiv ny\) modulo \(3\).
For \(B_6=k_1k_4k_q\tilde{k}_{pr}\), 
\(M_6\)
leads to
\(z\).
So, \(r_6=3\) iff \(z\ne 0\).
For \(B_7=k_1k_2k_r\tilde{k}_{pq}\), 
\(M_7\)
leads to \(y\) in the middle column.
So, \(r_7=3\) iff \(y\ne 0\).

Since \(\mathfrak{q}\) is principal in \(k_q\),
\(\lbrack\mathfrak{q}\rbrack\) capitulates in
\(B_6=k_1k_4k_q\tilde{k}_{pr}\) and
\(B_9=k_2k_3k_qk_{pr}\), and \\
since \(\mathfrak{r}\) is principal in \(k_r\), 
\(\lbrack\mathfrak{r}\rbrack\) capitulates in
\(B_7=k_1k_2k_r\tilde{k}_{pq}\) and
\(B_{10}=k_3k_4k_rk_{pq}\),
by Corollary
\ref{cor:Capitulation}.

Since \(\mathfrak{qr}^n\) is principal in \(k_{qr}\),
\(\lbrack\mathfrak{qr}^n\rbrack\) capitulates in
\(B_3=k_3\tilde{k}_{pq}\tilde{k}_{pr}k_{qr}\),
\(B_8=k_2k_4k_pk_{qr}\), and \\
since \(\mathfrak{q}^2\mathfrak{r}^n\) is principal in \(\tilde{k}_{qr}\),
\(\lbrack\mathfrak{q}^2\mathfrak{r}^n\rbrack\) capitulates in
\(B_2=k_2k_{pq}\tilde{k}_{pr}\tilde{k}_{qr}\),
\(B_4=k_4\tilde{k}_{pq}k_{pr}\tilde{k}_{qr}\), and
\(B_5=k_1k_3k_p\tilde{k}_{qr}\),
by Proposition
\ref{prp:Bicyc}.

In each case,
the minimal subfield unit index
\((U_j:V_j)=3\) for \(r_j=3\)
corresponds to
the maximal unit norm index \((U(k_\mu):N_{B_j/k_\mu}(U_j))=3\),
associated to a \textit{total} transfer kernel
\(\#\ker(T_{B_j/k_\mu})=9\),
by Lemma
\ref{lem:UnitIndices}.

The minimal unit norm index \((U(k_\mu):N_{B_8/k_\mu}(U_8))=1\),
associated to the partial transfer kernel
\(\ker(T_{B_8/k_\mu})=\langle\lbrack\mathfrak{qr}\rbrack\rangle\),
for \(\mu=2,4\),
corresponds to the tame maximal subfield unit index
\(h_3(B_8)=(U_8:V_8)=27\),
giving rise to type invariants
\(\mathrm{Cl}_3(B_8)\simeq (9,3)\).
\end{proof}

\noindent
Using Corollary
\ref{cor:Three},
Proposition
\ref{prp:Cat1Gph2}
and parts of its proof are now summarized in Table
\ref{tbl:UniCat1Gph2}
for the minimal transfer kernel type (mTKT) and \(n=1\),
with transposition in \textbf{boldface} font.


\renewcommand{\arraystretch}{1.1}

\begin{table}[ht]
\caption{Norm class groups and minimal transfer kernels with \(n=1\) for Graph I.2}
\label{tbl:UniCat1Gph2}
\begin{center}
{\normalsize
\begin{tabular}{|c||c|c|c|c||c|c|c|c||c|c|c|c|}
\hline
 Base         & \multicolumn{4}{c||}{\(k_2\)} & \multicolumn{4}{c||}{\(k_3\)} & \multicolumn{4}{c|}{\(k_4\)} \\
\hline
 Ext          & \(B_2\) & \(B_7\) & \(B_8\) & \(B_9\) & \(B_3\) & \(B_5\) & \(B_9\) & \(B_{10}\) & \(B_4\) & \(B_6\) & \(B_8\) & \(B_{10}\) \\
\hline
 NCG          & \(\mathfrak{r}\)    & \(\mathfrak{qr}^2\) & \(\mathfrak{qr}\) & \(\mathfrak{q}\)    
              & \(\mathfrak{qr}^2\) & \(\mathfrak{qr}\)   & \(\mathfrak{q}\)  & \(\mathfrak{r}\) & \(\mathfrak{q}\)    & \(\mathfrak{qr}^2\) & \(\mathfrak{qr}\) & \(\mathfrak{r}\) \\
 TK           & \(\mathfrak{qr}^2\) & \(\mathfrak{r}\)    & \(\mathfrak{qr}\) & \(\mathfrak{q}\) 
              & \(\mathfrak{qr}\)   & \(\mathfrak{qr}^2\) & \(\mathfrak{q}\)  & \(\mathfrak{r}\) & \(\mathfrak{qr}^2\) & \(\mathfrak{q}\)    & \(\mathfrak{qr}\) & \(\mathfrak{r}\) \\
\(\varkappa\) &  \(\mathbf{2}\) & \(\mathbf{1}\) & \(3\) & \(4\) &  \(\mathbf{2}\) & \(\mathbf{1}\) & \(3\) & \(4\) & \(\mathbf{2}\) & \(\mathbf{1}\) & \(3\) & \(4\) \\
\hline
\end{tabular}
}
\end{center}
\end{table}


\begin{theorem}
\label{thm:Cat1Gph2}
\textbf{(Second \(3\)-class group for \(\mathrm{I}.2\).)}
To identify the second \(3\)-class group
\(\mathfrak{M}=\mathrm{Gal}(\mathrm{F}_3^2(k_\mu)/k_\mu)\), \(2\le\mu\le 4\),
let the \textbf{principal factors} of
\(k_1\), and \(k_{qr}\), respectively \(\tilde{k}_{qr}\), be
\(A(k_1)=p^xq^yr^z\), \(x,y,z\in\lbrace 0,1,2\rbrace\), and
\(A(k_{qr})=qr^n\), respectively \(A(\tilde{k}_{qr})=q^2r^n\),
\(n\in\lbrace 1,2\rbrace\),
and additionally assume the \textbf{regular} situation where
\(\mathrm{Cl}_3(k_1)\simeq (3,3,3)\).

Then the
\textbf{minimal transfer kernel type} (mTKT) \(\varkappa_0\) of \(k_\mu\), \(1\le\mu\le 4\),
and other possible capitulation types in ascending partial order
\(\varkappa_0<\varkappa<\varkappa^\prime,\varkappa^{\prime\prime}\),
ending in two non-comparable types,
are \(\varkappa_0=(2134)\), type \(\mathrm{G}.16\),
\(\varkappa=(0134)\), type \(\mathrm{c}.21\),
\(\varkappa^\prime=(0004)\), type \(\mathrm{a}.2\),
\(\varkappa^{\prime\prime}=(0100)\), type \(\mathrm{a}.3\),
and the second \(3\)-class group is \(\mathfrak{M}\simeq\)

\begin{equation}
\label{eqn:Cat1Gph2}
\begin{cases}
\langle 81,8\rangle,\ \alpha=\lbrack 11,21,11,11\rbrack,\ \varkappa=(0100) & \text{ once if } y\ne 0,\ z\ne 0,\ \mathcal{N}=1, \\
\langle 81,10\rangle,\ \alpha=\lbrack 11,11,11,21\rbrack,\ \varkappa=(0004) & \text{ twice if } y\ne 0,\ z\ne 0,\ \mathcal{N}=1, \\
\langle 243,8\rangle,\ \alpha=\lbrack 21,21,21,21\rbrack,\ \varkappa=(0134) & \text{ if } y=z=0,\ \mathcal{N}=3, \\
\langle 729,52\rangle,\ \alpha=\lbrack 22,21,21,21\rbrack,\ \varkappa=(2134) & \text{ if } y=z=0,\ \mathcal{N}=4, \\
\end{cases}
\end{equation}
where \(\mathcal{N}:=\#\lbrace 1\le j\le 10\mid k_\mu<B_j,\ I_j=27\rbrace\).
Only for the leading three rows, the \(3\)-class field tower has certainly the group
\(\mathfrak{G}=\mathrm{Gal}(\mathrm{F}_3^\infty(k_\mu)/k_\mu)\simeq\mathfrak{M}\)
and length \(\ell_3(k_\mu)=2\),
otherwise length \(\ell_3(k_\mu)\ge 3\) cannot be excluded
although \(d_2(\mathfrak{M})\le 4\).
\end{theorem}

\begin{proof}
Let \(\mu\in\lbrace 2,3,4\rbrace\).

The first scenario,
\(y\ne 0\), \(z\ne 0\), and
\(-z\not\equiv ny\) modulo \(3\)
is equivalent to
\(\mathcal{N}=1\),
since \(r_j=3\), \((U_j:V_j)=3\),
\(h_3(B_j)=\frac{1}{3}h_3(k_1)=9\),
for the wild \(j=5,6,7\), and
\(h_3(B_j)=(U_j:V_j)=9\), for the tame \(j=2,4,9,10\),
whereas the distinguished tame \(j=3,8\) have
\(h_3(B_j)=(U_j:V_j)=27\).
This gives rise to Artin pattern
either \(\alpha=\lbrack 11,11,11,21\rbrack\)
\(\varkappa=(0004)\), for \(j=8\), \(\mu=2,4\)
(twice with fixed point),
characteristic for \(\langle 81,10\rangle\),
or \(\alpha=\lbrack 11,21,11,11\rbrack\)
\(\varkappa=(0100)\), for \(j=\mu=3\)
(only once with non-fixed point,
due to a hidden transposition),
characteristic for \(\langle 81,8\rangle\).

The other two scenarios share
\(y=z=0\), and thus also \(-z=ny\),
independently of \(n\),
which implies \(r_j=2\), \((U_j:V_j)\in\lbrace 9,27\rbrace\),
for \(j=5,6,7\), and
\(h_3(B_j)=(U_j:V_j)=27\), for the tame \(j=2,4,9,10\),
producing two fixed points at \(B_9\) and \(B_{10}\).

The second scenario with \(\mathcal{N}=3\)
is supplemented by \((U_j:V_j)=9\),
\(h_3(B_j)=h_3(k_1)=27\), for \(j=5,6,7\), and
\textit{total} capitulation,
\(\#\ker(T_{B_j/k_\mu})=9\),
for \(\mu=3,4,2\).
This gives rise to \(\alpha=\lbrack 21,21,21,21\rbrack\),
\(\varkappa=(0134)\),
characteristic for \(\langle 243,8\rangle\)
with annihilator ideal \(\mathfrak{L}\)
in the sense of Scholz and Taussky
\cite{SoTa1934}.

The third scenario with \(\mathcal{N}=4\)
is supplemented by \((U_j:V_j)=27\),
\(h_3(B_j)=3h_3(k_1)=81\), for \(j=5,6,7\) and
\textit{partial} non-fixed point capitulation.
This gives rise to \(\alpha=\lbrack 22,21,21,21\rbrack\),
\(\varkappa=(2134)\),
characteristic for \(\langle 729,52\rangle\)
with non-metabelian descendants.
Here, the hidden \textit{transposition} becomes explicit,
between either \(B_2\), \(B_7\)
or \(B_3\), \(B_5\)
or \(B_4\), \(B_6\).
\end{proof}


\begin{corollary}
\label{cor:UniCat1Gph2}
\textbf{(Non-uniformity of the sub-triplet for \(\mathrm{I}.2\).)}
The components of the sub-triplet with \(3\)-rank two
share a common capitulation type
\(\varkappa(k_\mu)\),
common abelian type invariants
\(\alpha(k_\mu)\),
and a common second \(3\)-class group
\(\mathrm{Gal}(\mathrm{F}_3^2(k_\mu)/k_\mu)\),
for \(\mu=2,3,4\),
only if \(y=z=0\), \(\mathcal{N}=3,4\).
For \(y\ne 0\), \(z\ne 0\), \(\mathcal{N}=1\), however,
only two fields \(k_2\) and \(k_4\) share common invariants,
whereas \(k_3\) has different \(\varkappa(k_3)\)
and different \(\mathrm{Gal}(\mathrm{F}_3^2(k_3)/k_3)\).
\end{corollary}

\begin{proof}
This follows immediately from Theorem
\ref{thm:Cat1Gph2},
whereas Table
\ref{tbl:UniCat1Gph2}
with minimal transfer kernel type \(\varkappa_0=(2134)\)
only shows the uniform situation,
which can become non-uniform by superposition
with total transfer kernels,
when \(\mathcal{N}=1\).
\end{proof}


\begin{example}
\label{exm:Cat1Gph2}
Prototypes for Graph I.2, i.e.,
minimal conductors for each scenario in Theorem
\ref{thm:Cat1Gph2}
have been found for each \(\mathcal{N}\in\lbrace 1,3,4\rbrace\).

Some are \textbf{regular}:
\(c=8\,001\) with symbol
\(\lbrace 9\leftarrow 127\rightarrow 7\rbrace\)
and, non-uniformly, once \(\mathfrak{G}=\mathfrak{M}=\langle 81,8\rangle\)
but twice \(\langle 81,10\rangle^2\);
\(c=21\,049\) with symbol
\(\lbrace 7\leftarrow 97\rightarrow 31\rbrace\)
and uniformly three times \(\mathfrak{G}=\mathfrak{M}=\langle 243,8\rangle^3\);
and \(c=59\,031\) with symbol
\(\lbrace 9\leftarrow 937\rightarrow 7\rbrace\)
and \(\mathfrak{M}=\langle 729,52\rangle^3\).

Others are \textbf{singular}:
\(c=7\,657\) with symbol
\(\lbrace 13\leftarrow 31\rightarrow 19\rbrace\)
and \(\mathfrak{G}=\mathfrak{M}=\langle 243,8\rangle^3\);
and \(c=48\,393\) with symbol
\(\lbrace 9\leftarrow 19\rightarrow 283\rbrace\)
and \(\mathfrak{M}=\langle 2187,301\vert 305\rangle^3\).

The groups of order \(\ge 729\)
with transfer kernel type \(\mathrm{G}.16\)
have non-metabelian extensions.
\end{example}


In Table
\ref{tbl:ProtoCat1Gph2},
we summarize the prototypes of Graph \(\mathrm{I}.2\)
in the same way as in Table
\ref{tbl:ProtoCat1Gph1},
except that two
critical exponents \(y,z\) in principal factor \(A(k_1)=p^xq^yr^z\)
and \(n\) in \(A(k_{qr})=qr^n\) are sufficient.

\renewcommand{\arraystretch}{1.1}

\begin{table}[ht]
\caption{Prototypes for Graph I.2}
\label{tbl:ProtoCat1Gph2}
\begin{center}
{\scriptsize
\begin{tabular}{|c|c||c|c|c|c|c|c|}
\hline
 No. & \(c\)        & \(q\leftarrow p\rightarrow r\)   & \(v\) & \(y,z;n\) & capitulation type & \(\mathfrak{M}\)                   & \(\ell_3(k)\) \\
\hline
   1 &  \(7\,657\) & \(13\leftarrow 31\rightarrow 19\) & \(4\) & \(1,1;2\) & \(\mathrm{c}.21\) & \(\langle 243,8\rangle^3\) & \(=2\) \\
   2 &  \(8\,001\) & \(9\leftarrow 127\rightarrow 7\)  & \(3\) & \(1,2;1\) & \(\mathrm{a}.3,\mathrm{a}.2\) & \(\langle 81,8\rangle,\langle 81,10\rangle^2\) & \(=2\) \\
  12 & \(21\,049\) & \(7\leftarrow 97\rightarrow 31\)  & \(3\) & \(0,0;1\) & \(\mathrm{c}.21\) & \(\langle 243,8\rangle^3\) & \(=2\) \\
  27 & \(48\,393\) & \(9\leftarrow 19\rightarrow 283\) & \(4\) & \(0,0;2\) & \(\mathrm{G}.16\) & \(\langle 2187,301\vert 305\rangle^3\) & \(\ge 2\) \\
  33 & \(59\,031\) & \(9\leftarrow 937\rightarrow 7\)  & \(3\) & \(0,0;1\) & \(\mathrm{G}.16\) & \(\langle 729,52\rangle^3\) & \(\ge 2\) \\
\hline
\end{tabular}
}
\end{center}
\end{table}


\subsection{Category II, Graph 1}
\label{ss:Cat2Gph1}

\noindent
Let \((k_1,\ldots,k_4)\) be a quartet of cyclic cubic number fields
sharing the common conductor \(c=pqr\),
belonging to Graph \(1\) of Category \(\mathrm{II}\)
with combined cubic residue symbol
\(\lbrack p,q,r\rbrack_3=\lbrace p\rightarrow q\leftarrow r\rbrace\).


\begin{lemma}
\label{lem:Cat2Gph1}
\textbf{(3-class ranks of components for \(\mathrm{II}.1\).)}
Under the normalizing assumption that \(q\) splits in \(\tilde{k}_{pr}\),
precisely the two components \(k_2\) and \(k_3\) of the quartet
have elementary bicyclic \(3\)-class group
\(\mathrm{Cl}_3(k_2)\simeq \mathrm{Cl}_3(k_3)\simeq (3,3)\)
of rank \(2\), whereas the other two components
have \(3\)-class rank \(\varrho_3(k_1)=\varrho_3(k_4)=3\).
Thus, the \textbf{tame} condition
\(9\mid h_3(B_j)=(U_j:V_j)\in\lbrace 9,27\rbrace\), \(r_j=2\),
is only satisfied
for the bicyclic bicubic fields \(B_j\) with \(j\in\lbrace 2,3,9\rbrace\).
\end{lemma}

\begin{proof}
Since \(p\rightarrow q\), \(p\) splits in \(k_q\).
Since \(q\leftarrow r\), \(r\) splits in \(k_q\),
and also splits in \(B_9=k_2k_3k_qk_{pr}\).
By the normalizing assumption that \(q\) splits in \(\tilde{k}_{pr}\),
it also splits in
\(B_2=k_2k_{pq}\tilde{k}_{pr}\tilde{k}_{qr}\) and
\(B_3=k_3\tilde{k}_{pq}\tilde{k}_{pr}k_{qr}\).
The primes \(p,q,r\) share the common decomposition type
\((e,f,g)=(3,1,3)\) in the bicyclic bicubic field
\(B_6=k_1k_4k_q\tilde{k}_{pr}\),
which implies that \(\varrho_3(k_1)=\varrho_3(k_4)=3\),
according to
\cite[Prop. 4.4, pp. 43--44]{Ay1995}.
Finally, only
\(B_2=k_2k_{pq}\tilde{k}_{pr}\tilde{k}_{qr}\),
\(B_3=k_3\tilde{k}_{pq}\tilde{k}_{pr}k_{qr}\),
\(B_9=k_2k_3k_qk_{pr}\)
do not contain \(k_1\), \(k_4\).
\end{proof}


\begin{proposition}
\label{prp:Cat2Gph1}
\textbf{(Sub-doublet with \(3\)-rank two for \(\mathrm{II}.1\).)}
For fixed \(\mu\in\lbrace 2,3\rbrace\),
let \(\mathfrak{p},\mathfrak{q},\mathfrak{r}\) be the prime ideals of \(k_\mu\)
over \(p,q,r\), that is
\(p\mathcal{O}_{k_\mu}=\mathfrak{p}^3\),
\(q\mathcal{O}_{k_\mu}=\mathfrak{q}^3\),
\(r\mathcal{O}_{k_\mu}=\mathfrak{r}^3\),
then the \(3\)-class group of \(k_\mu\) is generated by
the non-trivial classes \(\lbrack\mathfrak{q}\rbrack,\lbrack\mathfrak{r}\rbrack\),
that is,
\begin{equation}
\label{eqn:GenCat2Gph1}
\mathrm{Cl}_3(k_\mu)=
\langle\lbrack\mathfrak{q}\rbrack,\lbrack\mathfrak{r}\rbrack\rangle\simeq (3,3).
\end{equation}

The unramified cyclic cubic relative extensions of \(k_\mu\)
are among the absolutely bicyclic bicubic fields \(B_i\), \(1\le i\le 10\).
The unique \(B_\mu\), \(\mu\in\lbrace 2,3\rbrace\), which only contains \(k_\mu\),
has norm class group
\(N_{B_\mu/k_\mu}(\mathrm{Cl}_3(B_\mu))=\langle\lbrack\mathfrak{q}\rbrack\rangle\),
transfer kernel
\[
\ker(T_{B_\mu/k_\mu})\ge\langle\lbrack\mathfrak{r}\rbrack\rangle
\]
and \(3\)-class group
\(\mathrm{Cl}_3(B_\mu)=
\langle\lbrack\mathfrak{Q}_1\rbrack,\lbrack\mathfrak{Q}_2\rbrack,\lbrack\mathfrak{Q}_3\rbrack\rangle\ge (3,3)\),
generated by the classes of the prime ideals of \(B_\mu\) over
\(\mathfrak{q}\mathcal{O}_{B_\mu}=\mathfrak{Q}_1\mathfrak{Q}_2\mathfrak{Q}_3\).
The unique \(B_9=k_2k_3k_qk_{pr}\), which contains \(k_2\) and \(k_3\),
has norm class group
\(N_{B_9/k_\mu}(\mathrm{Cl}_3(B_9))=\langle\lbrack\mathfrak{r}\rbrack\rangle\),
\textbf{cyclic} transfer kernel
\begin{equation}
\label{eqn:PrtCat2Gph1}
\ker(T_{B_9/k_\mu})=\langle\lbrack\mathfrak{q}\rbrack\rangle
\end{equation}
of order \(3\),
and \textbf{elementary tricyclic} \(3\)-class group
\(\mathrm{Cl}_3(B_9)=
\langle\lbrack\mathfrak{R}_1\rbrack,\lbrack\mathfrak{R}_2\rbrack,\lbrack\mathfrak{R}_3\rbrack\rangle\simeq (3,3,3)\),
generated by the classes of the prime ideals of \(B_9\) over
\(\mathfrak{r}\mathcal{O}_{B_9}=\mathfrak{R}_1\mathfrak{R}_2\mathfrak{R}_3\).
The remaining two \(B_i>k_\mu\), \(i\in\lbrace 5,7,8,10\rbrace\),
more precisely, \(i\in\lbrace 7,8\rbrace\) for \(\mu=2\),
and \(i\in\lbrace 5,10\rbrace\) for \(\mu=3\),
have norm class group 
\(\langle\lbrack\mathfrak{qr}\rbrack\rangle\) respectively
\(\langle\lbrack\mathfrak{qr}^2\rbrack\rangle\),
and transfer kernel
\[
\ker(T_{B_i/k_\mu})\ge\langle\lbrack\mathfrak{r}\rbrack\rangle.
\]
In terms of \textbf{decisive principal factors}
\(A(k_\nu)=p^{x_\nu}q^{y_\nu}r^{z_\nu}\) for \(\nu\in\lbrace 1,4\rbrace\),
the ranks of principal factor matrices
\(M_i\), \(i\in\lbrace 1,4,5,7,8,10\rbrace\),
of \textbf{wild} bicyclic bicubic fields are
\begin{equation}
\label{eqn:Cat2Gph1Ranks}
r_1=r_5=r_7=3 \text{ iff } y_1\ne 0 \text{ iff } q\mid A(k_1)
\text{ and }
r_4=r_8=r_{10}=3 \text{ iff } y_4\ne 0 \text{ iff } q\mid A(k_4).
\end{equation}
\end{proposition}

\begin{proof}
Since \(p\rightarrow q\),
two principal factors are \(A(k_{pq})=A(\tilde{k}_{pq})=p\);
since \(q\leftarrow r\),
two further principal factors are \(A(k_{qr})=A(\tilde{k}_{qr})=r\);
both by Proposition
\ref{prp:Principal2}.

Since the \textit{tame} condition \(9\mid h_3(B_j)=(U_j:V_j)\) is satisfied
for \(j\in\lbrace 2,3,9\rbrace\),
the rank of the corresponding principal factor matrix \(M_j\)
must be \(r_2=r_3=r_9=2\).
We propose principal factors \(A(k_\mu)=p^{x_\mu}q^{y_\mu}r^{z_\mu}\),
for all \(1\le\mu\le 4\),
and \(A(k_{pr})=pr^\ell\), \(A(\tilde{k}_{pr})=p^2r^\ell\)
with \(\ell\in\lbrace 1,2\rbrace\).

For each bicyclic bicubic field \(B_j\),
the rank \(r_j\) is calculated with row operations
on the associated principal factor matrices \(M_j\): \\
\(M_2=
\begin{pmatrix}
x_2 & y_2 & z_2 \\
1 & 0 & 0 \\
2 & 0 & \ell \\
0 & 0 & 1 \\
\end{pmatrix}\),
\(M_3=
\begin{pmatrix}
x_3 & y_3 & z_3 \\
1 & 0 & 0 \\
2 & 0 & \ell \\
0 & 0 & 1 \\
\end{pmatrix}\),
\(M_9=
\begin{pmatrix}
x_2 & y_2 & z_2 \\
x_3 & y_3 & z_3 \\
0 & 1 & 0 \\
1 & 0 & \ell \\
\end{pmatrix}\).

For \(B_2=k_2k_{pq}\tilde{k}_{pr}\tilde{k}_{qr}\),
\(M_2\)
leads to the decisive pivot element
\(y_2\)
in the middle column,
and similarly,
for \(B_3=k_3\tilde{k}_{pq}\tilde{k}_{pr}k_{qr}\), 
\(M_3\)
leads to
\(y_3\).
So, \(r_2=r_3=2\) enforces \(y_2=y_3=0\),
i.e., \(q\nmid A(k_2)\), \(q\nmid A(k_3)\).
However,
for \(B_9=k_2k_3k_qk_{pr}\), 
\(M_9\)
leads to \(z_2-\ell x_2\) and \(z_3-\ell x_3\).
So, \(r_9=2\) enforces
\(z_2\equiv\ell x_2\) and \(z_3\equiv\ell x_3\) modulo \(3\),
i.e., \(A(k_2)=A(k_3)=pr^\ell\).

For every \textit{wild} bicyclic bicubic field \(B_j\),
\(j\in\lbrace 1,4,5,6,7,8,10\rbrace\),
the rank \(r_j\) is calculated by row operations
on the matrices \(M_j\), using \(A(k_2)=A(k_3)=pr^\ell\): \\
\(M_1=
\begin{pmatrix}
x_1 & y_1 & z_1 \\
1 & 0 & 0 \\
1 & 0 & \ell \\
0 & 0 & 1 \\
\end{pmatrix}\),
\(M_5=
\begin{pmatrix}
x_1 & y_1 & z_1 \\
1 & 0 & \ell \\
1 & 0 & 0 \\
0 & 0 & 1 \\
\end{pmatrix}\),
\(M_7=
\begin{pmatrix}
x_1 & y_1 & z_1 \\
1 & 0 & \ell \\
0 & 0 & 1 \\
1 & 0 & 0 \\
\end{pmatrix}\).

For \(B_1=k_1k_{pq}k_{pr}k_{qr}\),
\(M_1\)
leads to the decisive pivot element
\(y_1\)
in the middle column,
similarly,
for \(B_5=k_1k_3k_p\tilde{k}_{qr}\),
\(M_5\)
leads to \(y_1\),
and similarly,
for \(B_7=k_1k_2k_r\tilde{k}_{pq}\),
\(M_7\)
leads to \(y_1\).
So, \(r_1=r_5=r_7=3\)
iff \(y_1\ne 0\)
iff \(q\mid A(k_1)\).
Next we consider: \\
\(M_4=
\begin{pmatrix}
x_4 & y_4 & z_4 \\
1 & 0 & 0 \\
1 & 0 & \ell \\
0 & 0 & 1 \\
\end{pmatrix}\),
\(M_8=
\begin{pmatrix}
1 & 0 & \ell \\
x_4 & y_4 & z_4 \\
1 & 0 & 0 \\
0 & 0 & 1 \\
\end{pmatrix}\),
\(M_{10}=
\begin{pmatrix}
1 & 0 & \ell \\
x_4 & y_4 & z_4 \\
0 & 0 & 1 \\
1 & 0 & 0 \\
\end{pmatrix}\).

For \(B_4=k_4\tilde{k}_{pq}k_{pr}\tilde{k}_{qr}\),
\(M_4\)
leads to the decisive pivot element
\(y_4\)
in the middle column,
similarly,
for \(B_8=k_2k_4k_pk_{qr}\),
\(M_8\)
leads to \(y_4\),
and similarly,
for \(B_{10}=k_3k_4k_rk_{pq}\),
\(M_{10}\)
leads to \(y_4\).
So, \(r_4=r_8=r_{10}=3\)
iff \(y_4\ne 0\)
iff \(q\mid A(k_4)\).


By Lemma
\ref{lem:UnitIndices},
the minimal subfield unit index
\((U_j:V_j)=3\) for \(r_j=3\)
corresponds to
the maximal unit norm index \((U(k_\mu):N_{B_j/k_\mu}(U_j))=3\),
associated to a \textit{total} transfer kernel
\(\#\ker(T_{B_j/k_\mu})=9\).

Since \(q\) splits in \(\tilde{k}_{pr}\),
it also splits in 
\(B_2=k_2k_{pq}\tilde{k}_{pr}\tilde{k}_{qr}\),
\(B_3=k_3\tilde{k}_{pq}\tilde{k}_{pr}k_{qr}\),
\(\mathfrak{q}\mathcal{O}_{B_\mu}=\mathfrak{Q}_1\mathfrak{Q}_2\mathfrak{Q}_3\).

Since \(r\) splits in \(k_q\),
it also splits in \(B_9=k_2k_3k_qk_{pr}\),
\(\mathfrak{r}\mathcal{O}_{B_9}=\mathfrak{R}_1\mathfrak{R}_2\mathfrak{R}_3\). \\.

Since \(\mathfrak{r}\) is principal in \(k_r\), \(k_{qr}\), \(\tilde{k}_{qr}\),
\(\lbrack\mathfrak{r}\rbrack\) capitulates in
\(B_2=k_2k_{pq}\tilde{k}_{pr}\tilde{k}_{qr}\),
\(B_3=k_3\tilde{k}_{pq}\tilde{k}_{pr}k_{qr}\),
\(B_5=k_1k_3k_p\tilde{k}_{qr}\),
\(B_7=k_1k_2k_r\tilde{k}_{pq}\),
\(B_8=k_2k_4k_pk_{qr}\),
\(B_{10}=k_3k_4k_rk_{pq}\);
since \(\mathfrak{q}\) is principal in \(k_q\),
\(\lbrack\mathfrak{q}\rbrack\) capitulates in
\(B_9=k_2k_3k_qk_{pr}\)
(Proposition
\ref{prp:Bicyc}).
This gives a \textit{transposition}, either \((2,9)\) or \((3,9)\).

The minimal unit norm index \((U(k_\mu):N_{B_9/k_\mu}(U_9))=1\),
associated to the partial transfer kernel
\(\ker(T_{B_9/k_\mu})=\langle\lbrack\mathfrak{q}\rbrack\rangle\),
corresponds to the maximal subfield unit index
\(h_3(B_9)=(U_9:V_9)=27\),
giving rise to the \textit{elementary tricyclic} type invariants
\(\mathrm{Cl}_3(B_9)=
\langle\lbrack\mathfrak{R}_1\rbrack,\lbrack\mathfrak{R}_2\rbrack,\lbrack\mathfrak{R}_3\rbrack\rangle\simeq (3,3,3)\).
\end{proof}

\noindent
Using Corollary
\ref{cor:Three},
Proposition
\ref{prp:Cat2Gph1}
and parts of its proof are now summarized in Table
\ref{tbl:UniCat2Gph1}
with transposition in \textbf{bold} font.


\renewcommand{\arraystretch}{1.1}

\begin{table}[ht]
\caption{Norm class groups and minimal transfer kernels for Graph II.1}
\label{tbl:UniCat2Gph1}
\begin{center}
\begin{tabular}{|c||c|c|c|c||c|c|c|c|}
\hline
 Base         & \multicolumn{4}{c||}{\(k_2\)} & \multicolumn{4}{c||}{\(k_3\)} \\
\hline
 Ext          & \(B_2\) & \(B_7\) & \(B_8\) & \(B_9\) & \(B_3\) & \(B_5\) & \(B_9\) & \(B_{10}\) \\
\hline
 NCG          & \(\mathfrak{q}\) & \(\mathfrak{qr}\) & \(\mathfrak{qr}^2\) & \(\mathfrak{r}\) & \(\mathfrak{q}\) & \(\mathfrak{qr}^2\) & \(\mathfrak{r}\) & \(\mathfrak{qr}\) \\
 TK           & \(\mathfrak{r}\) & \(\mathfrak{r}\)  & \(\mathfrak{r}\)    & \(\mathfrak{q}\) & \(\mathfrak{r}\)    & \(\mathfrak{r}\)  & \(\mathfrak{q}\) & \(\mathfrak{r}\) \\
\(\varkappa\) & \(\mathbf{4}\) & \(4\) & \(4\) & \(\mathbf{1}\) & \(\mathbf{3}\) & \(3\) & \(\mathbf{1}\) & \(3\) \\
\hline
\end{tabular}
\end{center}
\end{table}


\begin{theorem}
\label{thm:Cat2Gph1}
\textbf{(Second \(3\)-class group for \(\mathrm{II}.1\).)}
Let \((k_1,\ldots,k_4)\) be a quartet of cyclic cubic number fields
sharing the common conductor \(c=pqr\),
belonging to Graph \(1\) of Category \(\mathrm{II}\)
with combined cubic residue symbol
\(\lbrack p,q,r\rbrack_3=\lbrace p\rightarrow q\leftarrow r\rbrace\).
Without loss of generality, suppose that \(q\) splits in \(\tilde{k}_{pr}\), and thus
\(\mathrm{Cl}_3(k_2)\simeq\mathrm{Cl}_3(k_3)\simeq (3,3)\),
and \(\varrho_3(k_1)=\varrho_3(k_4)=3\).

Then the \textbf{minimal transfer kernel type} (mTKT) of \(k_\mu\), \(2\le\mu\le 3\),
is \(\varkappa_0=(2111)\), type \(\mathrm{H}.4\),
and the other possible capitulation types in ascending order
\(\varkappa_0<\varkappa^\prime<\varkappa^{\prime\prime}<\varkappa^{\prime\prime\prime}\)
are
\(\varkappa^\prime=(2110)\), type \(\mathrm{d}.19\),
\(\varkappa^{\prime\prime}=(2100)\), type \(\mathrm{b}.10\), and
\(\varkappa^{\prime\prime\prime}=(2000)\), type \(\mathrm{a}.3^\ast\).

To identify the second \(3\)-class group
\(\mathfrak{M}=\mathrm{Gal}(\mathrm{F}_3^2(k_\mu)/k_\mu)\), \(2\le\mu\le 3\),
let the \textbf{decisive principal factors} of \(k_\nu\), \(\nu\in\lbrace 1,4\rbrace\), be
\(A(k_\nu)=p^{x_\nu}q^{y_\nu}r^{z_\nu}\), and additionally assume
the \textbf{regular situation} where both
\(\mathrm{Cl}_3(k_1)\simeq\mathrm{Cl}_3(k_4)\simeq (3,3,3)\)
are elementary tricyclic. Then

\begin{equation}
\label{eqn:Cat2Gph1}
\mathfrak{M}\simeq
\begin{cases}
\langle 81,7\rangle,\ \alpha=\lbrack 111,11,11,11\rbrack,\ \varkappa=(2000) & \text{ if } y_1\ne 0,\ y_4\ne 0,\ \mathcal{N}=1, \\
\langle 729,34..39\rangle,\ \alpha=\lbrack 111,111,21,21\rbrack,\ \varkappa=(2100) & \text{ if } y_1=y_4=0,\ \mathcal{N}=2, \\
\langle 729,41\rangle,\ \alpha=\lbrack 111,111,22,21\rbrack,\ \varkappa=(2110) & \text{ if } y_1=y_4=0,\ \mathcal{N}=3, \\
\langle 6561,714..719\vert 738..743\rangle \text{ or } & \\
\langle 2187,65\vert 67\rangle,\ \alpha=\lbrack 111,111,22,22\rbrack,\ \varkappa=(2111) & \text{ if } y_1=y_4=0,\ \mathcal{N}=4,
\end{cases}
\end{equation}
where \(\mathcal{N}:=\#\lbrace 1\le j\le 10\mid k_\mu<B_j,\ I_j=27\rbrace\).
Only in the leading row, the \(3\)-class field tower has warranted group
\(\mathfrak{G}=\mathrm{Gal}(\mathrm{F}_3^\infty(k_\mu)/k_\mu)\simeq\mathfrak{M}\),
with length \(\ell_3(k_\mu)=2\).
Otherwise, although the relation rank \(d_2(\mathfrak{M})\le 4\) is always admissible,
tower length \(\ell_3(k_\mu)\ge 3\) cannot be excluded.
\end{theorem}

\begin{proof}
We give the proof for \(k_3\)
with unramified cyclic cubic extensions
\(B_3,B_5,B_9,B_{10}\), 
(The proof for \(k_2\)
with unramified cyclic cubic extensions
\(B_2,B_7,B_8,B_9\)
is similar.)
We know that the tame ranks are
\(r_2=r_3=r_9=2\),
and thus \(I_2,I_3,I_9\in\lbrace 9,27\rbrace\),
in particular, \(I_9=27\), whence certainly \(\mathcal{N}\ge 1\).
Further, the wild ranks are
\(r_1=r_5=r_7=3\)
iff \(y_1\ne 0\),
and
\(r_4=r_8=r_{10}=3\)
iff \(y_4\ne 0\).

In the \textbf{regular situation} where the
\(3\)-class groups of \(k_1\) and \(k_4\)
are elementary tricyclic,
tight bounds arise for the abelian quotient invariants \(\alpha\) of
the group \(\mathfrak{M}\):

The first scenario,
\(y_1\ne 0\), \(y_4\ne 0\),
is equivalent to
\(\mathcal{N}=1\),
\(h_3(B_5)=h_3(B_7)=\frac{1}{3}h_3(k_1)=9\),
\(h_3(B_8)=h_3(B_{10})=\frac{1}{3}h_3(k_4)=9\),
\(h_3(B_2)=I_2=h_3(B_3)=I_3=9\),
\(h_3(B_9)=I_9=27\),
that is \(\alpha=\lbrack 111,11,11,11\rbrack\)
and consequently \(\varkappa=(2000)\),
since \(\langle 81,7\rangle\) is unique
with this \(\alpha\).

The other three scenarios share
\(y_1=y_4=0\),
and an explicit transposition
between \(B_3\) and \(B_9\),
giving rise to \(\varkappa=(21\ast\ast)\),
and common \(h_3(B_3)=I_3=27\), \(\alpha=\lbrack 111,111,\ast,\ast\rbrack\).

The second scenario with \(\mathcal{N}=2\)
is supplemented by 
\(h_3(B_5)=h_3(k_1)=27\),
\(h_3(B_{10})=h_3(k_4)=27\),
giving rise to \(\alpha=\lbrack 111,111,21,21\rbrack\),
\(\varkappa=(2100)\),
characteristic for \(\langle 729,34..39\rangle\)
(Cor.
\ref{cor:ClassGroup33}).

The third scenario with \(\mathcal{N}=3\)
is supplemented by 
\(h_3(B_5)=3h_3(k_1)=81\),
\(h_3(B_{10})=h_3(k_4)=27\),
giving rise to \(\alpha=\lbrack 111,111,22,21\rbrack\),
\(\varkappa=(2110)\),
characteristic for \(\langle 729,41\rangle\).

The fourth scenario with \(\mathcal{N}=4\)
is supplemented by 
\(h_3(B_5)=3h_3(k_1)=81\),
\(h_3(B_{10})=3h_3(k_4)=81\),
giving rise to \(\alpha=\lbrack 111,111,22,22\rbrack\),
\(\varkappa=(2111)\),
characteristic for either \(\langle 2187,65\vert 67\rangle\) or
\(\langle 6561,714..719\vert 738..743\rangle\)
with coclass \(\mathrm{cc}=3\).
If \(d_2(\mathfrak{M})=5\),
then tower length must be \(\ell_3(k_\mu)\ge 3\).
For this minimal capitulation type H.4, \(\varkappa=(2111)\),
all transfer kernels are cyclic of order \(3\),
and the minimal unit norm indices
correspond to maximal subfield unit indices.
\end{proof}


\begin{corollary}
\label{cor:UniCat2Gph1}
\textbf{(Uniformity of the sub-doublet for \(\mathrm{II}.1\).)}
The components of the sub-doublet with \(3\)-rank two
share a common capitulation type
\(\varkappa(k_2)\sim\varkappa(k_3)\),
common abelian type invariants
\(\alpha(k_2)\sim\alpha(k_3)\),
and a common second \(3\)-class group
\(\mathrm{Gal}(\mathrm{F}_3^2(k_2)/k_2)\simeq\mathrm{Gal}(\mathrm{F}_3^2(k_3)/k_3)\).
\end{corollary}

\begin{proof}
This is an immediate consequence of Theorem
\ref{thm:Cat2Gph1}
and Table
\ref{tbl:UniCat2Gph1}.
\end{proof}


\begin{example}
\label{exm:Cat2Gph1}
Prototypes for Graph \(\mathrm{II}.1\), i.e.,
minimal conductors for each scenario in Theorem
\ref{thm:Cat2Gph1}
have been detected for all \(N\in\lbrace 1,2,3,4\rbrace\).

There are
\textbf{regular} cases:
\(c=3\,913\) with symbol
\(\lbrace 13\rightarrow 7\leftarrow 43\rbrace\)
and \(\mathfrak{G}=\mathfrak{M}=\langle 81,7\rangle\);
\(c=22\,581\) with symbol
\(\lbrace 9\rightarrow 193\leftarrow 13\rbrace\)
and \(\mathfrak{M}=\langle 729,41\rangle\);
\(c=25\,929\) with symbol
\(\lbrace 9\rightarrow 67\leftarrow 43\rbrace\)
and \(\mathfrak{M}=\langle 729,34..36\rangle\) (Corollary
\ref{cor:ClassGroup33});
\(c=74\,043\) with symbol
\(\lbrace 19\rightarrow 9\leftarrow 433\rbrace\)
and either \(\mathfrak{M}=\langle 2187,65\vert 67\rangle\)
with \(d_2(\mathfrak{M})=5\)
or \(\mathfrak{M}=\langle 6561,714..719\vert 738..743\rangle\)
with \(d_2(\mathfrak{M})=4\); and
\(c=82\,327\) with symbol
\(\lbrace 7\rightarrow 19\leftarrow 619\rbrace\)
and \(\mathfrak{M}=\langle 729,37..39\rangle\) (Corollary
\ref{cor:ClassGroup33}).

We also have
\textbf{singular} cases:
\(c=30\,457\) with symbol
\(\lbrace 7\rightarrow 19\leftarrow 229\rbrace\)
and \(\mathfrak{M}=\langle 729,37..39\rangle\) (Corollary
\ref{cor:ClassGroup33});
\(c=34\,029\) with symbol
\(\lbrace 19\rightarrow 9\leftarrow 199\rbrace\)
and \(\mathfrak{M}=\langle 2187,248\vert 249\rangle\);
\(c=41\,839\) with symbol
\(\lbrace 43\rightarrow 7\leftarrow 139\rbrace\)
and \(\mathfrak{M}=\langle 6561,693..698\rangle\).

Finally, there is the \textbf{super-singular} 
\(c=83\,817\) with symbol
\(\lbrace 9\rightarrow 67\leftarrow 139\rbrace\)
and \(\mathfrak{M}=\langle 6561,693..698\rangle\).

With exception of \(\langle 81,7\rangle\),
all groups have non-metabelian descendants,
respectively extensions.
\end{example}


In Table
\ref{tbl:ProtoCat2Gph1},
we summarize the prototypes of Graph \(\mathrm{II}.1\)
in the same manner as in Table
\ref{tbl:ProtoCat1Gph1},
except that
regularity, resp. (super-)singularity, is expressed by
\(3\)-valuations \(v_\nu=v_3(\#\mathrm{Cl}(k_\nu))\)
of class numbers of critical fields \(k_\nu\), \(\nu=1,4\),
and critical exponents are \(y_\nu\) in principal factors
\(A(k_\nu)=p^{x_\nu}q^{y_\nu}r^{z_\nu}\),
\(\nu=1,4\).

\renewcommand{\arraystretch}{1.1}

\begin{table}[ht]
\caption{Prototypes for Graph II.1}
\label{tbl:ProtoCat2Gph1}
\begin{center}
{\scriptsize
\begin{tabular}{|c|c||c|c|c|c|c|c|}
\hline
 No. & \(c\)        & \(p\rightarrow q\leftarrow r\)    & \(v1,v4\) & \(y_1,y_4\) & capitulation type & \(\mathfrak{M}\)                   & \(\ell_3(k)\) \\
\hline
   1 &   \(3\,913\) & \(13\rightarrow 7\leftarrow 43\)  & \(3,3\) & \(1,1\) & \(\mathrm{a}.3^\ast\) & \(\langle 81,7\rangle^2\) & \(=2\) \\
   9 &  \(22\,581\) & \(9\rightarrow 193\leftarrow 13\) & \(3,3\) & \(0,0\) & \(\mathrm{d}.19\) & \(\langle 729,41\rangle^2\) & \(\ge 2\) \\
  11 &  \(25\,929\) & \(9\rightarrow 67\leftarrow 43\)  & \(3,3\) & \(0,0\) & \(\mathrm{b}.10\) & \(\langle 729,34..36\rangle^2\) & \(\ge 2\) \\
  15 &  \(30\,457\) & \(7\rightarrow 19\leftarrow 229\) & \(4,4\) & \(1,1\) & \(\mathrm{b}.10\) & \(\langle 729,37..39\rangle^2\) & \(\ge 2\) \\
  18 &  \(34\,029\) & \(19\rightarrow 9\leftarrow 199\) & \(4,4\) & \(1,0\) & \(\mathrm{d}.19\) & \(\langle 2187,248\vert 249\rangle^2\) & \(\ge 2\) \\
  23 &  \(41\,839\) & \(43\rightarrow 7\leftarrow 139\) & \(4,4\) & \(0,0\) & \(\mathrm{b}.10\) & \(\langle 6561,693..698\rangle^2\) & \(\ge 2\) \\
  35 &  \(74\,043\) & \(19\rightarrow 9\leftarrow 433\) & \(3,3\) & \(0,0\) & \(\mathrm{H}.4\)  & \(\langle 2187,65\vert 67\rangle^2\) & \(\ge 3\) \\
     &              &                                   &         &         & or                & \(\langle 6561,714..719\vert 738..743\rangle\) & \(\ge 2\) \\
  39 &  \(82\,327\) & \(7\rightarrow 19\leftarrow 619\) & \(3,3\) & \(0,0\) & \(\mathrm{b}.10\) & \(\langle 729,37..39\rangle^2\) & \(\ge 2\) \\
  42 &  \(83\,817\) & \(9\rightarrow 67\leftarrow 139\) & \(5,4\) & \(0,1\) & \(\mathrm{b}.10\) & \(\langle 6561,693..698\rangle^2\) & \(\ge 2\) \\
\hline
\end{tabular}
}
\end{center}
\end{table}


\subsection{Category II, Graph 2}
\label{ss:Cat2Gph2}

\noindent
Let \((k_1,\ldots,k_4)\) be a quartet of cyclic cubic number fields
sharing the common conductor \(c=pqr\),
belonging to Graph \(2\) of Category \(\mathrm{II}\)
with combined cubic residue symbol
\(\lbrack p,q,r\rbrack_3=\lbrace q\leftarrow p\rightarrow r\leftarrow q\rbrace\).

\begin{lemma}
\label{lem:Cat2Gph2}
\textbf{(3-class ranks of components for \(\mathrm{II}.2\).)}
Under the normalizing assumption that \(r\) splits in \(k_{pq}\),
precisely the two components \(k_1\) and \(k_2\) of the quartet
have elementary bicyclic \(3\)-class group
\(\mathrm{Cl}_3(k_1)\simeq\mathrm{Cl}_3(k_2)\simeq (3,3)\)
of rank \(2\),
whereas the other two components
have \(3\)-class rank \(\varrho_3(k_3)=\varrho_3(k_4)=3\).
Thus, the \textbf{tame} condition
\(9\mid h_3(B_j)=(U_j:V_j)\in\lbrace 9,27\rbrace\), \(r_j=2\),
is only satisfied for the bicyclic bicubic fields \(B_j\)
with \(j\in\lbrace 1,2,7\rbrace\).
\end{lemma}

\begin{proof}
Since \(p\rightarrow r\), \(p\) splits in \(k_r\).
Since \(r\leftarrow q\), \(q\) splits in \(k_r\),
and also splits in \(B_7=k_1k_2k_r\tilde{k}_{pq}\).
By the normalizing assumption that \(r\) splits in \(k_{pq}\),
it also splits in
\(B_1=k_1k_{pq}k_{pr}k_{qr}\) and \(B_2=k_2k_{pq}\tilde{k}_{pr}\tilde{k}_{qr}\).
The primes \(p,q,r\) share the common decomposition type
\((e,f,g)=(3,1,3)\) in the bicyclic bicubic field
\(B_{10}=k_3k_4k_rk_{pq}\),
which implies that \(\varrho_3(k_3)=\varrho_3(k_4)=3\),
according to
\cite[Prop. 4.4, pp. 43--44]{Ay1995}.
Finally, only
\(B_1=k_1k_{pq}k_{pr}k_{qr}\),
\(B_2=k_2k_{pq}\tilde{k}_{pr}\tilde{k}_{qr}\),
\(B_7=k_1k_2k_r\tilde{k}_{pq}\)
do not contain \(k_3\), \(k_4\).
\end{proof}


\begin{proposition}
\label{prp:Cat2Gph2}
\textbf{(Sub-doublet with \(3\)-rank two for \(\mathrm{II}.2\).)}
For fixed \(\mu\in\lbrace 1,2\rbrace\),
let \(\mathfrak{p},\mathfrak{q},\mathfrak{r}\) be the prime ideals of \(k_\mu\)
over \(p,q,r\), that is
\(p\mathcal{O}_{k_\mu}=\mathfrak{p}^3\),
\(q\mathcal{O}_{k_\mu}=\mathfrak{q}^3\),
\(r\mathcal{O}_{k_\mu}=\mathfrak{r}^3\),
then the \textbf{principal factor} of \(k_\mu\) is
\(A(k_\mu)=p\), with \(\lbrack\mathfrak{p}\rbrack=1\),
and the \(3\)-class group of \(k_\mu\) is
\begin{equation}
\label{eqn:GenCat2Gph2}
\mathrm{Cl}_3(k_\mu)=
\langle\lbrack\mathfrak{q}\rbrack,\lbrack\mathfrak{r}\rbrack\rangle\simeq (3,3).
\end{equation}

The unramified cyclic cubic relative extensions of \(k_\mu\)
are among the absolutely bicyclic bicubic fields \(B_i\), \(1\le i\le 10\).
The unique \(B_\mu\), \(\mu\in\lbrace 1,2\rbrace\), which only contains \(k_\mu\),
has norm class group
\(N_{B_\mu/k_\mu}(\mathrm{Cl}_3(B_\mu))=\langle\lbrack\mathfrak{r}\rbrack\rangle\),
transfer kernel
\[
\ker(T_{B_\mu/k_\mu})\ge\langle\lbrack\mathfrak{q}\rbrack\rangle
\]
and \(3\)-class group
\(\mathrm{Cl}_3(B_\mu)=
\langle\lbrack\mathfrak{R}_1\rbrack,\lbrack\mathfrak{R}_2\rbrack,\lbrack\mathfrak{R}_3\rbrack\rangle\ge (3,3)\),
generated by the classes of the prime ideals of \(B_\mu\) over
\(\mathfrak{r}\mathcal{O}_{B_\mu}=\mathfrak{R}_1\mathfrak{R}_2\mathfrak{R}_3\).
The unique \(B_7=k_1k_2k_r\tilde{k}_{pq}\), which contains \(k_1\) and \(k_2\),
has norm class group
\(N_{B_7/k_\mu}(\mathrm{Cl}_3(B_7))=\langle\lbrack\mathfrak{q}\rbrack\rangle\),
\textbf{cyclic} transfer kernel
\begin{equation}
\label{eqn:PrtCat2Gph2}
\ker(T_{B_7/k_\mu})=\langle\lbrack\mathfrak{r}\rbrack\rangle
\end{equation}
of order \(3\),
and \textbf{elementary tricyclic} \(3\)-class group
\(\mathrm{Cl}_3(B_7)=
\langle\lbrack\mathfrak{Q}_1\rbrack,\lbrack\mathfrak{Q}_2\rbrack,\lbrack\mathfrak{Q}_3\rbrack\rangle\simeq (3,3,3)\),
generated by the classes of the prime ideals of \(B_7\) over
\(\mathfrak{q}\mathcal{O}_{B_7}=\mathfrak{Q}_1\mathfrak{Q}_2\mathfrak{Q}_3\).
The remaining two \(B_i>k_\mu\), \(i\in\lbrace 5,6,8,9\rbrace\),
more precisely, \(i\in\lbrace 5,6\rbrace\) for \(\mu=1\),
and \(i\in\lbrace 8,9\rbrace\) for \(\mu=2\),
have norm class group 
\(\langle\lbrack\mathfrak{qr}\rbrack\rangle\) respectively
\(\langle\lbrack\mathfrak{qr}^2\rbrack\rangle\),
and transfer kernel
\[
\ker(T_{B_i/k_\mu})\ge\langle\lbrack\mathfrak{q}\rbrack\rangle.
\]
In terms of \textbf{decisive principal factors}
\(A(k_\nu)=p^{x_\nu}q^{y_\nu}r^{z_\nu}\) for \(\nu\in\lbrace 3,4\rbrace\),
the ranks of principal factor matrices
\(M_i\), \(i\in\lbrace 3,4,5,6,8,9\rbrace\),
of \textbf{wild} bicyclic bicubic fields are
\begin{equation}
\label{eqn:Cat2Gph2Ranks}
r_3=r_5=r_9=3 \text{ iff } z_3\ne 0 \text{ iff } r\mid A(k_3)
\text{ and } 
r_4=r_6=r_8=3 \text{ iff } z_4\ne 0 \text{ iff } r\mid A(k_4).
\end{equation}
\end{proposition}

\begin{proof}
Since \(q\leftarrow p\),
two principal factors are \(A(k_{pq})=A(\tilde{k}_{pq})=p\);
since \(p\rightarrow r\),
two further principal factors are \(A(k_{pr})=A(\tilde{k}_{pr})=p\);
since \(r\leftarrow q\),
two further principal factors are \(A(k_{qr})=A(\tilde{k}_{qr})=q\);
each by Proposition
\ref{prp:Principal2}.
Since \(q\leftarrow p\rightarrow r\) is universally repelling, we have
\(A(k_1)=A(k_2)=p\),
by Proposition
\ref{prp:Principal3}.

Thus \(\mathfrak{p}=\alpha\mathcal{O}_{k_\mu}\) is a principal ideal
with trivial class \(\lbrack\mathfrak{p}\rbrack=1\), for \(\mu\in\lbrace 1,2\rbrace\),
whereas the classes
\(\lbrack\mathfrak{q}\rbrack,\lbrack\mathfrak{r}\rbrack\)
are non-trivial.
We propose \(A(k_\nu)=p^{x_\nu}q^{y_\nu}r^{z_\nu}\) for \(\nu\in\lbrace 3,4\rbrace\).

Since the \textit{tame} condition \(9\mid h_3(B_j)=(U_j:V_j)\) is satisfied
for \(j\in\lbrace 1,2,7\rbrace\),
the rank of the corresponding principal factor matrix \(M_j\)
must be \(r_1=r_2=r_7=2\).
Due to the principal factors \(A(k_1)=A(k_2)=p\),
this also follows by direct calculation,
but has no further consequences.
For every \textit{wild} bicyclic bicubic field \(B_j\),
\(j\in\lbrace 3,4,5,6,8,9,10\rbrace\),
the rank \(r_j\) is calculated with row operations
on the associated principal factor matrices \(M_j\): \\
\(M_3=
\begin{pmatrix}
x_3 & y_3 & z_3 \\
1 & 0 & 0 \\
1 & 0 & 0 \\
0 & 1 & 0 \\
\end{pmatrix}\),
\(M_5=
\begin{pmatrix}
1 & 0 & 0 \\
x_3 & y_3 & z_3 \\
1 & 0 & 0 \\
0 & 1 & 0 \\
\end{pmatrix}\),
\(M_9=
\begin{pmatrix}
1 & 0 & 0 \\
x_3 & y_3 & z_3 \\
0 & 1 & 0 \\
1 & 0 & 0 \\
\end{pmatrix}\).

For \(B_3=k_3\tilde{k}_{pq}\tilde{k}_{pr}k_{qr}\), 
\(M_3\)
leads to the decisive pivot element
\(z_3\)
in the last column,
similarly,
for \(B_5=k_1k_3k_p\tilde{k}_{qr}\),
\(M_5\)
leads to \(z_3\),
and similarly,
for \(B_9=k_2k_3k_qk_{pr}\), 
\(M_9\)
leads to \(z_3\).
So, \(r_3=r_5=r_9=3\)
iff \(z_3\ne 0\)
iff \(r\mid A(k_3)\).
Next we consider: \\
\(M_4=
\begin{pmatrix}
x_4 & y_4 & z_4 \\
1 & 0 & 0 \\
1 & 0 & 0 \\
0 & 1 & 0 \\
\end{pmatrix}\),
\(M_6=
\begin{pmatrix}
1 & 0 & 0 \\
x_4 & y_4 & z_4 \\
0 & 1 & 0 \\
1 & 0 & 0 \\
\end{pmatrix}\),
\(M_8=
\begin{pmatrix}
1 & 0 & 0 \\
x_4 & y_4 & z_4 \\
1 & 0 & 0 \\
0 & 1 & 0 \\
\end{pmatrix}\).

For \(B_4=k_4\tilde{k}_{pq}k_{pr}\tilde{k}_{qr}\),
\(M_4\)
leads to the decisive pivot element
\(z_4\)
in the last column,
similarly,
for \(B_6=k_1k_4k_q\tilde{k}_{pr}\),
\(M_6\)
leads to \(z_4\),
and similarly,
for \(B_8=k_2k_4k_pk_{qr}\),
\(M_8\)
leads to \(z_4\).
So, \(r_4=r_6=r_8=3\)
iff \(z_4\ne 0\)
iff \(r\mid A(k_4)\).


By Lemma
\ref{lem:UnitIndices},
the minimal subfield unit index
\((U_j:V_j)=3\) for \(r_j=3\)
corresponds to
the maximal unit norm index \((U(k_\mu):N_{B_j/k_\mu}(U_j))=3\),
associated to a \textit{total} transfer kernel
\(\#\ker(T_{B_j/k_\mu})=9\).

As mentioned in the proof of Lemma
\ref{lem:Cat2Gph2}
already: \\
Since \(r\) splits in \(k_{pq}\),
it also splits in \(B_1=k_1k_{pq}k_{pr}k_{qr}\),
\(B_2=k_2k_{pq}\tilde{k}_{pr}\tilde{k}_{qr}\), i.e.,
\(\mathfrak{r}\mathcal{O}_{B_\mu}=\mathfrak{R}_1\mathfrak{R}_2\mathfrak{R}_3\). \\
Since \(q\) splits in \(k_r\),
it also splits in \(B_7=k_1k_2k_r\tilde{k}_{pq}\), i.e.,
\(\mathfrak{q}\mathcal{O}_{B_7}=\mathfrak{Q}_1\mathfrak{Q}_2\mathfrak{Q}_3\).

Since \(\mathfrak{q}\) is principal in \(k_q\), \(k_{qr}\), \(\tilde{k}_{qr}\),
\(\lbrack\mathfrak{q}\rbrack\) capitulates in
\(B_1=k_1k_{pq}k_{pr}k_{qr}\),
\(B_2=k_2k_{pq}\tilde{k}_{pr}\tilde{k}_{qr}\),
\(B_5=k_1k_3k_p\tilde{k}_{qr}\),
\(B_6=k_1k_4k_q\tilde{k}_{pr}\),
\(B_8=k_2k_4k_pk_{qr}\),
\(B_9=k_2k_3k_qk_{pr}\);
since \(\mathfrak{r}\) is principal in \(k_r\),
\(\lbrack\mathfrak{r}\rbrack\) capitulates in
\(B_7=k_1k_2k_r\tilde{k}_{pq}\)
(Proposition
\ref{prp:Bicyc}).
This gives a \textit{transposition}, either \((1,7)\) or \((2,7)\).

The minimal unit norm index \((U(k_\mu):N_{B_7/k_\mu}(U_7))=1\),
associated to the \textit{partial} transfer kernel
\(\ker(T_{B_7/k_\mu})=\langle\lbrack\mathfrak{r}\rbrack\rangle\),
corresponds to the maximal subfield unit index
\(h_3(B_7)=(U_7:V_7)=27\),
giving rise to the \textit{elementary tricyclic} type invariants
\(\mathrm{Cl}_3(B_7)=
\langle\lbrack\mathfrak{Q}_1\rbrack,\lbrack\mathfrak{Q}_2\rbrack,\lbrack\mathfrak{Q}_3\rbrack\rangle\simeq (3,3,3)\).
\end{proof}

\noindent
In terms of capitulation targets in Corollary
\ref{cor:Three},
Proposition
\ref{prp:Cat2Gph2}
and parts of its proof
are now summarized in Table
\ref{tbl:UniCat2Gph2}
with transposition in \textbf{bold} font.


\renewcommand{\arraystretch}{1.1}

\begin{table}[ht]
\caption{Norm class groups and minimal transfer kernels for Graph II.2}
\label{tbl:UniCat2Gph2}
\begin{center}
\begin{tabular}{|c||c|c|c|c||c|c|c|c|}
\hline
 Base         & \multicolumn{4}{c||}{\(k_1\)} & \multicolumn{4}{c||}{\(k_2\)} \\
\hline
 Ext          & \(B_1\) & \(B_5\) & \(B_6\) & \(B_7\) & \(B_2\) & \(B_7\) & \(B_8\) & \(B_9\) \\
\hline
 NCG          & \(\mathfrak{r}\) & \(\mathfrak{qr}\) & \(\mathfrak{qr}^2\) & \(\mathfrak{q}\) & \(\mathfrak{r}\) & \(\mathfrak{q}\) & \(\mathfrak{qr}\) & \(\mathfrak{qr}^2\) \\
 TK           & \(\mathfrak{q}\) & \(\mathfrak{q}\)  & \(\mathfrak{q}\)    & \(\mathfrak{r}\) & \(\mathfrak{q}\) & \(\mathfrak{r}\) & \(\mathfrak{q}\)  & \(\mathfrak{q}\) \\
\(\varkappa\) & \(\mathbf{4}\) & \(4\) & \(4\) & \(\mathbf{1}\) & \(\mathbf{2}\) & \(\mathbf{1}\) & \(2\) & \(2\) \\
\hline
\end{tabular}
\end{center}
\end{table}


\begin{theorem}
\label{thm:Cat2Gph2}
\textbf{(Second \(3\)-class group for \(\mathrm{II}.2\).)}
Let \((k_1,\ldots,k_4)\) be a quartet of cyclic cubic number fields
sharing the common conductor \(c=pqr\),
belonging to Graph \(2\) of Category \(\mathrm{II}\)
with combined cubic residue symbol
\(\lbrack p,q,r\rbrack_3=\lbrace q\leftarrow p\rightarrow r\leftarrow q\rbrace\).
Without loss of generality, suppose that \(r\) splits in \(k_{pq}\), and thus
\(\mathrm{Cl}_3(k_1)\simeq\mathrm{Cl}_3(k_2)\simeq (3,3)\),
and \(\varrho_3(k_3)=\varrho_3(k_4)=3\).

Then the \textbf{minimal transfer kernel type} (mTKT) of \(k_\mu\), \(1\le\mu\le 2\),
is \(\varkappa_0=(2111)\), type \(\mathrm{H}.4\),
and the other possible capitulation types in ascending order
\(\varkappa_0<\varkappa^\prime<\varkappa^{\prime\prime}<\varkappa^{\prime\prime\prime}\)
are
\(\varkappa^\prime=(2110)\), type \(\mathrm{d}.19\),
\(\varkappa^{\prime\prime}=(2100)\), type \(\mathrm{b}.10\), and
\(\varkappa^{\prime\prime\prime}=(2000)\), type \(\mathrm{a}.3^\ast\).

To identify the second \(3\)-class group
\(\mathfrak{M}=\mathrm{Gal}(\mathrm{F}_3^2(k_\mu)/k_\mu)\), \(1\le\mu\le 2\),
let the decisive \textbf{principal factors} of \(k_\nu\), \(3\le\nu\le 4\), be
\(A(k_\nu)=p^{x_\nu}q^{y_\nu}r^{z_\nu}\), and additionally assume
the \textbf{regular} situation where both
\(\mathrm{Cl}_3(k_3)\simeq\mathrm{Cl}_3(k_4)\simeq (3,3,3)\)
are elementary tricyclic. Then

\begin{equation}
\label{eqn:Cat2Gph2}
\mathfrak{M}\simeq
\begin{cases}
\langle 81,7\rangle,\ \alpha=\lbrack 111,11,11,11\rbrack,\ \varkappa=(2000) & \text{ if } z_3\ne 0,\ z_4\ne 0,\ \mathcal{N}=1, \\
\langle 729,34..39\rangle,\ \alpha=\lbrack 111,111,21,21\rbrack,\ \varkappa=(2100) & \text{ if } z_3=z_4=0,\ \mathcal{N}=2, \\
\langle 729,41\rangle,\ \alpha=\lbrack 111,111,22,21\rbrack,\ \varkappa=(2110) & \text{ if } z_3=z_4=0,\ \mathcal{N}=3, \\
\langle 6561,714..719\vert 738..743\rangle \text{ or } & \\
\langle 2187,65\vert 67\rangle,\ \alpha=\lbrack 111,111,22,22\rbrack,\ \varkappa=(2111) & \text{ if } z_3=z_4=0,\ \mathcal{N}=4,
\end{cases}
\end{equation}
where \(\mathcal{N}:=\#\lbrace 1\le j\le 10\mid k_\mu<B_j,\ I_j=27\rbrace\).
Only in the leading row, the \(3\)-class field tower has warranted group
\(\mathfrak{G}=\mathrm{Gal}(\mathrm{F}_3^\infty(k_\mu)/k_\mu)\simeq\mathfrak{M}\),
with length \(\ell_3(k_\mu)=2\). Otherwise,
even if the relation rank \(d_2(\mathfrak{M})\le 4\) is admissible,
tower length \(\ell_3(k_\mu)\ge 3\) cannot be excluded.
\end{theorem}

\begin{proof}
We give the proof for \(k_1\)
with unramified cyclic cubic extensions
\(B_1,B_5,B_6,B_7\), 
(The proof for \(k_2\)
with unramified cyclic cubic extensions
\(B_2,B_7,B_8,B_9\)
is similar.)
We know that the tame ranks are
\(r_1=r_2=r_7=2\),
and thus \(I_1,I_2,I_7\in\lbrace 9,27\rbrace\),
in particular, \(I_7=27\), whence certainly \(\mathcal{N}\ge 1\).
Further, the wild ranks are
\(r_4=r_6=r_8=3\)
iff \(z_4\ne 0\),
and
\(r_3=r_5=r_9=3\)
iff \(z_3\ne 0\).

In the \textbf{regular situation} where the
\(3\)-class groups of \(k_3\) and \(k_4\)
are elementary tricyclic,
tight bounds arise for the abelian quotient invariants \(\alpha\) of
the group \(\mathfrak{M}\):

The first scenario,
\(z_3\ne 0\), \(z_4\ne 0\),
is equivalent to
\(\mathcal{N}=1\),
\(h_3(B_5)=h_3(B_9)=\frac{1}{3}h_3(k_3)=9\),
\(h_3(B_6)=h_3(B_8)=\frac{1}{3}h_3(k_4)=9\),
\(h_3(B_1)=I_1=h_3(B_2)=I_2=9\),
\(h_3(B_7)=I_7=27\),
that is \(\alpha=\lbrack 111,11,11,11\rbrack\)
and consequently \(\varkappa=(2000)\),
since \(\langle 81,7\rangle\) is unique
with this \(\alpha\).

The other three scenarios share
\(z_3=z_4=0\),
and an explicit transposition
between \(B_1\) and \(B_7\),
giving rise to \(\varkappa=(21\ast\ast)\),
and common \(h_3(B_1)=I_1=27\), \(\alpha=\lbrack 111,111,\ast,\ast\rbrack\).

The second scenario with \(\mathcal{N}=2\)
is supplemented by 
\(h_3(B_5)=h_3(k_3)=27\),
\(h_3(B_6)=h_3(k_4)=27\),
giving rise to \(\alpha=\lbrack 111,111,21,21\rbrack\),
\(\varkappa=(2100)\),
characteristic for \(\langle 729,34..39\rangle\)
(Cor.
\ref{cor:ClassGroup33}).

The third scenario with \(\mathcal{N}=3\)
is supplemented by 
\(h_3(B_5)=3h_3(k_3)=81\),
\(h_3(B_6)=h_3(k_4)=27\),
giving rise to \(\alpha=\lbrack 111,111,22,21\rbrack\),
\(\varkappa=(2110)\),
characteristic for \(\langle 729,41\rangle\).

The fourth scenario with \(\mathcal{N}=4\)
is supplemented by 
\(h_3(B_5)=3h_3(k_3)=81\),
\(h_3(B_6)=3h_3(k_4)=81\),
giving rise to \(\alpha=\lbrack 111,111,22,22\rbrack\),
\(\varkappa=(2111)\),
characteristic for \(\langle 2187,65\vert 67\rangle\) or
\(\langle 6561,714..719\vert 738..743\rangle\)
with coclass \(\mathrm{cc}=3\).
If \(d_2(\mathfrak{M})=5\),
then \(\ell_3(k_\mu)\ge 3\).
\end{proof}


\begin{corollary}
\label{cor:UniCat2Gph2}
\textbf{(Uniformity of the sub-doublet for \(\mathrm{II}.2\).)}
The components of the sub-doublet with \(3\)-rank two
share a common capitulation type
\(\varkappa(k_1)\sim\varkappa(k_2)\),
common abelian type invariants
\(\alpha(k_1)\sim\alpha(k_2)\),
and a common second \(3\)-class group
\(\mathrm{Gal}(\mathrm{F}_3^2(k_1)/k_1)\simeq\mathrm{Gal}(\mathrm{F}_3^2(k_2)/k_2)\).
\end{corollary}

\begin{proof}
This follows immediately from Theorem
\ref{thm:Cat2Gph2}
and Table
\ref{tbl:UniCat2Gph2}.
\end{proof}


\begin{example}
\label{exm:Cat2Gph2}
Prototypes for Graph \(\mathrm{II}.2\), i.e.,
minimal conductors for each scenario in Theorem
\ref{thm:Cat2Gph2}
have been found for each \(\mathcal{N}\in\lbrace 1,2,3,4\rbrace\).

There are \textbf{regular} cases:
\(c=6\,327\) with symbol
\(\lbrace 19\rightarrow 9\leftarrow 37\rightarrow 19\rbrace\)
and \(\mathfrak{G}=\mathfrak{M}=\langle 81,7\rangle\);
\(c=41\,629\) with symbol
\(\lbrace 19\rightarrow 313\leftarrow 7\rightarrow 19\rbrace\)
and \(\mathfrak{M}=\langle 729,34..36\rangle\) (Corollary
\ref{cor:ClassGroup33});
\(c=56\,547\) with symbol
\(\lbrace 61\rightarrow 103\leftarrow 9\rightarrow 61\rbrace\)
and \(\mathfrak{M}=\langle 729,41\rangle\); and,
with \textbf{considerable statistic delay},
\(c=389\,329\) with ordinal number \(207\), symbol
\(\lbrace 19\rightarrow 661\leftarrow 31\rightarrow 19\rbrace\)
and either \(\mathfrak{M}=\langle 2187,65\vert 67\rangle\)
with \(d_2(\mathfrak{M})=5\)
or \(\mathfrak{M}=\langle 6561,714..719\vert 738..743\rangle\)
with \(d_2(\mathfrak{M})=4\).

Further, there are \textbf{singular} cases:
\(c=27\,873\) with symbol
\(\lbrace 19\rightarrow 9\leftarrow 163\rightarrow 19\rbrace\)
and \(\mathfrak{M}=\langle 729,34..36\rangle\) (Corollary
\ref{cor:ClassGroup33});
\(c=29\,197\) with symbol
\(\lbrace 43\rightarrow 7\leftarrow 97\rightarrow 43\rbrace\)
and \(\mathfrak{M}=\langle 2187,253\rangle\); and
\(c=63\,511\) with symbol
\(\lbrace 43\rightarrow 7\leftarrow 211\rightarrow 43\rbrace\)
and \(\mathfrak{M}=\langle 729,37..39\rangle\) (Corollary
\ref{cor:ClassGroup33}).

Finally, there is the \textbf{super-singular}
\(c=66\,157\) with symbol
\(\lbrace 13\rightarrow 7\leftarrow 727\rightarrow 13\rbrace\)
and \(\mathfrak{M}=\langle 6561,1989\rangle\).

With exception of \(\langle 81,7\rangle\),
all groups have non-metabelian descendants,
respectively extensions.

\end{example}


In Table
\ref{tbl:ProtoCat2Gph2},
we summarize the prototypes of Graph \(\mathrm{II}.2\)
in the same way as in Table
\ref{tbl:ProtoCat1Gph1},
except that
regularity, resp. (super-)singularity, is expressed by
\(3\)-valuations \(v_\nu=v_3(\#\mathrm{Cl}(k_\nu))\)
of class numbers of critical fields \(k_\nu\), \(\nu=3,4\),
and critical exponents are \(z_\nu\) in principal factors
\(A(k_\nu)=p^{x_\nu}q^{y_\nu}r^{z_\nu}\),
\(\nu=3,4\).

\renewcommand{\arraystretch}{1.1}

\begin{table}[ht]
\caption{Prototypes for Graph II.2}
\label{tbl:ProtoCat2Gph2}
\begin{center}
{\scriptsize
\begin{tabular}{|c|c||c|c|c|c|c|c|}
\hline
 No. & \(c\)        & \(q\rightarrow r\leftarrow p\rightarrow q\)      & \(v3,v4\) & \(z_3,z_4\) & capitulation type & \(\mathfrak{M}\)                   & \(\ell_3(k)\) \\
\hline
   1 &   \(6\,327\) & \(19\rightarrow 9\leftarrow 37\rightarrow 19\)   & \(3,3\) & \(1,1\) & \(\mathrm{a}.3^\ast\) & \(\langle 81,7\rangle^2\) & \(=2\) \\
   8 &  \(27\,873\) & \(19\rightarrow 9\leftarrow 163\rightarrow 19\)  & \(4,4\) & \(1,1\) & \(\mathrm{b}.10\) & \(\langle 729,34..36\rangle^2\) & \(\ge 2\) \\
  10 &  \(29\,197\) & \(43\rightarrow 7\leftarrow 97\rightarrow 43\)   & \(4,4\) & \(0,1\) & \(\mathrm{b}.10\) & \(\langle 2187,253\rangle^2\) & \(\ge 3\) \\
  14 &  \(41\,629\) & \(19\rightarrow 313\leftarrow 7\rightarrow 19\)  & \(3,3\) & \(0,0\) & \(\mathrm{b}.10\) & \(\langle 729,34..36\rangle^2\) & \(\ge 2\) \\
  23 &  \(56\,547\) & \(61\rightarrow 103\leftarrow 9\rightarrow 61\)  & \(3,3\) & \(0,0\) & \(\mathrm{d}.19\) & \(\langle 729,41\rangle^2\) & \(\ge 2\) \\
  28 &  \(63\,511\) & \(43\rightarrow 7\leftarrow 211\rightarrow 43\)  & \(4,4\) & \(1,1\) & \(\mathrm{b}.10\) & \(\langle 729,37..39\rangle^2\) & \(\ge 2\) \\
  31 &  \(66\,157\) & \(13\rightarrow 7\leftarrow 727\rightarrow 13\)  & \(5,4\) & \(1,0\) & \(\mathrm{d}.19\) & \(\langle 6561,1989\rangle^2\) & \(\ge 2\) \\
 207 & \(389\,329\) & \(19\rightarrow 661\leftarrow 31\rightarrow 19\) & \(3,3\) & \(0,0\) & \(\mathrm{H}.4\)  & \(\langle 2187,65\vert 67\rangle^2\) & \(\ge 3\) \\
     &              &                                                  &         &         & or                & \(\langle 6561,714..719\vert 738..743\rangle\) & \(\ge 2\) \\
\hline
\end{tabular}
}
\end{center}
\end{table}


\section{Category III, Graphs 1--4}
\label{s:Cat3Gph1To4}

\noindent
Let the combined cubic residue symbol
\(\lbrack p,q,r\rbrack_3\)
of three prime(power)s dividing the conductor \(c=pqr\)
be either
\(\lbrace p,q,r;\delta\not\equiv 0\,(\mathrm{mod}\,3)\rbrace\)
or
\(\lbrace p\rightarrow q;r\rbrace\)
or
\(\lbrace p\rightarrow q\rightarrow r\rbrace\)
or
\(\lbrace p\rightarrow q\rightarrow r\rightarrow p\rbrace\).
The symbol does not contain any mutual cubic residues.
We verify a conjecture in
\cite[Cnj. 1, p. 48]{Ma2022}.

\begin{theorem}
\label{thm:Cat3Gph1To4}
A cyclic cubic field \(k\) with conductor \(c=pqr\),
divisible by exactly three prime(power)s \(p,q,r\),
has an abelian \(3\)-class field tower with group
\(\mathfrak{G}=\mathrm{Gal}(\mathrm{F}_3^\infty(k)/k)\simeq\langle 9,2\rangle\),
\(\alpha=\lbrack 1,1,1,1\rbrack\), \(\varkappa=(0000)\),
if and only if the primes \(p,q,r\) form one of the four Graphs \(1\)--\(4\)
of Category \(\mathrm{III}\).
\end{theorem}

\begin{proof}
Ayadi
\cite[Thm. 4.1, pp. 76--77]{Ay1995}
has proved the sufficiency of the condition.
He does not claim explicitly that the condition is also necessary.
However, his techniques are able to prove both directions.
Recall that for both Graphs \(1\)--\(2\) of the Categories \(\mathrm{I}\) and \(\mathrm{II}\),
there is at least one component of the quartet \((k_\mu)_{\mu=1}^4\)
with \(3\)-class rank \(\varrho(k_\mu)=3\),
and that for all Graphs \(5\)--\(9\) of Category \(\mathrm{III}\), two primes
\(p\leftrightarrow q\) are mutual cubic residues, 
according to Theorem
\ref{thm:3ClassRank}.
In contrast,
precisely for the Graphs \(1\)--\(4\) of Category \(\mathrm{III}\),
the symbol \(\lbrack p,q,r\rbrack_3\) does not contain any mutual cubic residues,
and all four components
have \(3\)-class rank \(\varrho(k_\mu)=2\)
and elementary bicyclic \(3\)-class group \(\mathrm{Cl}_3(k_\mu)\simeq (3,3)\),
whence
these are the only cases where all bicyclic bicubic fields \(B_j\), \(1\le j\le 10\),
satisfy the \textit{tame} relation \(h_3(B_j)=(U_j:V_j)=3\) with matrix rank \(r_j=3\).
This is equivalent with abelian type invariants \(\alpha(k_\mu)=\lbrack 1,1,1,1\rbrack\)
for all \(1\le\mu\le 4\).
By the strategy of pattern recognition
\cite{Ma2020},
this enforces the abelian group \(\mathfrak{G}\simeq\langle 9,2\rangle\simeq (3,3)\),
which is the unique \(3\)-group \(G\) with \(G/G^\prime\simeq (3,3)\)
and abelian quotient invariants \(\alpha(G)=\lbrack 1,1,1,1\rbrack\).
\end{proof}

\noindent
For the prototypes of Graphs \(1,\ldots,4\) of Category \(\mathrm{III}\) see
\cite[Tbl. 6.3, p. 48]{Ma2022}.
Systematic tables have been presented at
\texttt{http://www.algebra.at/ResearchFrontier2013ThreeByThree.htm}
in sections \S\S\ 2.1--2.2.


\section{Category III, Graphs 5--9}
\label{s:Cat3Gph5To9}

\noindent
In this section,
the combined cubic residue symbol
\(\lbrack p,q,r\rbrack_3\)
of three prime(power)s dividing the conductor \(c=pqr\)
contains a unique pair \(p\leftrightarrow q\) of mutual cubic residues.

Consequently, the \textbf{decisive principal factors}
\begin{equation}
\label{eqn:PFCat3Gph5To9}
A(k_{pq})=p^mq^n, \quad A(\tilde{k}_{pq})=p^{\tilde{m}}q^{\tilde{n}}
\end{equation}
must be assumed with variable exponents in \(\lbrace 0,1,2\rbrace\),
such that \((m,n)\ne (0,0)\)
and \((\tilde{m},\tilde{n})\ne (0,0)\).
Concerning \(3\)-class groups of cyclic cubic subfields \(k<k^\ast\) with \(t=2\),
an elementary cyclic group \(\mathrm{Cl}_3(k)\simeq (3)\) is warranted for
\(k\in\lbrace k_{pr},\tilde{k}_{pr},k_{qr},\tilde{k}_{qr}\rbrace\).
For the critical fields \(k\in\lbrace k_{pq},\tilde{k}_{pq}\rbrace\),
however, we must distinguish
the \textbf{regular} situation \(\mathrm{Cl}_3(k_f^\ast)\simeq (3,3)\)
in terms of the sub-genus field \(k_f^\ast=k_{pq}\cdot\tilde{k}_{pq}\)
with partial conductor \(f=pq\) which divides \(c=pqr\),
where \(\mathrm{Cl}_3(k_{pq})\simeq\mathrm{Cl}_3(\tilde{k}_{pq})\simeq (3,3)\)
and equality \((m,n)=(\tilde{m},\tilde{n})\) is warranted,
as opposed to
the \textbf{singular} situation \(\mathrm{Cl}_3(k_f^\ast)\simeq (3,3,3)\),
and the \textbf{super-singular} situation \(81\mid h_3(k_f^\ast)\),
where usually
\(\mathrm{Cl}_3(k_{pq})\simeq\mathrm{Cl}_3(\tilde{k}_{pq})\simeq (9,3)\).

For doublets \((k_{pq},\tilde{k}_{pq})\)
with conductor \(f=pq\)
and non-elementary bicyclic \(3\)-class group,
a distinction arises from the \(3\)-valuation \(v^\ast:=v_3(h(k_f^\ast))\) of the
class number of the \(3\)-genus field \(k_f^\ast\):


\begin{definition}
\label{dfn:SingularDoublet}
A quartet \((k_\mu)_{1\le\mu\le 4}\) with conductor \(c=pqr\) and
its sub-doublet \((k_{pq},\tilde{k}_{pq})\) of cyclic cubic fields
with common partial conductor \(f=pq\) is called
\begin{equation}
\label{eqn:SingularDoublet}
\begin{cases}
\textbf{regular}        & \text{ if } v^\ast\in\lbrace 0,1,2\rbrace, \\
\textbf{singular}       & \text{ if } v^\ast=3,                      \\
\textbf{super-singular} & \text{ if } v^\ast\in\lbrace 4,5,6,\ldots\rbrace.
\end{cases}
\end{equation}
\end{definition}


Let \((k_1,\ldots,k_4)\) be the quartet of cyclic cubic number fields
sharing the common discriminant \(d=c^2\) with conductor \(c=pqr\),
divisible by exactly three primes \(\equiv 1\,(\mathrm{mod}\,3)\)
(one among them may be the prime power \(3^2\)),
and belonging to one of the Graphs \(5\)--\(9\) of Category \(\mathrm{III}\).
According to Theorem
\ref{thm:3ClassRank},
\(\mathrm{Cl}_3(k_\mu)\simeq (3,3)\)
and thus \(h_3(k_\mu)=9\), for \(1\le\mu\le 4\).

Due to these facts, the class number relation
\(243\cdot h_3(B_j)=(U_j:V_j)\cdot 9\cdot 9\cdot 1\cdot 3\)
for \(j\in\lbrace 5,6,8,9\rbrace\)
implies that there are precisely four \textbf{tame} bicyclic bicubic fields,
\(B_5=k_1k_3k_p\tilde{k}_{qr}\),
\(B_6=k_1k_4k_q\tilde{k}_{pr}\),
\(B_8=k_2k_4k_pk_{qr}\),
\(B_9=k_2k_3k_qk_{pr}\),
satisfying \(9\mid h_3(B_j)=(U_j:V_j)\),
for each \(j\in\lbrace 5,6,8,9\rbrace\),
and so we must have the matrix ranks \(r_5=r_6=r_8=r_9=2\)
with indices \((U_j:V_j)\in\lbrace 9,27\rbrace\).

In contrast, each of the six \textbf{wild} bicyclic bicubic fields,
\(B_1=k_1k_{pq}k_{pr}k_{qr}\),
\(B_2=k_2k_{pq}\tilde{k}_{pr}\tilde{k}_{qr}\),
\(B_{10}=k_3k_4k_rk_{pq}\),
\(B_3=k_3\tilde{k}_{pq}\tilde{k}_{pr}k_{qr}\),
\(B_4=k_4\tilde{k}_{pq}k_{pr}\tilde{k}_{qr}\),
\(B_7=k_1k_2k_r\tilde{k}_{pq}\),
with \(h_3(B_j)>(U_j:V_j)\),
either contains \(k_{pq}\) or \(\tilde{k}_{pq}\).
The class number relation
\eqref{eqn:ParryFormula}
implies
\[
243\cdot h_3(B_j)=(U_j:V_j)\cdot
\begin{cases}
9\cdot h_3(k_{pq})\cdot 3\cdot 3         & \text{ for } j=1,2, \\
9\cdot 9\cdot 1\cdot h_3(k_{pq})         & \text{ for } j=10, \\
9\cdot h_3(\tilde{k}_{pq})\cdot 3\cdot 3 & \text{ for } j=3,4, \\
9\cdot 9\cdot 1\cdot h_3(\tilde{k}_{pq}) & \text{ for } j=7.
\end{cases}
\]
Summarized, in dependence on the index \(I_j:=(U_j:V_j)\) of subfield units and the rank \(r_j\),
\begin{equation}
\label{eqn:Cat3Gph5To9}
h_3(B_j)=
\begin{cases}
h_3(k_{pq})                & \text{ for } j=1,2,10,\ I_j=3,\ r_j=3, \\
3\cdot h_3(k_{pq})         & \text{ for } j=1,2,10,\ I_j=9,\ r_j=2, \\
9\cdot h_3(k_{pq})         & \text{ for } j=1,2,10,\ I_j=27,\ r_j=2, \\
h_3(\tilde{k}_{pq})        & \text{ for } j=3,4,7,\ I_j=3,\ r_j=3, \\
3\cdot h_3(\tilde{k}_{pq}) & \text{ for } j=3,4,7,\ I_j=9,\ r_j=2, \\
9\cdot h_3(\tilde{k}_{pq}) & \text{ for } j=3,4,7,\ I_j=27,\ r_j=2, \\
\end{cases}
\end{equation}
with \(h_3(k_{pq})=h_3(\tilde{k}_{pq})=9\) in the regular situation, and
\(h_3(k_{pq}),h_3(\tilde{k}_{pq})\ge 27\) in the singular or super-singular situation.
Formula
\eqref{eqn:Cat3Gph5To9}
supplements Corollary
\ref{cor:ClassNumber}
in the case \(p\leftrightarrow q\).


\begin{lemma}
\label{lem:Cat3Gph5To9}
\textbf{(3-class ranks of components.)}
All four components \(k_\mu\), \(1\le\mu\le 4\), of the quartet
have elementary bicyclic \(3\)-class group
\(\mathrm{Cl}_3(k_\mu)\simeq (3,3)\).
The condition \(9\mid h_3(B_j)=(U_j:V_j)\in\lbrace 9,27\rbrace\), \(r_j=2\),
is satisfied for \(j\in\lbrace 5,6,8,9\rbrace\),
the so-called \textbf{tame} extensions.
\end{lemma}

\begin{proof}
This is a consequence of the definition
of Graph 5--9 in Category III
and the rank distribution in Theorem
\ref{thm:3ClassRank}.
The fields \(B_j\) with \(j\in\lbrace 5,6,8,9\rbrace\)
neither contain \(k_{pq}\) nor \(\tilde{k}_{pq}\).
\end{proof}

\noindent
All computations for examples in the following subsections were performed with Magma
\cite{BCP1997,BCFS2023,MAGMA2023}.


\subsection{Category III, Graph 5}
\label{ss:Cat3Gph5}

\noindent
In this section,
the combined cubic residue symbol
of three prime(power)s dividing the conductor \(c=pqr\)
is assumed to be
\(\lbrack p,q,r\rbrack_3=\lbrace p\leftrightarrow q; r\rbrace\).


Since there are no trivial cubic residue symbols
among the three prime(power) divisors \(p,q,r\)
of the conductor \(c=pqr\),
except \(p\leftrightarrow q\) with overall assumption
\eqref{eqn:PFCat3Gph5To9},
the principal factors of the subfields
\(k\in\lbrace k_{pr},\tilde{k}_{pr},k_{qr},\tilde{k}_{qr}\rbrace\)
with \(t=2\)
of the absolute genus field \(k^\ast\)
must be divisible by both relevant primes,
and we can use the general approach
\begin{equation}
\label{eqn:ElEsCat3Gph5}
\begin{aligned}
A(k_{pr})=p^\ell r, & \quad A(\tilde{k}_{pr})=p^{-\ell}r, \text{ and } \\
A(k_{qr})=q^sr, & \quad A(\tilde{k}_{qr})=q^{-s}r,
\end{aligned}
\end{equation}
with \(\ell,s\in\lbrace -1,1\rbrace\),
identifying \(-1\equiv 2\,(\mathrm{mod}\,3)\),
since it is easier to manage:
\(\ell^2=s^2=1\).


\begin{lemma}
\label{lem:Cat3Gph5}
In dependence on the 
\textbf{decisive principal factors} in Equation
\eqref{eqn:ElEsCat3Gph5},
the \textbf{principal factors} of the quartet \((k_\mu)_{\mu=1}^4\)
sharing common conductor \(c=pqr\) with Graph \(\mathrm{III}.5\) are given by
\begin{equation}
\label{eqn:PFCat3Gph5}
\begin{aligned}
A(k_1)=pqr^2,\quad A(k_2)=pqr,\quad A(k_3)=pq^2r,\quad A(k_4)=p^2qr \quad & \text{ if } (\ell,s)=(1,1), \\
A(k_1)=p^2qr,\quad A(k_2)=pq^2r,\quad A(k_3)=pqr,\quad A(k_4)=pqr^2 \quad & \text{ if } (\ell,s)=(1,2), \\
A(k_1)=pq^2r,\quad A(k_2)=p^2qr,\quad A(k_3)=pqr^2,\quad A(k_4)=pqr \quad & \text{ if } (\ell,s)=(2,1), \\
A(k_1)=pqr,\quad A(k_2)=pqr^2,\quad A(k_3)=p^2qr,\quad A(k_4)=pq^2r \quad & \text{ if } (\ell,s)=(2,2), \\
A(k_1)=p^\ell q^sr^{-1},\quad A(k_2)=p^\ell q^sr,\quad A(k_3)=p^\ell q^{-s}r,\quad A(k_4)=p^{-\ell}q^sr \quad & \text{ generally}.
\end{aligned}
\end{equation}
\end{lemma}

\begin{proof}
We implement the general approach
\eqref{eqn:ElEsCat3Gph5}.
From the ranks \(r_j=2\) for \(j=5,6,8,9\),
there arise constraints for the exponents in the proposal
\(A(k_\mu)=p^{x_\mu}q^{y_\mu}r^{z_\mu}\), \(1\le\mu\le 4\),
with the aid of principal factor matrices.
For these \textit{tame} bicyclic bicubic fields \(B_j\),
\(j\in\lbrace 5,6,8,9\rbrace\),
the rank \(r_j\) is calculated with row operations
on the associated matrix \(M_j\): \\
\(M_5=
\begin{pmatrix}
x_1 & y_1 & z_1 \\
x_3 & y_3 & z_3 \\
1 & 0 & 0 \\
0 & -s & 1 \\
\end{pmatrix}\),
\(M_6=
\begin{pmatrix}
x_1 & y_1 & z_1 \\
x_4 & y_4 & z_4 \\
0 & 1 & 0 \\
-\ell & 0 & 1 \\
\end{pmatrix}\),
\(M_8=
\begin{pmatrix}
x_2 & y_2 & z_2 \\
x_4 & y_4 & z_4 \\
1 & 0 & 0 \\
0 & s & 1 \\
\end{pmatrix}\),
\(M_9=
\begin{pmatrix}
x_2 & y_2 & z_2 \\
x_3 & y_3 & z_3 \\
0 & 1 & 0 \\
\ell & 0 & 1 \\
\end{pmatrix}\).

For \(B_5=k_1k_3k_p\tilde{k}_{qr}\),
\(M_5\)
leads to the decisive pivot elements \(z_1+sy_1\) and \(z_3+sy_3\)
in the last column,
similarly,
for \(B_6=k_1k_4k_q\tilde{k}_{pr}\),
\(M_6\)
leads to \(z_1+\ell x_1\) and \(z_4+\ell x_4\),
similarly,
for \(B_8=k_2k_4k_pk_{qr}\),
\(M_8\)
leads to \(z_2-sy_2\) and \(z_4-sy_4\),
and similarly,
for \(B_9=k_2k_3k_qk_{pr}\), 
\(M_9\)
leads to \(z_2-\ell x_2\) and \(z_3-\ell x_3\).
So, \(r_5=r_6=r_8=r_9=2\)
implies 
\(\ell x_1\equiv sy_1\equiv -z_1\),
\(\ell x_2\equiv sy_2\equiv z_2\),
\(\ell x_3\equiv -sy_3\equiv z_3\),
\(-\ell x_4\equiv sy_4\equiv z_4\),
and consequently
\eqref{eqn:PFCat3Gph5}.
\end{proof}


\begin{proposition}
\label{prp:Cat3Gph5}
\textbf{(Quartet with \(3\)-rank two for \(\mathrm{III}.5\).)}
Let \((k_\mu)_{\mu=1}^4\) be a quartet
with common conductor \(c=pqr\),
whose combined cubic residue symbol belongs to Graph \(5\) of Category \(\mathrm{III}\).
Then the ranks of principal factor matrices of \textbf{tame} bicyclic bicubic fields are
\(r_j=2\) for \(j=5,6,8,9\).
In terms of exponents of primes in four variable principal factors,
\(A(k_{pq})=p^mq^n\), \(A(\tilde{k}_{pq})=p^{\tilde{m}}q^{\tilde{n}}\), from
\eqref{eqn:PFCat3Gph5To9}, and
\(A(k_{pr})=p^\ell r\), \(A(k_{qr})=q^sr\), from 
\eqref{eqn:ElEsCat3Gph5},
the ranks of principal factor matrices of \textbf{wild} bicyclic bicubic fields are given by
\begin{equation}
\label{eqn:Cat3Gph5RanksFake}
r_1=r_2=r_{10}=3 \text{ iff } \ell m\not\equiv -sn\,(\mathrm{mod}\,3)
\quad
\text{ and }
\quad
r_3=r_4=r_7=3 \text{ iff } \ell\tilde{m}\not\equiv s\tilde{n}\,(\mathrm{mod}\,3).
\end{equation}
\end{proposition}

\begin{proof}
Up to this point, the parameters
\(m,n,\tilde{m},\tilde{n}\) did not come into the play yet.
They decide about the rank \(r_j\) of
the associated principal factor matrices \(M_j\)
of the \textit{wild} bicyclic bicubic fields \(B_j\),
\(j\in\lbrace 1,2,3,4,7,10\rbrace\).
Hence, we perform row operations
on these matrices: \\
\(M_1=
\begin{pmatrix}
\ell & s & -1 \\
m & n & 0 \\
\ell & 0 & 1 \\
0 & s & 1 \\
\end{pmatrix}\),
\(M_2=
\begin{pmatrix}
\ell & s & 1 \\
m & n & 0 \\
-\ell & 0 & 1 \\
0 & -s & 1 \\
\end{pmatrix}\),
\(M_{10}=
\begin{pmatrix}
\ell & -s & 1 \\
-\ell & s & 1 \\
0 & 0 & 1 \\
m & n & 0 \\
\end{pmatrix}\).

For \(B_1=k_1k_{pq}k_{pr}k_{qr}\),
\(M_1\)
leads to the decisive pivot element \(-\ell m-sn\)
in the last column,
similarly,
for \(B_2=k_2k_{pq}\tilde{k}_{pr}\tilde{k}_{qr}\),
\(M_2\)
leads to \(\ell m+sn\),
and similarly,
for \(B_{10}=k_3k_4k_rk_{pq}\),
\(M_{10}\)
leads to \(n+\ell sm\)
in the middle column.
So, \(r_1=r_2=r_{10}=3\) iff \(\ell m\not\equiv -sn\),
by viewing the pivot elements modulo \(3\).
Next we consider: \\
\(M_3=
\begin{pmatrix}
\ell & -s & 1 \\
\tilde{m} & \tilde{n} & 0 \\
-\ell & 0 & 1 \\
0 & s & 2 \\
\end{pmatrix}\),
\(M_4=
\begin{pmatrix}
-\ell & s & 1 \\
\tilde{m} & \tilde{n} & 0 \\
\ell & 0 & 1 \\
0 & -s & 1 \\
\end{pmatrix}\),
\(M_7=
\begin{pmatrix}
\ell & s & -1 \\
\ell & s & 1 \\
0 & 0 & 1 \\
\tilde{m} & \tilde{n} & 0 \\
\end{pmatrix}\).

For \(B_3=k_3\tilde{k}_{pq}\tilde{k}_{pr}k_{qr}\), 
\(M_3\)
leads to \(\ell\tilde{m}-s\tilde{n}\),
similarly,
for \(B_4=k_4\tilde{k}_{pq}k_{pr}\tilde{k}_{qr}\),
\(M_4\)
leads to \(-\ell\tilde{m}+s\tilde{n}\),
and similarly,
for \(B_7=k_1k_2k_r\tilde{k}_{pq}\),
\(M_7\)
leads to \(\tilde{n}-\ell s\tilde{m}\)
in the middle column.
So, \(r_3=r_4=r_7=3\) iff \(\ell\tilde{m}\not\equiv s\tilde{n}\).
\end{proof}

\noindent
In Ayadi's Thesis
\cite[p. 80]{Ay1995},
only the special case \(\ell=s=-1\) is elaborated.
As mentioned above already, the condition
\((m,n)=(\tilde{m},\tilde{n})\) is warranted
in the \textit{regular} situation \(9\Vert h(k_{pq})\).
In any situation,
at least one of the following two rank equations,
which imply a \textit{total} transfer kernel,
is satisfied --- in many cases even both simultaneously:
\begin{equation}
\label{eqn:WildRank3}
\begin{aligned}
r_1=r_2=r_{10}=3 \text{ for } & (m,n)\in\lbrace (0,1),(1,0)\rbrace, \\
r_3=r_4=r_7=3    \text{ for } & (\tilde{m},\tilde{n})\in\lbrace (0,1),(1,0)\rbrace.
\end{aligned}
\end{equation}


\begin{theorem}
\label{thm:Cat3Gph5}
\textbf{(Second \(3\)-class groups for \(\mathrm{III}.5\))}
There are several
\textbf{minimal transfer kernel types} (mTKT) \(\varkappa_0\) of \(k_\mu\), \(1\le\mu\le 4\),
and other possible capitulation types in ascending order
\(\varkappa_0<\varkappa<\varkappa^\prime<\varkappa^{\prime\prime}<\varkappa^{\prime\prime\prime}\),
either
\(\varkappa_0=(2134)\), type \(\mathrm{G}.16\),
\(\varkappa=(2130)\), type \(\mathrm{d}.23\),
or
\(\varkappa_0=(2143)\), type \(\mathrm{G}.19\),
\(\varkappa=(2140)\), type \(\mathrm{d}.25\),
ending in
\(\varkappa^\prime=(2100)\), type \(\mathrm{b}.10\),
\(\varkappa^{\prime\prime}=(2000)\), type \(\mathrm{a}.3^\ast\) or \(\mathrm{a}.3\),
or \(\varkappa^{\prime\prime}=(0004)\), type \(\mathrm{a}.2\),
and the maximal \(\varkappa^{\prime\prime\prime}=(0000)\), type \(\mathrm{a}.1\).

In terms of the counter
\(\mathcal{N}^\ast:=\#\lbrace 1\le j\le 10\mid (U_j:V_j)=27\rbrace\),
of \textbf{maximal indices of subfield units}
\(I_j=(U_j:V_j)\)
for all ten bicyclic bicubic fields \(B_j<k^\ast\) with conductor \(c=pqr\),
the second \(3\)-class groups
\(\mathfrak{M}=\mathrm{Gal}(\mathrm{F}_3^2(k_\mu)/k_\mu)\)
are given in the following way
as uniform or non-uniform quartets,
with abbreviation
\(P_7:=\langle 2187,64\rangle\): 
\begin{equation}
\label{eqn:Cat3Gph5}
\mathfrak{M}=
\begin{cases}
\langle 243,28..30\rangle^4,\ \alpha=\lbrack 21,11,11,11\rbrack,\ \varkappa=(0000) & \text{ if } \mathcal{N}^\ast=0, \\
\langle 243,27\rangle,\ \alpha=\lbrack 11,11,11,22\rbrack,\ \varkappa=(0004)       & \text{ once if } \mathcal{N}^\ast=1, \\
\langle 243,28..30\rangle^3,\ \alpha=\lbrack 21,11,11,11\rbrack,\ \varkappa=(0000) & \text{ thrice if } \mathcal{N}^\ast=1, \\
\langle 81,7\rangle^4,\ \alpha=\lbrack 111,11,11,11\rbrack,\ \varkappa=(2000)      & \text{ if } \mathcal{N}^\ast=2, \\
\langle 243,27\rangle^2,\ \alpha=\lbrack 11,11,11,22\rbrack,\ \varkappa=(0004)     & \text{ twice if } \mathcal{N}^\ast=3, \\
\langle 243,25\rangle^2,\ \alpha=\lbrack 22,11,11,11\rbrack,\ \varkappa=(2000)     & \text{ twice if } \mathcal{N}^\ast=3, \\
\langle 729,34..39\rangle^4,\ \alpha=\lbrack 111,111,21,21\rbrack,\ \varkappa=(2100) &\text{ if } \mathcal{N}^\ast=4, \\
\langle 2187,250\rangle^2,\ \alpha=\lbrack 111,111,32,21\rbrack,\ \varkappa=(2130)   & \text{ twice if } \mathcal{N}^\ast=7, \\
\langle 2187,251\vert 252\rangle^2,\ \alpha=\lbrack 111,111,32,21\rbrack,\ \varkappa=(2140) & \text{ twice if } \mathcal{N}^\ast=7, \\
(P_7-\#2;40\vert 48)^2,\ \alpha=\lbrack 111,111,32,32\rbrack,\ \varkappa=(2134)     & \text{ twice if } \mathcal{N}^\ast=10, \\
(P_7-\#2;42\vert 45\vert 49)^2,\ \alpha=\lbrack 111,111,32,32\rbrack,\ \varkappa=(2143) & \text{ twice if } \mathcal{N}^\ast=10.
\end{cases}
\end{equation}
The leading six rows concern the \textbf{regular} situation \(\mathrm{Cl}_3(k_{pq})\simeq (3,3)\).
In particular, the condition \((m,n)\in\lbrace (0,1),(1,0)\rbrace\) for \(\langle 81,7\rangle^4\) is equivalent to 
the extra special group \(\mathrm{Gal}(\mathrm{F}_3^2(k_{pq})/k_{pq})\simeq\langle 27,4\rangle\), whereas
\(\mathrm{Gal}(\mathrm{F}_3^2(k_{pq})/k_{pq})\simeq\langle 9,2\rangle\)
is abelian for all other pairs \((m,n)\).
The trailing rows concern the \textbf{(super-)singular} situation with
\(h_3(k_{pq})=h_3(\tilde{k}_{pq})=27\).
With exception of the trailing rows,
the \(3\)-class field tower has length \(\ell_3(K_\mu)=2\) and group
\(\mathfrak{G}=\mathrm{Gal}(\mathrm{F}_3^\infty(k_\mu)/k_\mu)\simeq\mathfrak{M}\).
\end{theorem}

\begin{proof}
Let \(\mathfrak{p},\mathfrak{q},\mathfrak{r}\)
be the prime ideals of \(k_\mu\)
over \(p,q,r\).

Since \(p\) splits in \(k_q\), it also splits in
\(B_6=k_1k_4k_q\tilde{k}_{pr}\) and
\(B_9=k_2k_3k_qk_{pr}\).

Since \(q\) splits in \(k_p\), it also splits in
\(B_5=k_1k_3k_p\tilde{k}_{qr}\) and
\(B_8=k_2k_4k_pk_{qr}\). By Corollary
\ref{cor:Capitulation},

since \(\mathfrak{q}\) is principal ideal in \(k_q\),
the class \(\lbrack\mathfrak{q}\rbrack\) capitulates in
\(B_6=k_1k_4k_q\tilde{k}_{pr}\) and
\(B_9=k_2k_3k_qk_{pr}\);

since \(\mathfrak{r}\) is principal ideal in \(k_r\),
the class \(\lbrack\mathfrak{r}\rbrack\) capitulates in
\(B_7=k_1k_2k_r\tilde{k}_{pq}\) and
\(B_{10}=k_3k_4k_rk_{pq}\).

Since \(\mathfrak{qr}\) is principal ideal in \(\tilde{k}_{qr}\),
the class \(\lbrack\mathfrak{qr}\rbrack\) capitulates in
\(B_2=k_2k_{pq}\tilde{k}_{pr}\tilde{k}_{qr}\),
\(B_4=k_4\tilde{k}_{pq}k_{pr}\tilde{k}_{qr}\), and
\(B_5=k_1k_3k_p\tilde{k}_{qr}\),
by Proposition
\ref{prp:Bicyc}.
Since \(A(k_1)=pqr\) and \(A(k_3)=p^2qr\), 
\(\lbrack\mathfrak{qr}\rbrack\) generates the same subgroup as \(\lbrack\mathfrak{p}\rbrack\)
in \(\ker(T_{B_5/k_\mu})\), \(\mu=1,3\).

Since \(\mathfrak{qr}^2\) is principal ideal in \(k_{qr}\),
the class \(\lbrack\mathfrak{qr}^2\rbrack\) capitulates in
\(B_1=k_1k_{pq}k_{pr}k_{qr}\),
\(B_3=k_3\tilde{k}_{pq}\tilde{k}_{pr}k_{qr}\), and
\(B_8=k_2k_4k_pk_{qr}\),
by Proposition
\ref{prp:Bicyc}.
Since \(A(k_2)=pqr^2\) and \(A(k_4)=pq^2r\),
\(\lbrack\mathfrak{qr}^2\rbrack\) generates the same subgroup as \(\lbrack\mathfrak{p}\rbrack\)
in \(\ker(T_{B_8/k_\mu})\), \(\mu=2,4\).

The \(3\)-class group of \(k_\mu\) is always
\(\mathrm{Cl}_3=\langle\lbrack\mathfrak{q}\rbrack,\lbrack\mathfrak{r}\rbrack\rangle\).
It contains the \textbf{norm class groups} of \(B_j>k_\mu\) as subgroups of index \(3\):
\(N_{B_j/k_\mu}\mathrm{Cl}_3(B_j)\) is always generated by \(\lbrack\mathfrak{q}\rbrack\) for \(j=5,8\),
due to the above mentioned splitting of \(q\).
See also Table
\ref{tbl:UniCat3Gph5}.

We recall that equality \((m,n)=(\tilde{m},\tilde{n})\)
is warranted for the \textbf{regular} situation \(\mathrm{Cl}_3(k_{pq})\simeq (3,3)\),
and there is an equivalence involving the counter \(\mathcal{P}\) in Theorem
\ref{thm:TwoPrimeCond}:
\(\mathrm{Gal}(\mathrm{F}_3^2(k_{pq})/k_{pq})\simeq\langle 27,4\rangle\) iff
\(\mathcal{P}=1\) iff (either \(m=0\) or \(n=0\)) iff \((m,n)\in\lbrace (0,1),(1,0)\rbrace\).
Let \(I_j:=(U_j:V_j)\).

\(\mathcal{N}^\ast=0\) implies wild ranks
\(r_j=2\) and \(I_j=9\), \(h_3(B_j)=3\cdot h_3(k_{pq})=3\cdot 9=27\) for \(j\in\lbrace 1,2,10\rbrace\),
but \(r_j=3\), \(I_j=3\), \(h_3(B_j)=h_3(\tilde{k}_{pq})=9\) for \(j\in\lbrace 3,4,7\rbrace\),
according to Equation
\eqref{eqn:Cat3Gph5To9},
and tame indices \(I_j=9\) for all \(j\in\lbrace 5,6,8,9\rbrace\).
The uniform minimal indices \(I_j\) of subgroup units correspond to
maximal norm unit indices \((U(k_\mu):N_{B_j/k_\mu}U(B_j))=3\)
and thus to total capitulations whenever \(k_\mu<B_j\) is a subfield
for \(1\le j\le 10\), \(1\le\mu\le 4\).
According to Theorem
\ref{thm:ClassGroup33}
and Corollary
\ref{cor:ClassGroup33},
the resulting abelian type invariants \(\alpha=\lbrack 21,11,11,11\rbrack\)
and transfer kernel type \(\varkappa=(0000)\),
that is the \textbf{Artin pattern} \((\alpha,\varkappa)\),
identify three possible groups \(\mathfrak{M}\simeq\langle 243,28..30\rangle\),
since \(\langle 81,9\rangle\) must be cancelled, due to wrong second layer \(\alpha_2\).

An exception arises for \(\mathcal{N}^\ast=1\),
which causes non-uniformity with \(I_1=27\), \(h_3(B_1)=9\cdot h_3(k_{pq})=9\cdot 9=81\),
as opposed to the remaining \(I_2=I_{10}=9\).
(Everything else is like \(\mathcal{N}^\ast=0\).)
Thus \((U(k_1):N_{B_1/k_1}U(B_1))=1\), and here we have a fixed point capitulation,
\(\ker(T_{B_1/k_1})=\langle\lbrack\mathfrak{qr}^2\rbrack\rangle\).
The corresponding abelian type invariants \(\alpha=\lbrack 11,11,11,22\rbrack\)
and transfer kernel type \(\varkappa=(0004)\)
uniquely identify the group \(\langle 243,27\rangle\) for \(\mu=1\).
The remaining three groups are \(\langle 243,28..30\rangle^3\).

For \(\mathcal{N}^\ast=2\) and \((m,n)\in\lbrace (0,1),(1,0)\rbrace\), the relations
\(m\not\equiv -n\) and \(m\not\equiv n\) imply
\(r_j=3\), \(I_j=3\), \(h_3(B_j)=h_3(\tilde{k}_{pq})=9\)
for all wild bicyclic cicubic fields \(B_j\),
\(j\in\lbrace 1,2,3,4,7,10\rbrace\).
For \(j\in\lbrace 5,8\rbrace\), we have \(I_j=9\),
but for \(j\in\lbrace 6,9\rbrace\), the maximal index \(I_j=27\)
is attained and enables an elementary tricyclic \(3\)-class group
\(\mathrm{Cl}_3(B_j)=\langle\mathfrak{P}_1,\mathfrak{P}_2,\mathfrak{P}_3\rangle\simeq (3,3,3)\)
generated by the prime ideals lying over
\(\mathfrak{p}\mathcal{O}_{B_j}=\mathfrak{P}_1\cdot\mathfrak{P}_2\cdot\mathfrak{P}_3\).
Here, we have a non-fixed point capitulation
\(\ker(T_{B_j/k_\mu})=\langle\lbrack\mathfrak{q}\rbrack\rangle\).
The transposition is hidden by total capitulation in \(B_5\) and \(B_8\)
with norm class group generated by \(\lbrack\mathfrak{q}\rbrack\).
The abelian type invariants \(\alpha=\lbrack 111,11,11,11\rbrack\)
and transfer kernel type \(\varkappa=(2000)\)
uniquely identify the group \(\langle 81,7\rangle\simeq\mathrm{Syl}_3(A_9)\) for \(\mu=1\).

\(\mathcal{N}^\ast=3\) implies \(r_j=3\) and thus wild indices \(I_j=3\), \(h_3(B_j)=h_3(\tilde{k}_{pq})=9\) for \(j\in\lbrace 1,2,10\rbrace\).
We also have tame indices \(I_j=9\) for \(j\in\lbrace 5,6,8,9\rbrace\).
However,
we have  \(r_j=2\) and remaining wild indices \(I_j=27\), \(h_3(B_j)=9\cdot h_3(k_{pq})=9\cdot 9=81\) for \(j\in\lbrace 3,4,7\rbrace\),
according to Equation
\eqref{eqn:Cat3Gph5To9}.
There arises a fixed point capitulation,
\(\ker(T_{B_7/k_\mu})=\langle\lbrack\mathfrak{r}\rbrack\rangle\)
and, non-uniformly, a non-fixed point capitulation,
\(\ker(T_{B_j/k_\mu})=\langle\lbrack\mathfrak{p}\rbrack\rangle\)
for \(j=3,4\) with norm class groups also generated by \(\lbrack\mathfrak{r}\rbrack\).
The corresponding abelian type invariants \(\alpha=\lbrack 11,11,11,22\rbrack\)
and transfer kernel type \(\varkappa=(0004)\), respectively \(\varkappa=(0003)\),
uniquely identify the two groups \(\langle 243,27\rangle^2\),
respectively the remaining two groups \(\langle 243,25\rangle^2\).

For \(\mathcal{N}^\ast=4\) and the simplest singular or super-singular situation
with \(h_3(k_{pq})=h_3(\tilde{k}_{pq})=27\),
\(\mathcal{P}=1\) implies \(r_j=3\), \(I_j=3\), \(h_3(B_j)=27\) for all wild \(j\in\lbrace 1,2,3,4,7,10\rbrace\),
and uniformly \(h_3(B_j)=I_j=27\) for all tame \(j\in\lbrace 5,6,8,9\rbrace\).
The latter correspond to elementary tricyclic \(3\)-class groups
\(\mathrm{Cl}_3(B_j)=\langle\mathfrak{Q}_1,\mathfrak{Q}_2,\mathfrak{Q}_3\rangle\simeq (3,3,3)\)
generated by the prime ideals lying over
\(\mathfrak{q}\mathcal{O}_{B_j}=\mathfrak{Q}_1\cdot\mathfrak{Q}_2\cdot\mathfrak{Q}_3\)
for \(j=5,8\), and
\(\mathrm{Cl}_3(B_j)=\langle\mathfrak{P}_1,\mathfrak{P}_2,\mathfrak{P}_3\rangle\simeq (3,3,3)\)
generated by the prime ideals lying over
\(\mathfrak{p}\mathcal{O}_{B_j}=\mathfrak{P}_1\cdot\mathfrak{P}_2\cdot\mathfrak{P}_3\)
for \(j=6,9\).
Here, we have a non-fixed point capitulation
\(\ker(T_{B_j/k_\mu})=\langle\lbrack\mathfrak{p}\rbrack\rangle\)
for \(j=5,8\), and
\(\ker(T_{B_j/k_\mu})=\langle\lbrack\mathfrak{q}\rbrack\rangle\)
for \(j=6,9\).
The transposition is not hidden by total capitulation
and characteristic for uniform transfer kernel type \(\mathrm{b}.10\).
According to Theorem
\ref{thm:ClassGroup33}
and Corollary
\ref{cor:ClassGroup33},
the abelian type invariants \(\alpha=\lbrack 111,111,21,21\rbrack\)
and the transfer kernel type \(\varkappa=(2100)\),
identify six possible groups \(\mathfrak{M}\simeq\langle 729,34..39\rangle\).

For \(\mathcal{N}^\ast=7\),
only three wild indices \(r_j=3\), \(I_j=3\), \(h_3(B_j)=27\) for \(j=3,4,7\)
are not maximal. The TKTs are not uniform,
\(\varkappa=(2130)\), type \(\mathrm{d}.23\), twice
and
\(\varkappa=(2140)\), type \(\mathrm{d}.25\), twice.

For \(\mathcal{N}^\ast=10\),
all tame and wild indices are maximal \(I_j=27\) for \(1\le j\le 10\),
which implies non-uniform minimal TKTs
\(\varkappa_0=(2134)\), type \(\mathrm{G}.16\), twice
and
\(\varkappa_0=(2143)\), type \(\mathrm{G}.19\), twice.
\end{proof}


\begin{corollary}
\label{cor:UniCat3Gph5}
\textbf{(Non-uniformity of the quartet for \(\mathrm{III}.5\).)}
For \(\mathcal{N}^\ast=1\),
only a sub-triplet of the quartet
shares a common capitulation type
\(\varkappa(k_\mu)\),
abelian type invariants
\(\alpha(k_\mu)\),
and second \(3\)-class group
\(\mathfrak{M}=\mathrm{Gal}(\mathrm{F}_3^2(k_\mu)/k_\mu)\).
The invariants of the fourth component
\textbf{differ}.
For \(\mathcal{N}^\ast\in\lbrace 3,7,10\rbrace\),
only two pairs of components of the quartet
share a common capitulation type
\(\varkappa(k_\mu)\),
and second \(3\)-class group
\(\mathfrak{M}=\mathrm{Gal}(\mathrm{F}_3^2(k_\mu)/k_\mu)\),
whereas the abelian type invariants
\(\alpha(k_\mu)\) are uniform.
However, the four components \textbf{agree} in all situations
with even \(\mathcal{N}^\ast\in\lbrace 0,2,4\rbrace\).
\end{corollary}

\begin{proof}
This is an immediate consequence of Theorem
\ref{thm:Cat3Gph5}
and Table
\ref{tbl:UniCat3Gph5}.
\end{proof}


\noindent
In terms of capitulation targets in Corollary
\ref{cor:Three},
Theorem
\ref{thm:Cat3Gph5}
and parts of its proof
are now summarized in Table
\ref{tbl:UniCat3Gph5}
with transpositions in \textbf{bold} font.

\renewcommand{\arraystretch}{1.1}

\begin{table}[ht]
\caption{Norm class groups and minimal transfer kernels for Graph \(\mathrm{III}.5\)}
\label{tbl:UniCat3Gph5}
\begin{center}
{\normalsize
\begin{tabular}{|c||c|c|c|c||c|c|c|c||c|c|c|c||c|c|c|c|}
\hline
 Base         & \multicolumn{4}{c||}{\(k_1\)} & \multicolumn{4}{c||}{\(k_2\)} & \multicolumn{4}{c||}{\(k_3\)} & \multicolumn{4}{c|}{\(k_4\)} \\
\hline
 Ext          & \(B_1\) & \(B_5\) & \(B_6\) & \(B_7\) & \(B_2\) & \(B_7\) & \(B_8\) & \(B_9\) & \(B_3\) & \(B_5\) & \(B_9\) & \(B_{10}\) & \(B_4\) & \(B_6\) & \(B_8\) & \(B_{10}\) \\
\hline
 NCG          & \(\mathfrak{qr}^2\) & \(\mathfrak{q}\) & \(\mathfrak{qr}\) & \(\mathfrak{r}\)    & \(\mathfrak{qr}\) & \(\mathfrak{r}\)    & \(\mathfrak{q}\) & \(\mathfrak{qr}^2\)
              & \(\mathfrak{r}\)    & \(\mathfrak{q}\) & \(\mathfrak{qr}\) & \(\mathfrak{qr}^2\) & \(\mathfrak{r}\)  & \(\mathfrak{qr}^2\) & \(\mathfrak{q}\) & \(\mathfrak{qr}\) \\
 TK           & \(\mathfrak{qr}^2\) & \(\mathfrak{qr}\) & \(\mathfrak{q}\) & \(\mathfrak{r}\) & \(\mathfrak{qr}\) & \(\mathfrak{r}\) & \(\mathfrak{qr}^2\) & \(\mathfrak{q}\)
              & \(\mathfrak{qr}^2\) & \(\mathfrak{qr}\) & \(\mathfrak{q}\) & \(\mathfrak{r}\) & \(\mathfrak{qr}\) & \(\mathfrak{q}\) & \(\mathfrak{qr}^2\) & \(\mathfrak{r}\) \\
\(\varkappa\) & \(1\) & \(\mathbf{3}\) & \(\mathbf{2}\) & \(4\) & \(1\) & \(2\) & \(\mathbf{4}\) & \(\mathbf{3}\) & \(\mathbf{4}\) & \(\mathbf{3}\) & \(\mathbf{2}\) & \(\mathbf{1}\) & \(\mathbf{4}\) & \(\mathbf{3}\) & \(\mathbf{2}\) & \(\mathbf{1}\) \\
\hline
\end{tabular}
}
\end{center}
\end{table}


\begin{example}
\label{exm:Cat3Gph5}
The prototypes for Graph \(\mathrm{III}.5\), i.e., the cases in Theorem
are four \textbf{regular} situations,
\(c=14\,049\) with \(\lbrace 7\leftrightarrow 223;9\rbrace\) and \(\mathcal{N}^\ast=0\);
\(c=17\,073\) with \(\lbrace 9\leftrightarrow 271;7\rbrace\) and \(\mathcal{N}^\ast=2\);
\(c=20\,367\) with \(\lbrace 9\leftrightarrow 73;31\rbrace\) and \(\mathcal{N}^\ast=1\);
\(c=21\,231\) with \(\lbrace 7\leftrightarrow 337;9\rbrace\) and \(\mathcal{N}^\ast=3\);
and the \textbf{singular} situation
\(c=42\,399\) with \(\lbrace 7\leftrightarrow 673;9\rbrace\) and \(\mathcal{N}^\ast=4\).
Here, we have distinct \((m,n)=(0,1)\), but \((\tilde{m},\tilde{n})=(1,0)\). 
There is also a \textbf{super-singular} prototype
\(c=48\,447\) with \(\lbrace 7\leftrightarrow 769;9\rbrace\) and \(\mathcal{N}^\ast=4\),
phenomenologically completely identical with the singular prototype,
except that \((m,n)=(\tilde{m},\tilde{n})=(0,1)\).
With \textbf{considerable statistic delay}, there appeared \(\mathcal{N}^\ast\in\lbrace 7,10\rbrace\).
\end{example}

In Table
\ref{tbl:ProtoCat3Gph5},
we summarize the prototypes of graph \(\mathrm{III}.5\).
Data comprises
ordinal number No.,
conductor \(c\) of \(k\),
combined cubic residue symbol \(\lbrack p,q,r\rbrack_3\),
regularity, resp. (super-)singularity, expressed by
\(3\)-valuation \(v^\ast=v_3(\#\mathrm{Cl}(k^\ast))\)
of class number of absolute \(3\)-genus field \(k^\ast\),
\(3\)-valuation \(v=v_3(\#\mathrm{Cl}(k_{pq}))\),
respectively \(\tilde{v}=v_3(\#\mathrm{Cl}(\tilde{k}_{pq}))\),
of class number of critical field \(k_{pq}\), respectively \(\tilde{k}_{pq}\),
critical exponents \(m,n\) in principal factor \(A(k_{pq})=p^mq^n\),
resp. \(\tilde{m},\tilde{n}\) in \(A(\tilde{k}_{pq})=p^{\tilde{m}}q^{\tilde{n}}\),
resp. \(\ell\) in \(A(k_{pr})=p^\ell r\),
resp. \(s\) in \(A(k_{qr})=q^s r\),
capitulation type of \(k\),
second \(3\)-class group
\(\mathfrak{M}=\mathrm{Gal}(\mathrm{F}_3^2(k)/k)\) of \(k\),
and length \(\ell_3(k)\) of \(3\)-class field tower of \(k\).
For abbreviation we put \\
\(P_7:=\langle 2187,64\rangle\), \quad
\(R_4^4:=P_7-\#2;54\), \quad 
\(R_5^4:=P_7-\#2;57\), \quad 
\(R_6^4:=P_7-\#2;59\), \\
\(S_4^4:=R_4^4-\#1;8-\#1;3\vert 7\), \quad
\(U_5^4:=R_5^4-\#1;1-\#1;3\vert 6\), \quad
\(V_6^4:=R_6^4-\#1;6-\#1;2\vert 6\). \\
See the tables and tree diagrams in
\cite[\S\S\ 11.3--11.4, pp. 108--116, Tbl. 4--5, Fig. 9--11]{Ma2018}.

\renewcommand{\arraystretch}{1.1}

\begin{table}[ht]
\caption{Prototypes for Graph III.5}
\label{tbl:ProtoCat3Gph5}
\begin{center}
{\scriptsize
\begin{tabular}{|r|r||c|c|c|c|c|c|c|c|c|c|c|}
\hline
 No. & \(c\) & \(p\leftrightarrow q,r\) & \(v^\ast\) & \(v\) & \(\tilde{v}\) & \(m,n\) & \(\tilde{m},\tilde{n}\) & \(\ell\) & \(s\) & capitulation type & \(\mathfrak{M}\) & \(\ell_3(k)\) \\
\hline
   1 & \(14\,049\) & \(7\leftrightarrow 223,9\) & \(1\) & \(2\) & \(2\) & \(2,1\) & \(2,1\) & \(1\) & \(1\) & \(\mathrm{a}.1\)              & \(\langle 243,28..30\rangle^4\)                       & \(=2\) \\
   2 & \(17\,073\) & \(9\leftrightarrow 271,7\) & \(2\) & \(2\) & \(2\) & \(0,1\) & \(0,1\) & \(1\) & \(1\) & \(\mathrm{a}.3^\ast\)         & \(\langle 81,7\rangle^4\)                             & \(=2\) \\
   3 & \(20\,367\) & \(9\leftrightarrow 73,31\) & \(1\) & \(2\) & \(2\) & \(2,1\) & \(2,1\) & \(2\) & \(2\) & \(\mathrm{a}.2,\mathrm{a}.1\) & \(\langle 243,27\rangle,\langle 243,28..30\rangle^3\) & \(=2\) \\
   4 & \(21\,231\) & \(7\leftrightarrow 337,9\) & \(1\) & \(2\) & \(2\) & \(1,1\) & \(1,1\) & \(1\) & \(1\) & \(\mathrm{a}.2,\mathrm{a}.3\) & \(\langle 243,27\rangle^2,\langle 243,25\rangle^2\)   & \(=2\) \\
  13 & \(42\,399\) & \(7\leftrightarrow 673,9\) & \(3\) & \(3\) & \(3\) & \(0,1\) & \(1,0\) & \(1\) & \(2\) & \(\mathrm{b}.10\)             & \(\langle 729,37..39\rangle^4\)                    & \(\ge 2\) \\
  16 & \(48\,447\) & \(7\leftrightarrow 769,9\) & \(4\) & \(3\) & \(3\) & \(0,1\) & \(0,1\) & \(1\) & \(1\) & \(\mathrm{b}.10\)             & \(\langle 729,37..39\rangle^4\)                    & \(\ge 2\) \\
  39 &\(100\,503\) &\(13\leftrightarrow 859,9\) & \(3\) & \(3\) & \(3\) & \(1,0\) & \(1,1\) & \(2\) & \(1\) & \(\mathrm{b}.10\)             & \(\langle 729,34..36\rangle^4\)                    & \(\ge 2\) \\
  67 &\(145\,593\) &\(7\leftrightarrow 2311,9\) & \(4\) & \(3\) & \(3\) & \(2,1\) & \(2,1\) & \(1\) & \(1\) & \(\mathrm{d}.23,\mathrm{d}.25\) & \(\langle 2187,250\rangle^2,\langle 2187,251\vert 252\rangle^2\) & \(\ge 2\) \\
 128 &\(256\,669\) &\(37\leftrightarrow 991,7\) & \(6\) & \(5\) & \(3\) & \(1,1\) & \(2,1\) & \(2\) & \(1\) & \(\mathrm{G}.16,\mathrm{G}.19\) & \((S_4^4)^2,(U_5^4\vert V_6^4)^2\)               & \(\ge 3\) \\
\hline
\end{tabular}
}
\end{center}
\end{table}


\subsection{Category III, Graph 6}
\label{ss:Cat3Gph6}

\noindent
Let the combined cubic residue symbol
of three primes dividing the conductor \(c=pqr\) be
\(\lbrack p,q,r\rbrack_3=\lbrace r\leftarrow p\leftrightarrow q\rbrace\).


\begin{proposition}
\label{prp:Cat3Gph6}
\textbf{(Quartet with \(3\)-rank two for \(\mathrm{III}.6\).)}
For fixed \(\mu\in\lbrace 1,2,3,4\rbrace\),
let \(\mathfrak{p},\mathfrak{q},\mathfrak{r}\) be the prime ideals of \(k_\mu\)
over \(p,q,r\), that is
\(p\mathcal{O}_{k_\mu}=\mathfrak{p}^3\),
\(q\mathcal{O}_{k_\mu}=\mathfrak{q}^3\),
\(r\mathcal{O}_{k_\mu}=\mathfrak{r}^3\),
then the \textbf{principal factor} of \(k_\mu\) is
\(A(k_\mu)=p\),
and the \(3\)-class group of \(k_\mu\) is,
\begin{equation}
\label{eqn:GenCat3Gph6}
\mathrm{Cl}_3(k_\mu)=
\langle\lbrack\mathfrak{q}\rbrack,\lbrack\mathfrak{r}\rbrack\rangle\simeq (3,3).
\end{equation}
The unramified cyclic cubic relative extensions of \(k_\mu\)
are among the absolutely bicyclic bicubic fields \(B_i\), \(1\le i\le 10\).
The \textbf{tame} extensions
with \(9\mid h_3(B_i)=(U_i:V_i)\in\lbrace 9,27\rbrace\)
are \(B_i\) with \(i=5,6,8,9\),
since they neither contain \(k_{pq}\) nor \(\tilde{k}_{pq}\).
For each \(\mu\), there are two tame extensions \(B_j/k_\mu\), \(B_\ell/k_\mu\)
with the following properties.
The first, \(B_j\) with \(j\in\lbrace 6,9\rbrace\),
has norm class group
\(N_{B_j/k_\mu}(\mathrm{Cl}_3(B_j))=\langle\lbrack\mathfrak{qr}^s\rbrack\rangle\)
with \(s\in\lbrace 1,2\rbrace\),
\textbf{cyclic} transfer kernel
\begin{equation}
\label{eqn:Trans1Cat3Gph6}
\ker(T_{B_j/k_\mu})=\langle\lbrack\mathfrak{q}\rbrack\rangle
\end{equation}
of order \(3\),
and \textbf{elementary tricyclic} \(3\)-class group
\(\mathrm{Cl}_3(B_j)=
\langle\lbrack\mathfrak{QR}^s\mathfrak{P}_1\rbrack,
\lbrack\mathfrak{QR}^s\mathfrak{P}_2\rbrack,
\lbrack\mathfrak{QR}^s\mathfrak{P}_3\rbrack\rangle\simeq (3,3,3)\),
generated by the classes of the prime ideals of \(B_j\) over
\(\mathfrak{p}\mathcal{O}_{B_j}=\mathfrak{P}_1\mathfrak{P}_2\mathfrak{P}_3\),
\(\mathfrak{q}\mathcal{O}_{B_j}=\mathfrak{Q}\),
\(\mathfrak{r}\mathcal{O}_{B_j}=\mathfrak{R}\).
The second, \(B_\ell\) with \(\ell\in\lbrace 5,8\rbrace\),
has norm class group
\(N_{B_\ell/k_\mu}(\mathrm{Cl}_3(B_\ell))=\langle\lbrack\mathfrak{q}\rbrack\rangle\),
transfer kernel
\begin{equation}
\label{eqn:Trans2Cat3Gph6}
\ker(T_{B_\ell/k_\mu})\ge\langle\lbrack\mathfrak{qr}^s\rbrack\rangle,
\end{equation}
and \(3\)-class group
\(\mathrm{Cl}_3(B_\ell)=
\langle\lbrack\mathfrak{Q}_1\rbrack,\lbrack\mathfrak{Q}_2\rbrack,\lbrack\mathfrak{Q}_3\rbrack\rangle\ge (3,3)\),
generated by the classes of the prime ideals of \(B_\ell\) over
\(\mathfrak{q}\mathcal{O}_{B_\ell}=\mathfrak{Q}_1\mathfrak{Q}_2\mathfrak{Q}_3\).
The pair \((j,\ell)\) forms a hidden or actual \textbf{transposition}
of the transfer kernel type \(\varkappa(k_\mu)\).
The remaining two \(B_i>k_\mu\), \(i\ne j\), \(i\ne\ell\),
have norm class group 
\(\langle\lbrack\mathfrak{r}\rbrack\rangle\), respectively
\(\langle\lbrack\mathfrak{q}^2\mathfrak{r}^s\rbrack\rangle\),
and transfer kernel
\[
\ker(T_{B_i/k_\mu})\ge\langle\lbrack\mathfrak{r}\rbrack\rangle, \text{ or }
\ge\langle\lbrack\mathfrak{q}^2\mathfrak{r}^s\rbrack\rangle,
\]
providing the option of either two possible \textbf{fixed points}
or a further \textbf{transposition}
in the transfer kernel type \(\varkappa(k_\mu)\).
In terms of \(n\) and \(\tilde{n}\) in
\(A(k_{pq})=p^mq^n\) and \(A(\tilde{k}_{pq})=p^{\tilde{m}}q^{\tilde{n}}\),
the ranks of the \textbf{wild} extensions are
\begin{equation}
\label{eqn:RankCat3Gph6}
r_1=r_2=r_{10}=3 \text{ iff } n\ne 0 \text{ iff } q\mid A(k_{pq})
\text{ and }
r_3=r_4=r_7=3 \text{ iff } \tilde{n}\ne 0 \text{ iff } q\mid A(\tilde{k}_{pq}).
\end{equation}
\end{proposition}

\begin{proof}
By Proposition
\ref{prp:Principal2},
principal factors are
\(A(k_{pr})=A(\tilde{k}_{pr})=p\),
since \(r\leftarrow p\).
Further, by Proposition
\ref{prp:Principal3},
\(A(k_\mu)=p\), for all \(1\le\mu\le 4\),
since \(p\) is universally repelling
\(r\leftarrow p\rightarrow q\).
Since \(\mathfrak{p}=\alpha\mathcal{O}_{k_\mu}\) is a principal ideal,
its class \(\lbrack\mathfrak{p}\rbrack=1\) is trivial,
whereas the classes
\(\lbrack\mathfrak{q}\rbrack,\lbrack\mathfrak{r}\rbrack\)
are non-trivial.

Assume the principal factors
\(A(k_{qr})=qr^2\) and \(A(\tilde{k}_{qr})=qr\).
The parameters \(m,n,\tilde{m},\tilde{n}\),
proposed for all Graphs \(5\)--\(9\) of Category \(\mathrm{III}\),
decide about the rank \(r_j\) of
the associated principal factor matrices \(M_j\)
of the \textit{wild} bicyclic bicubic fields \(B_j\),
\(j\in\lbrace 1,2,3,4,7,10\rbrace\).
As usual, we perform row operations
on these matrices: \\
\(M_1=
\begin{pmatrix}
1 & 0 & 0 \\
m & n & 0 \\
1 & 0 & 0 \\
0 & 1 & 2 \\
\end{pmatrix}\),
\(M_2=
\begin{pmatrix}
1 & 0 & 0 \\
m & n & 0 \\
1 & 0 & 0 \\
0 & 1 & 1 \\
\end{pmatrix}\),
\(M_{10}=
\begin{pmatrix}
1 & 0 & 0 \\
1 & 0 & 0 \\
0 & 0 & 1 \\
m & n & 0 \\
\end{pmatrix}\).

For \(B_1=k_1k_{pq}k_{pr}k_{qr}\),
\(M_1\)
leads to the decisive pivot element \(-2n\)
in the last column,
similarly
for \(B_2=k_2k_{pq}\tilde{k}_{pr}\tilde{k}_{qr}\),
\(M_2\)
leads to \(-n\),
and similarly
for \(B_{10}=k_3k_4k_rk_{pq}\),
\(M_{10}\)
leads to \(n\)
in the middle column.
So, \(r_1=r_2=r_{10}=3\) iff \(n\ne 0\),
by viewing the pivot elements modulo \(3\).
Next: \\
\(M_3=
\begin{pmatrix}
1 & 0 & 0 \\
\tilde{m} & \tilde{n} & 0 \\
1 & 0 & 0 \\
0 & 1 & 2 \\
\end{pmatrix}\),
\(M_4=
\begin{pmatrix}
1 & 0 & 0 \\
\tilde{m} & \tilde{n} & 0 \\
1 & 0 & 0 \\
0 & 1 & 1 \\
\end{pmatrix}\),
\(M_7=
\begin{pmatrix}
1 & 0 & 0 \\
1 & 0 & 0 \\
0 & 0 & 1 \\
\tilde{m} & \tilde{n} & 0 \\
\end{pmatrix}\).

For \(B_3=k_3\tilde{k}_{pq}\tilde{k}_{pr}k_{qr}\), 
\(M_3\)
leads to \(-2\tilde{n}\),
similarly
for \(B_4=k_4\tilde{k}_{pq}k_{pr}\tilde{k}_{qr}\),
\(M_4\)
leads to \(-\tilde{n}\),
and similarly
for \(B_7=k_1k_2k_r\tilde{k}_{pq}\),
\(M_7\)
leads to \(\tilde{n}\)
in the middle column.
So, \(r_3=r_4=r_7=3\) iff \(\tilde{n}\ne 0\).

In the regular case \(h_3(k_{pq})=h_3(\tilde{k}_{pq})=9\),
where \((m,n)=(\tilde{m},\tilde{n})\),
the condition \(n\ne 0\), that is \(q\mid A(k_{pq})\),
is certainly satisfied when \(\mathcal{P}=2\) or equivalently
\(\mathrm{Gal}(\mathrm{F}_3^2(k_{pq})/k_{pq})\simeq\langle 9,2\rangle\),
according to Theorem
\ref{thm:TwoPrimeCond}.
However, when \(\mathcal{P}=1\) or equivalently
\(\mathrm{Gal}(\mathrm{F}_3^2(k_{pq})/k_{pq})\simeq\langle 27,4\rangle\),
then we may either have \(q\mid A(k_{pq})\) and still \(n\ne 0\),
or \(p\mid A(k_{pq})\), \(n=0\), with completely different consequence
\(r_1=r_2=r_{10}=r_3=r_4=r_7=2\).
In the singular and super-singular cases, both pairs of parameters,
\((m,n)\) and \((\tilde{m},\tilde{n})\),
more precisely only \(n\) and \(\tilde{n}\),
must be taken into consideration, separately.
See also the proof of Theorem
\ref{thm:Cat3Gph6}.
\end{proof}


\noindent
In terms of capitulation targets in Corollary
\ref{cor:Three},
Theorem
\ref{thm:Cat3Gph6}
and parts of its proof
are now summarized in Table
\ref{tbl:UniCat3Gph6}
with transpositions in \textbf{bold} font.
\renewcommand{\arraystretch}{1.1}

\begin{table}[ht]
\caption{Norm class groups and minimal transfer kernels for Graph III.6}
\label{tbl:UniCat3Gph6}
\begin{center}
{\normalsize
\begin{tabular}{|c||c|c|c|c||c|c|c|c||c|c|c|c||c|c|c|c|}
\hline
 Base         & \multicolumn{4}{c||}{\(k_1\)} & \multicolumn{4}{c||}{\(k_2\)} & \multicolumn{4}{c||}{\(k_3\)} & \multicolumn{4}{c|}{\(k_4\)} \\
\hline
 Ext          & \(B_1\) & \(B_5\) & \(B_6\) & \(B_7\) & \(B_2\) & \(B_7\) & \(B_8\) & \(B_9\) & \(B_3\) & \(B_5\) & \(B_9\) & \(B_{10}\) & \(B_4\) & \(B_6\) & \(B_8\) & \(B_{10}\) \\
\hline
 NCG          & \(\mathfrak{r}\)    & \(\mathfrak{q}\) & \(\mathfrak{qr}\) & \(\mathfrak{qr}^2\) & \(\mathfrak{r}\)  & \(\mathfrak{qr}\)   & \(\mathfrak{q}\) & \(\mathfrak{qr}^2\)    
              & \(\mathfrak{qr}^2\) & \(\mathfrak{q}\) & \(\mathfrak{qr}\) & \(\mathfrak{r}\)    & \(\mathfrak{qr}\) & \(\mathfrak{qr}^2\) & \(\mathfrak{q}\) & \(\mathfrak{r}\) \\
 TK           & \(\mathfrak{qr}^2\) & \(\mathfrak{qr}\) & \(\mathfrak{q}\) & \(\mathfrak{r}\) & \(\mathfrak{qr}\) & \(\mathfrak{r}\) & \(\mathfrak{qr}^2\) & \(\mathfrak{q}\) 
              & \(\mathfrak{qr}^2\) & \(\mathfrak{qr}\) & \(\mathfrak{q}\) & \(\mathfrak{r}\) & \(\mathfrak{qr}\) & \(\mathfrak{q}\) & \(\mathfrak{qr}^2\) & \(\mathfrak{r}\) \\
\(\varkappa\) & \(\mathbf{4}\) & \(\mathbf{3}\) & \(\mathbf{2}\) & \(\mathbf{1}\) & \(\mathbf{2}\) & \(\mathbf{1}\) & \(\mathbf{4}\) & \(\mathbf{3}\) & \(1\) & \(\mathbf{3}\) & \(\mathbf{2}\) & \(4\) & \(1\) & \(\mathbf{3}\) & \(\mathbf{2}\) & \(4\) \\
\hline
\end{tabular}
}
\end{center}
\end{table}


\begin{theorem}
\label{thm:Cat3Gph6}
\textbf{(Second \(3\)-class group for \(\mathrm{III}.6\).)}
To identify the second \(3\)-class group
\(\mathfrak{M}=\mathrm{Gal}(\mathrm{F}_3^2(k_\mu)/k_\mu)\), \(1\le\mu\le 4\),
let the \textbf{principal factor} of \(k_{pq}\), respectively \(\tilde{k}_{pq}\), be
\(A(k_{pq})=p^mq^n\), respectively \(A(\tilde{k}_{pq})=p^{\tilde{m}}q^{\tilde{n}}\), and additionally assume the \textbf{regular} situation where both
\(\mathrm{Cl}_3(k_{pq})\simeq\mathrm{Cl}_3(\tilde{k}_{pq})\simeq (3,3)\)
are elementary bicyclic, whence \((m,n)=(\tilde{m},\tilde{n})\).

Then there are several
\textbf{minimal transfer kernel types} (mTKT) \(\varkappa_0\) of \(k_\mu\), \(1\le\mu\le 4\),
and other possible capitulation types in ascending order
\(\varkappa_0<\varkappa<\varkappa^\prime<\varkappa^{\prime\prime}\),
ending in the mandatory \(\varkappa^{\prime\prime}=(2000)\), type \(\mathrm{a}.3^\ast\),
either
\(\varkappa_0=(2134)\), type \(\mathrm{G}.16\),
\(\varkappa=(2130)\), type \(\mathrm{d}.23\),
\(\varkappa^\prime=(2100)\), type \(\mathrm{b}.10\),
or
\(\varkappa_0=(2143)\), type \(\mathrm{G}.19\),
\(\varkappa=(2140)\), type \(\mathrm{d}.25\), and again
\(\varkappa^\prime=(2100)\), type \(\mathrm{b}.10\).

In the \textbf{regular} situation, the second \(3\)-class group is \(\mathfrak{M}\simeq\)
\begin{equation}
\label{eqn:RegCat3Gph6}
\begin{cases}
\langle 81,7\rangle^4,\ \alpha=\lbrack 111,11,11,11\rbrack^4,\ \varkappa=(2000)^4 & \text{ if } n\ne 0,\ \mathcal{N}=1, \\
\langle 729,34..39\rangle^4,\ \alpha=\lbrack 111,111,21,21\rbrack^4,\ \varkappa=(2100)^4 & \text{ if } n=0,\ \mathcal{N}=2, \\
\langle 729,43\rangle^2,\langle 729,42\rangle^2,\ \alpha=\lbrack 111,111,22,21\rbrack^4,\ \varkappa=(2140)^2,(2130)^2 & \text{ if } n=0,\ \mathcal{N}=3, \\
\langle 2187,71\rangle^2,\langle 2187,69\rangle^2\ \alpha=\lbrack 111,111,22,22\rbrack^4,\ \varkappa=(2143)^2,(2134)^2 & \text{ if } n=0,\ \mathcal{N}=4, \\
\end{cases}
\end{equation}
where
\(\mathcal{N}:=\#\lbrace 1\le w\le 10\mid k_\mu<B_w,\ I_w=27\rbrace\).
Only in the first case, the \(3\)-class field tower has certainly the group
\(\mathfrak{G}=\mathrm{Gal}(\mathrm{F}_3^\infty(k_\mu)/k_\mu)\simeq\mathfrak{M}\)
and length \(\ell_3(k_\mu)=2\),
otherwise \(\ell_3(k_\mu)\ge 3\) cannot be excluded,
even if \(d_2(\mathfrak{M})\le 4\).

In the \textbf{singular} situation, the second \(3\)-class group is \(\mathfrak{M}\simeq\)
\begin{equation}
\label{eqn:SngCat3Gph6}
\begin{cases}
\langle 2187,251\vert 252\rangle^2,\langle 2187,250\rangle^2,\ \alpha=\lbrack 111,111,32,21\rbrack^4,\ \varkappa=(2140)^2,(2130)^2 & \text{ if } n=\tilde{n}=0,\ \mathcal{N}=3.
\end{cases}
\end{equation}

In the \textbf{super-singular} situation, no statement is possible,
since the order  of \(\mathfrak{M}\) may increase unboundedly. 
\end{theorem}

\begin{proof}
In the \textit{regular} situation
\(\mathrm{Cl}_3(k_{pq})=\mathrm{Cl}_3(\tilde{k}_{pq})=(3,3)\),
exponents \((m,n)\) and \((\tilde{m},\tilde{n})\)
of principal factors
\(A(k_{pq})=p^mq^n\) and \(A(\tilde{k}_{pq})=p^{\tilde{m}}q^{\tilde{n}}\) 
are equal.
Let \(\mathcal{P}\) be the number of primes dividing \(A(k_{pq})\).
According to the proof of Proposition
\ref{prp:Cat3Gph6},
ranks \(r_w\) and indices \(I_w:=(U_w:V_w)\) of subfield units
for \textit{wild} extensions are given by
\(r_w=3\), \(I_w=3\) iff \(n\ne 0\), for \(w\in\lbrace 1,2,10\rbrace\), and
\(r_w=3\), \(I_w=3\) iff \(\tilde{n}\ne 0\), for \(w\in\lbrace 3,4,7\rbrace\),
in particular, certainly for \(\mathcal{P}=2\).

This implies \(3\)-class numbers
\(h_3(B_w)=h_3(k_{pq})=h_3(\tilde{k}_{pq})=9\)
and \(3\)-class groups
\(\mathrm{Cl}_3(B_w)\simeq (3,3)\),
for \(w\in\lbrace 1,2,3,4,7,10\rbrace\),
whenever \(n\ne 0\), i.e. \(q\mid A(k_{pq})\),
a remarkable distinction of the prime \(q\)
against the primes \(p,r\).
We point out that this can occur
not only for \(\mathcal{P}=2\),
but also for \(\mathcal{P}=1\), provided that
\(A(k_{pq})=q\), \(n=1\), and not \(A(k_{pq})=p\), \(m=1\).

Indices of \textit{tame} extensions with
\(9\mid h_3(B_w)=I_w\in\lbrace 9,27\rbrace\) and \(r_w=2\)
are non-uniform:
corresponding to a unique elementary tricyclic
\(\mathrm{Cl}_3(B_w)\simeq (3,3,3)\),
we must have \(I_w=27\) for \(w\in\lbrace 6,9\rbrace\)
with norm class group \(N_{B_w/k_\mu}\mathrm{Cl}_3(B_w)\)
either \(\langle\lbrack\mathfrak{qr}\rbrack\rangle\)
or \(\langle\lbrack\mathfrak{qr}^2\rbrack\rangle\),
but corresponding to the remaining bicyclic
\(\mathrm{Cl}_3(B_w)\simeq (3,3)\),
the index \(I_w=9\) takes the minimal value
for \(w\in\lbrace 5,8\rbrace\)
with norm class group
\(N_{B_w/k_\mu}\mathrm{Cl}_3(B_w)=\langle\lbrack\mathfrak{q}\rbrack\rangle\).
Thus \(\mathcal{N}=1\) and the resulting Artin pattern
\(\alpha=\lbrack 111,11,11,11\rbrack\) uniquely identifies the group
\(\mathfrak{G}=\mathfrak{M}\simeq\langle 81,7\rangle\) of maximal class.

Now we come to \(n=0\), whence necessarily \(\mathcal{P}=1\).
Then \(r_w=2\) and \(I_w\in\lbrace 9,27\rbrace\)
for the wild extensions \(w\in\lbrace 1,2,3,4,7,10\rbrace\).
Indices of tame extensions now become uniform,
corresponding to a pair of elementary tricyclic
\(\mathrm{Cl}_3(B_w)\simeq (3,3,3)\),
which enforces \(I_w=27\) for \(w\in\lbrace 5,6,8,9\rbrace\),
i.e., \(\mathcal{N}\ge 2\).
The number \(\mathcal{N}\) of maximal unit indices decides about the group \(\mathfrak{M}\):
If \(\mathcal{N}=2\), then
for all \(w\in\lbrace 1,2,3,4,7,10\rbrace\):
\(I_w=9\), \(h_3(B_w)=3\cdot h_3(k_{pq})=27\),
and \(\mathrm{Cl}_3\simeq (9,3)\),
according to the laws for \(3\)-groups of coclass \(\ge 2\)
\cite[pp. 289--292]{Ma2015b}.
The Artin pattern
\(\alpha=\lbrack 111,111,21,21\rbrack\) identifies the possible groups
\(\mathfrak{M}\simeq\langle 729,34..39\rangle\).
If \(\mathcal{N}=3\), then
\(I_w=27\) for \(w\in\lbrace 3,4,7\rbrace\), but
\(I_w=9\) for \(w\in\lbrace 1,2,10\rbrace\).
The Artin pattern
\(\alpha=\lbrack 111,111,22,21\rbrack\)
together with \(\varkappa=(2140)^2\), \(\varkappa=(2130)^2\),
according to Table
\ref{tbl:UniCat3Gph6},
identifies the possible groups
\(\mathfrak{M}\simeq\langle 2187,251\vert 252\rangle^2,\langle 2187,250\rangle^2\)
of coclass \(2\).
If \(\mathcal{N}=4\), then
for all \(w\in\lbrace 1,2,3,4,7,10\rbrace\):
\(I_w=27\), \(h_3(B_w)=9\cdot h_3(k_{pq})=81\),
and \(\mathrm{Cl}_3\simeq (9,9)\).
The Artin pattern
\(\alpha=\lbrack 111,111,22,22\rbrack\)
together with \(\varkappa=(2143)^2\), \(\varkappa=(2134)^2\),
according to Table
\ref{tbl:UniCat3Gph6},
identifies the possible groups
\(\mathfrak{M}\simeq\langle 2187,71\rangle^2,\langle 2187,69\rangle^2\)
of coclass \(3\).
\end{proof}


\begin{corollary}
\label{cor:UniCat3Gph6}
\textbf{(Non-uniformity of the quartet for \(\mathrm{III}.6\).)}
Only for \(\mathcal{N}\le 2\),
the components of the quartet, all with \(3\)-rank two,
share a common capitulation type
\(\varkappa(k_\mu)\),
common abelian type invariants
\(\alpha(k_\mu)\),
and a common second \(3\)-class group
\(\mathrm{Gal}(\mathrm{F}_3^2(k_\mu)/k_\mu)\),
for \(1\le\mu\le 4\).
For \(\mathcal{N}\ge 3\),
the quartet splits into two sub-doublets
and thus becomes non-uniform.
\end{corollary}

\begin{proof}
This is an immediate consequence of Theorem
\ref{thm:Cat3Gph6}
and Table
\ref{tbl:UniCat3Gph6}.
\end{proof}


\begin{example}
\label{exm:Cat3Gph6}
Prototypes for Graph \(\mathrm{III}.6\), that is,
minimal conductors for each scenario in Theorem
\ref{thm:Cat3Gph6}
are the following.

There are the \textbf{regular} cases
\(c=8\,541\) with symbol
\(\lbrace 9\leftrightarrow 73\rightarrow 13\rbrace\),
\((m,n)=(1,2)\);
\(c=9\,373\) with symbol
\(\lbrace 103\leftrightarrow 13\rightarrow 7\rbrace\),
\((m,n)=(1,1)\);
\(c=56\,329\) with symbol
\(\lbrace 619\leftrightarrow 13\rightarrow 7\rbrace\),
\((m,n)=(0,1)\),
all uniformly with \(\mathfrak{G}=\mathfrak{M}=\langle 81,7\rangle^4\),
in contrast to
\(c=142\,519\) with symbol
\(\lbrace 19\leftrightarrow 577\rightarrow 13\rbrace\),
\((m,n)=(1,0)\),
and uniform \(\mathfrak{M}=\langle 729,37..39\rangle^4\);
\(c=152\,893\) with symbol
\(\lbrace 13\leftrightarrow 619\rightarrow 19\rbrace\),
\((m,n)=(1,0)\),
and uniform \(\mathfrak{M}=\langle 729,34..36\rangle^4\);
\(c=163\,681\) with symbol
\(\lbrace 67\leftrightarrow 349\rightarrow 7\rbrace\),
\((m,n)=(1,0)\),
and non-uniform \(\mathfrak{M}=\langle 729,42\rangle^2,\langle 729,43\rangle^2\);
\(c=193\,059\) with symbol
\(\lbrace 1129\leftrightarrow 19\rightarrow 9\rbrace\),
\((m,n)=(1,0)\),
and non-uniform \(\mathfrak{M}=\langle 2187,69\rangle^2,\langle 2187,71\rangle^2\)
with two distinct minimal transfer kernel types.

Further, the \textbf{singular} cases
\(c=78\,169\) with symbol
\(\lbrace 859\leftrightarrow 13\rightarrow 7\rbrace\),
\((m,n)=(1,0)\), \((\tilde{m},\tilde{n})=(1,1)\),
and non-uniform
\(\mathfrak{M}=\langle 2187,250\rangle^2,\langle 2187,251\vert 252\rangle^2\);
\(c=142\,947\) with symbol
\(\lbrace 9\leftrightarrow 2269\rightarrow 7\rbrace\),
\((m,n)=(1,0)\), \((\tilde{m},\tilde{n})=(0,1)\),
and uniform
\(\mathfrak{M}=\langle 2187,253\rangle^4\).

Finally, the \textbf{super-singular} cases
\(c=102\,277\) with symbol
\(\lbrace 769\leftrightarrow 7\rightarrow 19\rbrace\),
\((m,n)=(\tilde{m},\tilde{n})=(0,1)\),
and uniform \(\mathfrak{M}=\langle 729,37..39\rangle^4\);
\(c=199\,171\) with symbol
\(\lbrace 7\leftrightarrow 769\rightarrow 37\rbrace\),
\((m,n)=(\tilde{m},\tilde{n})=(1,0)\),
and uniform \(\mathfrak{M}=\langle 6561,693..698\rangle^4\).
\end{example}


In Table
\ref{tbl:ProtoCat3Gph6},
we summarize the prototypes of Graph \(\mathrm{III}.6\)
in the same way as in Table
\ref{tbl:ProtoCat3Gph5}.

\renewcommand{\arraystretch}{1.1}

\begin{table}[ht]
\caption{Prototypes for Graph III.6}
\label{tbl:ProtoCat3Gph6}
\begin{center}
{\tiny
\begin{tabular}{|c|c||c|c|c|c|c|c|c|c|c|}
\hline
 No. & \(c\)        & \(r\leftarrow p\leftrightarrow q\)    & \(v^\ast\) & \(v\) & \(m,n\) & \(\tilde{v}\) & \(\tilde{m},\tilde{n}\) & capitulation type & \(\mathfrak{M}\) & \(\ell_3(k)\) \\
\hline
   1 &   \(8\,541\) & \(13\leftarrow 73\leftrightarrow 9\)  & \(1\)      & \(2\) & \(1,2\) & \(2\) & \(1,2\) & \(\mathrm{a}.3^\ast\) & \(\langle 81,7\rangle^4\) & \(=2\) \\
   2 &   \(9\,373\) & \(7\leftarrow 13\leftrightarrow 103\) & \(1\)      & \(2\) & \(1,1\) & \(2\) & \(1,1\) & \(\mathrm{a}.3^\ast\) & \(\langle 81,7\rangle^4\) & \(=2\) \\
  20 &  \(56\,329\) & \(7\leftarrow 13\leftrightarrow 619\) & \(2\)      & \(2\) & \(0,1\) & \(2\) & \(0,1\) & \(\mathrm{a}.3^\ast\) & \(\langle 81,7\rangle^4\) & \(=2\) \\
  29 &  \(78\,169\) & \(7\leftarrow 13\leftrightarrow 859\) & \(3\)      & \(3\) & \(1,0\) & \(3\) & \(1,1\) & \(\mathrm{d}.23\)     & \(\langle 2187,250\rangle^2\) & \(\ge 2\) \\
     &              &                                       &            &       &         &       &         & \(\mathrm{d}.25\)     & \(\langle 2187,251\vert 252\rangle^2\) & \(\ge 2\) \\ 
  34 & \(102\,277\) & \(19\leftarrow 7\leftrightarrow 769\) & \(4\)      & \(3\) & \(0,1\) & \(3\) & \(0,1\) & \(\mathrm{b}.10\)     & \(\langle 729,37..39\rangle^4\) & \(\ge 2\) \\
  52 & \(142\,519\) & \(13\leftarrow 577\leftrightarrow 19\)& \(2\)      & \(2\) & \(1,0\) & \(2\) & \(1,0\) & \(\mathrm{b}.10\)     & \(\langle 729,37..39\rangle^4\) & \(\ge 2\) \\
  54 & \(142\,947\) & \(7\leftarrow 2269\leftrightarrow 9\) & \(3\)      & \(3\) & \(1,0\) & \(3\) & \(0,1\) & \(\mathrm{b}.10\)     & \(\langle 2187,253\rangle^4\) & \(\ge 2\) \\
  56 & \(152\,893\) & \(19\leftarrow 619\leftrightarrow 13\)& \(2\)      & \(2\) & \(1,0\) & \(2\) & \(1,0\) & \(\mathrm{b}.10\)     & \(\langle 729,34..36\rangle^4\) & \(\ge2\) \\
  58 & \(163\,681\) & \(7\leftarrow 349\leftrightarrow 67\) & \(2\)      & \(2\) & \(1,0\) & \(2\) & \(1,0\) & \(\mathrm{d}.23\)     & \(\langle 729,42\rangle^2\) & \(\ge 2\) \\
     &              &                                       &            &       &         &       &         & \(\mathrm{d}.25\)     & \(\langle 729,43\rangle^2\) & \(\ge 2\) \\
  71 & \(193\,059\) & \(9\leftarrow 19\leftrightarrow 1129\)& \(2\)      & \(2\) & \(1,0\) & \(2\) & \(1,0\) & \(\mathrm{G}.16\)     & \(\langle 2187,69\rangle^2\) & \(\ge 2\) \\
     &              &                                       &            &       &         &       &         & \(\mathrm{G}.19\)     & \(\langle 2187,71\rangle^2\) & \(\ge 2\) \\ 
  75 & \(199\,171\) & \(37\leftarrow 769\leftrightarrow 7\) & \(4\)      & \(3\) & \(1,0\) & \(3\) & \(1,0\) & \(\mathrm{b}.10\)     & \(\langle 6561,693..698\rangle^4\) & \(\ge 3\) \\
\hline
\end{tabular}
}
\end{center}
\end{table}


\subsection{Category III, Graph 7}
\label{ss:Cat3Gph7}

\noindent
Let \((k_1,\ldots,k_4)\) be a quartet of cyclic cubic number fields
sharing the common conductor \(c=pqr\),
belonging to Graph \(7\) of Category \(\mathrm{III}\)
with combined cubic residue symbol
\(\lbrack p,q,r\rbrack_3=\lbrace r\rightarrow p\leftrightarrow q\rbrace\).


\begin{proposition}
\label{prp:Cat3Gph7}
\textbf{(Quartet with \(3\)-rank two for \(\mathrm{III}.7\).)}
For fixed \(\mu\in\lbrace 1,2,3,4\rbrace\),
let \(\mathfrak{p},\mathfrak{q},\mathfrak{r}\) be the prime ideals of \(k_\mu\)
over \(p,q,r\), that is,
\(p\mathcal{O}_{k_\mu}=\mathfrak{p}^3\),
\(q\mathcal{O}_{k_\mu}=\mathfrak{q}^3\),
\(r\mathcal{O}_{k_\mu}=\mathfrak{r}^3\).
Under the normalizing assumption
\(A(k_{qr})=qr^2\), \(A(\tilde{k}_{qr})=qr\),
the \textbf{principal factors} of \(k_\mu\) are
\begin{equation}
\label{eqn:PFCat3Gph7}
A(k_1)=A(k_3)=qr \quad \text{ and } \quad
A(k_2)=A(k_4)=qr^2,
\end{equation}
and the \(3\)-class group of \(k_\mu\) is
\begin{equation}
\label{eqn:GenCat3Gph7}
\mathrm{Cl}_3(k_\mu)=
\langle\lbrack\mathfrak{p}\rbrack,\lbrack\mathfrak{q}\rbrack\rangle=
\langle\lbrack\mathfrak{p}\rbrack,\lbrack\mathfrak{r}\rbrack\rangle\simeq (3,3).
\end{equation}
The unramified cyclic cubic relative extensions of \(k_\mu\)
are among the absolutely bicyclic bicubic fields \(B_i\), \(1\le i\le 10\).
The \textbf{tame} extensions
with \(9\mid h_3(B_i)=(U_i:V_i)\in\lbrace 9,27\rbrace\)
are \(B_i\) with \(i=5,6,8,9\),
since they neither contain \(k_{pq}\) nor \(\tilde{k}_{pq}\).
For each \(\mu\), there are two tame extensions \(B_j/k_\mu\), \(B_\ell/k_\mu\)
with the following properties.
The first, \(B_j\) with \(j\in\lbrace 6,9\rbrace\),
has norm class group
\(N_{B_j/k_\mu}(\mathrm{Cl}_3(B_j))=\langle\lbrack\mathfrak{p}\rbrack\rangle\),
transfer kernel
\begin{equation}
\label{eqn:Trans1Cat3Gph7}
\ker(T_{B_j/k_\mu})\ge\langle\lbrack\mathfrak{q}\rbrack\rangle,
\end{equation}
and \(3\)-class group
\(\mathrm{Cl}_3(B_j)=
\langle\lbrack\mathfrak{P}_1\rbrack,
\lbrack\mathfrak{P}_2\rbrack,
\lbrack\mathfrak{P}_3\rbrack\rangle\ge (3,3)\),
generated by the classes of the prime ideals of \(B_j\) over
\(\mathfrak{p}\mathcal{O}_{B_j}=\mathfrak{P}_1\mathfrak{P}_2\mathfrak{P}_3\).
The second, \(B_\ell\) with \(\ell\in\lbrace 5,8\rbrace\),
has norm class group
\(N_{B_\ell/k_\mu}(\mathrm{Cl}_3(B_\ell))=\langle\lbrack\mathfrak{q}\rbrack\rangle\),
\textbf{cyclic} transfer kernel
\begin{equation}
\label{eqn:Trans2Cat3Gph7}
\ker(T_{B_\ell/k_\mu})=\langle\lbrack\mathfrak{p}\rbrack\rangle
\end{equation}
of order \(3\),
and \textbf{elementary tricyclic} \(3\)-class group
\(\mathrm{Cl}_3(B_\ell)=
\langle\lbrack\mathfrak{Q}_1\rbrack,
\lbrack\mathfrak{Q}_2\rbrack,
\lbrack\mathfrak{Q}_3\rbrack\rangle\simeq (3,3,3)\),
generated by the classes of the prime ideals of \(B_\ell\) over
\(\mathfrak{q}\mathcal{O}_{B_\ell}=\mathfrak{Q}_1\mathfrak{Q}_2\mathfrak{Q}_3\).
The pair \((j,\ell)\) forms a hidden or actual \textbf{transposition}
of the transfer kernel type \(\varkappa(k_\mu)\).
The remaining two \(B_i>k_\mu\), \(i\ne j\), \(i\ne\ell\),
have norm class group 
\(\langle\lbrack\mathfrak{pq}\rbrack\rangle\), respectively
\(\langle\lbrack\mathfrak{pq}^2\rbrack\rangle\),
and transfer kernel
\[
\ker(T_{B_i/k_\mu})\ge\langle\lbrack\mathfrak{q}\rbrack\rangle, 
\]
providing the option of two possible \textbf{repetitions}
in the transfer kernel type \(\varkappa(k_\mu)\).\\
In terms of \(n\) and \(\tilde{n}\) in
\(A(k_{pq})=p^mq^n\) and \(A(\tilde{k}_{pq})=p^{\tilde{m}}q^{\tilde{n}}\),
the ranks of the \textbf{wild} extensions are
\begin{equation}
\label{RankCat3Gph7}
r_1=r_2=r_{10}=3 \text{ iff } m\ne 0 \text{ iff } p\mid A(k_{pq})
\text{ and }
r_3=r_4=r_7=3 \text{ iff } \tilde{m}\ne 0 \text{ iff } p\mid A(\tilde{k}_{pq}).
\end{equation}
\end{proposition}

\begin{proof}
By Proposition
\ref{prp:Principal2},
the symbol \(r\rightarrow p\) implies principal factors
\(A(k_{pr})=A(\tilde{k}_{pr})=r\).

We assume principal factors \(A(k_\mu)=p^{x_\mu}q^{y_\mu}r^{z_\mu}\),
for \(1\le\mu\le 4\),
and \(A(k_{qr})=qr^2\), \(A(\tilde{k}_{qr})=qr\).

We generally have the \textit{tame} matrix ranks \(r_5=r_6=r_8=r_9=2\)
and draw conclusions by explicit calculations.
For these bicyclic bicubic fields \(B_j\),
\(j\in\lbrace 5,6,8,9\rbrace\),
the rank \(r_j\) is calculated with row operations
on the associated principal factor matrices \(M_j\): \\
\(M_5=
\begin{pmatrix}
x_1 & y_1 & z_1 \\
x_3 & y_3 & z_3 \\
1 & 0 & 0 \\
0 & 1 & 1 \\
\end{pmatrix}\),
\(M_6=
\begin{pmatrix}
x_1 & y_1 & z_1 \\
x_4 & y_4 & z_4 \\
0 & 1 & 0 \\
0 & 0 & 1 \\
\end{pmatrix}\),
\(M_8=
\begin{pmatrix}
x_2 & y_2 & z_2 \\
x_4 & y_4 & z_4 \\
1 & 0 & 0 \\
0 & 1 & 2 \\
\end{pmatrix}\),
\(M_9=
\begin{pmatrix}
x_2 & y_2 & z_2 \\
x_3 & y_3 & z_3 \\
0 & 1 & 0 \\
0 & 0 & 1 \\
\end{pmatrix}\).

For \(B_5=k_1k_3k_p\tilde{k}_{qr}\),
\(M_5\)
leads to the decisive pivot elements \(z_1-y_1\) and \(z_3-y_3\),
similarly,
for \(B_6=k_1k_4k_q\tilde{k}_{pr}\),
\(M_6\)
leads to \(x_1\) and \(x_4\),
similarly,
for \(B_8=k_2k_4k_pk_{qr}\),
\(M_8\)
leads to \(z_2-2y_2\) and \(z_4-2y_4\),
and similarly,
for \(B_9=k_2k_3k_qk_{pr}\), 
\(M_9\)
leads to \(x_2\) and \(x_3\).
So, \(r_5=r_6=2\)
implies 
\(z_1=y_1\), \(z_3=y_3\),
\(x_1=x_4=0\),
and 
\(r_8=r_9=2\) implies
\(z_2=2y_2\), \(z_4=2y_4\),
\(x_2=x_3=0\),
i.e. \(A(k_1)=A(k_3)=qr\) and
\(A(k_2)=A(k_4)=qr^2\).

A consequence of these principal factors
is the coincidence of the subgroups of
\(\mathrm{Cl}_3(k_\mu)\)
generated by the classes
\(\lbrack\mathfrak{q}\rbrack\) and \(\lbrack\mathfrak{r}\rbrack\)
in \(k_\mu\), \(\mu=1,\ldots,4\).
By Corollary
\ref{cor:Capitulation}, \\
since \(\mathfrak{p}\) is principal ideal in \(k_p\),
the class \(\lbrack\mathfrak{p}\rbrack\) capitulates in 
\(B_5=k_1k_3k_p\tilde{k}_{qr}\) and
\(B_8=k_2k_4k_pk_{qr}\); \\
since \(\mathfrak{q}\) is principal ideal in \(k_q\),
the class \(\lbrack\mathfrak{q}\rbrack\) capitulates in
\(B_6=k_1k_4k_q\tilde{k}_{pr}\) and
\(B_9=k_2k_3k_qk_{pr}\); \\
since \(\mathfrak{r}\) is principal ideal in \(k_r\),
the class \(\lbrack\mathfrak{r}\rbrack\),
and thus \(\lbrack\mathfrak{q}\rbrack\),
capitulates in
\(B_7=k_1k_2k_r\tilde{k}_{pq}\) and
\(B_{10}=k_3k_4k_rk_{pq}\).

Moreover, since \(\mathfrak{r}\) is principal ideal in
\(k_{pr}\) and \(\tilde{k}_{pr}\),
the class \(\lbrack\mathfrak{r}\rbrack\),
and thus \(\lbrack\mathfrak{q}\rbrack\),
also capitulates in \\
\(B_1=k_1k_{pq}k_{pr}k_{qr}\),
\(B_2=k_2k_{pq}\tilde{k}_{pr}\tilde{k}_{qr}\),
\(B_3=k_3\tilde{k}_{pq}\tilde{k}_{pr}k_{qr}\),
\(B_4=k_4\tilde{k}_{pq}k_{pr}\tilde{k}_{qr}\), \\
\(B_6=k_1k_4k_q\tilde{k}_{pr}\), and
\(B_9=k_2k_3k_qk_{pr}\), by Proposition
\ref{prp:Bicyc}.

The parameters \(m,n,\tilde{m},\tilde{n}\),
proposed for all Graphs \(5\)--\(9\) of Category \(\mathrm{III}\),
decide about the rank \(r_j\) of
the associated principal factor matrices \(M_j\)
of the \textit{wild} bicyclic bicubic fields \(B_j\),
\(j\in\lbrace 1,2,3,4,7,10\rbrace\).
As usual, we perform row operations
on these matrices: \\
\(M_1=
\begin{pmatrix}
0 & 1 & 1 \\
m & n & 0 \\
0 & 0 & 1 \\
0 & 1 & 2 \\
\end{pmatrix}\),
\(M_2=
\begin{pmatrix}
0 & 1 & 2 \\
m & n & 0 \\
0 & 0 & 1 \\
0 & 1 & 1 \\
\end{pmatrix}\),
\(M_{10}=
\begin{pmatrix}
0 & 1 & 1 \\
0 & 1 & 2 \\
0 & 0 & 1 \\
m & n & 0 \\
\end{pmatrix}\).

For \(B_1=k_1k_{pq}k_{pr}k_{qr}\),
\(M_1\)
leads to the decisive pivot element \(m\)
in the first column,
similarly,
for \(B_2=k_2k_{pq}\tilde{k}_{pr}\tilde{k}_{qr}\),
\(M_2\)
leads to \(m\),
and similarly,
for \(B_{10}=k_3k_4k_rk_{pq}\),
\(M_{10}\)
leads to \(m\)
in the first column.
So, \(r_1=r_2=r_{10}=3\) iff \(m\ne 0\).
Next we consider: \\
\(M_3=
\begin{pmatrix}
0 & 1 & 1 \\
\tilde{m} & \tilde{n} & 0 \\
0 & 0 & 1 \\
0 & 1 & 2 \\
\end{pmatrix}\),
\(M_4=
\begin{pmatrix}
0 & 1 & 2 \\
\tilde{m} & \tilde{n} & 0 \\
0 & 0 & 1 \\
0 & 1 & 1 \\
\end{pmatrix}\),
\(M_7=
\begin{pmatrix}
0 & 1 & 1 \\
0 & 1 & 2 \\
0 & 0 & 1 \\
\tilde{m} & \tilde{n} & 0 \\
\end{pmatrix}\).

For \(B_3=k_3\tilde{k}_{pq}\tilde{k}_{pr}k_{qr}\), 
\(M_3\)
leads to \(\tilde{m}\),
similarly,
for \(B_4=k_4\tilde{k}_{pq}k_{pr}\tilde{k}_{qr}\),
\(M_4\)
leads to \(\tilde{m}\),
and similarly,
for \(B_7=k_1k_2k_r\tilde{k}_{pq}\),
\(M_7\)
leads to \(\tilde{m}\)
in the first column.
So, \(r_3=r_4=r_7=3\) iff \(\tilde{m}\ne 0\).

Since \(r\) splits in \(k_{p}\),
it also splits in
\(B_5=k_1k_3k_p\tilde{k}_{qr}\),
\(B_8=k_2k_4k_pk_{qr}\).

Since \(q\) splits in \(k_{p}\),
it also splits in
\(B_5=k_1k_3k_p\tilde{k}_{qr}\),
\(B_8=k_2k_4k_pk_{qr}\).

Since \(p\) splits in \(k_{q}\),
it also splits in
\(B_6=k_1k_4k_q\tilde{k}_{pr}\),
\(B_9=k_2k_3k_qk_{pr}\).
\end{proof}


\noindent
In terms of capitulation targets in Corollary
\ref{cor:Three},
Theorem
\ref{thm:Cat3Gph7}
and parts of its proof
are now summarized in Table
\ref{tbl:UniCat3Gph7}
with transpositions in \textbf{bold} font.

\renewcommand{\arraystretch}{1.1}

\begin{table}[ht]
\caption{Norm class groups and minimal transfer kernels for Graph III.7}
\label{tbl:UniCat3Gph7}
\begin{center}
\begin{tabular}{|c||c|c|c|c||c|c|c|c||c|c|c|c||c|c|c|c|}
\hline
 Base         & \multicolumn{4}{c||}{\(k_1\)} & \multicolumn{4}{c||}{\(k_2\)} & \multicolumn{4}{c||}{\(k_3\)} & \multicolumn{4}{c|}{\(k_4\)} \\
\hline
 Ext          & \(B_1\) & \(B_5\) & \(B_6\) & \(B_7\) & \(B_2\) & \(B_7\) & \(B_8\) & \(B_9\) & \(B_3\) & \(B_5\) & \(B_9\) & \(B_{10}\) & \(B_4\) & \(B_6\) & \(B_8\) & \(B_{10}\) \\
\hline
 NCG          & \(\mathfrak{pq}^2\) & \(\mathfrak{q}\) & \(\mathfrak{p}\) & \(\mathfrak{pq}\) & \(\mathfrak{pq}^2\) & \(\mathfrak{pq}\) & \(\mathfrak{q}\) & \(\mathfrak{p}\)
              & \(\mathfrak{pq}^2\) & \(\mathfrak{q}\) & \(\mathfrak{p}\) & \(\mathfrak{pq}\) & \(\mathfrak{pq}^2\) & \(\mathfrak{p}\)  & \(\mathfrak{q}\) & \(\mathfrak{pq}\) \\
 TK           & \(\mathfrak{q}\)  & \(\mathfrak{p}\) & \(\mathfrak{q}\) & \(\mathfrak{q}\)    & \(\mathfrak{q}\)    & \(\mathfrak{q}\)  & \(\mathfrak{p}\) & \(\mathfrak{q}\)
              & \(\mathfrak{q}\)  & \(\mathfrak{p}\) & \(\mathfrak{q}\) & \(\mathfrak{q}\)    & \(\mathfrak{q}\)    & \(\mathfrak{q}\)  & \(\mathfrak{p}\) & \(\mathfrak{q}\) \\
\(\varkappa\) & \(2\) & \(\mathbf{3}\) & \(\mathbf{2}\) & \(2\) & \(3\) & \(3\) & \(\mathbf{4}\) & \(\mathbf{3}\) & \(2\) & \(\mathbf{3}\) & \(\mathbf{2}\) & \(2\) & \(3\) & \(\mathbf{3}\) & \(\mathbf{2}\) & \(3\) \\
\hline
\end{tabular}
\end{center}
\end{table}


\begin{theorem}
\label{thm:Cat3Gph7}
\textbf{(Second \(3\)-class group for \(\mathrm{III}.7\).)}
Let \((k_1,\ldots,k_4)\) be a quartet of cyclic cubic number fields
sharing the common conductor \(c=pqr\),
belonging to Graph \(7\) of Category \(\mathrm{III}\)
with combined cubic residue symbol
\(\lbrack p,q,r\rbrack_3=\lbrace q\leftrightarrow p\rightarrow r\rbrace\).

Then the \textbf{minimal transfer kernel type} (mTKT) of \(k_\mu\), \(1\le\mu\le 4\),
is \(\varkappa_0=(2111)\), type \(\mathrm{H}.4\),
and the other possible capitulation types in ascending order
\(\varkappa_0<\varkappa^\prime<\varkappa^{\prime\prime}<\varkappa^{\prime\prime\prime}\)
are
\(\varkappa^\prime=(2110)\), type \(\mathrm{d}.19\),
\(\varkappa^{\prime\prime}=(2100)\), type \(\mathrm{b}.10\), and
\(\varkappa^{\prime\prime\prime}=(2000)\), type \(\mathrm{a}.3^\ast\).

To identify the second \(3\)-class group
\(\mathfrak{M}=\mathrm{Gal}(\mathrm{F}_3^2(k_\mu)/k_\mu)\), \(1\le\mu\le 4\),
let the \textbf{decisive principal factors} be
\(A(k_{pq})=p^mq^n\), \(A(\tilde{k}_{pq})=p^{\tilde{m}}q^{\tilde{n}}\), and additionally assume
the \textbf{regular situation} where both
\(\mathrm{Cl}_3(k_{pq})\simeq\mathrm{Cl}_3(\tilde{k}_{pq})\simeq (3,3)\)
are elementary bicyclic, whence \((m,n)=(\tilde{m},\tilde{n})\). Then

\begin{equation}
\label{eqn:Cat3Gph7}
\mathfrak{M}\simeq
\begin{cases}
\langle 81,7\rangle,\ \alpha=\lbrack 111,11,11,11\rbrack,\ \varkappa=(2000) & \text{ if } m\ne 0,\  \mathcal{N}=1, \\
\langle 729,34..39\rangle,\ \alpha=\lbrack 111,111,21,21\rbrack,\ \varkappa=(2100) & \text{ if } m=0,\ \mathcal{N}=2, \\
\langle 729,41\rangle,\ \alpha=\lbrack 111,111,22,21\rbrack,\ \varkappa=(2110) & \text{ if } m=0,\ \mathcal{N}=3, \\
\langle 2187,65\vert 67\rangle,\ \alpha=\lbrack 111,111,22,22\rbrack,\ \varkappa=(2111) & \text{ if } m=0,\ \mathcal{N}=4,
\end{cases}
\end{equation}
where \(\mathcal{N}:=\#\lbrace 1\le j\le 10\mid k_\mu<B_j,\ I_j=27\rbrace\).
Only in the leading row, the \(3\)-class field tower has warranted group
\(\mathfrak{G}=\mathrm{Gal}(\mathrm{F}_3^\infty(k_\mu)/k_\mu)\simeq\mathfrak{M}\),
with length \(\ell_3(k_\mu)=2\).
Otherwise,
tower length \(\ell_3(k_\mu)\ge 3\) cannot be excluded,
even if \(d_2(\mathfrak{M})\le 4\).

In \textbf{(super-)singular situations},
the group \(\mathfrak{M}\) must be of coclass \(\mathrm{cc}(\mathfrak{M})\ge 2\),
and capitulation of type \(\varkappa^{\prime\prime\prime}=(2000)\) is impossible.
\end{theorem}

\begin{proof}
We know that the tame ranks are
\(r_5=r_6=r_8=r_9=2\),
and thus \(I_5,I_6,I_8,I_9\in\lbrace 9,27\rbrace\),
in particular, \(I_5=I_8=27\), whence certainly \(\mathcal{N}\ge 1\).
Further, the wild ranks are
\(r_1=r_2=r_{10}=3\)
iff \(m\ne 0\),
and
\(r_3=r_4=r_7=3\)
iff \(\tilde{m}\ne 0\).

In the \textbf{regular situation} where the
\(3\)-class groups of \(k_{pq}\) and \(\tilde{k}_{pq}\)
are elementary bicyclic,
tight bounds arise for the abelian quotient invariants \(\alpha\) of
the group \(\mathfrak{M}\):

The first scenario,
\(m\ne 0\),
is equivalent to
\(\mathcal{N}=1\),
with wild ranks
\(h_3(B_j)=h_3(k_{pq})=9\), for \(j=1,2,10\),
\(h_3(B_j)=h_3(\tilde{k}_{pq})=9\), for \(j=3,4,7\),
and tame ranks
\(h_3(B_j)=I_j=9\), for \(j=6,9\),
\(h_3(B_j)=I_j=27\), for \(j=5,8\),
that is \(\alpha=\lbrack 111,11,11,11\rbrack\)
and consequently \(\varkappa=(2000)\),
since \(\langle 81,7\rangle\) is unique
with this \(\alpha\).

The other three scenarios share
\(m=0\),
and an explicit transposition
between \(B_5\), \(B_6\), respectively \(B_5\), \(B_9\),
and \(B_6\), \(B_8\), respectively \(B_8\), \(B_9\),
giving rise to \(\varkappa=(21\ast\ast)\),
and common \(h_3(B_j)=I_j=27\), for \(j=5,6,8,9\),
implying \(\alpha=\lbrack 111,111,\ast,\ast\rbrack\).

The second scenario with \(\mathcal{N}=2\)
is supplemented by 
\(I_j=9\), \(h_3(B_j)=3\cdot h_3(k_{pq})=27\), for \(j=1,2,10\).
\(I_j=9\), \(h_3(B_j)=3\cdot h_3(\tilde{k}_{pq})=27\), for \(j=3,4,7\),
giving rise to \(\alpha=\lbrack 111,111,21,21\rbrack\),
\(\varkappa=(2100)\),
characteristic for \(\langle 729,34..39\rangle\)
(Cor.
\ref{cor:ClassGroup33}).

The third scenario with \(\mathcal{N}=3\)
is supplemented by 
\(I_j=27\), \(h_3(B_j)=9\cdot h_3(k_{pq})=81\), for \(j=1,2,10\),
but still
\(I_j=9\), \(h_3(B_j)=3\cdot h_3(\tilde{k}_{pq})=27\), for \(j=3,4,7\),
giving rise to \(\alpha=\lbrack 111,111,22,21\rbrack\),
\(\varkappa=(2110)\),
characteristic for \(\langle 729,41\rangle\).

The fourth scenario with \(\mathcal{N}=4\)
is supplemented by 
\(I_j=27\), \(h_3(B_j)=9\cdot h_3(k_{pq})=81\), for \(j=1,2,10\),
\(I_j=27\), \(h_3(B_j)=9\cdot h_3(k_{pq})=81\), for \(j=3,4,7\),
giving rise to \(\alpha=\lbrack 111,111,22,22\rbrack\),
\(\varkappa=(2111)\),
characteristic for either \(\langle 2187,65\vert 67\rangle\) or
\(\langle 6561,714..719\vert 738..743\rangle\)
with coclass \(\mathrm{cc}=3\).
If \(d_2(\mathfrak{M})=5\),
then tower length must be \(\ell_3(k_\mu)\ge 3\)
For this minimal capitulation type H.4, \(\varkappa=(2111)\),
all transfer kernels are cyclic of order \(3\),
and the minimal unit norm indices
correspond to maximal subfield unit indices.

In \textbf{(super-)singular situations}, the
\(3\)-class groups of \(k_{pq}\) and \(\tilde{k}_{pq}\)
are non-elementary bicyclic,
and even in the simplest case \(m\ne 0\), \(\tilde{m}\ne 0\),
we have
\(I_j=3\), \(27\mid h_3(B_j)=h_3(k_{pq})\), for \(j=1,2,10\),
\(I_j=3\), \(27\mid h_3(B_j)=h_3(\tilde{k}_{pq})\), for \(j=3,4,7\),
which prohibits the occurrence of abelian type invariants \((11)\),
required for
\(3\)-groups of coclass \(\mathrm{cc}(\mathfrak{M})=1\) (maximal class).
\end{proof}


\begin{corollary}
\label{cor:UniCat3Gph7}
\textbf{(Uniformity of the quartet for \(\mathrm{III}.7\).)}
The components of the quartet, all with \(3\)-rank two,
share a common capitulation type
\(\varkappa(k_\mu)\),
common abelian type invariants
\(\alpha(k_\mu)\),
and a common second \(3\)-class group
\(\mathrm{Gal}(\mathrm{F}_3^2(k_\mu)/k_\mu)\),
for \(1\le\mu\le 4\).
\end{corollary}

\begin{proof}
This is a consequence of Theorem
\ref{thm:Cat3Gph7}
and Table
\ref{tbl:UniCat3Gph7}.
\end{proof}


\begin{example}
\label{exm:Cat3Gph7}
We have found prototypes for Graph \(\mathrm{III}.7\) in the form of
minimal conductors for each scenario in Theorem 
\ref{thm:Cat3Gph7}
as follows.
There are \textbf{regular} cases:
\(c=4\,599\) with symbol
\(\lbrace 9\leftrightarrow 73\leftarrow 7\rbrace\), \(v^\ast=1\),
and
\(\mathfrak{G}=\mathfrak{M}=\langle 81,7\rangle^4\);
\(c=31\,707\) with symbol
\(\lbrace 9\leftrightarrow 271\leftarrow 13\rbrace\), \(v^\ast=2\),
and
\(\mathfrak{G}=\mathfrak{M}=\langle 81,7\rangle^4\);
\(c=76\,741\) with symbol
\(\lbrace 577\leftrightarrow 19\leftarrow 7\rbrace\), \(v^\ast=2\),
and
\(\mathfrak{M}=\langle 2187,65\vert 67\rangle^4\)
of elevated coclass \(3\); and
\(c=90\,243\) with symbol
\(\lbrace 271\leftrightarrow 9\leftarrow 37\rbrace\), \(v^\ast=2\),
and
\(\mathfrak{M}=\langle 729,41\rangle^4\).
There is also a \textbf{singular} case \(c=61\,243\) with symbol
\(\lbrace 673\leftrightarrow 7\leftarrow 13\rbrace\), \(v^\ast=3\),
and
\(\mathfrak{M}=\langle 2187,253\rangle^4\);
and \textbf{super-singular} cases
\(c=69\,979\) with symbol
\(\lbrace 769\leftrightarrow 7\leftarrow 13\rbrace\), \(v^\ast=4\),
and
\(\mathfrak{M}=\langle 6561,676\vert 677\rangle^4\)
of elevated coclass \(3\); and
\(c=86\,821\) with symbol
\(\lbrace 79\leftrightarrow 157\leftarrow 7\rbrace\), \(v^\ast=4\),
and
\(\mathfrak{M}=\langle 729,37..39\rangle^4\).
\end{example}


In Table
\ref{tbl:ProtoCat3Gph7},
we summarize the prototypes of Graph \(\mathrm{III}.7\)
in the same way as in Table
\ref{tbl:ProtoCat3Gph5}.

\renewcommand{\arraystretch}{1.1}

\begin{table}[hb]
\caption{Prototypes for Graph III.7}
\label{tbl:ProtoCat3Gph7}
\begin{center}
{\tiny
\begin{tabular}{|c|c||c|c|c|c|c|c|c|c|c|}
\hline
 No. & \(c\)        & \(q\leftrightarrow p\leftarrow r\)    & \(v^\ast\) & \(v\) & \(m,n\) & \(\tilde{v}\) & \(\tilde{m},\tilde{n}\) & capitulation type & \(\mathfrak{M}\) & \(\ell_3(k)\) \\
\hline
   1 &   \(4\,599\) & \(9\leftrightarrow 73\leftarrow 7\)   & \(1\)      & \(2\) & \(1,2\) & \(2\) & \(1,2\) & \(\mathrm{a}.3^\ast\) & \(\langle 81,7\rangle\)            & \(=2\) \\
   2 &  \(12\,051\) & \(13\leftrightarrow 103\leftarrow 9\) & \(1\)      & \(2\) & \(1,1\) & \(2\) & \(1,1\) & \(\mathrm{a}.3^\ast\) & \(\langle 81,7\rangle\)            & \(=2\) \\
   6 &  \(31\,707\) & \(9\leftrightarrow 271\leftarrow 13\) & \(2\)      & \(2\) & \(1,0\) & \(2\) & \(1,0\) & \(\mathrm{a}.3^\ast\) & \(\langle 81,7\rangle\)            & \(=2\) \\
  21 &  \(76\,741\) & \(577\leftrightarrow 19\leftarrow 7\) & \(2\)      & \(2\) & \(0,1\) & \(2\) & \(0,1\) & \(\mathrm{H}.4\)      & \(\langle 2187,65\vert 67\rangle\) & \(\ge 3\) \\
  27 &  \(90\,243\) & \(271\leftrightarrow 9\leftarrow 37\) & \(2\)      & \(2\) & \(0,1\) & \(2\) & \(0,1\) & \(\mathrm{d}.19\)     & \(\langle 729,41\rangle\)          & \(\ge 2\) \\
\hline
  13 &  \(61\,243\) & \(673\leftrightarrow 7\leftarrow 13\) & \(3\)      & \(3\) & \(0,1\) & \(3\) & \(1,0\) & \(\mathrm{b}.10\)     & \(\langle 2187,253\rangle\)          & \(\ge 2\) \\
  17 &  \(69\,979\) & \(769\leftrightarrow 7\leftarrow 13\) & \(4\)      & \(3\) & \(0,1\) & \(3\) & \(0,1\) & \(\mathrm{d}.19\)     & \(\langle 6561,676\vert 677\rangle\) & \(\ge 3\) \\
  25 &  \(86\,821\) & \(79\leftrightarrow 157\leftarrow 7\) & \(4\)      & \(3\) & \(1,1\) & \(3\) & \(1,2\) & \(\mathrm{b}.10\)     & \(\langle 729,37..39\rangle\)        & \(\ge 2\) \\
\hline
\end{tabular}
}
\end{center}
\end{table}


\subsection{Category III, Graph 8}
\label{ss:Cat3Gph8}

\noindent
Let \((k_1,\ldots,k_4)\) be a quartet of cyclic cubic number fields
belonging to Graph \(8\) of Category \(\mathrm{III}\)
with combined cubic residue symbol
\(\lbrack p,q,r\rbrack_3=\lbrace r\rightarrow p\leftrightarrow q\leftarrow r\rbrace\)
of three prime(power)s dividing the conductor \(c=pqr\).


\begin{proposition}
\label{prp:Cat3Gph8}
\textbf{(Quartet with \(3\)-rank two for \(\mathrm{III}.8\).)}
For fixed \(\mu\in\lbrace 1,2,3,4\rbrace\),
let \(\mathfrak{p},\mathfrak{q},\mathfrak{r}\) be the prime ideals of \(k_\mu\)
over \(p,q,r\), that is,
\(p\mathcal{O}_{k_\mu}=\mathfrak{p}^3\),
\(q\mathcal{O}_{k_\mu}=\mathfrak{q}^3\),
\(r\mathcal{O}_{k_\mu}=\mathfrak{r}^3\),
then the \textbf{principal factor} of \(k_\mu\) is
\(A(k_\mu)=r\),
and the \(3\)-class group of \(k_\mu\) is
\begin{equation}
\label{eqn:GenCat3Gph8}
\mathrm{Cl}_3(k_\mu)=
\langle\lbrack\mathfrak{p}\rbrack,\lbrack\mathfrak{q}\rbrack\rangle\simeq (3,3).
\end{equation}
The unramified cyclic cubic relative extensions of \(k_\mu\)
are among the absolutely bicyclic bicubic fields \(B_i\), \(1\le i\le 10\).
The \textbf{wild} ranks for \(i=1,2,3,4,7,10\) are \(r_i=2\), independently of
\(m,n,\tilde{m},\tilde{n}\).
For each \(\mu\), there are two tame extensions \(B_j/k_\mu\), \(B_\ell/k_\mu\)
with the following properties. \\
The first, \(B_j\),
has norm class group
\(N_{B_j/k_\mu}(\mathrm{Cl}_3(B_j))=\langle\lbrack\mathfrak{p}\rbrack\rangle\),
\textbf{cyclic} transfer kernel
\begin{equation}
\label{eqn:Trans1Cat3Gph8}
\ker(T_{B_j/k_\mu})=\langle\lbrack\mathfrak{q}\rbrack\rangle
\end{equation}
of order \(3\),
and \textbf{elementary tricyclic} \(3\)-class group
\(\mathrm{Cl}_3(B_j)=
\langle\lbrack\mathfrak{P}_1\rbrack,\lbrack\mathfrak{P}_2\rbrack,\lbrack\mathfrak{P}_3\rbrack\rangle\simeq (3,3,3)\),
generated by the classes of the prime ideals of \(B_j\) over
\(\mathfrak{p}\mathcal{O}_{B_j}=\mathfrak{P}_1\mathfrak{P}_2\mathfrak{P}_3\). \\
The second, \(B_\ell\),
has norm class group
\(N_{B_\ell/k_\mu}(\mathrm{Cl}_3(B_\ell))=\langle\lbrack\mathfrak{q}\rbrack\rangle\),
\textbf{cyclic} transfer kernel
\begin{equation}
\label{eqn:Trans2Cat3Gph8}
\ker(T_{B_\ell/k_\mu})=\langle\lbrack\mathfrak{p}\rbrack\rangle
\end{equation}
of order \(3\),
and \textbf{elementary tricyclic} \(3\)-class group
\(\mathrm{Cl}_3(B_\ell)=
\langle\lbrack\mathfrak{Q}_1\rbrack,\lbrack\mathfrak{Q}_2\rbrack,\lbrack\mathfrak{Q}_3\rbrack\rangle\simeq (3,3,3)\),
generated by the classes of the prime ideals of \(B_\ell\) over
\(\mathfrak{q}\mathcal{O}_{B_\ell}=\mathfrak{Q}_1\mathfrak{Q}_2\mathfrak{Q}_3\). \\
The pair \((j,\ell)\) forms a mandatory \textbf{transposition}
of the transfer kernel type \(\varkappa(k_\mu)\).
The remaining two \(B_i>k_\mu\), \(i\ne j\), \(i\ne\ell\),
have norm class group 
\(\langle\lbrack\mathfrak{pq}\rbrack\rangle\), respectively
\(\langle\lbrack\mathfrak{pq}^2\rbrack\rangle\),
necessarily \textbf{non-elementary bicyclic} \(3\)-class group of order
\(27\mid h_3(B_i)=\frac{(U_i:V_i)}{3}h_3(k_{pq})\),
respectively \(\frac{(U_i:V_i)}{3}h_3(\tilde{k}_{pq})\),
and transfer kernel
\[
\ker(T_{B_i/k_\mu})\ge\langle\lbrack\mathfrak{p}^m\mathfrak{q}^n\rbrack\rangle, \text{ respectively }
\ge\langle\lbrack\mathfrak{p}^{\tilde{m}}\mathfrak{q}^{\tilde{n}}\rbrack\rangle,
\]
providing the option of a possible \textbf{fixed point}
in the transfer kernel type \(\varkappa(k_\mu)\).
\end{proposition}

\begin{proof}
Since \(r\rightarrow p\),
there are two principal factors \(A(k_{pr})=A(\tilde{k}_{pr})=r\).
Since \(q\leftarrow r\),
there are two further principal factors \(A(k_{qr})=A(\tilde{k}_{qr})=r\).
Since \(q\leftarrow r\rightarrow p\) is universally repelling,
we also have four other principal factors \(A(k_\mu)=r\), for all \(1\le\mu\le 4\)
according to
\cite[Prop. 4.6, p. 49]{Ay1995}.
Since \(\mathfrak{r}=\alpha\mathcal{O}_{k_\mu}\) is a principal ideal,
its class \(\lbrack\mathfrak{r}\rbrack=1\) is trivial,
whereas the classes \(\lbrack\mathfrak{p}\rbrack,\lbrack\mathfrak{q}\rbrack\)
are non-trivial and generate \(\mathrm{Cl}_3(k_\mu)\).

There are four \textit{tame} bicyclic bicubic fields,
\(B_5=k_1k_3k_p\tilde{k}_{qr}\),
\(B_6=k_1k_4k_q\tilde{k}_{pr}\),
\(B_8=k_2k_4k_pk_{qr}\),
\(B_9=k_2k_3k_qk_{pr}\),
satisfying \(9\mid h_3(B_i)=(U_i:V_i)\), for \(i\in\lbrace 5,6,8,9\rbrace\).
Consequently, we must have the indices \(I_i=(U_i:V_i)\in\lbrace 9,27\rbrace\),
and thus the matrix ranks \(r_5=r_6=r_8=r_9=2\).

On the other hand, there are six \textit{wild} bicyclic bicubic fields,
\(B_1=k_1k_{pq}k_{pr}k_{qr}\),
\(B_2=k_2k_{pq}\tilde{k}_{pr}\tilde{k}_{qr}\),
\(B_3=k_3\tilde{k}_{pq}\tilde{k}_{pr}k_{qr}\),
\(B_4=k_4\tilde{k}_{pq}k_{pr}\tilde{k}_{qr}\),
\(B_7=k_1k_2k_r\tilde{k}_{pq}\),
\(B_{10}=k_3k_4k_rk_{pq}\),
with \(h_3(B_i)>(U_i:V_i)\).

For these bicyclic bicubic fields \(B_i\),
\(i\in\lbrace 1,2,3,4,7,10\rbrace\),
the rank \(r_i\) is calculated with row operations
on the associated principal factor matrices \(M_i\): \\
\(M_1=M_2=
\begin{pmatrix}
0 & 0 & 1 \\
m & n & 0 \\
0 & 0 & 1 \\
0 & 0 & 1 \\
\end{pmatrix}\),
\(M_{10}=
\begin{pmatrix}
0 & 0 & 1 \\
0 & 0 & 1 \\
0 & 0 & 1 \\
m & n & 0 \\
\end{pmatrix}\),
\(M_3=M_4=
\begin{pmatrix}
0 & 0 & 1 \\
\tilde{m} & \tilde{n} & 0 \\
0 & 0 & 1 \\
0 & 0 & 1 \\
\end{pmatrix}\),
\(M_7=
\begin{pmatrix}
0 & 0 & 1 \\
0 & 0 & 1 \\
0 & 0 & 1 \\
\tilde{m} & \tilde{n} & 0 \\
\end{pmatrix}\).

For \(B_1=k_1k_{pq}k_{pr}k_{qr}\)
and \(B_2=k_2k_{pq}\tilde{k}_{pr}\tilde{k}_{qr}\),
\(M_1=M_2\)
immediately leads to rank \(r_1=r_2=2\),
since \((m,n)\ne (0,0)\),
and similarly,
for \(B_{10}=k_3k_4k_rk_{pq}\),
\(M_{10}\)
leads to rank \(r_{10}=2\).

For \(B_3=k_3\tilde{k}_{pq}\tilde{k}_{pr}k_{qr}\), 
and \(B_4=k_4\tilde{k}_{pq}k_{pr}\tilde{k}_{qr}\),
\(M_3=M_4\)
immediately leads to rank \(r_3=r_4=2\),
since \((\tilde{m},\tilde{n})\ne (0,0)\),
and similarly,
for \(B_7=k_1k_2k_r\tilde{k}_{pq}\),
\(M_7\)
leads to rank \(r_7=2\). \\
So Graph 8 of Category III
is the \textbf{unique} situation where
\(r_i=2\), for all \(1\le i\le 10\),
without any conditions, and thus
\(h_3(B_i)=\frac{I_i}{3}h_3(k_{pq})\), for \(i\in\lbrace 1,2,10\rbrace\), and
\(h_3(B_i)=\frac{I_i}{3}h_3(\tilde{k}_{pq})\), for \(i\in\lbrace 3,4,7\rbrace\),
where \(I_i=(U_i:V_i)\in\lbrace 9,27\rbrace\), and
\(9\mid h_3(k_{pq})\), \(9\mid h_3(\tilde{k}_{pq})\).

In each case,
the minimal subfield unit index
\((U_i:V_i)=9\)
corresponds to
the maximal unit norm index \((U(k_\mu):N_{B_i/k_\mu}(U_i))=3\),
associated with a \textit{total} transfer kernel
\(\#\ker(T_{B_i/k_\mu})=9\),
whenever \(k_\mu<B_i\) \(1\le\mu\le 4\), \(1\le i\le 10\).

According to Theorem
\ref{thm:Three},
the unramified cyclic cubic relative extensions 
of \(k_\mu\)
among the absolutely bicyclic bicubic subfields
of the \(3\)-genus field \(k^\ast=k_pk_qk_r\) are
\(B_1,B_5,B_6,B_7\), for \(\mu=1\),
\(B_2,B_7,B_8,B_9\), for \(\mu=2\),
\(B_3,B_5,B_9,B_{10}\), for \(\mu=3\), and
\(B_4,B_6,B_8,B_{10}\), for \(\mu=4\).

Since \(p\) splits in \(k_{q}\),
it also splits in
\(B_6=k_1k_4k_q\tilde{k}_{pr}\) and
\(B_9=k_2k_3k_qk_{pr}\).

Since \(\mathfrak{q}\) is principal in \(k_q\),
\(\lbrack\mathfrak{q}\rbrack\) capitulates in
\(B_6=k_1k_4k_q\tilde{k}_{pr}\) and
\(B_9=k_2k_3k_qk_{pr}\).

For \((\mu,j)\in\lbrace(1,6),(4,6),(2,9),(3,9)\rbrace\),
the minimal unit norm index \((U(k_\mu):N_{B_j/k_\mu}(U_j))=1\),
associated to the \textit{partial} transfer kernel
\(\ker(T_{B_j/k_\mu})=\langle\lbrack\mathfrak{q}\rbrack\rangle\),
corresponds to the maximal subfield unit index
\(h_3(B_j)=(U_j:V_j)=27\),
giving rise to the characteristic abelian type invariants
\(\mathrm{Cl}_3(B_j)=
\langle\lbrack\mathfrak{P}_1\rbrack,\lbrack\mathfrak{P}_2\rbrack,\lbrack\mathfrak{P}_3\rbrack\rangle\simeq (3,3,3)\)
generated by the classes of the prime ideals of \(B_j\) over
\(\mathfrak{p}\mathcal{O}_{B_j}=\mathfrak{P}_1\mathfrak{P}_2\mathfrak{P}_3\).
The field \(B_j\), which contains \(k_\mu\),
has norm class group
\(N_{B_j/k_\mu}(\mathrm{Cl}_3(B_j))=\langle\lbrack\mathfrak{p}\rbrack\rangle\).

Since \(q\) splits in \(k_p\),
it also splits in 
\(B_5=k_1k_3k_p\tilde{k}_{qr}\) and
\(B_8=k_2k_4k_pk_{qr}\).

Since \(\mathfrak{p}\) is principal in \(k_p\),
\(\lbrack\mathfrak{p}\rbrack\) capitulates in
\(B_5=k_1k_3k_p\tilde{k}_{qr}\) and
\(B_8=k_2k_4k_pk_{qr}\).

For \((\mu,\ell)\in\lbrace(1,5),(3,5),(2,8),(4,8)\rbrace\),
the minimal unit norm index \((U(k_\mu):N_{B_\ell/k_\mu}(U_\ell))=1\),
associated to the \textit{partial} transfer kernel
\(\ker(T_{B_\ell/k_\mu})=\langle\lbrack\mathfrak{p}\rbrack\rangle\),
corresponds to the maximal subfield unit index
\(h_3(B_\ell)=(U_\ell:V_\ell)=27\),
giving rise to the characteristic abelian type invariants
\(\mathrm{Cl}_3(B_\ell)=
\langle\lbrack\mathfrak{Q}_1\rbrack,\lbrack\mathfrak{Q}_2\rbrack,\lbrack\mathfrak{Q}_3\rbrack\rangle\simeq (3,3,3)\)
generated by the classes of the prime ideals of \(B_\ell\) over
\(\mathfrak{q}\mathcal{O}_{B_\ell}=\mathfrak{Q}_1\mathfrak{Q}_2\mathfrak{Q}_3\).
The field \(B_\ell\), which contains \(k_\mu\),
has norm class group
\(N_{B_\ell/k_\mu}(\mathrm{Cl}_3(B_\ell))=\langle\lbrack\mathfrak{q}\rbrack\rangle\).

Since \(\mathfrak{p}^m\mathfrak{q}^n\) is principal in \(k_{pq}\),
\(\lbrack\mathfrak{p}^m\mathfrak{q}^n\rbrack\) capitulates in
\(B_1=k_1k_{pq}k_{pr}k_{qr}\),
\(B_2=k_2k_{pq}\tilde{k}_{pr}\tilde{k}_{qr}\), and
\(B_{10}=k_3k_4k_rk_{pq}\).

Since \(\mathfrak{p}^{\tilde{m}}\mathfrak{q}^{\tilde{n}}\) is principal in \(\tilde{k}_{pq}\),
\(\lbrack\mathfrak{p}^{\tilde{m}}\mathfrak{q}^{\tilde{n}}\rbrack\) capitulates in
\(B_3=k_3\tilde{k}_{pq}\tilde{k}_{pr}k_{qr}\),
\(B_4=k_4\tilde{k}_{pq}k_{pr}\tilde{k}_{qr}\), and
\(B_7=k_1k_2k_r\tilde{k}_{pq}\).

The remaining two \(B_i>k_\mu\), \(i\in\lbrace 1,2,3,4,7,10\rbrace\),
more precisely, \(i\in\lbrace 1,7\rbrace\) for \(\mu=1\),
and \(i\in\lbrace 2,7\rbrace\) for \(\mu=2\),
and \(i\in\lbrace 3,10\rbrace\) for \(\mu=3\),
and \(i\in\lbrace 4,10\rbrace\) for \(\mu=4\),
have norm class group 
\(\langle\lbrack\mathfrak{pq}\rbrack\rangle\), respectively
\(\langle\lbrack\mathfrak{pq}^2\rbrack\rangle\),
and \textit{minimal} transfer kernel
\(\ker(T_{B_i/k_\mu})\ge\langle\lbrack\mathfrak{p}^m\mathfrak{q}^n\rbrack\rangle\), respectively
\(\ker(T_{B_i/k_\mu})\ge\langle\lbrack\mathfrak{p}^{\tilde{m}}\mathfrak{q}^{\tilde{n}}\rbrack\rangle\).
\end{proof}


Proposition
\ref{prp:Cat3Gph8}
and parts of its proof are now summarized in Table
\ref{tbl:UniCat3Gph8},
with transposition in \textbf{boldface} font,
based on Corollary
\ref{cor:Three}.
In this table,
we give the norm class group (NCG)
\(N_{B_i/k_\mu}(\mathrm{Cl}_3(B_i))\)
and the transfer kernel (TK)
\(\ker(T_{B_i/k_\mu})\),
also in the symbolic form \(\varkappa\)
with place holders \(1\le x,\tilde{x},y,\tilde{y},z,\tilde{z},w,\tilde{w}\le 4\),
for each collection of four unramified cyclic cubic relative extensions
\(B_j\), \(j=1,\ldots,10\),
of each base field \(k_\mu\), \(\mu=1,\ldots,4\), of the quartet.


\renewcommand{\arraystretch}{1.1}

\begin{table}[ht]
\caption{Norm class groups and minimal transfer kernels for Graph III.8}
\label{tbl:UniCat3Gph8}
\begin{center}
{\scriptsize
\begin{tabular}{|c||c|c|c|c||c|c|c|c||c|c|c|c||c|c|c|c|}
\hline
 Base         & \multicolumn{4}{c||}{\(k_1\)} & \multicolumn{4}{c||}{\(k_2\)} & \multicolumn{4}{c||}{\(k_3\)} & \multicolumn{4}{c|}{\(k_4\)} \\
\hline
 Ext          & \(B_1\) & \(B_5\) & \(B_6\) & \(B_7\) & \(B_2\) & \(B_7\) & \(B_8\) & \(B_9\) & \(B_3\) & \(B_5\) & \(B_9\) & \(B_{10}\) & \(B_4\) & \(B_6\) & \(B_8\) & \(B_{10}\) \\
\hline
 NCG          & \(\mathfrak{pq}\)   & \(\mathfrak{q}\) & \(\mathfrak{p}\) & \(\mathfrak{pq}^2\) & \(\mathfrak{pq}\) & \(\mathfrak{pq}^2\) & \(\mathfrak{q}\) & \(\mathfrak{p}\)    
              & \(\mathfrak{pq}\) & \(\mathfrak{q}\) & \(\mathfrak{p}\) & \(\mathfrak{pq}^2\)   & \(\mathfrak{pq}\) & \(\mathfrak{p}\)    & \(\mathfrak{q}\) & \(\mathfrak{pq}^2\) \\
 TK           & \(\mathfrak{p}^m\mathfrak{q}^n\) & \(\mathfrak{p}\) & \(\mathfrak{q}\) & \(\mathfrak{p}^{\tilde{m}}\mathfrak{q}^{\tilde{n}}\) & \(\mathfrak{p}^m\mathfrak{q}^n\) & \(\mathfrak{p}^{\tilde{m}}\mathfrak{q}^{\tilde{n}}\)    & \(\mathfrak{p}\) & \(\mathfrak{q}\) 
              & \(\mathfrak{p}^{\tilde{m}}\mathfrak{q}^{\tilde{n}}\) & \(\mathfrak{p}\) & \(\mathfrak{q}\) & \(\mathfrak{p}^m\mathfrak{q}^n\) & \(\mathfrak{p}^{\tilde{m}}\mathfrak{q}^{\tilde{n}}\) & \(\mathfrak{q}\) & \(\mathfrak{p}\) & \(\mathfrak{p}^m\mathfrak{q}^n\)    \\
\(\varkappa\) & \(x\) & \(\mathbf{3}\) & \(\mathbf{2}\) & \(\tilde{x}\) & \(y\) & \(\tilde{y}\) & \(\mathbf{4}\) & \(\mathbf{3}\) & \(\tilde{z}\) & \(\mathbf{3}\) & \(\mathbf{2}\) & \(z\) & \(\tilde{w}\) & \(\mathbf{3}\) & \(\mathbf{2}\) & \(w\) \\
\hline
\end{tabular}
}
\end{center}
\end{table}


\begin{theorem}
\label{thm:Cat3Gph8}
\textbf{(Second \(3\)-class group for \(\mathrm{III}.8\).)}
Let \((k_1,\ldots,k_4)\) be the quartet of cyclic cubic number fields
sharing the common conductor \(c=pqr\),
belonging to Graph \(8\) of Category \(\mathrm{III}\)
with combined cubic residue symbol
\(\lbrack p,q,r\rbrack_3=\lbrace r\rightarrow p\leftrightarrow q\leftarrow r\rbrace\).

To identify the second \(3\)-class group
\(\mathfrak{M}=\mathrm{Gal}(\mathrm{F}_3^2(k_\mu)/k_\mu)\), \(1\le\mu\le 4\),
let the \textbf{principal factor} of \(k_{pq}\), respectively \(\tilde{k}_{pq}\), be
\(A(k_{pq})=p^mq^n\), respectively \(A(\tilde{k}_{pq})=p^{\tilde{m}}q^{\tilde{n}}\), and additionally assume the \textbf{regular} situation where both
\(\mathrm{Cl}_3(k_{pq})\simeq\mathrm{Cl}_3(\tilde{k}_{pq})\simeq (3,3)\)
are elementary bicyclic, whence \((m,n)=(\tilde{m},\tilde{n})\).

Then there are several
\textbf{minimal transfer kernel types} (mTKT) \(\varkappa_0\) of \(k_\mu\), \(1\le\mu\le 4\),
and the other possible capitulation types in ascending order
\(\varkappa_0<\varkappa^\prime<\varkappa^{\prime\prime}\),
ending in the mandatory \(\varkappa^{\prime\prime}=(2100)\), type \(\mathrm{b}.10\): \\
either \(\varkappa_0=(2111)\), type \(\mathrm{H}.4\),
\(\varkappa^\prime=(2110)\), type \(\mathrm{d}.19\),
for \(\mathcal{P}=1\),
or \(\varkappa_0=(2133)\), type \(\mathrm{F}.11\),
\(\varkappa^\prime=(2130)\), type \(\mathrm{d}.23\), or \((2103)\), type \(\mathrm{d}.25\),
for \(\mathcal{P}=2\),
and the second \(3\)-class group is \(\mathfrak{M}\simeq\)

\begin{equation}
\label{eqn:RegCat3Gph8}
\begin{cases}
\langle 729,34..39\rangle,\ \alpha=\lbrack 111,111,21,21\rbrack,\ \varkappa=(2100) & \text{ if } \mathcal{N}=2, \\
\langle 729,41\rangle,\ \alpha=\lbrack 111,111,22,21\rbrack,\ \varkappa=(2110) & \text{ if } \mathcal{P}=1,\ \mathcal{N}=3, \\
\langle 729,42\vert 43\rangle,\ \alpha=\lbrack 111,111,22,21\rbrack,\ \varkappa=(2130)\vert(2140) & \text{ if } \mathcal{P}=2,\ \mathcal{N}=3, \\
\langle 2187,65\vert 67\rangle,\ \alpha=\lbrack 111,111,22,22\rbrack,\ \varkappa=(2111) & \text{ if } \mathcal{P}=1,\ \mathcal{N}=4, \\
\langle 2187,66\vert 73\rangle,\ \alpha=\lbrack 111,111,22,22\rbrack,\ \varkappa=(2133) & \text{ if } \mathcal{P}=2,\ \mathcal{N}=4,
\end{cases}
\end{equation}
where \(\mathcal{N}:=\#\lbrace 1\le j\le 10\mid k_\mu<B_j,\ I_j=27\rbrace\)
and \(\mathcal{P}\) is the number of prime divisors of \(p^mq^n\).
In any case, the \(3\)-class field tower may have a group
\(\mathfrak{G}=\mathrm{Gal}(\mathrm{F}_3^\infty(k_\mu)/k_\mu)\)
bigger than \(\mathfrak{M}\)
although \(d_2(\mathfrak{M})\le 4\).
Further, \(\mathrm{III}.8\) is the unique graph where the second \(3\)-class group
\(\mathfrak{M}=\mathrm{Gal}(\mathrm{F}_3^2(k_\mu)/k_\mu)\)
\textbf{cannot be of maximal class}.

Since the group order cannot be specified
in the \textbf{(super-)singular} situation,
only the capitulation type can be given.
Additionally, the number \(\tilde{\mathcal{P}}\)
of prime divisors of \(p^{\tilde{m}}q^{\tilde{n}}\) is used,
and two cases are separated. \\
If \((m,n)=(\tilde{m},\tilde{n})\), then \(\mathcal{P}=\tilde{\mathcal{P}}\),
and all four types \(\varkappa(k_\mu)\), \(\mu=1,\ldots,4\), coincide:
\begin{equation}
\label{eqn:Sng1Cat3Gph8}
\begin{cases}
\varkappa=(2100),\ \mathrm{b}.10 & \text{ if } \mathcal{N}=2, \\
\varkappa=(2111),\ \mathrm{H}.4 & \text{ if } \mathcal{P}=1,\ \mathcal{N}=4, \\
\varkappa=(2133),\ \mathrm{F}.11 & \text{ if } \mathcal{P}=2,\ \mathcal{N}=4.
\end{cases}
\end{equation}
If \((m,n)\ne (\tilde{m},\tilde{n})\), then
the number of coinciding types must be indicated by formal exponents:
\begin{equation}
\label{eqn:Sng2Cat3Gph8}
\begin{cases}
\varkappa=(2100)^4,\ \mathrm{b}.10 & \text{ if } \mathcal{N}=2, \\
\varkappa=(2130)^2,\ \mathrm{d}.23,\
\varkappa=(2140)^2,\ \mathrm{d}.25 & \text{ if } \mathcal{N}=3, \\
\varkappa=(2112)^4,\ \mathrm{F}.7  & \text{ if } \mathcal{P}=\tilde{\mathcal{P}}=1,\ \mathcal{N}=4, \\
\varkappa=(2134)^2,\ \mathrm{G}.16,\
\varkappa=(2143)^2,\ \mathrm{G}.19 & \text{ if } \mathcal{P}=\tilde{\mathcal{P}}=2,\ \mathcal{N}=4, \\
\varkappa=(2131)^2,\ \mathrm{F}.12,\
\varkappa=(2113)^2,\ \mathrm{F}.13 & \text{ if } \mathcal{P}\ne\tilde{\mathcal{P}},\ \mathcal{N}=4. \\
\end{cases}
\end{equation}
\end{theorem}

\begin{proof}
We normalize the transpositions in Table
\ref{tbl:UniCat3Gph8}
by the following convention
\(\varkappa=(21\ast\ast)\sim (\ast 32\ast)\sim (\ast\ast 43)\),
taking the leading type of three equivalent types.

Since the transfer kernels \(\ker(T_{B_i/k_\mu})\)
for the \textit{tame} extensions with \(i\in\lbrace 5,6,8,9\rbrace\) are \textit{partial},
the corresponding indices \(I_i=(U_i:V_i)=27\) of subfield units
must be maximal, whence necessarily \(\mathcal{N}\ge 2\).
The associated \(3\)-class numbers \(h_3(B_i)=27\) are consistent
with occurrence of \textit{two elementary tricyclic} \(3\)-class groups
\(\mathrm{Cl}_3(B_i)\simeq (3,3,3)\),
connected with a \textit{transposition} in the capitulation type
\(\varkappa(k_\mu)=(21\ast\ast)\),
according to Proposition
\ref{prp:Cat3Gph8}.

In the proof of this proposition,
it was also derived that,
due to \(r_i=2\),
the \(3\)-class numbers of the \textit{wild} extensions are given by
\(h_3(B_i)=\frac{I_i}{3}h_3(k_{pq})\ge 27\), for \(i\in\lbrace 1,2,10\rbrace\), and
\(h_3(B_i)=\frac{I_i}{3}h_3(\tilde{k}_{pq})\ge 27\), for \(i\in\lbrace 3,4,7\rbrace\),
where \(I_i=(U_i:V_i)\in\lbrace 9,27\rbrace\), and
\(9\mid h_3(k_{pq})\), \(9\mid h_3(\tilde{k}_{pq})\).

It follows that maximal class \(\mathrm{cc}(\mathfrak{M})=1\)
is prohibited for two reasons,
firstly by the Artin pattern
\(\varkappa\sim (21\varkappa_3\varkappa_4)\),
\(\alpha\sim\lbrack 111,111,\alpha_3,\alpha_4\rbrack\),
and secondly by the bi-polarization of order at least \(27\),
which implies that \(\alpha_3\) and \(\alpha_4\)
are bicyclic equal to \((21)\) or bigger
\cite[pp. 289--292]{Ma2015b}.

In fact, even \(\mathrm{cc}(\mathfrak{M})=2\) is very restricted,
because the candidates for \(\mathfrak{M}\) must be descendants
of the group \(\langle 243,3\rangle\).
The other two groups with two or three components \((111)\)
in the abelian type invariants \(\alpha\) are discouraged,
since \(\varkappa\sim (1133)\), \(\alpha\sim\lbrack 21,111,21,111\rbrack\)
for \(\langle 243,7\rangle\)
does not contain a transposition,
and in \(\varkappa\sim (2111)\), \(\alpha\sim\lbrack 111,21,111,111\rbrack\)
for \(\langle 243,4\rangle\),
the transposition in \(\varkappa\) is not associated with
two elementary tricyclic components of \(\alpha\). 

If \(\mathcal{N}=2\),
then the Artin pattern \(\alpha=\lbrack 111,111,21,21\rbrack\), \(\varkappa=(2100)\)
identifies one of the six groups
\(\langle 729,34..39\rangle,\ \alpha=\lbrack 111,111,21,21\rbrack\),
since \(\langle 243,3\rangle\) is forbidden by Corollary
\ref{cor:ClassGroup33}.

If \(\mathcal{N}=3\),
then generally \(\alpha=\lbrack 111,111,22,21\rbrack\),
with bipolarization consisting of copolarization \((21)\), i.e. coclass \(2\),
and polarization \((22)\), i.e., class \(4\).
Now, if \(\mathcal{P}=1\), then the capitulation type \(\varkappa=(2110)\sim (2120)\)
contains a repetition, which identifies the group \(\langle 729,41\rangle\).
On the other hand, if \(\mathcal{P}=2\), then the capitulation type \(\varkappa=(2130)\)
either contains a fixed point, which gives \(\langle 729,42\rangle\),
or \(\varkappa=(2140)\) neither contains a repetition nor a fixed point,
which gives \(\langle 729,43\rangle\).

If \(\mathcal{N}=4\),
then generally \(\alpha=\lbrack 111,111,22,22\rbrack\),
but a finer distinction is provided by \(\mathcal{P}\).
If \(\mathcal{P}=1\), then the capitulation type \(\varkappa=(2111)\sim (2122)\)
contains two repetitions and becomes nearly constant,
which identifies the groups \(\langle 2187,65\vert 67\rangle\) of coclass \(3\).
However, if \(\mathcal{P}=2\), then a fixed point and its repetition occurs
in the capitulation type \(\varkappa=(2133)\sim (2144)\),
which leads to the groups
\(\langle 2187,66\vert 73\rangle,\ \alpha=\lbrack 111,111,22,22\rbrack\).

\noindent
Concerning the \textbf{(super-)singular} situation, two cases are distinguished.
If \((m,n)=(\tilde{m},\tilde{n})\), then
\(\mathcal{P}=1\), i.e., \((m,n)\in\lbrace (0,1),(1,0)\rbrace\),
implies two identical repetitions in
\(\varkappa_0\sim (2111)\sim (2122)\), \(\mathrm{H}.4\); but 
\(\mathcal{P}=2\), i.e., \((m,n)\in\lbrace (1,1),(1,2),(2,1)\rbrace\),
produces a single fixed point in
\(\varkappa_0\sim (2133)\sim (2144)\), \(\mathrm{F}.11\).
These two minimal transfer kernel types for \(\mathcal{N}=4\)
both expand to \(\varkappa^{\prime\prime}\sim (2100)\)
for \(\mathcal{N}=2\).
All three cases are uniform.
If \(\mathcal{N}=3\) were possible, then
\(\mathcal{P}=1\) would lead to
type \(\varkappa\sim (2110)\sim (2120)\), \(\mathrm{d}.19\),
and \(\mathcal{P}=2\) would either imply type
\(\varkappa\sim (2130)\), \(\mathrm{d}.23\),
or \(\varkappa\sim (2140)\), \(\mathrm{d}.25\).
The latter case would be non-uniform,
but \(\mathcal{N}=3\) does not seem to occur at all.

If \((m,n)\ne (\tilde{m},\tilde{n})\), then
nevertheless \(\mathcal{P}=\tilde{\mathcal{P}}\) is possible,
and then \(\mathcal{P}=1\)
implies two distinct repetitions in
\(\varkappa\sim (2112)\sim (2121)\), \(\mathrm{F}.7\), uniformly,
whereas \(\mathcal{P}=2\) leads to either two fixed points in
\(\varkappa\sim (2134)\), \(\mathrm{G}.16\)
or a second transposition in
\(\varkappa\sim (2143)\), \(\mathrm{G}.19\).
These permutation types would be
non-uniform in two sub-doublets,
but they are obviously forbidden,
for an unknown reason.
Finally, \(\mathcal{P}\ne\tilde{\mathcal{P}}\)
admits several distinct realizations
with identical result:
it always leads to a repetition,
and additionally
either to a fixed point in \(\varkappa\sim (2131)\), \(\mathrm{F}.12\),
or a non-fixed point in \(\varkappa\sim (2131)\), \(\mathrm{F}.13\),
non-uniformly in two sub-doublets.
\end{proof}


\begin{corollary}
\label{cor:UniCat3Gph8}
\textbf{(Non-uniformity of the quartet for \(\mathrm{III}.8\).)}
If \((m,n)=(\tilde{m},\tilde{n})\), in particular always
in the \textbf{regular} situation,
the components of the quartet, all with \(3\)-rank two,
\textbf{uniformly} share a common capitulation type
\(\varkappa(k_\mu)\),
common abelian type invariants
\(\alpha(k_\mu)\),
and a common second \(3\)-class group
\(\mathrm{Gal}(\mathrm{F}_3^2(k_\mu)/k_\mu)\),
for \(1\le\mu\le 4\).
\textbf{Otherwise}, the invariants may be \textbf{non-uniform},
divided in two sub-doublets.
\end{corollary}

\begin{proof}
In the \textit{regular} case, we must have \((m,n)=(\tilde{m},\tilde{n})\).
All TKTs are either equivalent to 
F.11, with mandatory fixed point,
if \(\mathfrak{p}^m\mathfrak{q}^n\in\lbrace\mathfrak{p}\mathfrak{q},\mathfrak{p}\mathfrak{q}^2\rbrace\), or to
H.4, if \(\mathfrak{p}^m\mathfrak{q}^n\in\lbrace\mathfrak{p},\mathfrak{q}\rbrace\),
according to Table
\ref{tbl:UniCat3Gph8}.
The potential non-uniformity was proved in Theorem
\ref{thm:Cat3Gph8}.
\end{proof}


In Table
\ref{tbl:ProtoCat3Gph8},
we summarize the prototypes of Graph \(\mathrm{III}.8\)
in the same way as in Table
\ref{tbl:ProtoCat3Gph5}.
The group with multifurcation of order four is abbreviated by
\(P_7:=\langle 2187,64\rangle\).
See the table and tree diagram
\cite[\S\ 11, pp. 96--100, Tbl. 1, Fig. 5]{Ma2018}.

\renewcommand{\arraystretch}{1.1}

\begin{table}[ht]
\caption{Prototypes for Graph III.8}
\label{tbl:ProtoCat3Gph8}
\begin{center}
{\tiny
\begin{tabular}{|c|c||c|c|c|c|c|c|c|c|c|}
\hline
 No. & \(c\)         & \(r\rightarrow p\leftrightarrow q\leftarrow r\)      & \(v^\ast\) & \(v\) & \(m,n\) & \(\tilde{v}\) & \(\tilde{m},\tilde{n}\) & capitulation type & \(\mathfrak{M}\)                   & \(\ell_3(k)\) \\
\hline
   1 &  \(20\,293\)  & \(13\rightarrow 7\leftrightarrow 223\leftarrow 13\)  & \(1\) & \(2\) & \(2,1\) &  &  & \(\mathrm{b}.10\) & \(\langle 729,37..39\rangle\) & \(\ge 2\) \\
   2 &  \(41\,509\)  & \(31\rightarrow 13\leftrightarrow 103\leftarrow 31\) & \(1\) & \(2\) & \(1,1\) &  &  & \(\mathrm{b}.10\) & \(\langle 729,37..39\rangle\) & \(\ge 2\) \\
   3 &  \(46\,341\)  & \(19\rightarrow 9\leftrightarrow 271\leftarrow 19\)  & \(2\) & \(2\) & \(0,1\) &  &  & \(\mathrm{b}.10\) & \(\langle 729,37..39\rangle\) & \(\ge 2\) \\
   5 &  \(52\,497\)  & \(19\rightarrow 9\leftrightarrow 307\leftarrow 19\)  & \(1\) & \(2\) & \(2,1\) &  &  & \(\mathrm{F}.11\) & \(\langle 2187,66\vert 73\rangle\) & \(\ge 2\) \\
   7 &  \(92\,911\)  & \(13\rightarrow 7\leftrightarrow 1021\leftarrow 13\) & \(1\) & \(2\) & \(1,1\) &  &  & \(\mathrm{F}.11\) & \(\langle 2187,66\vert 73\rangle\) & \(\ge 2\) \\
  18 & \(191\,007\)  & \(19\rightarrow 9\leftrightarrow 1117\leftarrow 19\) & \(2\) & \(2\) & \(1,0\) &  &  & \(\mathrm{b}.10\) & \(\langle 729,37..39\rangle\) & \(\ge 2\) \\
\hline
  26 & \(231\,469\)  & \(43\rightarrow 7\leftrightarrow 769\leftarrow 43\)  & \(4\) & \(3\) & \(0,1\) & \(3\) & \(0,1\) & \(\mathrm{H}.4\)  & \(P_7-\#2;34\vert 35\) & \(\ge 3\) \\
  40 & \(387\,729\)  & \(9\rightarrow 67\leftrightarrow 643\leftarrow 9\)   & \(3\) & \(3\) & \(2,1\) & \(3\) & \(2,1\) & \(\mathrm{F}.11\) & \(P_7-\#2;36\vert 38\) & \(\ge 2\) \\
  92 & \(756\,499\)  & \(43\rightarrow 73\leftrightarrow 241\leftarrow 43\) & \(3\) & \(3\) & \(0,1\) & \(3\) & \(1,1\) & \(\mathrm{F}.12\) & \(P_7-\#2;43\vert 46\vert 51\vert 53\) & \(\ge 2\) \\
     &               &                                                      &       &       &         &       &         & \(\mathrm{F}.13\) & \(P_7-\#2;41\vert 47\vert 50\vert 52\) & \(\ge 3\) \\
  93 & \(758\,233\)  & \(7\rightarrow 19\leftrightarrow 5701\leftarrow 7\)  & \(3\) & \(3\) & \(1,1\) & \(3\) & \(1,1\) & \(\mathrm{F}.11\) & \(P_7-\#2;36\vert 38\) & \(\ge 2\) \\
 101 & \(806\,869\)  & \(7\rightarrow 73\leftrightarrow 1579\leftarrow 7\)  & \(5\) & \(3\) & \(1,1\) & \(5\) & \(1,1\) & \(\mathrm{F}.11\) & \(\) & \(\ge 2\) \\
 105 & \(831\,001\)  & \(67\rightarrow 79\leftrightarrow 157\leftarrow 67\) & \(4\) & \(3\) & \(1,1\) & \(3\) & \(2,1\) & \(\mathrm{d}.23\) & \(\langle 6561,678\rangle\) & \(\ge 3\) \\
     &               &                                                      &       &       &         &       &         & \(\mathrm{d}.25\) & \(\langle 6561,679\vert 680\rangle\) & \(\ge 3\) \\
 102 & \(945\,117\)  & \(19\rightarrow 9\leftrightarrow 5527\leftarrow 19\) & \(3\) & \(3\) & \(1,0\) & \(3\) & \(1,1\) & \(\mathrm{F}.12\) & \(P_7-\#2;43\vert 46\vert 51\vert 53\) & \(\ge 2\) \\
     &               &                                                      &       &       &         &       &         & \(\mathrm{F}.13\) & \(P_7-\#2;41\vert 47\vert 50\vert 52\) & \(\ge 3\) \\
 162 &\(1\,301\,287\)& \(31\rightarrow 13\leftrightarrow 3229\leftarrow 31\)& \(4\) & \(4\) & \(0,1\) & \(3\) & \(2,1\) & \(\mathrm{d}.23\) & \(\) & \(\ge 2\) \\
     &               &                                                      &       &       &         &       &         & \(\mathrm{d}.25\) & \(\) & \(\ge 2\) \\
 164 &\(1\,305\,937\)& \(31\rightarrow 103\leftrightarrow 409\leftarrow 31\)& \(6\) & \(4\) & \(1,1\) & \(5\) & \(1,1\) & \(\mathrm{F}.11\) & \(\) & \(\ge 2\) \\
 183 &\(1\,463\,917\)& \(13\rightarrow 7\leftrightarrow 16087\leftarrow 13\)& \(3\) & \(3\) & \(0,1\) & \(3\) & \(1,0\) & \(\mathrm{F}.7\) & \(P_7-\#2;55\vert 56\vert 58\) & \(\ge 2\) \\
 185 &\(1\,483\,767\)& \(19\rightarrow 9\leftrightarrow 8677\leftarrow 19\) & \(3\) & \(3\) & \(2,1\) & \(3\) & \(0,1\) & \(\mathrm{F}.12\) & \(P_7-\#2;43\vert 46\vert 51\vert 53\) & \(\ge 2\) \\
     &               &                                                      &       &       &         &       &         & \(\mathrm{F}.13\) & \(P_7-\#2;41\vert 47\vert 50\vert 52\) & \(\ge 2\) \\
 253 &\(2\,068\,587\)& \(19\rightarrow 9\leftrightarrow 12097\leftarrow 19\)& \(5\) & \(3\) & \(1,1\) & \(5\) & \(0,1\) & \(\mathrm{F}.12\) & \(\) & \(\ge 2\) \\
     &               &                                                      &       &       &         &       &         & \(\mathrm{F}.13\) & \(\) & \(\ge 2\) \\
 385 &\(2\,991\,987\)& \(19\rightarrow 9\leftrightarrow 17497\leftarrow 19\)& \(3\) & \(3\) & \(1,0\) & \(3\) & \(2,1\) & \(\mathrm{F}.12\) & \(P_7-\#2;43\vert 46\vert 51\vert 53\) & \(\ge 2\) \\
     &               &                                                      &       &       &         &       &         & \(\mathrm{F}.13\) & \(P_7-\#2;41\vert 47\vert 50\vert 52\) & \(\ge 2\) \\
 468 &\(3\,556\,699\)& \(97\rightarrow 37\leftrightarrow 991\leftarrow 97\) & \(6\) & \(5\) & \(1,1\) & \(3\) & \(2,1\) & \(\mathrm{d}.23\) & \(\) & \(\ge 2\) \\
     &               &                                                      &       &       &         &       &         & \(\mathrm{d}.25\) & \(\) & \(\ge 2\) \\
 651 &\(4\,686\,019\)&\(109\rightarrow 13\leftrightarrow 3307\leftarrow 109\)& \(4\)& \(4\) & \(0,1\) & \(3\) & \(1,0\) & \(\mathrm{F}.7\) & \(P_7-\#2;55\vert 56\vert 58\) & \(\ge 2\) \\
\hline
\end{tabular}
}
\end{center}
\end{table}


\begin{example}
\label{exm:Cat3Gph8}
Since \(\mathrm{III}.8\) is the graph with most sparse population by far, Ayadi
\cite[pp. 89--90]{Ay1995}
was unable to give any examples.
We found many, but not all, prototypes.
These are the minimal conductors for each scenario in Theorem
\ref{thm:Cat3Gph8}.
In the \textbf{regular} case,
they have been found for \(\mathcal{N}\in\lbrace 2,4\rbrace\),
but not for \(\mathcal{N}=3\).
There are some \textbf{regular} prototypes:
\(c=20\,293\) with symbol
\(\lbrace 13\rightarrow 7\leftrightarrow 223\leftarrow 13\rbrace\), \(v=1\),
and \(\mathfrak{M}=\langle 729,37..39\rangle\);
\(c=46\,341\) with symbol
\(\lbrace 19\rightarrow 9\leftrightarrow 271\leftarrow 19\rbrace\), \(v=2\),
and \(\mathfrak{M}=\langle 729,37..39\rangle\);
\(c=52\,497\) with symbol
\(\lbrace 19\rightarrow 9\leftrightarrow 307\leftarrow 19\rbrace\), \(v=1\),
and \(\mathfrak{M}=\langle 2187,66\vert 73\rangle\).
Furthermore, there is a \textbf{super-singular} prototype
\(c=231\,469\) with symbol
\(\lbrace 43\rightarrow 7\leftrightarrow 769\leftarrow 43\rbrace\),
type \(\mathrm{H}.4\),
and \(\mathfrak{M}=\langle 2187,64\rangle-\#2;i\), \(i\in\lbrace 34,35\rbrace\),
with \(d_2(\mathfrak{M})=5\), outside of the library
\cite{BEO2005},
not treated by Theorem
\ref{thm:Cat3Gph8}.
\end{example}


\subsection{Category III, Graph 9}
\label{ss:Cat3Gph9}

\noindent
Let \((k_1,\ldots,k_4)\) be a quartet of cyclic cubic number fields
sharing the common conductor \(c=pqr\),
belonging to Graph \(9\) of Category \(\mathrm{III}\)
with combined cubic residue symbol
\(\lbrack p,q,r\rbrack_3=\lbrace r\leftarrow p\leftrightarrow q\leftarrow r\rbrace\).


\begin{proposition}
\label{prp:Cat3Gph9}
\textbf{(Quartet with \(3\)-rank two for \(\mathrm{III}.9\).)}
For fixed \(\mu\in\lbrace 1,2,3,4\rbrace\),
let \(\mathfrak{p},\mathfrak{q},\mathfrak{r}\) be the prime ideals of \(k_\mu\)
over \(p,q,r\), that is,
\(p\mathcal{O}_{k_\mu}=\mathfrak{p}^3\),
\(q\mathcal{O}_{k_\mu}=\mathfrak{q}^3\),
\(r\mathcal{O}_{k_\mu}=\mathfrak{r}^3\),
then the \textbf{principal factor} of \(k_\mu\) is
\(A(k_\mu)=p\),
and the \(3\)-class group of \(k_\mu\) is
\begin{equation}
\label{eqn:GenCat3Gph9}
\mathrm{Cl}_3(k_\mu)=
\langle\lbrack\mathfrak{q}\rbrack,\lbrack\mathfrak{r}\rbrack\rangle\simeq (3,3).
\end{equation}
In terms of \(n\) and \(\tilde{n}\) in
\(A(k_{pq})=p^mq^n\) and \(A(\tilde{k}_{pq})=p^{\tilde{m}}q^{\tilde{n}}\),
the ranks of the \textbf{wild} extensions are
\begin{equation}
\label{eqn:RankCat3Gph9}
r_1=r_2=r_{10}=3 \text{ iff } n\ne 0 \text{ iff } q\mid A(k_{pq})
\text{ and }
r_3=r_4=r_7=3 \text{ iff } \tilde{n}\ne 0 \text{ iff } q\mid A(\tilde{k}_{pq}).
\end{equation}
\end{proposition}

\begin{proof}
By Proposition
\ref{prp:Principal2},
principal factors are
\(A(k_{pr})=A(\tilde{k}_{pr})=p\),
since \(r\leftarrow p\), and
\(A(k_{qr})=A(\tilde{k}_{qr})=r\),
since \(q\leftarrow r\).
Further, by Proposition
\ref{prp:Principal3},
\(A(k_\mu)=p\), for all \(1\le\mu\le 4\),
since \(p\) is universally repelling
\(r\leftarrow p\rightarrow q\).
Since \(\mathfrak{p}=\alpha\mathcal{O}_{k_\mu}\) is a principal ideal,
its class \(\lbrack\mathfrak{p}\rbrack=1\) is trivial,
whereas the classes
\(\lbrack\mathfrak{q}\rbrack,\lbrack\mathfrak{r}\rbrack\)
are non-trivial.
By Corollary
\ref{cor:Capitulation}, \\
since \(\mathfrak{q}\) is principal ideal in \(k_q\),
the class \(\lbrack\mathfrak{q}\rbrack\) capitulates in
\(B_6=k_1k_4k_q\tilde{k}_{pr}\) and
\(B_9=k_2k_3k_qk_{pr}\); \\
since \(\mathfrak{r}\) is principal ideal in \(k_r\),
the class \(\lbrack\mathfrak{r}\rbrack\) capitulates in
\(B_7=k_1k_2k_r\tilde{k}_{pq}\) and
\(B_{10}=k_3k_4k_rk_{pq}\). \\
However, since \(\mathfrak{r}\) is principal ideal in
\(k_{qr}\) and \(\tilde{k}_{qr}\),
the class \(\lbrack\mathfrak{r}\rbrack\) also capitulates in
\(B_1=k_1k_{pq}k_{pr}k_{qr}\),
\(B_2=k_2k_{pq}\tilde{k}_{pr}\tilde{k}_{qr}\),
\(B_3=k_3\tilde{k}_{pq}\tilde{k}_{pr}k_{qr}\),
\(B_4=k_4\tilde{k}_{pq}k_{pr}\tilde{k}_{qr}\), 
\(B_5=k_1k_3k_p\tilde{k}_{qr}\), and
\(B_8=k_2k_4k_pk_{qr}\).

For the wild bicyclic bicubic fields \(B_j\),
\(j\in\lbrace 1,2,3,4,7,10\rbrace\),
the rank \(r_j\) is calculated with row operations
on the associated principal factor matrices \(M_j\): \\
\(M_1=M_2=
\begin{pmatrix}
1 & 0 & 0 \\
m & n & 0 \\
1 & 0 & 0 \\
0 & 0 & 1 \\
\end{pmatrix}\),
\(M_{10}=
\begin{pmatrix}
1 & 0 & 0 \\
1 & 0 & 0 \\
0 & 0 & 1 \\
m & n & 0 \\
\end{pmatrix}\),
\(M_3=M_4=
\begin{pmatrix}
1 & 0 & 0 \\
\tilde{m} & \tilde{n} & 0 \\
1 & 0 & 0 \\
0 & 0 & 1 \\
\end{pmatrix}\),
\(M_7=
\begin{pmatrix}
1 & 0 & 0 \\
1 & 0 & 0 \\
0 & 0 & 1 \\
\tilde{m} & \tilde{n} & 0 \\
\end{pmatrix}\).

For \(B_1=k_1k_{pq}k_{pr}k_{qr}\)
and \(B_2=k_2k_{pq}\tilde{k}_{pr}\tilde{k}_{qr}\),
\(M_1=M_2\)
leads to the decisive pivot element \(n\)
in the middle column,
for \(B_{10}=k_3k_4k_rk_{pq}\),
\(M_{10}\)
also leads to \(n\).
So rank \(r_1=r_2=r_{10}=3\) iff \(n\ne 0\).

For \(B_3=k_3\tilde{k}_{pq}\tilde{k}_{pr}k_{qr}\)
and \(B_4=k_4\tilde{k}_{pq}k_{pr}\tilde{k}_{qr}\),
\(M_3=M_4\)
leads to the decisive pivot element \(\tilde{n}\)
in the middle column,
for \(B_7=k_1k_2k_r\tilde{k}_{pq}\),
\(M_7\)
also leads to \(\tilde{n}\).
So rank \(r_3=r_4=r_7=3\) iff \(\tilde{n}\ne 0\).
\end{proof}


In terms of capitulation targets in Corollary
\ref{cor:Three},
Proposition
\ref{prp:Cat3Gph9}
and parts of its proof
are now summarized in Table
\ref{tbl:UniCat3Gph9}
with transpositions in \textbf{bold} font.

\renewcommand{\arraystretch}{1.1}

\begin{table}[hb]
\caption{Norm class groups and minimal transfer kernels for Graph III.9}
\label{tbl:UniCat3Gph9}
\begin{center}
{\normalsize
\begin{tabular}{|c||c|c|c|c||c|c|c|c||c|c|c|c||c|c|c|c|}
\hline
 Base         & \multicolumn{4}{c||}{\(k_1\)} & \multicolumn{4}{c||}{\(k_2\)} & \multicolumn{4}{c||}{\(k_3\)} & \multicolumn{4}{c|}{\(k_4\)} \\
\hline
 Ext          & \(B_1\) & \(B_5\) & \(B_6\) & \(B_7\) & \(B_2\) & \(B_7\) & \(B_8\) & \(B_9\) & \(B_3\) & \(B_5\) & \(B_9\) & \(B_{10}\) & \(B_4\) & \(B_6\) & \(B_8\) & \(B_{10}\) \\
\hline
 NCG          & \(\mathfrak{qr}\)   & \(\mathfrak{q}\) & \(\mathfrak{r}\) & \(\mathfrak{qr}^2\) & \(\mathfrak{qr}^2\) & \(\mathfrak{qr}\) & \(\mathfrak{q}\) & \(\mathfrak{r}\)    
              & \(\mathfrak{qr}^2\) & \(\mathfrak{q}\) & \(\mathfrak{r}\) & \(\mathfrak{qr}\)   & \(\mathfrak{qr}\)   & \(\mathfrak{r}\)  & \(\mathfrak{q}\) & \(\mathfrak{qr}^2\) \\
 TK           & \(\mathfrak{r}\)    & \(\mathfrak{r}\) & \(\mathfrak{q}\) & \(\mathfrak{r}\)    & \(\mathfrak{r}\)    & \(\mathfrak{r}\)  & \(\mathfrak{r}\) & \(\mathfrak{q}\) 
              & \(\mathfrak{r}\)    & \(\mathfrak{r}\) & \(\mathfrak{q}\) & \(\mathfrak{r}\)    & \(\mathfrak{r}\)    & \(\mathfrak{q}\)  & \(\mathfrak{r}\) & \(\mathfrak{r}\) \\
\(\varkappa\) & \(3\) & \(\mathbf{3}\) & \(\mathbf{2}\) & \(3\) & \(4\) & \(4\) & \(\mathbf{4}\) & \(\mathbf{3}\) & \(3\) & \(\mathbf{3}\) & \(\mathbf{2}\) & \(3\) & \(2\) & \(\mathbf{3}\) & \(\mathbf{2}\) & \(2\) \\
\hline
\end{tabular}
}
\end{center}
\end{table}


\begin{theorem}
\label{thm:Cat3Gph9}
\textbf{(Second \(3\)-class group for \(\mathrm{III}.9\).)}
To identify the second \(3\)-class group
\(\mathfrak{M}=\mathrm{Gal}(\mathrm{F}_3^2(k_\mu)/k_\mu)\), \(1\le\mu\le 4\),
let the \textbf{principal factor} of \(k_{pq}\), respectively \(\tilde{k}_{pq}\), be
\(A(k_{pq})=p^mq^n\), respectively \(A(\tilde{k}_{pq})=p^{\tilde{m}}q^{\tilde{n}}\), and additionally assume the \textbf{regular} situation where both
\(\mathrm{Cl}_3(k_{pq})\simeq\mathrm{Cl}_3(\tilde{k}_{pq})\simeq (3,3)\)
are elementary bicyclic, whence \((m,n)=(\tilde{m},\tilde{n})\).

Then the
\textbf{minimal transfer kernel type} (mTKT) \(\varkappa_0\) of \(k_\mu\), \(1\le\mu\le 4\),
and other possible capitulation types in ascending order
\(\varkappa_0<\varkappa^\prime<\varkappa^{\prime\prime}<\varkappa^{\prime\prime\prime}\),
ending in the mandatory \(\varkappa^{\prime\prime\prime}=(2000)\), type \(\mathrm{a}.3^\ast\),
are \(\varkappa_0=(2111)\), type \(\mathrm{H}.4\),
\(\varkappa^\prime=(2110)\), type \(\mathrm{d}.19\),
\(\varkappa^{\prime\prime}=(2100)\), type \(\mathrm{b}.10\),
and the second \(3\)-class group is \(\mathfrak{M}\simeq\)

\begin{equation}
\label{eqn:Cat3Gph9}
\begin{cases}
\langle 81,7\rangle,\ \alpha=\lbrack 111,11,11,11\rbrack,\ \varkappa=(2000) & \text{ if } n\ne 0,\ \mathcal{N}=1, \\
\langle 729,34..39\rangle,\ \alpha=\lbrack 111,111,21,21\rbrack,\ \varkappa=(2100) & \text{ if } n=0,\ \mathcal{N}=2, \\
\langle 729,41\rangle,\ \alpha=\lbrack 111,111,22,21\rbrack,\ \varkappa=(2110) & \text{ if } n=0,\ \mathcal{N}=3, \\
\langle 2187,65\vert 67\rangle,\ \alpha=\lbrack 111,111,22,22\rbrack,\ \varkappa=(2111) & \text{ if } n=0,\ \mathcal{N}=4, \\
\end{cases}
\end{equation}
where \(\mathcal{N}:=\#\lbrace 1\le j\le 10\mid k_\mu<B_j,\ I_j=27\rbrace\).
Only in the first case, the \(3\)-class field tower has certainly the group
\(\mathfrak{G}=\mathrm{Gal}(\mathrm{F}_3^\infty(k_\mu)/k_\mu)\simeq\mathfrak{M}\)
and length \(\ell_3(k_\mu)=2\),
otherwise \(\ell_3(k_\mu)\ge 3\) cannot be excluded,
even if \(d_2(\mathfrak{M})\le 4\).
\end{theorem}

\begin{proof}
The essence of the proof is
a systematic evaluation of the facts proved in Proposition
\ref{prp:Cat3Gph9}
and illustrated by Table
\ref{tbl:UniCat3Gph9},
ordered by increasing indices \(I_j:=(U_j:V_j)\) of subfield units
and, accordingly, by Lemma
\ref{lem:UnitIndices},
shrinking transfer kernels \(\ker(T_{B_j/k_\mu})\),
with \(1\le j\le 10\), \(1\le\mu\le 4\).
\begin{enumerate}
\item
For the maximal TKT, \(\varkappa^{\prime\prime\prime}\sim (2000)\),
called \(\mathrm{a}.3^\ast\) in conjunction with ATI \(\lbrack 111,11,11,11\rbrack\),
we must have \(n\ne 0\), \(\tilde{n}\ne 0\) and by
\eqref{eqn:RankCat3Gph9}
\textit{wild} ranks \(r_j=3\) and indices \(I_j=3\) for all \(j=1,2,3,4,7,10\),
causing eight (because \(7\) is used twice over \(k_1\) and \(k_2\)
and \(10\) is used twice over \(k_3\) and \(k_4\))
minimal \(3\)-class numbers
\(h_3(B_j)=h_3(k_{pq})=h_3(\tilde{k}_{pq})=9\), by
\eqref{eqn:Cat3Gph5To9},
and ATI \(\mathrm{Cl}_3(B_j)\simeq (11)\),
characterisic for a group of coclass \(\mathrm{cc}(\mathfrak{M})=1\),
i.e. maximal class,
namely \(\mathfrak{M}\simeq\langle 81,7\rangle\).
However, the elementary tricyclic component \((111)\) of the ATI
requires \textit{tame} indices \(I_j=27\) for \(j=5,8\),
and thus \(\mathcal{N}=1\) for each \(1\le\mu\le 4\)
(because \(5\) is used twice over \(k_1\) and \(k_3\)
and \(8\) is used twice over \(k_2\) and \(k_4\)).
\item
Next, one of the total TK shrinks to a transposition,
\(\varkappa^{\prime\prime}\sim (2100)\), \(\mathrm{b}.10\),
which requires a group of coclass \(\mathrm{cc}(\mathfrak{M})\ge 2\),
implying, firstly, \textit{tame} indices \(I_j=27\) also for \(j=6,9\),
and thus \(\mathcal{N}=2\) for each \(1\le\mu\le 4\)
(because \(6\) is used twice over \(k_1\) and \(k_4\)
and \(9\) is used twice over \(k_2\) and \(k_3\)),
and, secondly, (from now on) necessarily both \(n=\tilde{n}=0\), implying
\textit{wild} ranks \(r_j=2\) indices \(I_j\in\lbrace 9,27\rbrace\) for all \(j=1,2,3,4,7,10\),
here \(I_j=9\), \(3\)-class numbers
\(h_3(B_j)=3\cdot h_3(k_{pq})=3\cdot h_3(\tilde{k}_{pq})=3\cdot 9=27\),
and thus ATI \(\alpha\sim\lbrack 111,111,21,21\rbrack\),
leading to \(\mathfrak{M}\simeq\langle 729,34..39\rangle\),
in view of Corollary
\ref{cor:ClassGroup33}.
\item
Now another total TK shrinks to a repetition,
\(\varkappa^\prime\sim (2110)\), \(\mathrm{d}.19\),
the first three wild indices for \(j=1,2,10\) become maximal \(I_j=27\),
causing \(\mathcal{N}=3\) and four
(because \(10\) is used twice over \(k_3\) and \(k_4\))
maximal new \(3\)-class numbers
\(h_3(B_j)=9\cdot h_3(\tilde{k}_{pq})=9\cdot 9=81\),
and thus ATI \(\alpha\sim\lbrack 111,111,22,21\rbrack\),
uniquely identifying \(\mathfrak{M}\simeq\langle 729,41\rangle\).
\item
Finally, for the minimal TKT, \(\varkappa_0\sim (2111)\), \(\mathrm{H}.4\),
the remaining three wild indices for \(j=3,4,7\) become maximal \(I_j=27\),
causing \(\mathcal{N}=4\) and four
(because \(7\) is used twice over \(k_1\) and \(k_2\))
maximal new \(3\)-class numbers
\(h_3(B_j)=9\cdot h_3(\tilde{k}_{pq})=9\cdot 9=81\),
and thus ATI \(\alpha\sim\lbrack 111,111,22,22\rbrack\),
enforcing a group of coclass \(\mathrm{cc}(\mathfrak{M})=3\)
namely \(\mathfrak{M}\simeq\langle 2187,65\vert 67\rangle\).
\end{enumerate}
\end{proof}


\begin{corollary}
\label{cor:UniCat3Gph9}
\textbf{(Uniformity of the quartet for \(\mathrm{III}.9\).)}
The components of the quartet, all with \(3\)-rank two,
share a common capitulation type
\(\varkappa(k_\mu)\),
common abelian type invariants
\(\alpha(k_\mu)\),
and a common second \(3\)-class group
\(\mathrm{Gal}(\mathrm{F}_3^2(k_\mu)/k_\mu)\),
for \(1\le\mu\le 4\).
\end{corollary}

\begin{proof}
This follows immediately from Theorem
\ref{thm:Cat3Gph9}.
\end{proof}

\begin{example}
\label{exm:Cat3Gph9}
Prototypes for Graph \(\mathrm{III}.9\) are
the minimal conductors for each scenario in Theorem
\ref{thm:Cat3Gph9}.
They have been found for all \(\mathcal{N}\in\lbrace 1,2,3,4\rbrace\).

There are
\textbf{regular} cases: 
\(c=16\,471\) with symbol
\(\lbrace 13\leftarrow 181\leftrightarrow 7\leftarrow 13\rbrace\), \(v^\ast=1\),
and \(\mathfrak{G}=\mathfrak{M}=\langle 81,7\rangle\);
\(c=89\,487\) with symbol
\(\lbrace 9\leftarrow 163\leftrightarrow 61\leftarrow 9\rbrace\), \(v^\ast=2\),
and \(\mathfrak{M}=\langle 729,41\rangle\);
\(c=109\,291\) with symbol
\(\lbrace 7\leftarrow 13\leftrightarrow 1201\leftarrow 7\rbrace\), \(v^\ast=2\),
and \(\mathfrak{M}=\langle 729,34..36\rangle\);
\(c=193\,921\) with symbol
\(\lbrace 7\leftarrow 13\leftrightarrow 2131\leftarrow 7\rbrace\), \(v^\ast=2\),
and \(\mathfrak{M}=\langle 729,37..39\rangle\);
and, with \textbf{extreme statistic delay},
\(c=707\,517\) with ordinal number \(145\), symbol
\(\lbrace 9\leftarrow 127\leftrightarrow 619\leftarrow 9\rbrace\), \(v^\ast=2\),
and \(\mathfrak{M}=\langle 2187,65\vert 67\rangle\)
with \(d_2(\mathfrak{M})=5\).

Only one \textbf{super-singular} case for \(c<2\cdot 10^5\): It is
\(c=197\,239\) with symbol
\(7\leftarrow 1483\leftrightarrow 19\leftarrow 7\), \(v^\ast=4\),
and \(\mathfrak{M}=\langle 729,37..39\rangle\).
Astonishingly, no bigger order and coclass of \(\mathfrak{M}\),
due to \(n\ne 0\). 
\end{example}


In Table
\ref{tbl:ProtoCat3Gph9},
we summarize the prototypes of Graph \(\mathrm{III}.9\)
in the same way as in Table
\ref{tbl:ProtoCat3Gph5}.

\renewcommand{\arraystretch}{1.1}

\begin{table}[ht]
\caption{Prototypes for Graph III.9}
\label{tbl:ProtoCat3Gph9}
\begin{center}
{\tiny
\begin{tabular}{|r|r||c|c|c|c|c|c|c|c|c|}
\hline
 No. & \(c\) & \(r\leftarrow p\leftrightarrow q\leftarrow r\) & \(v^\ast\) & \(v\) & \(m,n\) & \(\tilde{v}\) & \(\tilde{m},\tilde{n}\) & capitulation type & \(\mathfrak{M}\) & \(\ell_3(k)\) \\
\hline
   1 &  \(16\,471\) & \(13\leftarrow 181\leftrightarrow 7\leftarrow 13\) & \(1\) & \(2\) & \(1,1\) & \(2\) & \(1,1\) & \(\mathrm{a}.3^\ast\) & \(\langle 81,7\rangle\) & \(=2\) \\
  15 &  \(89\,487\) & \(9\leftarrow 163\leftrightarrow 61\leftarrow 9\)  & \(2\) & \(2\) & \(1,0\) & \(2\) & \(1,0\) & \(\mathrm{d}.19\) & \(\langle 729,41\rangle\) & \(\ge 2\) \\
  19 & \(109\,291\) & \(7\leftarrow 13\leftrightarrow 1201\leftarrow 7\) & \(2\) & \(2\) & \(1,0\) & \(2\) & \(1,0\) & \(\mathrm{b}.10\) & \(\langle 729,34..36\rangle\) & \(\ge 2\) \\
  28 & \(193\,921\) & \(7\leftarrow 13\leftrightarrow 2131\leftarrow 7\) & \(2\) & \(2\) & \(1,0\) & \(2\) & \(1,0\) & \(\mathrm{b}.10\) & \(\langle 729,37..39\rangle\) & \(\ge 2\) \\
  31 & \(197\,239\) & \(7\leftarrow 1483\leftrightarrow 19\leftarrow 7\) & \(4\) & \(3\) & \(0,1\) & \(3\) & \(0,1\) & \(\mathrm{b}.10\) & \(\langle 729,37..39\rangle\) & \(\ge 2\) \\
 145 & \(707\,517\) & \(9\leftarrow 127\leftrightarrow 619\leftarrow 9\) & \(2\) & \(2\) & \(1,0\) & \(2\) & \(1,0\) & \(\mathrm{H}.4\) & \(\langle 2187,65\vert 67\rangle\) & \(\ge 3\) \\
\hline
\end{tabular}
}
\end{center}
\end{table}


\section{Conclusions}
\label{s:Conclusion}

\noindent
In this work, we have seen that
order and structure of the \textbf{second} \(3\)-class group
\(\mathfrak{M}=\mathrm{Gal}(\mathrm{F}_3^2(k)/k)\)
of a cyclic cubic number field \(k\) with
conductor \(c=pqr\) divisible by
three prime(power)s \(p,q,r\) and
elementary bicyclic \(3\)-class group
\(\mathrm{Cl}_3(k)\simeq (3,3)\)
depends on arithmetical invariants of
other cyclic cubic \textbf{auxiliary fields},
associated with \(k\).
The field \(k\) is component of a quartet
\((k_1,\ldots,k_4)\)
of cyclic cubic fields sharing the common conductor \(c\).
The graph \(\lbrack p,q,r\rbrack_3\) which is combined by
the cubic residue symbols
\(\left(\frac{p}{q}\right)_3,\left(\frac{q}{p}\right)_3,
\left(\frac{p}{r}\right)_3,\left(\frac{r}{p}\right)_3,
\left(\frac{q}{r}\right)_3,\left(\frac{r}{q}\right)_3\)
decides whether one, or two, or no, component(s) of the quartet
have a \(3\)-class group of rank \(\varrho(k_\mu)=3\),
and accordingly the conductor \(c=pqr\) is called of Category
\(\mathrm{I}\), or \(\mathrm{II}\), or \(\mathrm{III}\).
For Category \(\mathrm{I}\),
the order of the \(3\)-class group
of the unique component with \(\varrho(k_{\mu_0})=3\)
is crucial.
For Category \(\mathrm{II}\),
the orders of both \(3\)-class groups
of the two components with
\(\varrho(k_{\mu_1})=\varrho(k_{\mu_2})=3\)
exert an impact.
For Category \(\mathrm{III}\),
the behavior is uniform with abelian \(\mathfrak{M}\simeq (3,3)\),
if \(\lbrack p,q,r\rbrack_3\) does not contain
mutual cubic residues (Graphs \(1\)--\(4\)),
otherwise there is exactly one pair \(p\leftrightarrow q\)
of mutual cubic residues (Graphs \(5\)--\(9\)),
and the auxiliary fields
with decisive \(3\)-class groups
are the two subfields
\(k_{pq}\) and \(\tilde{k}_{pq}\)
of the absolute genus field \(k^\ast\) of \(k\),
having the partial conductor \(pq\).
In each case, the \textbf{principal factors}
(norms of ambiguous principal ideals)
determine the fine structure in form of
uniform or non-uniform second \(3\)-class groups
\(\mathfrak{M}=\mathrm{Gal}(\mathrm{F}_3^2(k_\mu)/k_\mu)\).
Explicit numerical investigations indicate that
there is no upper bound for the orders
of the \(3\)-class groups
\(\mathrm{Cl}_3(k_{\mu_0})\),
respectively
\(\mathrm{Cl}_3(k_{\mu_1})\) and \(\mathrm{Cl}_3(k_{\mu_2})\),
respectively
\(\mathrm{Cl}_3(k_{pq})\) and \(\mathrm{Cl}_3(\tilde{k}_{pq})\).
In the regular situation, these orders are \(27\), respectively \(9\),
in the singular situation, they are \(81\), respectively \(27\),
but in the super-singular situation,
they are at least \(243\), respectively \(27\),
and the \textbf{orders may increase unboundedly}.
Concrete numerical examples are known with orders up to \(729\).

Bicyclic bicubic fields \(B_j\), \(j=1,\ldots,10\),
constitute the \textbf{capitulation targets}
of the cyclic cubic fields \(k_\mu\), \(\mu=1,\ldots,4\).
The introduction of important new concepts,
the \textbf{minimal and maximal capitulation type} (mTKT),
\(\varkappa_0\) and \(\varkappa_\infty\),
permitted recognition of common patterns
for several Graphs, partially in distinct Categories.

The four Graphs \(\mathrm{II}.1\), \(\mathrm{II}.2\),
\(\mathrm{III}.7\), \(\mathrm{III}.9\)
share the same ordered sequence of TKTs,
\(\varkappa_0\sim (2111)<(2110)<(2100)<(2000)\sim\varkappa_\infty\),
called \(\mathrm{H}.4\), \(\mathrm{d}.19\),
\(\mathrm{b}.10\), \(\mathrm{a}.3^\ast\),
although the proofs and details are quite different.
In terms of splitting prime ideals
\(\mathfrak{q}\mathcal{O}_{B_j}=\mathfrak{Q}_1\mathfrak{Q}_2\mathfrak{Q}_3\),
\(\mathfrak{r}\mathcal{O}_{B_\ell}=\mathfrak{R}_1\mathfrak{R}_2\mathfrak{R}_3\),
all these TKTs contain a crucial \textbf{transposition},
due to elementary tricyclic \(3\)-class groups
\(\mathrm{Cl}_3(B_j)=\langle\lbrack\mathfrak{Q}_1,\rbrack\lbrack\mathfrak{Q}_2,\rbrack\lbrack\mathfrak{Q}_3\rbrack\rangle\),
\(\mathrm{Cl}_3(B_\ell)=\langle\lbrack\mathfrak{R}_1,\rbrack\lbrack\mathfrak{R}_2,\rbrack\lbrack\mathfrak{R}_3\rbrack\rangle\),
and twisted capitulation kernels
\(\ker(T_{B_j/k_\mu})=\langle\lbrack\mathfrak{r}\rbrack\rangle\),
\(\ker(T_{B_\ell/k_\mu})=\langle\lbrack\mathfrak{q}\rbrack\rangle\),
which restricts the group \(\mathfrak{M}\) to descendants of \(\langle 243,3\rangle\)
(except \(\langle 81,7\rangle\), where a total transfer kernel hides the transposition).

Similarly, the two graphs \(\mathrm{I}.1\), \(\mathrm{I}.2\)
admit another characteristic ordered sequence of TKTs,
\(\varkappa_0\sim (4231)<(0231)<(0200),(0001)<(0000)\sim\varkappa_\infty\),
called \(\mathrm{G}.16\), \(\mathrm{c}.21\),
\(\mathrm{a}.2\), \(\mathrm{a}.3\), \(\mathrm{a}.1\),
with \textbf{two fixed points},
which restrict the group \(\mathfrak{M}\) to descendants of \(\langle 243,8\rangle\)
(except \(\langle 81,8\rangle\), \(\langle 81,10\rangle\),
\(\langle 243,25\rangle\), \(\langle 243,27\rangle\),
where total transfer kernels partially or completely hide the fixed points).

A remarkable outsider is Graph \(\mathrm{III}.8\)
with a veritable wealth of exotic capitulation types,
but restricted to the unusual maximal TKT
\(\varkappa_\infty\sim (2100)\), \(\mathrm{b}.10\),
forced by mandatory transposition.

Due to the lack of cubic residue conditions between the prime divisors
of the conductor \(c=pqr\),
two Graphs \(\mathrm{I}.1\), \(\mathrm{III}.5\)
admit the absolute maximum of all TKTs \(\varkappa=(0000)\)
(non-abelian!). 

It might be worth one's while to point out
that a glance at \(\alpha_2\) in Tables
\ref{tbl:Metabelian33}
and
\ref{tbl:NonMetabelian33}
reveals that the commutator subgroup of all encountered
second \(3\)-class groups \(\mathfrak{M}\), respectively
\(3\)-class tower groups \(\mathfrak{G}\),
has order \(\#(\mathfrak{M}^\prime)\ge 9\),
respectively \(\#(\mathfrak{G}^\prime)\ge 9\),
which means that the class number of the 
Hilbert \(3\)-class field \(\mathrm{F}_3^1(k)\)
is divisible by \(9\),
for all cyclic cubic fields \(k\),
with the exception of \(t=1\), the regular cases for \(t=2\),
and the Graphs \(1,\ldots,4\) of Category \(\mathrm{III}\)
for \(t=3\).

For Category \(\mathrm{I}\) and \(\mathrm{II}\),
we expect a rather rigid impact
of the groups \(\mathfrak{M}=\mathrm{Gal}(\mathrm{F}_3^2(k_\mu)/k_\mu)\)
for \(\mathrm{Cl}_3(k_\mu)\simeq (3,3)\)
on the groups \(\mathrm{Gal}(\mathrm{F}_3^2(k_\nu)/k_\nu)\)
for \(\mathrm{Cl}_3(k_\nu)\simeq (3,3,3)\),
as suggested by the numerous tables in  
\cite{Ma2022}.
This research line will be pursued further in a forthcoming paper.


\section{Acknowledgements}
\label{s:Acknowledgements}

\noindent
Both authors gratefully acknowledge computational aid
by Bill Allombert, who used the algorithm of Aurel Page
\cite{Pg2020}
to determine the \(3\)-class groups of Hilbert \(3\)-class fields
and thus \(\alpha_2\) in Corollary
\ref{cor:ClassGroup33}.
The second author acknowledges that his research was supported by
the Austrian Science Fund (FWF): projects J0497-PHY and P26008-N25,
and by the Research Executive Agency of the European Union (EUREA):
project Horizon Europe 2021--2027.



\end{document}